\documentclass{amsart}%
\usepackage{amsmath}
\usepackage{amsfonts}
\usepackage{amssymb}
\usepackage{graphicx}%
\tolerance=1600 \hbadness=10000 \theoremstyle{plain}
\newtheorem{theorem}{Theorem}[section]
\newtheorem{corollary}[theorem]{Corollary}
\newtheorem{lemma}[theorem]{Lemma}

\newtheorem{proposition}[theorem]{Proposition}

\newtheorem{solution}[theorem]{Solution}

\renewenvironment{proof}[1][Proof]{\textbf{#1.} }{\hfill \rule{0.5em}{0.5em}}

\theoremstyle{definition}
\newtheorem{definition}[theorem]{Definition}
\newtheorem{notation}[theorem]{Notation}

\newtheorem{example}[theorem]{Example}
\newtheorem{exercise}[theorem]{Exercise}

\theoremstyle{remark}
\newtheorem{remark}[theorem]{Remark}

\newtheorem{fact}[theorem]{Fact}
\theoremstyle{plain}

\theoremstyle{definition}
\newtheorem{assumption}[]{Assumption}
\numberwithin{equation}{section}

\begin{document}

\title{Curved Wiener Space Analysis}

%\title{A Primer on Riemannian Geometry and Stochastic Analysis on Path Spaces}
%\title{Riemannian Geometry, Stochastic Analysis and Path Spaces}

\author{Bruce K. Driver}

\address{Department of Mathematics, 0112, University
of California at San Diego, La Jolla, CA, 92093-0112}

%  Current address

%\curraddr{Department of Mathematics, 0112, University of California at San Diego,
%La Jolla, CA, 92093-0112}

\email{driver@math.ucsd.edu}
\thanks{This research was partially supported by NSF Grants
DMS 96-12651, DMS 99-71036 and DMS 0202939. This article will
appear in ``Real and Stochastic Analysis: New Perspectives.''}

%\date{September 5, 1995.}

%\date{\emph{File:\jobname{.tex}}\qquad Last revised: \today}
\maketitle

\setcounter{tocdepth}{1}

\tableofcontents

\section{Introduction\label{s.1}}

These notes represent a much expanded and updated version of the
\textquotedblleft mini course\textquotedblright\ that the author gave at the
ETH (Z\"{u}rich) and the University of Z\"{u}rich in February of 1995. The
purpose of these notes is to first provide some basic background to Riemannian
geometry and stochastic calculus on manifolds and then to cover some of the
more recent developments pertaining to analysis on \textquotedblleft curved
Wiener spaces.\textquotedblright\ Essentially no differential geometry is
assumed. However, it is assumed that the reader is comfortable with stochastic
calculus and differential equations on Euclidean spaces. Here is a brief
description of what will be covered in the text below.

Section \ref{s.2} is a basic introduction to differential geometry through
imbedded submanifolds. Section \ref{s.3} is an introduction to the Riemannian
geometry that will be needed in the sequel. Section \ref{s.4} records a number
of results pertaining to flows of vector fields and \textquotedblleft Cartan's
rolling map.\textquotedblright\ The stochastic version of these results will
be important tools in the sequel. Section \ref{s.5} is a rapid introduction to
stochastic calculus on manifolds and related geometric constructions. Section
\ref{s.6} briefly gives applications of stochastic calculus on manifolds to
representation formulas for derivatives of heat kernels. Section \ref{s.7} is
devoted to the study of the calculus and integral geometry associated with the
path space of a Riemannian manifold equipped with \textquotedblleft Wiener
measure.\textquotedblright\ In particular, quasi-invariance, Poincar\'{e} and
logarithmic Sobolev inequalities are developed for the Wiener measure on path
spaces in this section. Section \ref{s.8} is a short introduction to
Malliavin's probabilistic methods for dealing with hypoelliptic diffusions.
The appendix in section \ref{s.9} records some basic martingale and stochastic
differential equation estimates which are mostly used in section \ref{s.8}.

Although the majority of these notes form a survey of known results, many
proofs have been cleaned up and some proofs are new. Moreover, Section
\ref{s.8} is written using the geometric language introduced in these notes
which is not completely standard in the literature. I have also tried (without
complete success) to give an overview of many of the major techniques which
have been used to date in this subject. Although numerous references are given
to the literature, the list is far from complete. I apologize in advance to
anyone who feels cheated by not being included in the references. However, I
do hope the list of references is sufficiently rich that the interested reader
will be able to find additional information by looking at the related articles
and the references that they contain.

\emph{Acknowledgement:\/} It is pleasure to thank Professor A. Sznitman and
the ETH for their hospitality and support and the opportunity to give the
talks which started these notes. I also would like to thank Professor E.
Bolthausen for his hospitality and his role in arranging the first lecture to
be held at University of Z\"{u}rich.

\section{Manifold Primer\label{s.2}}

\textbf{Conventions:}

\begin{enumerate}
\item If $A,B$ are linear operators on some vector space, then $\left[
A,B\right]  :=AB-BA$ is the \textbf{commutator }of $A$ and $B.$

\item If $X$ is a topological space we will write $A\subset_{o}X,$ $A\sqsubset
X$ and $A\sqsubset\sqsubset X$ to mean $A$ is an open, closed, and
respectively a compact subset of $X.$

\item Given two sets $A$ and $B,$ the notation $f:A\rightarrow B$ will mean
that $f$ is a function from a subset $\mathcal{D}(f)\subset A$ to $B.$ (We
will allow $\mathcal{D}(f)$ to be the empty set.) The set $\mathcal{D}%
(f)\subset A$ is called the domain of $f$ and the subset $\mathcal{R}%
(f):=f(\mathcal{D}(f))\subset B$ is called the range of $f.$ If $f$ is
injective, let $f^{-1}:B\rightarrow A$ denote the inverse function with domain
$\mathcal{D}(f^{-1})=\mathcal{R}(f)$ and range $\mathcal{R}(f^{-1}%
)=\mathcal{D}(f).$ If $f:A\rightarrow B$ and $g:B\rightarrow C,$ then $g\circ
f$ denotes the composite function from $A$ to $C$ with domain $\mathcal{D}%
(g\circ f):=f^{-1}(\mathcal{D}(g))$ and range $\mathcal{R}(g\circ f):=g\circ
f(\mathcal{D}(g\circ f))=g(\mathcal{R}(f)\cap\mathcal{D}(g)).$ \medskip
\end{enumerate}

\begin{notation}
\label{n.2.1}Throughout these notes, let $E$ and $V$ denote finite dimensional
vector spaces. A function $F:E\rightarrow V$ is said to be smooth if
$\mathcal{D}(F)$ is open in $E$ ($\mathcal{D}(F)=\emptyset$ is allowed) and
$F:\mathcal{D}(F)\rightarrow V$ is infinitely differentiable. Given a
\textbf{smooth} function $F:E\rightarrow V,$ let $F^{\prime}(x)$ denote the
differential of $F$ at $x\in\mathcal{D}(F).$ Explicitly, $F^{\prime
}(x)=DF\left(  x\right)  $ denotes the linear map from $E$ to $V$ determined
by
\begin{equation}
DF\left(  x\right)  a=F^{\prime}(x)a:=\frac{d}{dt}|_{0}F(x+ta)~\forall~a\in E.
\label{e.2.1}%
\end{equation}
We also let
\begin{equation}
F^{\prime\prime}\left(  x\right)  \left(  v,w\right)  =F^{\prime\prime}\left(
x\right)  \left(  v,w\right)  :=\left(  \partial_{v}\partial_{w}F\right)
\left(  x\right)  =\frac{d}{dt}|_{0}\frac{d}{ds}|_{0}F\left(  x+tv+sw\right)
. \label{e.2.2}%
\end{equation}

\end{notation}

\subsection{Imbedded Submanifolds}

Rather than describe the most abstract setting for Riemannian geometry, for
simplicity we choose to restrict our attention to imbedded submanifolds of a
Euclidean space $E=\mathbb{R}^{N}.$ \footnote{Because of the Whitney imbedding
theorem (see for example Theorem 6-3 in Auslander and MacKenzie \cite{AMa}),
this is actually not a restriction.} We will equip $\mathbb{R}^{N}$ with the
standard inner product,%
\[
\langle a,b\rangle=\langle a,b\rangle_{\mathbb{R}^{N}}:=\sum_{i=1}^{N}%
a_{i}b_{i}.
\]
In general, we will denote inner products in these notes by $\langle
\cdot,\cdot\rangle.$

\begin{definition}
\label{d.2.2} A subset $M$ of $E$ (see Figure \ref{f.1}) is a $d$\textbf{ --
dimensional imbedded submanifold }(without boundary) of $E$ iff for all $m\in
M,$ there is a function $z:E\rightarrow\mathbb{R}^{N}$ such that:

\begin{enumerate}
\item $\mathcal{D}(z)$ is an open neighborhood of $E$ containing $m,$

\item $\mathcal{R}(z)$ is an open subset of $\mathbb{R}^{N}$,

\item $z:\mathcal{D}(z)\rightarrow\mathcal{R}(z)$ is a diffeomorphism (a
smooth invertible map with smooth inverse), and

\item $z(M\cap\mathcal{D}(z))=\mathcal{R}(z)\cap(\mathbb{R}^{d}\times
\{0\})\subset\mathbb{R}^{N}.$
\end{enumerate}

(We write $M^{d}$ if we wish to emphasize that $M$ is a $d$ -- dimensional manifold.)%

%TCIMACRO{\FRAME{ftbphFU}{3.2624in}{2.0349in}{0pt}{\Qcb{An imbedded one
%dimensional submanifold in $\mathbb{R}^{2}.$}}{\Qlb{f.1}}{f1.eps}%
%width 3.2624in;  height 2.0349in;  depth 0pt;  original-width 5.715in;
%original-height 4.4857in;  cropleft "0";  croptop "1";  cropright "1";
%cropbottom "0";  filename 'F1.eps';file-properties "XNPEU";}%
%}}%
%BeginExpansion
\begin{figure}
[ptbh]
\begin{center}
\includegraphics[
height=2.0349in,
width=3.2624in
]%
{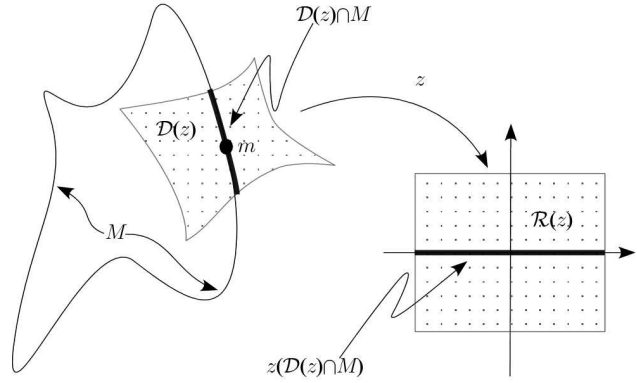}%
\caption{An imbedded one dimensional submanifold in $\mathbb{R}^{2}.$}%
\label{f.1}%
\end{center}
\end{figure}
%EndExpansion

\end{definition}

\begin{notation}
\label{n.2.3}Given an imbedded submanifold and diffeomorphism $z$ as in the
above definition, we will write $z=(z_{<},z_{>})$ where $z_{<}$ is the first
$d$ components of $z$ and $z_{>}$ consists of the last $N-d$ components of
$z.$ Also let $x:M\rightarrow\mathbb{R}^{d}$ denote the function defined by
$\mathcal{D}(x):=M\cap\mathcal{D}(z)$ and $x:=z_{<}|_{\mathcal{D}(x)}.$ Notice
that $\mathcal{R}(x):=x(\mathcal{D}(x))$ is an open subset of $\mathbb{R}^{d}$
and that $x^{-1}:\mathcal{R}(x)\rightarrow\mathcal{D}(x),$ thought of as a
function taking values in $E,$ is smooth. The bijection $x:\mathcal{D}%
(x)\rightarrow\mathcal{R}(x)$ is called a \textbf{chart} on $M.$ Let
$\mathcal{\ A}=\mathcal{A}(M)$ denote the collection of charts on $M.$ The
collection of charts $\mathcal{A}=\mathcal{A}(M)$ is often referred to as an
\textbf{atlas} for $M.$
\end{notation}

\begin{remark}
\label{r.2.4}The imbedded submanifold $M$ is made into a topological space
using the induced topology from $E.$ With this topology, each chart
$x\in\mathcal{A}(M)$ is a homeomorphism from $\mathcal{D}(x)\subset_{o}M$ to
$\mathcal{R}(x)\subset_{o}\mathbb{R}^{d}.$
\end{remark}

\begin{theorem}
[A Basic Construction of Manifolds]\label{t.2.5} Let $F:E\rightarrow
\mathbb{R}^{N-d}$ be a smooth function and $M:=F^{-1}(\{0\})\subset E$ which
we assume to be non-empty. Suppose that $F^{\prime}(m):E\rightarrow
\mathbb{R}^{N-d}$ is surjective for all $m\in M.$ Then $M$ is a $d$ --
dimensional imbedded submanifold of $E.$
\end{theorem}

\begin{proof}
Let $m\in M,$ we will begin by constructing a smooth function $G:E\rightarrow
\mathbb{R}^{d}$ such that $(G,F)^{\prime}(m):E\rightarrow\mathbb{R}%
^{N}=\mathbb{R}^{d}\times\mathbb{R}^{N-d}$ is invertible. To do this, let
$X=\operatorname*{Nul}(F^{\prime}(m))$ and $Y$ be a complementary subspace so
that $E=X\oplus Y$ and let $P:E\rightarrow X$ be the associated projection
map, see Figure \ref{fig.2}. Notice that $F^{\prime}(m):Y\rightarrow
\mathbb{R}^{N-d}$ is a linear isomorphism of vector spaces and hence
\[
\dim(X)=\dim(E)-\dim(Y)=N-(N-d)=d.
\]
In particular, $X$ and $\mathbb{R}^{d}$ are isomorphic as vector spaces. Set
$G(m)=APm$ where $A:X\rightarrow\mathbb{R}^{d}$ is an arbitrary but fixed
linear isomorphism of vector spaces. Then for $x\in X$ and $y\in Y,$
\begin{align*}
(G,F)^{\prime}(m)(x+y)  &  =(G^{\prime}(m)(x+y),F^{\prime}(m)(x+y))\\
&  =(AP(x+y),F^{\prime}(m)y)=(Ax,F^{\prime}(m)y)\in\mathbb{R}^{d}%
\times\mathbb{R}^{N-d}%
\end{align*}
from which it follows that $(G,F)^{\prime}(m)$ is an isomorphism.%

%TCIMACRO{\FRAME{ftbphFU}{3.9487in}{1.8028in}{0pt}{\Qcb{Constructing charts for
%$M$ using the inverse function theorem. For simplicity of the drawing, $m\in
%M$ is assumed to be the origin of $E=X\oplus Y.$}}{\Qlb{fig.2}}{ff.eps}%
%width 3.9487in;  height 1.8028in;  depth 0pt;  original-width 5.4964in;
%original-height 2.4941in;  cropleft "0";  croptop "1";  cropright "1";
%cropbottom "0";  filename 'FF.eps';file-properties "XNPEU";}%
%}}%
%BeginExpansion
\begin{figure}
[ptbh]
\begin{center}
\includegraphics[
height=1.8028in,
width=3.9487in
]%
{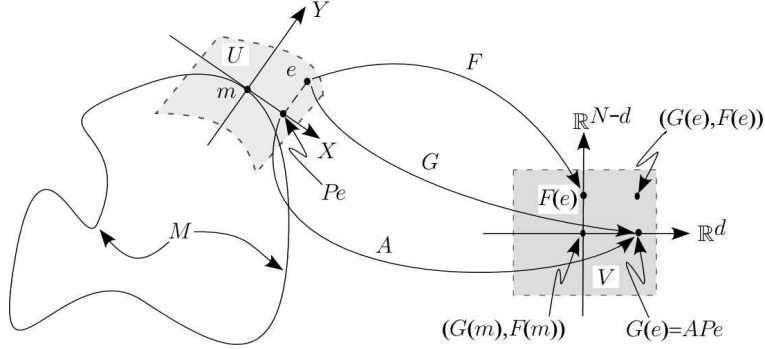}%
\caption{Constructing charts for $M$ using the inverse function theorem. For
simplicity of the drawing, $m\in M$ is assumed to be the origin of $E=X\oplus
Y.$}%
\label{fig.2}%
\end{center}
\end{figure}
%EndExpansion

By the inverse function theorem, there exists a neighborhood $U\subset_{o}E$
of $m$ such that $V:=(G,F)(U)\subset_{o}\mathbb{R}^{N}$ and
$(G,F):U\rightarrow V$ is a diffeomorphism. Let $z=(G,F)$ with $\mathcal{D}%
(z)=U$ and $\mathcal{R}(z)=V.$ Then $z$ is a chart of $E$ about $m$ satisfying
the conditions of Definition \ref{d.2.2}. Indeed, items 1) -- 3) are clear by
construction. If $p\in M\cap\mathcal{D}(z)$ then $z(p)=(G(p),F(p))=(G(p),0)\in
\mathcal{R}(z)\cap(\mathbb{R}^{d}\times\{0\}).$ Conversely, if $p\in
\mathcal{D}(z)$ is a point such that $z(p)=(G(p),F(p))\in\mathcal{R}%
(z)\cap(\mathbb{R}^{d}\times\{0\}),$ then $F(p)=0$ and hence $p\in
M\cap\mathcal{D}(z);$ so item 4) of Definition \ref{d.2.2} is verified.
\end{proof}

\begin{example}
\label{ex.2.6}Let $gl(n,\mathbb{R})$ denote the set of all $n\times n$ real
matrices. The following are examples of imbedded submanifolds.

\begin{enumerate}
\item Any open subset $M$ of $E$.

\item The graph,
\[
\Gamma\left(  f\right)  :=\left\{  \left(  x,f(x)\right)  \in\mathbb{R}%
^{d}\times\mathbb{R}^{N-d}:x\in\mathcal{D}\left(  f\right)  \right\}
\subset\mathcal{D}\left(  f\right)  \times\mathbb{R}^{N-d}\subset
\mathbb{R}^{N},
\]
of any smooth function $f:\mathbb{R}^{d}\rightarrow\mathbb{R}^{N-d}$ as can be
seen by applying Theorem \ref{t.2.5} with $F\left(  x,y\right)  :=y-f\left(
x\right)  .$ In this case it would be a good idea for the reader to produce an
explicit chart $z$ as in Definition \ref{d.2.2} such that $\mathcal{D}\left(
z\right)  =\mathcal{R}\left(  z\right)  =\mathcal{D}\left(  f\right)
\times\mathbb{R}^{N-d}.$

\item The unit sphere, $S^{N-1}:=\{x\in\mathbb{R}^{N}:\langle x,x\rangle
_{\mathbb{R}^{N}}=1\},$ as is seen by applying Theorem \ref{t.2.5} with
$E=\mathbb{R}^{\mathbb{N}}$ and $F(x):=\langle x,x\rangle_{\mathbb{R}^{N}}-1.$
Alternatively, express $S^{N-1}$ locally as the graph of smooth functions and
then use item 2.

\item $GL(n,\mathbb{R}):=\{g\in gl(n,\mathbb{R})|\det(g)\neq0\},$ see item 1.

\item $SL(n,\mathbb{R}):=\{g\in gl(n,\mathbb{R})|\det(g)=1\}$ as is seen by
taking $E=gl(n,\mathbb{R})$ and $F(g):=\det(g)$ and then applying Theorem
\ref{t.2.5} with the aid of Lemma \ref{l.2.7} below.

\item $O(n):=\{g\in gl(n,\mathbb{R})|g^{\mathrm{tr}}g=I\}$ where
$g^{\mathrm{tr}}$ denotes the transpose of $g.$ In this case take
$F(g):=g^{\mathrm{tr}}g-I$ thought of as a function from $E=gl(n,\mathbb{R})$
to$\mathcal{\ S}(n),$ where
\[
\mathcal{S}(n):=\left\{  A\in gl(n,\mathbb{R}):A^{\mathrm{tr}}=A\right\}
\]
is the subspace of symmetric matrices. To show $F^{\prime}(g)$ is surjective,
show
\[
F^{\prime}(g)(gB)=B+B^{\mathrm{tr}}\text{ for all }g\in O(n)\text{ and }B\in
gl(n,\mathbb{R}).
\]

\item $SO(n):=\{g\in O(n)|\det(g)=1\}$, an open subset of $O(n).$

\item $M\times N\subset E\times V,$ where $M$ and $N$ are imbedded
submanifolds of $E$ and $V$ respectively. The reader should verify this by
constructing appropriate charts for $E\times V$ by taking \textquotedblleft
tensor\textquotedblright\ products of the charts for $E$ and $V$ associated to
$M$ and $N$ respectively.

\item The $n$ -- dimensional torus,%
\[
T^{n}:=\{z\in\mathbb{C}^{n}:|z_{i}|=1\text{\textrm{\ for }}i=1,2,\ldots
,n\}=(S^{1})^{n},
\]
where $z=\left(  z_{1},\dots,z_{n}\right)  $ and $\left\vert z_{i}\right\vert
=\sqrt{z_{i}\bar{z}_{i}}.$ This follows by induction using items 3. and 8.
Alternatively apply Theorem \ref{t.2.5} with $F\left(  z\right)  :=\left(
\left\vert z_{1}\right\vert ^{2}-1,\dots,\left\vert z_{n}\right\vert
^{2}-1\right)  .$
\end{enumerate}
\end{example}

\begin{lemma}
\label{l.2.7}Suppose $g\in GL(n,\mathbb{R})$ and $A\in gl(n,\mathbb{R}),$ then%
\begin{equation}
\det{}^{\prime}(g)A=\det(g)\text{$\operatorname*{tr}$}(g^{-1}A). \label{e.2.3}%
\end{equation}

\end{lemma}

\begin{proof}
By definition we have
\[
\det\,^{\prime}(g)A=\frac{d}{dt}|_{0}\det(g+tA)=\det(g)\frac{d}{dt}|_{0}%
\det(I+tg^{-1}A).
\]
So it suffices to prove $\frac{d}{dt}|_{0}\det
(I+tB)=\mathrm{\operatorname*{tr}}(B)$ for all matrices $B.$ If $B$ is upper
triangular, then $\det(I+tB)=\prod_{i=1}^{n}(1+tB_{ii})$ and hence by the
product rule,
\[
\frac{d}{dt}|_{0}\det(I+tB)=\sum_{i=1}^{n}B_{ii}=\operatorname*{tr}(B).
\]
This completes the proof because; 1) every matrix can be put into upper
triangular form by a similarity transformation, and 2) \textquotedblleft$\det
$\textquotedblright\ and \textquotedblleft$\operatorname*{tr}$%
\textquotedblright\ are invariant under similarity transformations.
\end{proof}

\begin{definition}
\label{d.2.8} Let $E$ and $V$ be two finite dimensional vector spaces and
$M^{d}\subset E$ and $N^{k}\subset V$ be two imbedded submanifolds. A function
$f:M\rightarrow N$ is said to be \textbf{smooth} if for all charts
$x\in\mathcal{A}(M)$ and $y\in\mathcal{A}(N)$ the function $y\circ f\circ
x^{-1}:\mathbb{R}^{d}\rightarrow\mathbb{R}^{k}$ is smooth.
\end{definition}

\begin{exercise}
\ \label{exr.2.9} Let $M^{d}\subset E$ and $N^{k}\subset V$ be two imbedded
submanifolds as in Definition \ref{d.2.8}.

\begin{enumerate}
\item Show that a function $f:\mathbb{R}^{k}\rightarrow M$ is smooth iff $f$
is smooth when thought of as a function from $\mathbb{R}^{k}$ to $E.$

\item If $F:E\rightarrow V$ is a smooth function such that $F(M\cap
\mathcal{D}(F))\subset N,$ show that $f:=F|_{M}:M\rightarrow N$ is smooth.

\item Show the composition of smooth maps between imbedded submanifolds is smooth.
\end{enumerate}
\end{exercise}

\begin{proposition}
\label{p.2.10}Assuming the notation in Definition \ref{d.2.8}, a function
$f:M\rightarrow N$ is smooth iff there is a smooth function $F:E\rightarrow V$
such that $f=F|_{M}.$
\end{proposition}

\begin{proof}
(Sketch.) Suppose that $f:M\rightarrow N$ is smooth, $m\in M$ and $n=f(m).$
Let $z$ be as in Definition \ref{d.2.2} and $w$ be a chart on $N$ such that
$n\in\mathcal{D}(w).$ By shrinking the domain of $z$ if necessary, we may
assume that $\mathcal{R}(z)=U\times W$ where $U\subset_{o}\mathbb{R}^{d}$ and
$W\subset_{o}\mathbb{R}^{N-d}$ in which case $z(M\cap\mathcal{D}%
(z))=U\times\left\{  0\right\}  .$ For $\xi\in\mathcal{D}(z),$ let
$F(\xi):=f(z^{-1}(z_{<}(\xi),0))$ with $z=\left(  z_{<},z_{>}\right)  $ as in
Notation \ref{n.2.3}. Then $F:\mathcal{D}(z)\rightarrow N$ is a smooth
function such that $F|_{M\cap\mathcal{D}(z)}=f|_{M\cap\mathcal{D}(z)}.$ The
function $F$ is smooth. Indeed, letting $x=z_{<}|_{\mathcal{D}(z)\cap M},$%
\[
w_{<}\circ F=w_{<}\circ f(z^{-1}(z_{<}(\xi),0))=w_{<}\circ f\circ x^{-1}%
\circ(z_{<}(\cdot),0)
\]
which, being the composition of the smooth maps $w_{<}\circ f\circ x^{-1}$
(smooth by assumption) and $\xi\rightarrow(z_{<}(\xi),0),$ is smooth as well.
Hence by definition, $F$ is smooth as claimed. Using a standard partition of
unity argument (which we omit), it is possible to piece this local argument
together to construct a globally defined smooth function $F:E\rightarrow V$
such that $f=F|_{M}.$
\end{proof}

\begin{definition}
\label{d.2.11}A function $f:M\rightarrow N$ is a \textbf{diffeomorphism} if
$f$ is smooth and has a smooth inverse. The set of diffeomorphisms
$f:M\rightarrow M$ is a group under composition which will be denoted by
$\mathrm{Diff}(M).$
\end{definition}

\subsection{Tangent Planes and Spaces}

\begin{figure}
[ptbh]
\begin{center}
\includegraphics[
height=1.4723in,
width=2.154in
]%
{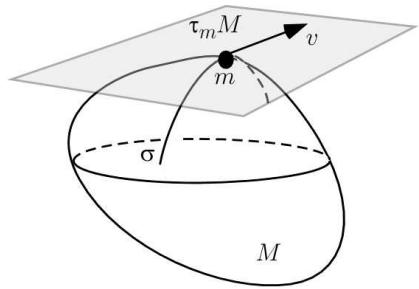}%
\caption{Tangent plane, $\tau_{m}M,$ to $M$ at $m$ and a vector, $v,$ in
$\tau_{m}M.$}%
\label{fig.3}%
\end{center}
\end{figure}

\begin{definition}
\label{d.2.12}Given an imbedded submanifold $M\subset E$ and $m\in M,$ let
$\tau_{m}M\subset E$ denote the collection of all vectors $v\in E$ such there
exists a smooth path $\sigma:(-\varepsilon,\varepsilon)\rightarrow M$ with
$\sigma(0)=m$ and $v=\frac{d}{ds}|_{0}\sigma(s).$ The subset $\tau_{m}M$ is
called the \textbf{tangent plane} to $M$ at $m$ and $v\in\tau_{m}M$ is called
a \textbf{tangent vector}, see Figure \ref{fig.3}.
\end{definition}

%

%TCIMACRO{\FRAME{ftbphFU}{2.154in}{1.4723in}{0pt}{\Qcb{Tangent plane, $\tau
%_{m}M,$ to $M$ at $m$ and a vector, $v,$ in $\tau_{m}M.$}}{\Qlb{fig.3}%
%}{f2.eps}{\special{ language "Scientific Word";  type "GRAPHIC";
%width 2.154in;  height 1.4723in;  depth 0pt;  original-width 2.7474in;
%original-height 1.8705in;  cropleft "0";  croptop "1";  cropright "1";
%cropbottom "0";  filename 'f2.eps';file-properties "XNPEU";}%
%}}%
%BeginExpansion
%EndExpansion

\begin{theorem}
\label{t.2.13}For each $m\in M,$ $\tau_{m}M$ is a $d$ -- dimensional subspace
of $E$. If $z:E\rightarrow\mathbb{R}^{N}$ is as in Definition \ref{d.2.2},
then $\tau_{m}M=\operatorname*{Nul}(z_{>}^{\prime}(m)).$ If $x$ is a chart on
$M$ such that $m\in\mathcal{D}(x),$ then
\[
\{\frac{d}{ds}|_{0}x^{-1}(x(m)+se_{i})\}_{i=1}^{d}%
\]
is a basis for $\tau_{m}M,$ where $\{e_{i}\}_{i=1}^{d}$ is the standard basis
for $\mathbb{R}^{d}.$
\end{theorem}

\begin{proof}
Let $\sigma:(-\varepsilon,\varepsilon)\rightarrow M$ be a smooth path with
$\sigma(0)=m$ and $v=\frac{d}{ds}|_{0}\sigma(s)$ and $z$ be a chart (for $E)$
around $m$ as in Definition \ref{d.2.2} such that $x=z_{<}.$ Then
$z_{>}(\sigma(s))=0$ for all $s$ and therefore,
\[
0=\frac{d}{ds}|_{0}z_{>}(\sigma(s))=z_{>}^{\prime}(m)v
\]
which shows that $v\in\operatorname*{Nul}(z_{>}^{\prime}(m)),$ i.e. $\tau
_{m}M\subset\operatorname*{Nul}(z_{>}^{\prime}(m)).$

Conversely, suppose that $v\in\operatorname*{Nul}(z_{>}^{\prime}(m)).$ Let
$w=z_{<}^{\prime}(m)v\in\mathbb{R}^{d}\ $and $\sigma(s):=x^{-1}(z_{<}%
(m)+sw)\in M$ -- defined for $s$ near $0.$ Differentiating the identity
$z^{-1}\circ z=id$ at $m$ shows
\[
\left(  z^{-1}\right)  ^{\prime}(z(m))z^{\prime}(m)=I.
\]
Therefore,%
\begin{align*}
\sigma^{\prime}(0)  &  =\frac{d}{ds}|_{0}x^{-1}(z_{<}(m)+sw)=\frac{d}{ds}%
|_{0}z^{-1}(z_{<}(m)+sw,0)\\
&  =\left(  z^{-1}\right)  ^{\prime}((z_{<}(m),0))(z_{<}^{\prime}(m)v,0)\\
&  =\left(  z^{-1}\right)  ^{\prime}((z_{<}(m),0))(z_{<}^{\prime}%
(m)v,z_{>}^{\prime}(m)v)\\
&  =\left(  z^{-1}\right)  ^{\prime}(z(m))z^{\prime}(m)v=v,
\end{align*}
and so by definition $v=\sigma^{\prime}(0)\in\tau_{m}M.$ We have now shown
$\operatorname*{Nul}(z_{>}^{\prime}(m))\subset\tau_{m}M$ which completes the
proof that $\tau_{m}M=\operatorname*{Nul}(z_{>}^{\prime}(m)).$

Since $z_{<}^{\prime}(m):\tau_{m}M\rightarrow\mathbb{R}^{d}$ is a linear
isomorphism, the above argument also shows%
\[
\frac{d}{ds}|_{0}x^{-1}(x(m)+sw)=\left(  z_{<}^{\prime}(m)|_{\tau_{m}%
M}\right)  ^{-1}w\in\tau_{m}M\ \forall~w\in\mathbb{R}^{d}.
\]
In particular it follows that
\[
\{\frac{d}{ds}|_{0}x^{-1}(x(m)+se_{i})\}_{i=1}^{d}=\{\left(  z_{<}^{\prime
}(m)|_{\tau_{m}M}\right)  ^{-1}e_{i}\}_{i=1}^{d}%
\]
is a basis for $\tau_{m}M,$ see Figure \ref{fig.4} below.
\end{proof}

The following proposition is an easy consequence of Theorem \ref{t.2.13} and
the proof of Theorem \ref{t.2.5}.

\begin{proposition}
\label{p.2.14}Suppose that $M$ is an imbedded submanifold constructed as in
Theorem \ref{t.2.5}. Then $\tau_{m}M=\operatorname*{Nul}$$\left(  F^{\prime
}(m)\right)  .$
\end{proposition}

\begin{exercise}
\label{exr.2.15}Show:

\begin{enumerate}
\item $\tau_{m}M=E,$ if $M$ is an open subset of $E.$

\item $\tau_{g}GL(n,\mathbb{R})=gl(n,\!\mathbb{R}),$ for all $g\in
GL(n,\mathbb{R}).$

\item $\tau_{m}S^{N-1}=\{m\}^{\perp}$ for all $m\in S^{N-1}.$

\item Let $sl(n,\mathbb{R})$ be the traceless matrices,%
\begin{equation}
sl(n,\mathbb{R}):=\{A\in gl(n,\mathbb{R})|\operatorname*{tr}(A)=0\}.
\label{e.2.4}%
\end{equation}
Then%
\[
\tau_{g}SL(n,\mathbb{R})=\{A\in gl(n,\mathbb{R})|g^{-1}A\in sl(n,\mathbb{R}%
)\}
\]
and in particular $\tau_{I}SL(n,\mathbb{R})=sl(n,\mathbb{R}).$

\item Let $so\left(  n,\mathbb{R}\right)  $ be the skew symmetric matrices,%
\[
so\left(  n,\mathbb{R}\right)  :=\{A\in gl(n,\mathbb{R})|A=-A^{\mathrm{tr}%
}\}.
\]
Then
\[
\tau_{g}O(n)=\{A\in gl(n,\mathbb{R})|g^{-1}A\in so\left(  n,\mathbb{R}\right)
\}
\]
and in particular $\tau_{I}O\left(  n\right)  =so\left(  n,\mathbb{R}\right)
.$ \textbf{Hint:} $g^{-1}=g^{\mathrm{tr}}$ for all $g\in O(n).$

\item If $M\subset E$ and $N\subset V$ are imbedded submanifolds then
\[
\tau_{(m,n)}(M\times N)=\tau_{m}M\times\tau_{n}N\subset E\times V.
\]

\end{enumerate}
\end{exercise}

It is quite possible that $\tau_{m}M=\tau_{m^{\prime}}M$ for some $m\neq
m^{\prime},$ with $m$ and $m^{\prime}$ in $M$ (think of the sphere). Because
of this, it is helpful to label each of the tangent planes with their base point.

\begin{definition}
\label{d.2.16}The \textbf{tangent space} $(T_{m}M)$ to $M$ at $m$ is given by
\[
T_{m}M:=\{m\}\times\tau_{m}M\subset M\times E.
\]
Let
\[
TM:=\cup_{m\in M}T_{m}M,
\]
and call $TM$ the \textbf{tangent space (or tangent bundle)} of $M.$ A
\textbf{tangent vector} is a point $v_{m}:=(m,v)\in TM$ and we let
$\pi:TM\rightarrow M$ denote the \textbf{canonical projection} defined by
$\pi(v_{m})=m.$ Each tangent space is made into a vector space with the vector
space operations being defined by: $c(v_{m}):=(cv)_{m}$ and $v_{m}%
+w_{m}:=(v+w)_{m}.$
\end{definition}

\begin{exercise}
\ \label{exr.2.17}Prove that $TM$ is an imbedded submanifold of $E\times E.$
\textbf{Hint:} suppose that $z:E\rightarrow\mathbb{R}^{N}$ is a function as in
the Definition \ref{d.2.2}. Define $\mathcal{D}(Z):=\mathcal{D}(z)\times E$
and $Z:\mathcal{D}(Z)\rightarrow\mathbb{R}^{N}\times\mathbb{R}^{N}$ by
$Z(x,a):=(z(x),z^{\prime}(x)a).$ Use $Z$'s of this type to check $TM$
satisfies Definition \ref{d.2.2}.
\end{exercise}

\begin{notation}
\label{n.2.18}In the sequel, given a smooth path $\sigma:(-\varepsilon
,\varepsilon)\rightarrow M,$ we will abuse notation and write $\sigma^{\prime
}\left(  0\right)  $ for either%
\[
\frac{d}{ds}|_{0}\sigma(s)\in\tau_{\sigma\left(  0\right)  }M
\]
or for%
\[
(\sigma(0),\frac{d}{ds}|_{0}\sigma(s))\in T_{\sigma(0)}M=\left\{
\sigma(0)\right\}  \times\tau_{\sigma(0)}M.
\]
Also given a chart $x=(x^{1},x^{2},\ldots,x^{d})$ on $M$ and $m\in
\mathcal{D}(x),$ let $\partial/\partial x^{i}|_{m}$ denote the element
$T_{m}M$ determined by $\partial/\partial x^{i}|_{m}=\sigma^{\prime}(0),$
where $\sigma(s):=x^{-1}(x(m)+se_{i}),$ i.e.
\begin{equation}
\frac{\partial}{\partial x^{i}}|_{m}=(m,\frac{d}{ds}|_{0}x^{-1}(x(m)+se_{i})),
\label{e.2.5}%
\end{equation}
see Figure \ref{fig.4}.
\end{notation}

%

%TCIMACRO{\FRAME{ftbphFU}{3.9359in}{1.9595in}{0pt}{\Qcb{Forming a basis of
%tangent vectors.}}{\Qlb{fig.4}}{f3.eps}{\special{ language "Scientific Word";
%type "GRAPHIC";  maintain-aspect-ratio TRUE;  display "USEDEF";
%valid_file "F";  width 3.9359in;  height 1.9595in;  depth 0pt;
%original-width 5.0601in;  original-height 2.5068in;  cropleft "0";
%filename 'F3.eps';file-properties "XNPEU";}}}%
%BeginExpansion
\begin{figure}
[ptbh]
\begin{center}
\includegraphics[
height=1.9595in,
width=3.9359in
]%
{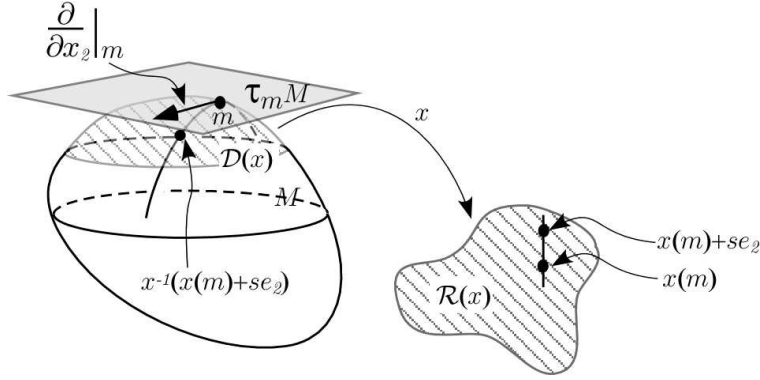}%
\caption{Forming a basis of tangent vectors.}%
\label{fig.4}%
\end{center}
\end{figure}
%EndExpansion

The reason for the strange notation in Eq. (\ref{e.2.5}) will be explained
after Notation \ref{n.2.20}. By definition, every element of $T_{m}M$ is of
the form $\sigma^{\prime}\left(  0\right)  $ where $\sigma$ is a smooth path
into $M$ such that $\sigma\left(  0\right)  =m.$ Moreover by Theorem
\ref{t.2.13}, $\{\partial/\partial x^{i}|_{m}\}_{i=1}^{d}$ is a basis for
$T_{m}M.$

\begin{definition}
\label{d.2.19}Suppose that $f:M\rightarrow V$ is a smooth function,
$m\in\mathcal{D}(f)$ and $v_{m}\in T_{m}M.$ Write
\[
v_{m}f=df(v_{m}):=\frac{d}{ds}|_{0}f(\sigma(s)),
\]
where $\sigma$ is any smooth path in $M$ such that $\sigma^{\prime}(0)=v_{m}.$
The function $df:TM\rightarrow V$ will be called the \textbf{differential of
}$f.$
\end{definition}

\begin{notation}
\label{n.2.20}If $M$ and $N$ are two manifolds $f:M\times N\rightarrow V$ is a
smooth function, we will write $d_{M}f\left(  \cdot,n\right)  $ to indicate
that we are computing the differential of the function $m\in M\rightarrow
f(m,n)\in V$ for fixed $n\in N.$
\end{notation}

To understand the notation in \text{(\ref{e.2.5})}, suppose that $f=F\circ
x=F(x^{1},x^{2},\ldots,x^{d})$ where $F:\mathbb{R}^{d}\rightarrow\mathbb{R}$
is a smooth function and $x$ is a chart on $M.$ Then
\[
\frac{\partial f(m)}{\partial x^{i}}:=\frac{\partial}{\partial x^{i}}%
|_{m}f=(D_{i}F)(x(m)),
\]
where $D_{i}$ denotes the $i^{\text{th}}$ -- partial derivative of $F.$ Also
notice that $dx^{j}\left(  \frac{\partial}{\partial x^{i}}|_{m}\right)
=\delta_{ij}$ so that $\left\{  dx^{i}|_{T_{m}M}\right\}  _{i=1}^{d}$ is the
dual basis of $\{\partial/\partial x^{i}|_{m}\}_{i=1}^{d}$ and therefore if
$v_{m}\in T_{m}M$ then
\begin{equation}
v_{m}=\sum_{i=1}^{d}dx^{i}\left(  v_{m}\right)  \frac{\partial}{\partial
x^{i}}|_{m}. \label{e.2.6}%
\end{equation}
This explicitly exhibits $v_{m}$ as a first order differential operator acting
on \textquotedblleft germs\textquotedblright\ of smooth functions defined near
$m\in M.$

\begin{remark}
[Product Rule]\label{r.2.21}Suppose that $f:M\rightarrow V$ and
$g:M\rightarrow\operatorname*{End}(V)$ are smooth functions, then
\[
v_{m}(gf)=\frac{d}{ds}|_{0}\left[  g(\sigma(s))f(\sigma(s))\right]
=v_{m}g\cdot f(m)+g(m)v_{m}f
\]
or equivalently
\[
d(gf)(v_{m})=dg(v_{m})f(m)+g(m)df(v_{m}).
\]
This last equation will be abbreviated as $d(gf)=dg\cdot f+gdf.$
\end{remark}

\begin{definition}
\label{d.2.22} Let $f:M\rightarrow N$ be a smooth map of imbedded
submanifolds. Define the \textbf{differential, }$f_{\ast},$ of $f$ by
\[
f_{\ast}v_{m}=(f\circ\sigma)^{\prime}(0)\in T_{f(m)}N,
\]
where $v_{m}=\sigma^{\prime}(0)\in T_{m}M,$ and $m\in\mathcal{D}(f).$
\end{definition}

%

%TCIMACRO{\FRAME{ftbphFU}{3.4166in}{1.8087in}{0pt}{\Qcb{The differential of
%$f.$}}{\Qlb{fig.5}}{f4.eps}{\special{ language "Scientific Word";
%type "GRAPHIC";  maintain-aspect-ratio TRUE;  display "USEDEF";
%valid_file "F";  width 3.4166in;  height 1.8087in;  depth 0pt;
%original-width 4.7069in;  original-height 2.4797in;  cropleft "0";
%filename 'F4.eps';file-properties "XNPEU";}}}%
%BeginExpansion
\begin{figure}
[ptbh]
\begin{center}
\includegraphics[
height=1.8087in,
width=3.4166in
]%
{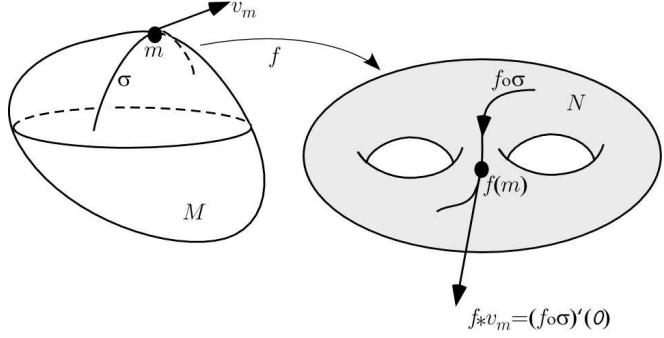}%
\caption{The differential of $f.$}%
\label{fig.5}%
\end{center}
\end{figure}
%EndExpansion

\begin{lemma}
\label{l.2.23} The differentials defined in Definitions \ref{d.2.19} and
\ref{d.2.22} are well defined linear maps on $T_{m}M$ for each $m\in
\mathcal{D}(f).$
\end{lemma}

\begin{proof}
I will only prove that $f_{\ast}$ is well defined, since the case of $df$ is
similar. By Proposition \ref{p.2.10}, there is a smooth function
$F:E\rightarrow V,$ such that $f=F|_{M}.$ Therefore by the chain rule
\begin{equation}
f_{\ast}v_{m}=(f\circ\sigma)^{\prime}(0):=\left[  \frac{d}{ds}|_{0}%
f(\sigma(s))\right]  _{f(\sigma(0))}=\left[  F^{\prime}(m)v\right]  _{f\left(
m\right)  }, \label{e.2.7}%
\end{equation}
where $\sigma$ is a smooth path in $M$ such that $\sigma^{\prime}(0)=v_{m}.$
It follows from \text{(\ref{e.2.7})} that $f_{\ast}v_{m}$ does not depend on
the choice of the path $\sigma.$ It is also clear from \text{(\ref{e.2.7})},
that $f_{\ast}$ is linear on $T_{m}M.$
\end{proof}

\begin{remark}
\label{r.2.24} Suppose that $F:E\rightarrow V$ is a smooth function and that
$f:=F|_{M}.$ Then as in the proof of Lemma \ref{l.2.23},%
\begin{equation}
df(v_{m})=F^{\prime}(m)v \label{e.2.8}%
\end{equation}
for all $v_{m}\in T_{m}M$, and $m\in\mathcal{D}(f).$ Incidentally, since the
left hand sides of \text{(\ref{e.2.7})} and \text{(\ref{e.2.8})} are defined
\textquotedblleft intrinsically,\textquotedblright\ the right members of
\text{(\ref{e.2.7})} and \text{(\ref{e.2.8})} are independent of the possible
choices of functions $F$ which extend $f.$
\end{remark}

\begin{lemma}
[Chain Rules]\label{l.2.25} Suppose that $M,$ $N,$ and $P$ are imbedded
submanifolds and $V$ is a finite dimensional vector space. Let $f:M\rightarrow
N,$ $g:N\rightarrow P,$ and $h:N\rightarrow V$ be smooth functions. Then:
\begin{equation}
(g\circ f)_{\ast}v_{m}=g_{\ast}(f_{\ast}v_{m})\!,\,\quad\forall~v_{m}\in TM
\label{e.2.9}%
\end{equation}
and
\begin{equation}
d(h\circ f)(v_{m})=dh(f_{\ast}v_{m}),\,\quad\forall~v_{m}\in TM.
\label{e.2.10}%
\end{equation}
These equations will be written more concisely as $(g\circ f)_{\ast}=g_{\ast
}f_{\ast}$ and $d(h\circ f)=dhf_{\ast}$ respectively.
\end{lemma}

\begin{proof}
Let $\sigma$ be a smooth path in $M$ such that $v_{m}=\sigma^{\prime}(0).$
Then, see Figure \ref{fig.6},
\begin{align*}
(g\circ f)_{\ast}v_{m}  &  :=(g\circ f\circ\sigma)^{\prime}(0)=g_{\ast}%
(f\circ\sigma)^{\prime}(0)\\
&  =g_{\ast}f_{\ast}\sigma^{\prime}(0)=g_{\ast}f_{\ast}v_{m}.
\end{align*}
Similarly,
\begin{align*}
d(h\circ f)(v_{m})  &  :=\frac{d}{ds}|_{0}(h\circ f\circ\sigma)(s)=dh((f\circ
\sigma)^{\prime}(0))\\
&  =dh(f_{\ast}\sigma^{\prime}(0))=dh(f_{\ast}v_{m}).
\end{align*}%
%TCIMACRO{\FRAME{ftbphFU}{4.1825in}{2.345in}{0pt}{\Qcb{The chain rule.}%
%}{\Qlb{fig.6}}{f5.eps}{\special{ language "Scientific Word";  type "GRAPHIC";
%width 4.1825in;  height 2.345in;  depth 0pt;  original-width 6.0022in;
%original-height 3.3539in;  cropleft "0";  croptop "1";  cropright "1";
%cropbottom "0";  filename 'F5.eps';file-properties "XNPEU";}%
%}}%
%BeginExpansion
\begin{figure}
[ptbh]
\begin{center}
\includegraphics[
height=2.345in,
width=4.1825in
]%
{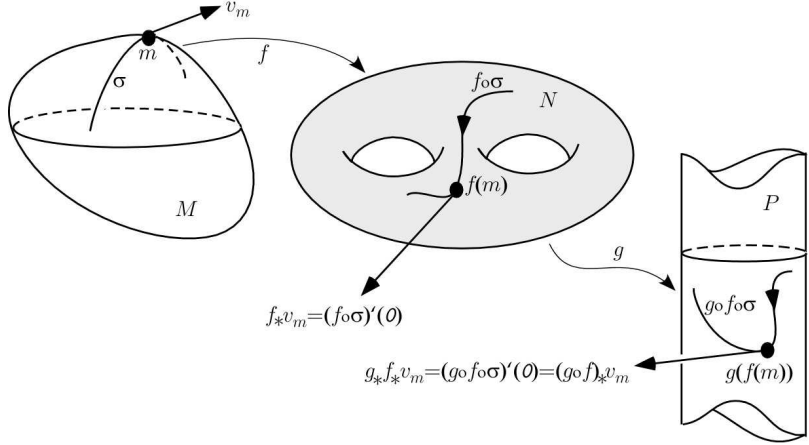}%
\caption{The chain rule.}%
\label{fig.6}%
\end{center}
\end{figure}
%EndExpansion

\end{proof}

If $f:M\rightarrow V$ is a smooth function, $x$ is a chart on $M$, and
$m\in\mathcal{D}(f)\cap\mathcal{D}(x),$ we will write $\partial f(m)/\partial
x^{i}$ for $df\left(  \partial/\partial x^{i}|_{m}\right)  .$ Combining this
notation with Eq. (\ref{e.2.6}) leads to the pleasing formula,%
\begin{equation}
df=\sum_{i=1}^{d}\frac{\partial f}{\partial x^{i}}dx^{i}, \label{e.2.11}%
\end{equation}
by which we mean%
\[
df(v_{m})=\sum_{i=1}^{d}\frac{\partial f(m)}{\partial x^{i}}dx^{i}(v_{m}).
\]

Suppose that $f:M^{d}\rightarrow N^{k}$ is a smooth map of imbedded
submanifolds, $m\in M,$ $x$ is a chart on $M$ such that $m\in\mathcal{D}(x),$
and $y$ is a chart on $N$ such that $f(m)\in\mathcal{D}(y).$ Then the matrix
of
\[
f_{\ast m}:=f_{\ast}|_{T_{m}M}:T_{m}M\rightarrow T_{f(m)}N
\]
relative to the bases $\{\partial/\partial x^{i}|_{m}\}_{i=1}^{d}$ of $T_{m}M$
and $\{\partial/\partial y^{j}|_{f(m)}\}_{j=1}^{k}$ of $T_{f(m)}N$ is
$(\partial(y^{j}\circ f)(m)/\partial x^{i}).$ Indeed, if $v_{m}=\sum_{i=1}%
^{d}v^{i}\partial/\partial x^{i}|_{m},$ then%

\begin{align*}
f_{\ast}v_{m}  &  =\sum_{j=1}^{k}dy^{j}(f_{\ast}v_{m})\partial/\partial
y^{j}|_{f(m)}\\
&  =\sum_{j=1}^{k}d(y^{j}\circ f)(v_{m})\partial/\partial y^{j}|_{f(m)}%
\qquad\qquad\qquad\qquad\text{(by\ Eq. (\ref{e.2.10}))}\\
&  =\sum_{j=1}^{k}\sum_{i=1}^{d}\frac{\partial(y^{j}\circ f)(m)}{\partial
x^{i}}\cdot dx^{i}(v_{m})\partial/\partial y^{j}|_{f(m)}\qquad\text{(by\ Eq.
(\ref{e.2.11}))}\\
&  =\sum_{j=1}^{k}\sum_{i=1}^{d}\frac{\partial(y^{j}\circ f)(m)}{\partial
x^{i}}v^{i}\partial/\partial y^{j}|_{f(m)}.
\end{align*}

\begin{example}
\label{ex.2.26} Let $M=O(n),$ $k\in O(n),$ and $f:O(n)\rightarrow O(n)$ be
defined by $f(g):=kg.$ Then $f$ is a smooth function on $O(n)$ because it is
the restriction of a smooth function on $gl(n,\mathbb{R}).$ Given $A_{g}\in
T_{g}O(n),$ by Eq. \text{(\ref{e.2.7})},
\[
f_{\ast}A_{g}=(kg,kA)=(kA)_{kg}%
\]
(In the future we denote $f$ by $L_{k};$ $L_{k}$ is \textbf{left translation}
by $k\in O(n).)$
\end{example}

\begin{definition}
\label{d.2.27}A \textbf{Lie group} is a manifold, $G,$ which is also a group
such that the group operations are smooth functions. The tangent space,
$\mathfrak{g}:=\mathrm{Lie}\left(  G\right)  :=T_{e}G,$ to $G$ at the identity
$e\in G$ is called the \textbf{Lie algebra} of $G.$
\end{definition}

\begin{exercise}
\label{exr.2.28}Verify that $GL(n,\mathbb{R}),$ $SL(n,\mathbb{R}),$ $O(n),$
$SO(n)$ and $T^{n}$ (see Example \ref{ex.2.6}) are all Lie groups and
\begin{align*}
\mathrm{Lie}\left(  GL(n,\mathbb{R})\right)   &  \cong gl(n,\mathbb{R}),\\
\mathrm{Lie}\left(  SL(n,\mathbb{R}))\right)   &  \cong sl(n,\mathbb{R})\\
\mathrm{Lie}\left(  O(n))\right)   &  =\mathrm{Lie}\left(  SO(n))\right)
\cong so(n,\mathbb{R})\text{ and}\\
\mathrm{Lie}\left(  T^{n})\right)   &  \cong\left(  i\mathbb{R}\right)
^{n}\subset\mathbb{C}^{n}.
\end{align*}
See Exercise \ref{exr.2.15} for the notation being used here.
\end{exercise}

\begin{exercise}
[Continuation of Exercise \ref{exr.2.17}]\label{exr.2.29} Show for each chart
$x$ on $M$ that the function
\[
\phi(v_{m}):=(x(m),dx(v_{m}))=x_{\ast}v_{m}%
\]
is a chart on $TM.$ Note that $\mathcal{D}(\phi):=\cup_{m\in\mathcal{D}%
(x)}T_{m}M.$
\end{exercise}

The following lemma gives an important example of a smooth function on $M$
which will be needed when we consider $M$ as a \textquotedblleft Riemannian
manifold.\textquotedblright

\begin{lemma}
\label{l.2.30}Suppose that $(E,\langle\cdot,\cdot\rangle)$ is an inner product
space and the $M\subset E$ is an imbedded submanifold. For each $m\in M,$ let
$P(m)$ denote the orthogonal projection of $E$ onto $\tau_{m}M$ and
$Q(m):=I-P(m)$ denote the orthogonal projection onto $\tau_{m}M^{\perp}.$ Then
$P$ and $Q$ are smooth functions from $M$ to $gl(E),$ where $gl(E)$ denotes
the vector space of linear maps from $E$ to $E.$
\end{lemma}

\begin{proof}
Let $z:E\rightarrow\mathbb{R}^{N}$ be as in Definition \ref{d.2.2}. To
simplify notation, let $F(p):=z_{>}(p)$ for all $p\in\mathcal{D}(z),$ so that
$\tau_{m}M=\operatorname*{Nul}\left(  F^{\prime}(m)\right)  $ for
$m\in\mathcal{D}(x)=\mathcal{D}(z)\cap M.$ Since $F^{\prime}(m):E\rightarrow
\mathbb{R}^{N-d}$ is surjective, an elementary exercise in linear algebra
shows
\[
(F^{\prime}(m)F^{\prime}(m)^{\ast}):\mathbb{R}^{N-d}\rightarrow\mathbb{R}%
^{N-d}%
\]
is invertible for all $m\in\mathcal{D}(x).$ The orthogonal projection
$Q\left(  m\right)  $ may be expressed as;
\begin{equation}
Q(m)=F^{\prime}(m)^{\ast}(F^{\prime}(m)F^{\prime}(m)^{\ast})^{-1}F^{\prime
}(m). \label{e.2.12}%
\end{equation}
Since being invertible is an open condition, $(F^{\prime}(\cdot)F^{\prime
}(\cdot)^{\ast})$ is invertible in an open neighborhood $\mathcal{\ N}\subset
E$ of $\mathcal{D}(x).$ Hence $Q$ has a smooth extension $\tilde{Q}$ to
$\mathcal{\ N}$ given by
\[
\tilde{Q}(x):=F^{\prime}(x)^{\ast}(F^{\prime}(x)F^{\prime}(x)^{\ast}%
)^{-1}F^{\prime}(x).
\]
Since $Q|_{\mathcal{D}(x)}=\tilde{Q}|_{\mathcal{D}(x)}$ and $\tilde{Q}$ is
smooth on $\mathcal{\ N},$ $Q|_{\mathcal{D}(x)}$ is also smooth. Since $z$ as
in Definition \ref{d.2.2} was arbitrary and smoothness is a local property, it
follows that $Q$ is smooth on $M.$ Clearly, $P:=I-Q$ is also a smooth function
on $M.$
\end{proof}

\begin{definition}
\label{d.2.31}A \textbf{local vector field} $Y$ on $M$ is a smooth function
$Y:M\rightarrow TM$ such that $Y(m)\in T_{m}M$ for all $m\in\mathcal{D}(Y),$
where $\mathcal{D}(Y)$ is assumed to be an open subset of $M.$ Let
$\Gamma(TM)$ denote the collection of globally defined (i.e. $\mathcal{D}%
(Y)=M)$ smooth vector-fields $Y$ on $M.$
\end{definition}

Note that $\partial/\partial x^{i}$ are local vector-fields on $M$ for each
chart $x\in\mathcal{A}(M)$ and $i=1,2,\ldots,d.$ The next exercise asserts
that these vector fields are smooth.

\begin{exercise}
\ \label{exr.2.32}Let $Y$ be a vector field on $M,$ $x\in\mathcal{A}(M)$ be a
chart on $M$ and $Y^{i}:=dx^{i}(Y).$ Then
\[
Y(m):=\sum_{i=1}^{d}Y^{i}\left(  m\right)  \partial/\partial x^{i}%
|_{m}~\forall~m\in\mathcal{D}\left(  x\right)  ,
\]
which we abbreviate as $Y=\sum_{i=1}^{d}Y^{i}\partial/\partial x^{i}.$ Show
the condition that $Y$ is smooth translates into the statement that each of
the functions $Y^{i}$ is smooth.
\end{exercise}

\begin{exercise}
\label{exr.2.33}Let $Y:M\rightarrow TM,$ be a vector field. Then
\[
Y(m)=(m,y(m))=y(m)_{m}%
\]
for some function $y:M\rightarrow E$ such that $y(m)\in\tau_{m}M$ for all
$m\in\mathcal{D}(Y)=\mathcal{D}(y).$ Show that $Y$ is smooth iff
$y:M\rightarrow E$ is smooth.
\end{exercise}

\begin{example}
\label{ex.2.34}Let $M=SL(n,\mathbb{R})$ and $A\in sl(n,\mathbb{R})=\tau
_{I}SL(n,\mathbb{R}),$ i.e. $A$ is a $n\times n$ real matrix such that
$\operatorname*{tr}\left(  A\right)  =0.$ Then $\tilde{A}(g):=L_{g\ast}%
A_{e}=(g,gA)$ for $g\in M$ is a smooth vector field on $M.$
\end{example}

\begin{example}
\label{ex.2.35} Keep the notation of Lemma \ref{l.2.30}. Let $y:M\rightarrow
E$ be any smooth function. Then $Y(m):=(m,P(m)y(m))$ for all $m\in M$ is a
smooth vector-field on $M.$
\end{example}

\begin{definition}
\label{d.2.36}Given $Y\in\Gamma(TM)$ and $f\in C^{\infty}(M),$ let $Yf\in
C^{\infty}(M)$ be defined by $(Yf)(m):=df(Y(m)),$ for all $m\in\mathcal{D}%
(f)\cap\mathcal{D}(Y).$ In this way the vector-field $Y$ may be viewed as a
first order differential operator on $C^{\infty}(M).$
\end{definition}

\begin{notation}
\label{n.2.37}The \textbf{Lie bracket }of two smooth vector fields, $Y$ and
$W,$ on $M$ is the vector field $[Y,W]$ which acts on $C^{\infty}(M)$ by the
formula
\begin{equation}
\lbrack Y,W]f:=Y(Wf)-W(Yf),\,\quad\forall~f\in C^{\infty}(M). \label{e.2.13}%
\end{equation}
(In general one might suspect that $[Y,W]$ is a second order differential
operator, however this is not the case, see Exercise \ref{exr.2.38}.)
Sometimes it will be convenient to write $L_{Y}W$ for $[Y,W].$
\end{notation}

\begin{exercise}
\label{exr.2.38}Show that $[Y,W]$ is again a first order differential operator
on $C^{\infty}(M)$ coming from a vector-field. In particular, if $x$ is a
chart on $M,$ $Y=\sum_{i=1}^{d}Y^{i}\partial/\partial x^{i}$ and $W=\sum
_{i=1}^{d}W^{i}\partial/\partial x^{i},$ then on $\mathcal{D}(x),$%
\begin{equation}
\lbrack Y,W]=\sum_{i=1}^{d}(YW^{i}-WY^{i})\partial/\partial x^{i}.
\label{e.2.14}%
\end{equation}

\end{exercise}

\begin{proposition}
\label{p.2.39}If $Y(m)=(m,y(m))$ and $W(m)=(m,w(m))$ and $y,w:M\rightarrow E$
are smooth functions such that $y(m),w(m)\in\tau_{m}M,$ then we may express
the Lie bracket, $[Y,W](m),$ as
\begin{equation}
\lbrack Y,W](m)=(m,(Yw-Wy)(m))=(m,dw(Y(m))-dy(W(m))). \label{e.2.15}%
\end{equation}

\end{proposition}

\begin{proof}
Let $f$ be a smooth function $M$ which we may take, by Proposition
\ref{p.2.10}, to be the restriction of a smooth function on $E.$ Similarly we
we may assume that $y$ and $w$ are smooth functions on $E$ such that
$y(m),w(m)\in\tau_{m}M$ for all $m\in M.$ Then
\begin{align}
(YW-WY)f  &  =Y\left[  f^{\prime}w\right]  -W\left[  f^{\prime}y\right]
\nonumber\\
&  =f^{\prime\prime}(y,w)-f^{\prime\prime}(w,y)+f^{\prime}\left(  Yw\right)
-f^{\prime}\left(  Wy\right) \nonumber\\
&  =f^{\prime}\left(  Yw-Wy\right)  \label{e.2.16}%
\end{align}
wherein the last equality we have use the fact that mixed partial derivatives
commute to conclude%
\[
f^{\prime\prime}(u,v)-f^{\prime\prime}(v,u):=\left(  \partial_{u}\partial
_{v}-\partial_{v}\partial_{u}\right)  f=0\text{ }\forall~u,v\in E.
\]
Taking $f=z_{>}$ in Eq. (\ref{e.2.16}) with $z=(z_{<},z_{>})$ being a chart on
$E$ as in Definition \ref{d.2.2}, shows
\[
0=(YW-WY)z_{>}\left(  m\right)  =z_{>}^{\prime}\left(
dw(Y(m))-dy(W(m))\right)
\]
and thus $(m,dw(Y(m))-dy(W(m)))\in T_{m}M.$ With this observation, we then
have
\[
f^{\prime}\left(  Yw-Wy\right)  =df\left(  (m,dw(Y(m))-dy(W(m)))\right)
\]
which combined with Eq. (\ref{e.2.16}) verifies Eq. (\ref{e.2.15}).
\end{proof}

\begin{exercise}
\label{exr.2.40}Let $M=SL(n,\mathbb{R})$ and $A,B\in sl(n,\mathbb{R})$ and
$\tilde{A}$ and $\tilde{B}$ be the associated left invariant vector fields on
$M$ as introduced in Example \ref{ex.2.34}. Show $\left[  \tilde{A},\tilde
{B}\right]  =\widetilde{\left[  A,B\right]  }$ where $\left[  A,B\right]
:=AB-BA$ is the matrix commutator of $A$ and $B.$
\end{exercise}

\subsection{More References}

The reader wishing to learn about manifolds is referred to
\cite{Ab2,AMa,BC,Davies90,DoC,Hi,Kl1,Kl2,Kl3,KN1,KN2,Spivak1}. The texts by
Kobayashi and Nomizu are very thorough while the books by Klingenberg give an
idea of why differential geometers are interested in loop spaces. There is a
vast literature on Lie groups and there representations. Here are just two
books which I have found very useful, \cite{BD,Wal1}.

\section{Riemannian Geometry Primer}

\label{s.3}

This section introduces the following objects: 1) Riemannian metrics, 2)
Riemannian volume forms, 3) gradients, 4) divergences, 5) Laplacians, 6)
covariant derivatives, 7) parallel translations, and 8) curvatures.

\subsection{Riemannian Metrics\label{s.3.1}}

\begin{definition}
\label{d.3.1}A \textbf{Riemannian metric}, $\langle\cdot,\cdot\rangle$ (also
denoted by $g),$ on $M$ is a smoothly varying choice of inner product,
$g_{m}=\langle\cdot,\cdot\rangle_{m},$ on each of the tangent spaces $T_{m}M,$
$m\in M.$ The smoothness condition is the requirement that the function $m\in
M\rightarrow\langle X(m),Y(m)\rangle_{m}\in\mathbb{R}$ is smooth for all
smooth vector fields $X$ and $Y$ on $M.$
\end{definition}

It is customary to write $ds^{2}$ for the function on $TM$ defined by
\begin{equation}
ds^{2}(v_{m}):=\langle v_{m},v_{m}\rangle_{m}=g_{m}\left(  v_{m},v_{m}\right)
. \label{e.3.1}%
\end{equation}
By polarization, the Riemannian metric $\langle\cdot,\cdot\rangle$ is uniquely
determined by the function $ds^{2}.$ Given a chart $x$ on $M$ and $v\in
T_{m}M,$ by Eqs. (\ref{e.3.1}) and (\ref{e.2.6}) we have
\begin{equation}
ds^{2}(v_{m})=\sum_{i,j=1}^{d}\langle\partial/\partial x^{i}|_{m}%
,\partial/\partial x^{j}|_{m}\rangle_{m}dx^{i}(v_{m})dx^{j}(v_{m}).
\label{e.3.2}%
\end{equation}
We will abbreviate this equation in the future by writing
\begin{equation}
ds^{2}=\sum_{i,j=1}^{d}g_{ij}^{x}dx^{i}dx^{j} \label{e.3.3}%
\end{equation}
where
\[
g_{i,j}^{x}(m):=\langle\partial/\partial x^{i}|_{m},\partial/\partial
x^{j}|_{m}\rangle_{m}=g\left(  \partial/\partial x^{i}|_{m},\partial/\partial
x^{j}|_{m}\right)  .
\]
Typically $g_{i,j}^{x}$ will be abbreviated by $g_{ij}$ if no confusion is
likely to arise.

\begin{example}
\label{ex.3.2} Let $M=\mathbb{R}^{N}$ and let $x=(x^{1},x^{2},\ldots,x^{N})$
denote the standard chart on $M,$ i.e. $x(m)=m$ for all $m\in M.$ The standard
Riemannian metric on $\mathbb{R}^{N}$ is determined by
\[
ds^{2}=\sum_{i=1}^{N}(dx^{i})^{2}=\sum_{i=1}^{N}dx^{i}\cdot dx^{i},
\]
and so $g^{x}$ is the identity matrix here. The general Riemannian metric on
$\mathbb{R}^{N}$ is determined by $ds^{2}=\sum_{i,j=1}^{N}g_{ij}dx^{i}dx^{j},$
where $g=(g_{ij})$ is a smooth $gl(N,\mathbb{R})$ -- valued function on
$\mathbb{R}^{N}$ such that $g(m)$ is positive definite matrix for all
$m\in\mathbb{R}^{N}.$
\end{example}

Let $M$ be an imbedded submanifold of a finite dimensional inner product space
$(E,\langle\cdot,\cdot\rangle).$ The manifold $M$ \textbf{inherits} a metric
from $E$ determined by
\[
ds^{2}(v_{m})=\langle v,v\rangle\ \forall~v_{m}\in TM.
\]
It is a well known deep fact that \textbf{all} finite dimensional Riemannian
manifolds may be constructed in this way, see Nash \cite{Nash} and Moser
\cite{Moser1,Moser3,Moser4}. To simplify the exposition, in the sequel we will
usually assume that $(E,\langle\cdot,\cdot\rangle)$ is an inner product space,
$M^{d}\subset E$ is an imbedded submanifold, and the Riemannian metric on $M$
is determined in this way, i.e.%
\[
\langle v_{m},w_{m}\rangle=\langle v,w\rangle_{\mathbb{R}^{N}},\,\quad
\forall~v_{m},w_{m}\in T_{m}M\text{\textrm{\ and }}m\in M.
\]
In this setting the components $g_{i,j}^{x}$ of the metric $ds^{2}$ relative
to a chart $x$ may be computed as $g_{i,j}^{x}(m)=(\phi_{;i}(x(m)),\phi
_{;j}(x(m))),$ where $\{e_{i}\}_{i=1}^{d}$ is the standard basis for
$\mathbb{R}^{d},$%
\[
\phi:=x^{-1}\text{ and }\phi_{;i}(a):=\frac{d}{dt}|_{0}\phi(a+te_{i}).
\]

\begin{example}
\ \label{ex.3.3}Let $M=G:=SL(n,\mathbb{R})$ and $A_{g}\in T_{g}M.$

\begin{enumerate}
\item Then%
\begin{equation}
ds^{2}(A_{g}):=\text{$\operatorname*{tr}$}(A^{\ast}A) \label{e.3.4}%
\end{equation}
defines a Riemannian metric on $G.$ This metric is the inherited metric from
the inner product space $E=gl(n,\mathbb{R})$ with inner product $\langle
A,B\rangle:=\operatorname*{tr}(A^{\ast}B).$

\item A more \textquotedblleft natural\textquotedblright\ choice of a metric
on $G$ is%
\begin{equation}
ds^{2}(A_{g}):=\text{$\operatorname*{tr}$}((g^{-1}A)^{\ast}g^{-1}A).
\label{e.3.5}%
\end{equation}
This metric is invariant under left translations, i.e. $ds^{2}(L_{k\ast}%
A_{g})=ds^{2}(A_{g}),$ for all $k\in G$ and $A_{g}\in TG.$ According to the
imbedding theorem of Nash and Moser, it would be possible to find another
imbedding of $G$ into a Euclidean space, $E,$ so that the metric in Eq.
(\ref{e.3.5}) is inherited from an inner product on $E.$
\end{enumerate}
\end{example}

\begin{example}
\label{ex.3.4}Let $M=\mathbb{R}^{3}$ be equipped with the standard Riemannian
metric and $(r,\varphi,\theta)$ be spherical coordinates on $M$, see Figure
\ref{fig.7}.%
%TCIMACRO{\FRAME{ftbphFU}{1.2623in}{1.7409in}{0pt}{\Qcb{Defining the spherical
%coordinates, $\left(  r,\theta,\phi\right)  $ on $\mathbb{R}^{3}.$}%
%}{\Qlb{fig.7}}{f6.eps}{\special{ language "Scientific Word";  type "GRAPHIC";
%width 1.2623in;  height 1.7409in;  depth 0pt;  original-width 2.5406in;
%original-height 1.7096in;  cropleft "0";  croptop "1";  cropright "1";
%cropbottom "0";  filename 'F6.eps';file-properties "XNPEU";}%
%}}%
%BeginExpansion
\begin{figure}
[ptbh]
\begin{center}
\includegraphics[
height=1.7409in,
width=1.2623in
]%
{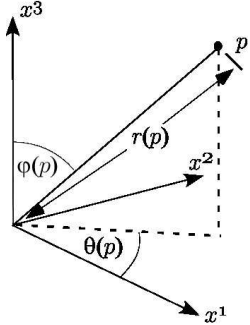}%
\caption{Defining the spherical coordinates, $\left(  r,\theta,\phi\right)  $
on $\mathbb{R}^{3}.$}%
\label{fig.7}%
\end{center}
\end{figure}
%EndExpansion
Here $r,$ $\varphi,$ and $\theta$ are taken to be functions on $\mathbb{R}%
^{3}\setminus\{p\in\mathbb{R}^{3}:p_{2}=0$\textrm{\ and }$p_{1}>0\}$ defined
by $r(p)=|p|,$ $\varphi(p)=\cos^{-1}(p_{3}/|p|)\in(0,\pi),$ and $\theta
(p)\in(0,2\pi)$ is given by $\theta(p)=\tan^{-1}(p_{2}/p_{1})$ if $p_{1}>0$
and $p_{_{2}}>0$ with similar formulas for $(p_{1},p_{2})$ in the other three
quadrants of $\mathbb{R}^{2}.$ Since $x^{1}=r\sin\varphi\cos\theta,$
$x^{2}=r\sin\varphi\sin\theta,$ and $x^{3}=r\cos\varphi,$ it follows using Eq.
(\ref{e.2.11}) that,%
\begin{align*}
dx^{1} &  =\frac{\partial x^{1}}{\partial r}dr+\frac{\partial x^{1}}%
{\partial\varphi}d\varphi+\frac{\partial x^{1}}{\partial\theta}d\theta\\
&  =\sin\varphi\cos\theta dr+r\cos\varphi\cos\theta d\varphi-r\sin\varphi
\sin\theta d\theta,
\end{align*}%
\[
dx^{2}=\sin\varphi\sin\theta dr+r\cos\varphi\sin\theta d\varphi+r\sin
\varphi\cos\theta d\theta,
\]
and
\[
dx^{3}=\cos\varphi dr-r\sin\varphi d\varphi.
\]
An elementary calculation now shows that
\begin{equation}
ds^{2}=\sum_{i=1}^{3}(dx^{i})^{2}=dr^{2}+r^{2}d\varphi^{2}+r^{2}\sin
^{2}\varphi d\theta^{2}.\label{e.3.6}%
\end{equation}
From this last equation, we see that
\begin{equation}
g^{(r,\varphi,\theta)}=\left[
\begin{array}
[c]{ccc}%
1 & 0 & 0\\
0 & r^{2} & 0\\
0 & 0 & r^{2}\sin^{2}\varphi
\end{array}
\right]  .\label{e.3.7}%
\end{equation}

\end{example}

\begin{exercise}
\label{exr.3.5}Let $M:=\{m\in\mathbb{R}^{3}:|m|^{2}=\rho^{2}\},$ so that $M$
is a sphere of radius $\rho$ in $\mathbb{R}^{3}.$ Since $r=\rho$ on $M$ and
$dr\left(  v\right)  =0$ for all $v\in T_{m}M,$ it follows from Eq.
(\ref{e.3.6}) that the induced metric $ds^{2}$ on $M$ is given by%
\begin{equation}
ds^{2}=\rho^{2}d\varphi^{2}+\rho^{2}\sin^{2}\varphi d\theta^{2},\label{e.3.8}%
\end{equation}
and hence%
\begin{equation}
g^{(\varphi,\theta)}=\left[
\begin{array}
[c]{cc}%
\rho^{2} & 0\\
0 & \rho^{2}\sin^{2}\varphi
\end{array}
\right]  .\label{e.3.9}%
\end{equation}

\end{exercise}

\subsection{Integration and the Volume Measure\label{s.3.2}}

\begin{definition}
\label{d.3.6} Let $f\in C_{c}^{\infty}(M)$ (the smooth functions on $M^{d}$
with compact support) and assume the support of $f$ is contained in
$\mathcal{D}(x),$ where $x$ is some chart on $M.$ Set
\[
\int_{M}fdx=\int_{\mathcal{R}(x)}f\circ x^{-1}(a)da,
\]
where $da$ denotes Lebesgue measure on $\mathbb{R}^{d}.$
\end{definition}

\begin{figure}
[ptbh]
\begin{center}
\includegraphics[
height=2.184in,
width=3.8665in
]%
{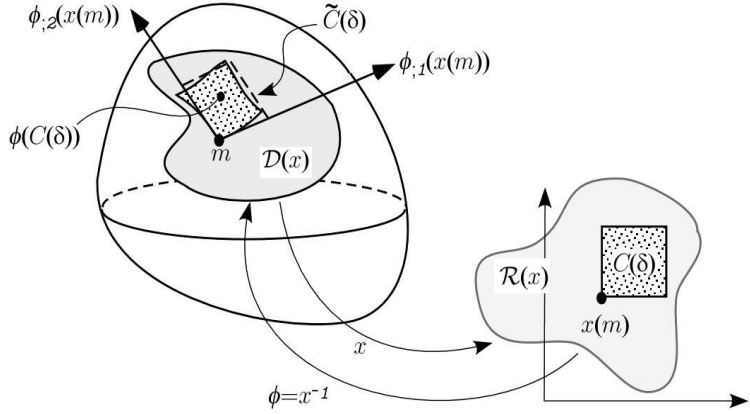}%
\caption{Defining the Riemannian \textquotedblleft volume
element.\textquotedblright}%
\label{fig.8}%
\end{center}
\end{figure}

The problem with this notion of integration is that (as the notation
indicates) $\int_{M}fdx$ depends on the choice of chart $x.$ To remedy this,
consider a small cube $C(\delta)$ of side $\delta$ contained in $\mathcal{R}%
(x)$, see Figure \ref{fig.8}. We wish to estimate \textquotedblleft the
volume\textquotedblright\ of $\phi(C(\delta))$ where $\phi:=x^{-1}%
:\mathcal{R}(x)\rightarrow\mathcal{D}(x).$ Heuristically, we expect the volume
of $\phi(C(\delta))$ to be approximately equal to the volume of the
parallelepiped, $\tilde{C}(\delta),$ in the tangent space $T_{m}M$ determined
by
\begin{equation}
\tilde{C}(\delta):=\left\{  \sum_{i=1}^{d}s_{i}\delta\cdot\phi_{;i}(x\left(
m\right)  )|0\leq s_{i}\leq1,\text{\textrm{\ for }}i=1,2,\ldots,d\right\}  ,
\label{e.3.10}%
\end{equation}
where we are using the notation proceeding Example \ref{ex.3.3}, see Figure
\ref{fig.8}.
%TCIMACRO{\FRAME{ftbphFU}{3.8665in}{2.184in}{0pt}{\Qcb{Defining the Riemannian
%\textquotedblleft volume element.\textquotedblright}}{\Qlb{fig.8}}%
%{f7.eps}{\special{ language "Scientific Word";  type "GRAPHIC";
%width 3.8665in;  height 2.184in;  depth 0pt;  original-width 5.4286in;
%original-height 3.054in;  cropleft "0";  croptop "1";  cropright "1";
%cropbottom "0";  filename 'F7.eps';file-properties "XNPEU";}%
%}}%
%BeginExpansion
%EndExpansion
Since $T_{m}M$ is an inner product space, the volume of $\tilde{C}(\delta)$ is
well defined. For example choose an isometry $\theta:T_{m}M\rightarrow
\mathbb{R}^{d}$ and define the volume of $\tilde{C}(\delta)$ to be $m\left(
\theta(\tilde{C}(\delta))\right)  $ where $m$ is Lebesgue measure on
$\mathbb{R}^{d}.$ The next elementary lemma will be used to give a formula for
the volume of $\tilde{C}\left(  \delta\right)  .$

\begin{lemma}
\label{l.3.7}If $V$ is a finite dimensional inner product space, $\left\{
v_{i}\right\}  _{i=1}^{\dim V}$ is any basis for $V$ and $A:V\rightarrow V$ is
a linear transformation, then
\begin{equation}
\det\left(  A\right)  =\frac{\det\left[  \langle Av_{i},v_{j}\rangle\right]
}{\det\left[  \langle v_{i},v_{j}\rangle\right]  }, \label{e.3.11}%
\end{equation}
where $\det\left[  \langle Av_{i},v_{j}\rangle\right]  $ is the determinant of
the matrix with $i$-$j^{\text{th}}$ -- entry being $\langle Av_{i}%
,v_{j}\rangle.$ Moreover if
\[
\tilde{C}(\delta):=\left\{  \sum_{i=1}^{d}\delta s_{i}\cdot v_{i}:0\leq
s_{i}\leq1,\text{\textrm{\ for }}i=1,2,\ldots,d\right\}
\]
then the volume of $\tilde{C}\left(  \delta\right)  $ is $\delta^{d}\sqrt
{\det\left[  \langle v_{i},v_{j}\rangle\right]  }.$
\end{lemma}

\begin{proof}
Let $\left\{  e_{i}\right\}  _{i=1}^{\dim V}$ be an orthonormal basis for $V,$
then%
\[
\langle Av_{i},v_{j}\rangle=\sum_{l,k}\langle v_{i},e_{l}\rangle\langle
Ae_{l},e_{k}\rangle\langle e_{k},v_{j}\rangle
\]
and therefore by the multiplicative property of the determinant,%
\begin{align}
\det\left[  \langle Av_{i},v_{j}\rangle\right]   &  =\det\left[  \langle
v_{i},e_{l}\rangle\right]  \det\left[  \langle Ae_{l},e_{k}\rangle\right]
\det\left[  \langle e_{k},v_{j}\rangle\right] \nonumber\\
&  =\det\left(  A\right)  \det\left[  \langle v_{i},e_{l}\rangle\right]
\cdot\det\left[  \langle e_{k},v_{j}\rangle\right]  . \label{e.3.12}%
\end{align}
Taking $A=I$ in this equation then shows
\begin{equation}
\det\left[  \langle v_{i},v_{j}\rangle\right]  =\det\left[  \langle
v_{i},e_{l}\rangle\right]  \cdot\det\left[  \langle e_{k},v_{j}\rangle\right]
. \label{e.3.13}%
\end{equation}
Dividing Eq. (\ref{e.3.13}) into Eq. (\ref{e.3.12}) proves Eq. (\ref{e.3.11}).

For the second assertion, it suffices to assume $V=\mathbb{R}^{d}$ with the
usual inner-product. Define $T:\mathbb{R}^{d}\rightarrow\mathbb{R}^{d}$ so
that $Te_{i}=v_{i}$ where $\left\{  e_{i}\right\}  _{i=1}^{d}$ is the standard
basis for $\mathbb{R}^{d},$ then $\tilde{C}\left(  \delta\right)  =T\left(
\left[  0,\delta\right]  ^{d}\right)  $ and hence%
\begin{align*}
m\left(  \tilde{C}\left(  \delta\right)  \right)   &  =\left\vert \det
T\right\vert m\left(  \left[  0,\delta\right]  ^{d}\right)  =\delta
^{d}\left\vert \det T\right\vert =\delta^{d}\sqrt{\det T^{\mathrm{tr}}T}\\
&  =\delta^{d}\sqrt{\det\left[  \langle T^{\mathrm{tr}}Te_{i},e_{j}%
\rangle\right]  }=\delta^{d}\sqrt{\det\left[  \left(  Te_{i},Te_{j}\right)
\right]  }=\delta^{d}\sqrt{\det\left[  \langle v_{i},v_{j}\rangle\right]  }.
\end{align*}

\end{proof}

Using the second assertion in Lemma \ref{l.3.7}, the volume of $\tilde
{C}(\delta)$ in Eq. (\ref{e.3.10}) is $\delta^{d}\sqrt{\det g^{x}(m)},$ where
$g_{ij}^{x}(m)=\langle\phi_{;i}(x(m)),\phi_{;j}(x(m))\rangle_{m}.$ Because of
the above computations, it is reasonable to try to define a new integral on
$\mathcal{D}\left(  x\right)  \subset M$ by
\[
\int_{\mathcal{D}(x)}f\,d\lambda_{\mathcal{D}\left(  x\right)  }%
:=\int_{\mathcal{D}(x)}f\sqrt{g^{x}}dx,
\]
i.e. let $\lambda_{\mathcal{D}\left(  x\right)  }$ be the measure satisfying%
\begin{equation}
d\lambda_{\mathcal{D}\left(  x\right)  }=\sqrt{g^{x}}dx. \label{e.3.14}%
\end{equation}

\begin{lemma}
\label{l.3.8}Suppose that $y$ and $x$ are two charts on $M,$ then
\begin{equation}
g_{l,k}^{y}=\sum_{i,j=1}^{d}g_{i,j}^{x}\frac{\partial x^{i}}{\partial y^{k}%
}\frac{\partial x^{j}}{\partial y^{l}}. \label{e.3.15}%
\end{equation}

\end{lemma}

\begin{proof}
Inserting the identities
\[
dx^{i}=\sum_{k=1}^{d}\frac{\partial x^{i}}{\partial y^{k}}dy^{k}\text{ and
}dx^{j}=\sum_{l=1}^{d}\frac{\partial x^{j}}{\partial y^{l}}dy^{l}%
\]
and into the formula $ds^{2}=\sum_{i,j=1}^{d}g_{i,j}^{x}dx^{i}dx^{j}$ gives
\[
ds^{2}=\sum_{i,j,k,l=1}^{d}g_{i,j}^{x}\frac{\partial x^{i}}{\partial y^{k}%
}\frac{\partial x^{j}}{\partial y^{l}}dy^{l}dy^{k}%
\]
from which \text{(\ref{e.3.15})} follows.
\end{proof}

\begin{exercise}
\label{exr.3.9}Suppose that $x$ and $y$ are two charts on $M$ and $f\in
C_{c}^{\infty}(M)$ such that the support of $f$ is contained in $\mathcal{D}%
(x)\cap\mathcal{D}(y).$ Using Lemma \ref{l.3.8} and the change of variable
formula show,
\[
\int_{\mathcal{D}(x)\cap\mathcal{D}(y)}f\sqrt{g^{x}}dx=\int_{\mathcal{D}%
(x)\cap\mathcal{D}(y)}f\sqrt{g^{y}}dy.
\]

\end{exercise}

\begin{theorem}
[Riemann Volume Measure]\label{t.3.10}There exists a unique measure,
$\lambda_{M}$ on the Borel $\sigma$ -- algebra of $M$ such that for any chart
$x$ on $M,$%
\begin{equation}
d\lambda_{M}\left(  x\right)  =d\lambda_{\mathcal{D}\left(  x\right)  }%
=\sqrt{g^{x}}dx\text{ on }\mathcal{D}\left(  x\right)  . \label{e.3.16}%
\end{equation}

\end{theorem}

\begin{proof}
Choose a countable collection of charts, $\{x_{i}\}_{i=1}^{\infty}$ such that
$M=\cup_{i=1}^{\infty}\mathcal{D}\left(  x_{i}\right)  $ and let
$U_{1}:=\mathcal{D}(x_{1})$ and $U_{i}:=\mathcal{D}(x_{i})\setminus(\cup
_{j=1}^{i-1}\mathcal{D}(x_{j}))$ for $i\geq1.$ Then if $B\subset X$ is a Borel
set, define the measure $\lambda_{M}\left(  B\right)  $ by%
\begin{equation}
\lambda_{M}\left(  B\right)  :=\sum_{i=1}^{\infty}\lambda_{\mathcal{D}\left(
x_{i}\right)  }\left(  B\cap U_{i}\right)  . \label{e.3.17}%
\end{equation}
If $x$ is any chart on $M$ and $B\subset\mathcal{D}\left(  x\right)  ,$ then
$B\cap U_{i}\subset\mathcal{D}\left(  x_{i}\right)  \cap\mathcal{D}\left(
x\right)  $ and so by Exercise \ref{exr.3.9}, $\lambda_{\mathcal{D}\left(
x_{i}\right)  }\left(  B\cap U_{i}\right)  =\lambda_{\mathcal{D}\left(
x\right)  }(B).$ Using this identity in Eq. (\ref{e.3.17}) implies%
\[
\lambda_{M}\left(  B\right)  :=\sum_{i=1}^{\infty}\lambda_{\mathcal{D}\left(
x\right)  }\left(  B\cap U_{i}\right)  =\lambda_{\mathcal{D}\left(  x\right)
}\left(  B\right)
\]
and hence we have proved the existence of $\lambda_{M}.$ The uniqueness
assertion is easy and will be left to the reader.
\end{proof}

\begin{example}
\label{ex.3.11}Let $M=\mathbb{R}^{3}$ with the standard Riemannian metric, and
let $x$ denote the standard coordinates on $M$ determined by $x(m)=m$ for all
$m\in M.$ Then $\lambda_{\mathbb{R}^{3}}$ is Lebesgue measure which in
spherical coordinates may be written as%
\[
d\lambda_{\mathbb{R}^{3}}=r^{2}\sin\varphi drd\varphi d\theta
\]
because $\sqrt{g^{(r,\varphi,\theta)}}=r^{2}\sin\varphi$ by Eq.
\text{(\ref{e.3.7})}. Similarly using Eq. \text{(\ref{e.3.9})},
\[
d\lambda_{M}=\rho^{2}\sin\varphi d\varphi d\theta
\]
when $M\subset\mathbb{R}^{3}$ is the sphere of radius $\rho$ centered at
$0\in\mathbb{R}^{3}.$
\end{example}

\begin{exercise}
\label{exr.3.12}Compute the \textquotedblleft volume
element,\textquotedblright\ $d\lambda_{\mathbb{R}^{3}}\mathrm{,}$ for
$\mathbb{R}^{3}$ in cylindrical coordinates.
\end{exercise}

\begin{theorem}
[Change of Variables Formula]\label{t.3.13}Let $\left(  M,\langle\cdot
,\cdot\rangle_{M}\right)  $ and $\left(  N,\langle\cdot,\cdot\rangle
_{N}\right)  $ be two Riemannian manifolds, $\psi:M\rightarrow N$ be a
diffeomorphism and $\rho\in C^{\infty}\left(  M,\left(  0,\infty\right)
\right)  $ be determined by the equation
\[
\rho\left(  m\right)  =\sqrt{\det\left[  \psi_{\ast m}^{\mathrm{tr}}\psi_{\ast
m}\right]  }\text{ for all }m\in M,
\]
where $\psi_{\ast m}^{\mathrm{tr}}$ denotes the adjoint of $\psi_{\ast m}$
relative to Riemannian inner products on $T_{m}M$ and $T_{\psi\left(
m\right)  }N.$ If $f:N\rightarrow\mathbb{R}_{+}$ is a positive Borel
measurable function, then%
\[
\int_{N}fd\lambda_{N}=\int_{M}\rho\cdot\left(  f\circ\psi\right)  d\lambda
_{M}.
\]
In particular if $\psi$ is an isometry, i.e. $\psi_{\ast m}:T_{m}M\rightarrow
T_{\psi\left(  m\right)  }N$ is orthogonal for all $m,$ then
\[
\int_{N}fd\lambda_{N}=\int_{M}f\circ\psi~d\lambda_{M}.
\]

\end{theorem}

\begin{proof}
By a partition of unity argument (see the proof of Theorem \ref{t.3.10}), it
suffices to consider the case where $f$ has \textquotedblleft
small\textquotedblright\ support, i.e. we may assume that the support of
$f\circ\psi$ is contained in $\mathcal{D}\left(  x\right)  $ for some chart
$x$ on $M.$ Letting $\phi:=x^{-1},$ by Eq. (\ref{e.3.11}) of Lemma
\ref{l.3.7},%
\begin{align*}
&  \frac{\det\left[  \langle\partial_{i}\left(  \psi\circ\phi\right)  \left(
t\right)  ,\partial_{j}\left(  \psi\circ\phi\right)  \left(  t\right)
\rangle_{N}\right]  }{\det\left[  \langle\partial_{i}\phi\left(  t\right)
,\partial_{j}\phi\left(  t\right)  \rangle_{M}\right]  }\\
&  \qquad=\frac{\det\left[  \langle\psi_{\ast}\partial_{i}\phi\left(
t\right)  ,\psi_{\ast}\partial_{j}\phi\left(  t\right)  \rangle_{N}\right]
}{\det\left[  \langle\partial_{i}\phi\left(  t\right)  ,\partial_{j}%
\phi\left(  t\right)  \rangle_{M}\right]  }=\frac{\det\left[  \langle
\psi_{\ast}^{\mathrm{tr}}\psi_{\ast}\partial_{i}\phi\left(  t\right)
,\partial_{j}\phi\left(  t\right)  \rangle_{M}\right]  }{\det\left[
\langle\partial_{i}\phi\left(  t\right)  ,\partial_{j}\phi\left(  t\right)
\rangle_{M}\right]  }\\
&  \qquad=\det\left[  \psi_{\ast\phi\left(  t\right)  }^{\mathrm{tr}}%
\psi_{\ast\phi\left(  t\right)  }\right]  =\rho^{2}\left(  \phi\left(
t\right)  \right)  .
\end{align*}
This implies%
\begin{align*}
\int_{N}fd\lambda_{N}  &  =\int_{\mathcal{R}\left(  x\right)  }f\circ\left(
\psi\circ\phi\right)  (t)\sqrt{\det\left[  \langle\partial_{i}\left(
\psi\circ\phi\right)  \left(  t\right)  ,\partial_{j}\left(  \psi\circ
\phi\right)  \left(  t\right)  \rangle_{N}\right]  }dt\\
&  =\int_{\mathcal{R}\left(  x\right)  \ }\left(  f\circ\psi\right)  \circ
\phi(t)\cdot\rho\left(  \phi\left(  t\right)  \right)  \sqrt{\det\left[
\langle\partial_{i}\phi\left(  t\right)  ,\partial_{j}\phi\left(  t\right)
\rangle_{M}\right]  }dt\\
&  =\int_{\mathcal{D}\left(  x\right)  \ }\left(  f\circ\psi\right)  \cdot
\rho\cdot\sqrt{g^{x}}dx=\int_{M}\rho\cdot f\circ\psi\ d\lambda_{M}.
\end{align*}

\end{proof}

\begin{example}
\label{ex.3.14}Let $M=SL(n,\mathbb{R})$ as in Example \ref{ex.3.3} and let
$\langle\cdot,\cdot\rangle_{M}$ be the metric given by Eq. (\ref{e.3.5}).
Because $L_{g}:M\rightarrow M$ is an isometry, Theorem \ref{t.3.13} implies%
\[
\int_{SL(n,\mathbb{R})}f\left(  gx\right)  d\lambda_{G}\left(  x\right)
=\int_{SL(n,\mathbb{R})}f\left(  x\right)  d\lambda_{G}\left(  x\right)
\text{ for all }g\in G.
\]
That is $\lambda_{G}$ is invariant under left translations by elements of $G$
and such an invariant left invariant measure is called a \textquotedblleft%
\textbf{left Haar}\textquotedblright\ measure on $G.$

Similarly if $G=O\left(  n\right)  $ with Riemannian metric determined by Eq.
(\ref{e.3.5}), then, since $g\in G$ is orthogonal, we have%
\[
ds^{2}(A_{g}):=\text{$\operatorname*{tr}$}((g^{-1}A)^{\ast}g^{-1}%
A)=\text{$\operatorname*{tr}$}((g^{\ast}A)^{\ast}g^{-1}%
A)=\text{$\operatorname*{tr}$}(A^{\ast}gg^{-1}A)=\text{$\operatorname*{tr}$%
}(A^{\ast}A)
\]
and%
\[
\text{$\operatorname*{tr}$}((Ag^{-1})^{\ast}Ag^{-1})=\text{$\operatorname*{tr}%
$}(gA^{\ast}Ag^{-1})=\text{$\operatorname*{tr}$}(A^{\ast}Ag^{-1}%
g)=\text{$\operatorname*{tr}$}(A^{\ast}A).
\]
Therefore, both left and right translations by element $g\in G$ are isometries
for this Riemannian metric on $O\left(  m\right)  $ and so by Theorem
\ref{t.3.13},%
\[
\int_{O\left(  n\right)  }f\left(  gx\right)  d\lambda_{G}\left(  x\right)
=\int_{O\left(  n\right)  }f\left(  x\right)  d\lambda_{G}\left(  x\right)
=\int_{O\left(  n\right)  }f\left(  xg\right)  d\lambda_{G}\left(  x\right)
\]
for all $g\in G.$
\end{example}

\subsection{Gradients, Divergence, and Laplacians\label{s.3.3}}

In the sequel, let $M$ be a Riemannian manifold, $x$ be a chart on $M,$
$g_{ij}:=\langle\partial/\partial x^{i},\partial/\partial x^{j}\rangle,$ and
$ds^{2}=\sum_{i,j=1}^{d}g_{ij}dx^{i}dx^{j}.$

\begin{definition}
\label{d.3.15} Let $g^{ij}$ denote the $i$-$j^{\text{th}}$ -- matrix element
for the inverse matrix to the matrix, $(g_{ij}).$
\end{definition}

Given $f\in C^{\infty}(M)$ and $m\in M,$ $df_{m}:=df|_{T_{m}M}$ is a linear
functional on $T_{m}M.$ Hence there is a unique vector $v_{m}\in T_{m}M$ such
that $df_{m}=\langle v_{m},\cdot\rangle_{m}.$

\begin{definition}
\label{d.3.16}The vector $v_{m}$ above is called the \textbf{gradient} of $f$
at $m$ and will be denoted by either $\operatorname*{grad}f(m)$ or
$\vec{\nabla}f\left(  m\right)  .$
\end{definition}

\begin{exercise}
\ \label{exr.3.17}If $x$ is a chart on $M$ and $m\in\mathcal{D}(x)$ then%
\begin{equation}
\vec{\nabla}f(m)=\text{$\operatorname*{grad}$}f(m)=\sum_{i,j=1}^{d}%
g^{ij}(m)\frac{\partial f(m)}{\partial x^{i}}\frac{\partial}{\partial x^{j}%
}|_{m}, \label{e.3.18}%
\end{equation}
where as usual, $g_{ij}=g_{ij}^{x}$ and $g^{ij}=\left(  g_{ij}\right)  ^{-1}.$
Notice from Eq. (\ref{e.3.18}) that $\vec{\nabla}f$ is a smooth vector field
on $M.$
\end{exercise}

\begin{exercise}
\label{exr.3.18}Suppose $M\subset\mathbb{R}^{N}$ is an imbedded submanifold
with the induced Riemannian structure. Let $F:\mathbb{R}^{N}\rightarrow
\mathbb{R}$ be a smooth function and set $f:=F|_{M}.$ Then
$\operatorname*{grad}f(m)=(P(m)\vec{\nabla}F(m))_{m},$ where $\vec{\nabla
}F(m)$ denotes the usual gradient on $\mathbb{R}^{N},$ and $P(m)$ denotes
orthogonal projection of $\mathbb{R}^{N}$ onto $\tau_{m}M.$
\end{exercise}

We now introduce the divergence of a vector field $Y$ on $M.$

\begin{lemma}
[Divergence]\label{l.3.19}To every smooth vector field $Y$ on $M$ there is a
unique smooth function, $\vec{\nabla}\cdot Y=\operatorname*{div}Y,$ on $M$
such that
\begin{equation}
\int_{M}Yf\,d\lambda_{M}=-\int_{M}\text{$\operatorname*{div}$}Y\cdot
f\,d\lambda_{M},\,\quad\forall~f\in C_{c}^{\infty}(M). \label{e.3.19}%
\end{equation}
(The function, $\vec{\nabla}\cdot Y=\operatorname*{div}Y,$ is called the
\textbf{divergence }of $Y.$) Moreover if $x$ is a chart on $M,$ then on its
domain, $\mathcal{D}(x),$
\begin{equation}
\vec{\nabla}\cdot Y=\operatorname*{div}Y=\sum_{i=1}^{d}{\frac{1}{\sqrt{g}}%
}\frac{\partial(\sqrt{g}Y^{i})}{\partial x^{i}}=\sum_{i=1}^{d}\{\frac{\partial
Y^{i}}{\partial x^{i}}+\frac{\partial\log\sqrt{g}}{\partial x^{i}}Y^{i}\}
\label{e.3.20}%
\end{equation}
where $Y^{i}:=dx^{i}(Y)$ and $\sqrt{g}=\sqrt{g^{x}}.$
\end{lemma}

\begin{proof}
(Sketch) Suppose that $f\in C_{c}^{\infty}(M)$ such that the support of $f$ is
contained in $\mathcal{D}(x).$ Because $Yf=\sum_{i=1}^{d}Y^{i}\partial
f/\partial x^{i},$
\begin{align*}
\int_{M}Yf\,d\lambda_{M}  &  =\int_{M}\sum_{i=1}^{d}Y^{i}\partial f/\partial
x^{i}\cdot\sqrt{g}dx=-\int_{M}\sum_{i=1}^{d}f\frac{\partial(\sqrt{g}\ Y^{i}%
)}{\partial x^{i}}dx\\
&  =-\int_{M}f\sum_{i=1}^{d}\frac{1}{\sqrt{g}}\frac{\partial(\sqrt{g}Y^{i}%
)}{\partial x^{i}}\,d\lambda_{M},
\end{align*}
where the second equality follows by an integration by parts. This shows that
if $\operatorname*{div}Y$ exists it must be given on $\mathcal{D}(x)$ by Eq.
\text{(\ref{e.3.20})}. This proves the uniqueness assertion. Using what we
have already proved, it is easy to conclude that the formula for
$\operatorname*{div}Y$ is chart independent. Hence we may define smooth
function $\operatorname*{div}Y$ on $M$ using Eq. \text{(\ref{e.3.20})} in each
coordinate chart $x$ on $M.$ It is then possible to show (again using a smooth
partition of unity argument) that this function satisfies Eq.
\text{(\ref{e.3.19})}.
\end{proof}

\begin{remark}
\label{r.3.20}We may write Eq. \text{(\ref{e.3.19})} as
\begin{equation}
\int_{M}\langle Y,\operatorname*{grad}f\rangle\,d\lambda_{M}=-\int
_{M}\text{$\operatorname*{div}$}Y\cdot f\,d\lambda_{M},\,\quad\forall~f\in
C_{c}^{\infty}(M), \label{e.3.21}%
\end{equation}
so that $\operatorname*{div}$ is the negative of the formal adjoint of
$\operatorname*{grad}.$
\end{remark}

\begin{exercise}
[Product Rule]\label{exr.3.21}If $f\in C^{\infty}\left(  M\right)  $ and
$Y\in\Gamma\left(  TM\right)  $ then%
\[
\vec{\nabla}\cdot\left(  fY\right)  =\langle\vec{\nabla}f,Y\rangle
+f\ \vec{\nabla}\cdot Y.
\]

\end{exercise}

\begin{lemma}
[Integration by Parts]\label{l.3.22} Suppose that $Y\in\Gamma(TM),$ $f\in
C_{c}^{\infty}(M),$ and $h\in C^{\infty}(M),$ then
\[
\int_{M}Yf\cdot h\,d\lambda_{M}=\int_{M}f\{-Yh-h\cdot
\text{$\operatorname*{div}$}Y\}\,d\lambda_{M}.
\]

\end{lemma}

\begin{proof}
By the definition of $\operatorname*{div}Y$ and the product rule,%
\[
\int_{M}fh\ \text{$\operatorname*{div}$}Y\,d\lambda_{M}=-\int_{M}%
Y(fh)\,d\lambda_{M}=-\int_{M}\{hYf+fYh\}\,d\lambda_{M}.
\]

\end{proof}

\begin{definition}
\label{d.3.23}The \textbf{Laplacian }on $M$ is the second order differential
operator, $\Delta:C^{\infty}(M)\rightarrow C^{\infty}(M),$ defined by
\begin{equation}
\Delta f:=\text{$\operatorname*{div}$}(\operatorname*{grad}f)=\vec{\nabla
}\cdot\vec{\nabla}f. \label{e.3.22}%
\end{equation}
In local coordinates,%
\begin{equation}
\Delta f={\frac{1}{\sqrt{g}}}\sum_{i,j=1}^{d}\partial_{i}\{\sqrt{g}%
g^{ij}\partial_{j}f\}, \label{e.3.23}%
\end{equation}
where $\partial_{i}=\partial/\partial x^{i},$ $g=g^{x},$ $\sqrt{g}=\sqrt{\det
g},$ and $(g^{ij})=(g_{ij}^{x})^{-1}.$
\end{definition}

\begin{remark}
\label{r.3.24}The Laplacian, $\Delta f,$ may be characterized by the
equation:
\[
\int_{M}\Delta f\cdot h\,d\lambda_{M}=-\int_{M}\langle\vec{\nabla}%
f,\vec{\nabla}h\rangle\,d\lambda_{M},
\]
which is to hold for all $f\in C^{\infty}(M)$ and $h\in C_{c}^{\infty}(M).$
\end{remark}

\begin{example}
\label{ex.3.25}Suppose that $M=\mathbb{R}^{N}$ with the standard Riemannian
metric $ds^{2}=\sum_{i=1}^{N}(dx^{i})^{2},$ then the standard formulas:
\[
\operatorname*{grad}f=\sum_{i=1}^{N}\partial f/\partial x^{i}\cdot
\partial/\partial x^{i},\text{ $\operatorname*{div}$}Y=\sum_{i=1}^{N}\partial
Y^{i}/\partial x^{i}\text{ and }\Delta f=\sum_{i=1}^{N}\frac{\partial^{2}%
f}{(\partial x^{i})^{2}}%
\]
are easily verified, where $f$ is a smooth function on $\mathbb{R}^{N}$ and
$Y=\sum_{i=1}^{N}Y^{i}\partial/\partial x^{i}$ is a smooth vector-field.
\end{example}

\begin{exercise}
\label{exr.3.26}Let $M=\mathbb{R}^{3},$ $(r,\varphi,\theta)$ be spherical
coordinates on $\mathbb{R}^{3},$ $\partial_{r}=\partial/\partial r,$
$\partial_{\varphi}=\partial/\partial\varphi,$ and $\partial_{\theta}%
=\partial/\partial_{\theta}.$ Given a smooth function $f$ and a vector-field
$Y=Y_{r}\partial_{r}+Y_{\varphi}\partial_{\varphi}+Y_{\theta}\partial_{\theta
}$ on $\mathbb{R}^{3}$ verify:
\[
\operatorname*{grad}f=(\partial_{r}f)\partial_{r}+\frac{1}{r^{2}}%
(\partial_{\varphi}f)\partial_{\varphi}+\frac{1}{r^{2}\sin^{2}\varphi
}(\partial_{\theta}f)\partial_{\theta},
\]%
\begin{align*}
\text{$\operatorname*{div}$}Y &  =\frac{1}{r^{2}\sin\varphi}\{\partial
_{r}(r^{2}\sin\varphi Y_{r})+\partial_{\varphi}(r^{2}\sin\varphi Y_{\varphi
})+r^{2}\sin\varphi\partial_{\theta}Y_{\theta}\}\\
&  =\frac{1}{r^{2}}\partial_{r}(r^{2}Y_{r})+\frac{1}{\sin\varphi}%
\partial_{\varphi}(\sin\varphi Y_{\varphi})+\partial_{\theta}Y_{\theta},
\end{align*}
and
\[
\Delta f=\frac{1}{r^{2}}\partial_{r}(r^{2}\partial_{r}f)+\frac{1}{r^{2}%
\sin\varphi}\partial_{\varphi}(\sin\varphi\partial_{\varphi}f)+\frac{1}%
{r^{2}\sin^{2}\varphi}\partial_{\theta}^{2}f.
\]

\end{exercise}

\begin{example}
\label{ex.3.27}Let $M=G=O\left(  n\right)  $ with Riemannian metric determined
by Eq. (\ref{e.3.5}) and for $A\in\mathfrak{g}:=T_{e}G$ let $\tilde{A}%
\in\Gamma\left(  TG\right)  $ be the left invariant vector field,%
\[
\tilde{A}\left(  x\right)  :=L_{x\ast}A=\frac{d}{dt}|_{0}xe^{tA}%
\]
as was done for $SL(n,\mathbb{R})$ in Example \ref{ex.2.34}. Using the
invariance of $d\lambda_{G}$ under right translations established in Example
\ref{ex.3.14}, we find for $f,h\in C^{1}\left(  G\right)  $ that%
\begin{align*}
\int_{G}\tilde{A}f\left(  x\right)  \cdot h\left(  x\right)  d\lambda
_{G}\left(  x\right)   &  =\int_{G}\frac{d}{dt}|_{0}f\left(  xe^{tA}\right)
\cdot h\left(  x\right)  d\lambda_{G}\left(  x\right) \\
&  =\frac{d}{dt}|_{0}\int_{G}f\left(  xe^{tA}\right)  \cdot h\left(  x\right)
d\lambda_{G}\left(  x\right) \\
&  =\frac{d}{dt}|_{0}\int_{G}f\left(  x\right)  \cdot h\left(  xe^{-tA}%
\right)  d\lambda_{G}\left(  x\right) \\
&  =\int_{G}f\left(  x\right)  \cdot\frac{d}{dt}|_{0}h\left(  xe^{-tA}\right)
d\lambda_{G}\left(  x\right) \\
&  =-\int_{G}f\left(  x\right)  \cdot\tilde{A}h\left(  x\right)  d\lambda
_{G}\left(  x\right)  .
\end{align*}
Taking $h\equiv1$ implies%
\begin{align*}
0  &  =\int_{G}\tilde{A}f\left(  x\right)  d\lambda_{G}\left(  x\right)
=\int_{G}\left\langle \tilde{A}\left(  x\right)  ,\vec{\nabla}f\left(
x\right)  \right\rangle d\lambda_{G}\left(  x\right) \\
&  =-\int_{G}\vec{\nabla}\cdot\tilde{A}\left(  x\right)  \cdot f\left(
x\right)  d\lambda_{G}\left(  x\right)
\end{align*}
from which we learn $\vec{\nabla}\cdot\tilde{A}=0.$

Now letting $S_{0}\subset\mathfrak{g}$ be an orthonormal basis for
$\mathfrak{g},$ because $L_{g\ast}$ is an isometry, $\{\tilde{A}\left(
g\right)  :A\in S_{0}\}$ is an orthonormal basis for $T_{g}G$ for all $g\in
G.$ Hence%
\[
\vec{\nabla}f\left(  g\right)  =\sum_{A\in S_{0}}\left\langle \vec{\nabla
}f\left(  g\right)  ,\tilde{A}\left(  g\right)  \right\rangle \tilde{A}\left(
g\right)  =\sum_{A\in S_{0}}\left(  \tilde{A}f\right)  \left(  g\right)
\tilde{A}\left(  g\right)  .
\]
and, by the product rule and $\vec{\nabla}\cdot\tilde{A}=0,$%
\[
\Delta f=\vec{\nabla}\cdot\vec{\nabla}f=\sum_{A\in S_{0}}\vec{\nabla}%
\cdot\left[  \left(  \tilde{A}f\right)  \tilde{A}\right]  =\sum_{A\in S_{0}%
}\left\langle \vec{\nabla}\tilde{A}f,\tilde{A}\right\rangle =\sum_{A\in S_{0}%
}\tilde{A}^{2}f.
\]

\end{example}

\subsection{Covariant Derivatives and Curvature\label{s.3.4}}

\begin{definition}
\label{d.3.28}We say a smooth path $s\rightarrow V(s)$ in $TM$ is a
\textbf{vector-field along a smooth path} $s\rightarrow\sigma(s)$ in $M$ if
$\pi\circ V(s)=\sigma(s),$ i.e. $V(s)\in T_{\sigma(s)}M$ for all $s.$ (Recall
that $\pi$ is the canonical projection defined in Definition \ref{d.2.16}.)
\end{definition}

Note: if $V$ is a smooth path in $TM$ then $V$ is a vector-field along
$\sigma:=\pi\circ V.$ This section is motivated by the desire to have the
notion of the derivative of a smooth path $V(s)\in TM.$ On one hand, since
$TM$ is a manifold, we may write $V^{\prime}(s)$ as an element of $TTM.$
However, this is not what we will want for later purposes. We would like the
derivative of $V$ to again be a path back in $TM,$ not in $TTM.$ In order to
define such a derivative, we will need to use more than just the manifold
structure of $M,$ see Definition \ref{d.3.31} below.

\begin{notation}
\label{n.3.29}In the sequel, we assume that $M^{d}$ is an imbedded submanifold
of an inner product space $(E=\mathbb{R}^{N},\langle\cdot,\cdot\rangle),$ and
that $M$ is equipped with the inherited Riemannian metric. Also let $P(m)$
denote orthogonal projection of $E$ onto $\tau_{m}M$ for all $m\in M$ and
$Q(m):=I-P(m)$ be orthogonal projection onto $(\tau_{m}M)^{\perp}.$
\end{notation}

The following elementary lemma will be used throughout the sequel.

\begin{lemma}
\label{l.3.30}The differentials of the orthogonal projection operators, $P$
and $Q,$ satisfy%
\begin{align*}
0  &  =dP+dQ,\\
PdQ  &  =-dPQ=dQQ\text{ and}\\
QdP  &  =-dQP=dPP.
\end{align*}
In particular,%
\[
QdPQ=QdQQ=PdPP=PdQP=0.
\]

\end{lemma}

\begin{proof}
The first equality comes from differentiating the identity, $I=P+Q,$ the
second from differentiating $0=PQ$ and the third from differentiating $0=QP.$
\end{proof}

\begin{figure}
[ptbh]
\begin{center}
\includegraphics[
height=2.0823in,
width=3.3226in
]%
{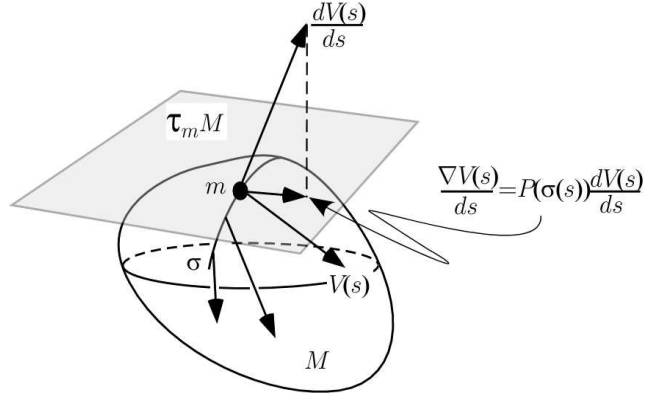}%
\caption{The Levi-Civita covariant derivative.}%
\label{fig.9}%
\end{center}
\end{figure}

\begin{definition}
[Levi-Civita Covariant Derivative]\label{d.3.31}Let $V(s)=(\sigma
(s),v(s))=v(s)_{\sigma(s)}$ be a smooth path in $TM$ (see Figure \ref{fig.9}),
then the \textbf{covariant derivative},\textbf{ }$\nabla V(s)/ds,$ is the
vector field along $\sigma$ defined by%
\begin{equation}
\frac{\nabla V(s)}{ds}:=(\sigma(s),P(\sigma(s))\frac{d}{ds}v(s)).
\label{e.3.24}%
\end{equation}

\end{definition}%

%TCIMACRO{\FRAME{ftbphFU}{3.3226in}{2.0823in}{0pt}{\Qcb{The Levi-Civita
%covariant derivative.}}{\Qlb{fig.9}}{f8.eps}%
%width 3.3226in;  height 2.0823in;  depth 0pt;  original-width 4.1045in;
%original-height 2.5627in;  cropleft "0";  croptop "1";  cropright "1";
%cropbottom "0";  filename 'f8.eps';file-properties "XNPEU";}%
%}}%
%BeginExpansion
%EndExpansion

\begin{proposition}
[Properties of $\nabla/ds$]\label{p.3.32}Let $W(s)=(\sigma(s),w(s))$ and
$V(s)=(\sigma(s),v(s))$ be two smooth vector fields along a path\ $\sigma$ in
$M.$ Then:

\begin{enumerate}
\item $\nabla W(s)/ds$ may be computed as:
\begin{equation}
\frac{\nabla W(s)}{ds}:=(\sigma(s),\frac{d}{ds}w(s)+(dQ(\sigma^{\prime
}(s)))w(s)). \label{e.3.25}%
\end{equation}

\item $\nabla$ is \textbf{metric compatible}, i.e.
\begin{equation}
\frac{d}{ds}\langle W(s),V(s)\rangle=\langle\frac{\nabla W(s)}{ds}%
,V(s)\rangle+\langle W(s),\frac{\nabla V(s)}{ds}\rangle. \label{e.3.26}%
\end{equation}

Now suppose that $(s,t)\rightarrow\sigma(s,t)$ is a smooth function into $M,$
$W(s,t)=(\sigma(s,t),w(s,t))$ is a smooth function into $TM,$ $\sigma^{\prime
}(s,t):=(\sigma(s,t),\frac{d}{ds}\sigma(s,t))$ and $\dot{\sigma}%
(s,t)=(\sigma(s,t),\frac{d}{dt}\sigma(s,t)).$ (Notice by assumption that
$w(s,t)\in T_{\sigma(s,t)}M$ for all $(s,t).$)

\item $\nabla$ has \textbf{zero torsion}, i.e.%
\begin{equation}
\frac{\nabla\sigma^{\prime}}{dt}=\frac{\nabla\dot{\sigma}}{ds}. \label{e.3.27}%
\end{equation}

\item If $R$ is the \textbf{curvature tensor} of $\nabla$ defined by
\begin{equation}
R(u_{m},v_{m})w_{m}=(m,[dQ(u_{m}),dQ(v_{m})]w), \label{e.3.28}%
\end{equation}
then
\begin{equation}
\left[  \frac{\nabla}{dt},\frac{\nabla}{ds}\right]  W:=(\frac{\nabla}{dt}%
\frac{\nabla}{ds}-\frac{\nabla}{ds}\frac{\nabla}{dt})W=R(\dot{\sigma}%
,\sigma^{\prime})W. \label{e.3.29}%
\end{equation}

\end{enumerate}
\end{proposition}

\begin{proof}
Differentiate the identity, $P(\sigma(s))w(s)=w(s),$ relative to $s$ implies
\[
(dP(\sigma^{\prime}(s)))w(s)+P(\sigma(s))\frac{d}{ds}w(s)=\frac{d}{ds}w(s)
\]
from which Eq. \text{(\ref{e.3.25}) follows.}

For Eq. \text{(\ref{e.3.26})} just compute:
\begin{align*}
\frac{d}{ds}\langle W(s),V(s)\rangle &  =\frac{d}{ds}\left\langle
w(s),v(s)\right\rangle \\
&  =\left\langle \frac{d}{ds}w(s),v(s)\right\rangle +\left\langle
w(s),\frac{d}{ds}v(s)\right\rangle \\
&  =\left\langle \frac{d}{ds}w(s),P(\sigma(s))v(s)\right\rangle +\left\langle
P(\sigma(s))w(s),\frac{d}{ds}v(s)\right\rangle \\
&  =\left\langle P(\sigma(s))\frac{d}{ds}w(s),v(s)\right\rangle +\left\langle
w(s),P(\sigma(s))\frac{d}{ds}v(s)\right\rangle \\
&  =\left\langle \frac{\nabla W(s)}{ds},V(s)\right\rangle +\left\langle
W(s),\frac{\nabla V(s)}{ds}\right\rangle ,
\end{align*}
where the third equality relies on $v(s)$ and $w(s)$ being in $\tau
_{\sigma(s)}M$ and the fourth equality relies on $P(\sigma(s))$ being an
orthogonal projection.

From the definitions of $\sigma^{\prime},$ $\dot{\sigma},$ $\nabla/dt,$
$\nabla/ds$ and the fact that mixed partial derivatives commute,%
\begin{align*}
\frac{\nabla\sigma^{\prime}(s,t)}{dt}  &  =\frac{\nabla}{dt}(\sigma
(t,s),\sigma^{\prime}(s,t))=(\sigma(t,s),P(\sigma(s,t))\frac{d}{dt}\frac
{d}{ds}\sigma(t,s))\\
&  =(\sigma(t,s),P(\sigma(s,t))\frac{d}{ds}\frac{d}{dt}\sigma(t,s))=\nabla
\dot{\sigma}(s,t)/ds,
\end{align*}
which proves Eq. \text{(\ref{e.3.27})}.

For Eq. \text{(\ref{e.3.29})} we observe,
\begin{align*}
\frac{\nabla}{dt}\frac{\nabla}{ds}W(s,t)  &  =\frac{\nabla}{dt}(\sigma
(s,t),\frac{d}{ds}w(s,t)+dQ(\sigma^{\prime}(s,t))w(s,t))\\
&  =(\sigma(s,t),\eta_{+}(s,t))
\end{align*}
where (with the arguments $(s,t)$ suppressed from the notation)
\begin{align*}
\eta_{+}  &  =\frac{d}{dt}\left[  \frac{d}{ds}w+dQ(\sigma^{\prime})w\right]
+dQ(\dot{\sigma})\left[  \frac{d}{ds}w+dQ(\sigma^{\prime})w\right] \\
&  =\frac{d}{dt}\frac{d}{ds}w+\left(  \frac{d}{dt}\left[  dQ(\sigma^{\prime
})\right]  \right)  w+dQ(\sigma^{\prime})\frac{d}{dt}w+dQ(\dot{\sigma}%
)\frac{d}{ds}w+dQ(\dot{\sigma})dQ(\sigma^{\prime})w.
\end{align*}
Therefore
\[
\left[  \frac{\nabla}{dt},\frac{\nabla}{ds}\right]  W=(\sigma,\eta_{+}%
-\eta_{-}),
\]
where $\eta_{-}$ is defined the same as $\eta_{+}$ with all $s$ and $t$
derivatives interchanged. Hence, it follows (using again $\frac{d}{dt}\frac
{d}{ds}w=\frac{d}{ds}\frac{d}{dt}w)$ that
\[
\left[  \frac{\nabla}{dt},\frac{\nabla}{ds}\right]  W=(\sigma,[\frac{d}%
{dt}(dQ(\sigma^{\prime}))]w-[\frac{d}{ds}(dQ(\dot{\sigma}))]w+[dQ(\dot{\sigma
}),dQ(\sigma^{\prime})]w).
\]
The proof of Eq. (\ref{e.3.28}) is finished because
\[
\frac{d}{dt}(dQ(\sigma^{\prime}))-\frac{d}{ds}(dQ(\dot{\sigma}))=\frac{d}%
{dt}\frac{d}{ds}(Q\circ\sigma)-\frac{d}{ds}\frac{d}{dt}(Q\circ\sigma)=0.
\]

\end{proof}

\begin{example}
\label{ex.3.33}Let $M=\{m\in\mathbb{R}^{N}:|m|=\rho\}$ be the sphere of radius
$\rho.$ In this case $Q(m)=${$\frac{1}{\rho^{2}}$}$mm^{\mathrm{tr}}$ for all
$m\in M.$ Therefore
\[
dQ(v_{m})=\frac{1}{\rho^{2}}\{vm^{\mathrm{tr}}+mv^{\mathrm{tr}}\}~\forall
~v_{m}\in T_{m}M
\]
and hence%
\begin{align*}
dQ(u_{m})dQ(v_{m})  &  =\frac{1}{\rho^{4}}\{um^{\mathrm{tr}}+mu^{\mathrm{tr}%
}\}\{vm^{\mathrm{tr}}+mv^{\mathrm{tr}}\}\\
&  =\frac{1}{\rho^{4}}\{\rho^{2}uv^{\mathrm{tr}}+\langle u,v\rangle Q(m)\}.
\end{align*}
So the curvature tensor is given by
\[
R(u_{m},v_{m})w_{m}=(m,{\frac{1}{\rho^{2}}}\{uv^{\mathrm{tr}}-vu^{\mathrm{tr}%
}\}w)=(m,{\frac{1}{\rho^{2}}}\{\langle v,w\rangle u-\langle u,w\rangle v\}).
\]

\end{example}

\begin{exercise}
\ \label{exr.3.34} Show the curvature tensor of the cylinder%
\[
M=\{(x,y,z)\in\mathbb{R}^{3}:x^{2}+y^{2}=1\}
\]
is zero.
\end{exercise}

\begin{definition}
[Covariant Derivative on $\Gamma(TM)$]\label{d.3.35}Suppose that $Y$ is a
vector field on $M$ and $v_{m}\in T_{m}M.$ Define $\nabla_{v_{m}}Y\in T_{m}M$
by
\[
\nabla_{v_{m}}Y:=\frac{\nabla Y(\sigma(s))}{ds}|_{s=0},
\]
where $\sigma$ is any smooth path in $M$ such that $\sigma^{\prime}(0)=v_{m}.$
\end{definition}

If $Y(m)=(m,y(m)),$ then
\[
\nabla_{v_{m}}Y=(m,P(m)dy(v_{m}))=(m,dy(v_{m})+dQ(v_{m})y(m)),
\]
from which it follows $\nabla_{v_{m}}Y$ is well defined, i.e. $\nabla_{v_{m}%
}Y$ is independent of the choice of $\sigma$ such that $\sigma^{\prime}\left(
0\right)  =v_{m}.$ The following proposition relates curvature and torsion to
the covariant derivative $\nabla$ on vector fields.

\begin{proposition}
\label{p.3.36}Let $m\in M,$ $v\in T_{m}M,$ $X,Y,Z\in\Gamma(TM),$ and $f\in
C^{\infty}(M),$ then the following relations hold.

\begin{description}
\item[1. Product Rule] $\nabla_{v}(f\cdot X)=df(v)\cdot X(m)+f(m)\cdot
\nabla_{v}X.$

\item[2. Zero Torsion] $\nabla_{X}Y-\nabla_{Y}X-[X,Y]=0.$

\item[3. Zero Torsion] For all $v_{m},w_{m}\in T_{m}M,$ $dQ(v_{m}%
)w_{m}=dQ(w_{m})v_{m}.$

\item[4. Curvature Tensor] $R(X,Y)Z=[\nabla_{X},\nabla_{Y}]Z-\nabla_{\lbrack
X,Y]}Z,$ where
\[
\lbrack\nabla_{X},\nabla_{Y}]Z:=\nabla_{X}(\nabla_{Y}Z)-\nabla_{Y}(\nabla
_{X}Z).
\]
Moreover if $u,v,w,z\in T_{m}M,$ then $R$ has the following symmetries

\begin{description}
\item[a] $R(u_{m},v_{m})=-R(v_{m},u_{m})$

\item[b] $\left[  R(u_{m},v_{m})\right]  ^{\mathrm{tr}}=-R(u_{m},v_{m})$ and

\item[c] if $z_{m}\in\tau_{m}M,$ then
\begin{equation}
\langle R(u_{m},v_{m})w_{m},z_{m}\rangle=\langle R(w_{m},z_{m})u_{m}%
,v_{m}\rangle. \label{e.3.30}%
\end{equation}

\end{description}

\item[5. Ricci Curvature Tensor] For each $m\in M,$ let
$\mathrm{\operatorname*{Ric}}_{m}:T_{m}M\rightarrow T_{m}M$ be defined by%
\begin{equation}
\operatorname{Ric}_{m}v_{m}:=\sum_{a\in S}R(v_{m},a)a, \label{e.3.31}%
\end{equation}
where $S\subset T_{m}M$ is an orthonormal basis. Then $\operatorname*{Ric}%
_{m}^{\mathrm{tr}}=\operatorname{Ric}_{m}$ and $\operatorname*{Ric}%
\nolimits_{m}$ may be computed as%
\begin{equation}
\langle\mathrm{\operatorname*{Ric}}_{m}u,v\rangle=\operatorname*{tr}%
(dQ(dQ(u)v)-dQ(v)dQ(u))\text{ for all }u,v\in T_{m}M. \label{e.3.32}%
\end{equation}

\end{description}
\end{proposition}

\begin{proof}
The product rule is easily checked and may be left to the reader. For the
second and third items, write $X(m)=(m,x(m)),$ $Y(m)=(m,y(m)),$ and
$Z(m)=(m,z(m))$ where $x,y,z:M\rightarrow\mathbb{R}^{N}$ are smooth functions
such that $x(m),$ $y(m),$ and $z(m)$ are in $\tau_{m}M$ for all $m\in M.$ Then
using Eq. \text{(\ref{e.2.15})}, we have
\begin{align}
(\nabla_{X}Y-\nabla_{Y}X)(m)  &  =(m,P(m)(dy(X(m))-dx(Y(m))))\nonumber\\
&  =(m,(dy(X(m))-dx(Y(m))))=[X,Y](m), \label{e.3.33}%
\end{align}
which proves the second item. Since $(\nabla_{X}Y)(m)$ may also be written as%
\[
(\nabla_{X}Y)(m)=(m,dy(X(m))+dQ(X(m))y(m)),
\]
Eq. (\ref{e.3.33}) may be expressed as $dQ(X(m))y(m)=dQ(Y(m))x(m)$ which
implies the third item.

Similarly for fourth item:
\begin{align*}
\nabla_{X}\nabla_{Y}Z  &  =\nabla_{X}(\cdot,Yz+(YQ)z)\\
&  =(\cdot,XYz+(XYQ)z+(YQ)Xz+(XQ)(Yz+(YQ)z)),
\end{align*}
where $YQ:=dQ(Y)$ and $Yz:=dz(Y).$ Interchanging $X$ and $Y$ in this last
expression and then subtracting gives:
\begin{align*}
\lbrack\nabla_{X},\nabla_{Y}]Z  &  =(\cdot,[X,Y]z+([X,Y]Q)z+[XQ,YQ]z)\\
&  =\nabla_{\lbrack X,Y]}Z+R(X,Y)Z.
\end{align*}
The anti-symmetry properties in items 4a) and 4b) follow easily from Eq.
(\ref{e.3.28}). For example for 4b), $dQ\left(  u_{m}\right)  $ and
$dQ(v_{m})$ are symmetric operators and hence%
\begin{align*}
\left[  R(u_{m},v_{m})\right]  ^{\mathrm{tr}}  &  =[dQ(u_{m}),dQ(v_{m}%
)]^{\mathrm{tr}}=[dQ(v_{m})^{\mathrm{tr}},dQ(u_{m})^{\mathrm{tr}}]\\
&  =[dQ(v_{m}),dQ(u_{m})]=-[dQ(u_{m}),dQ(v_{m})]=-R(u_{m},v_{m}).
\end{align*}
To prove Eq. (\ref{e.3.30}) we make use of the zero - torsion condition
$dQ(v_{m})w_{m}=dQ(w_{m})v_{m}$ and the fact that $dQ\left(  u_{m}\right)  $
is symmetric to learn%
\begin{align}
\langle R(u_{m},v_{m})w,z\rangle &  =\langle\lbrack dQ(u_{m}),dQ(v_{m}%
)]w,z\rangle\nonumber\\
&  =\langle\lbrack dQ(u_{m})dQ(v_{m})-dQ(v_{m})dQ(u_{m})]w,z\rangle\nonumber\\
&  =\langle dQ(v_{m})w,dQ(u_{m})z\rangle-\langle dQ(u_{m})w,dQ(v_{m}%
)z\rangle\nonumber\\
&  =\langle dQ(w)v,dQ(z)u\rangle-\langle dQ(w)u,dQ(z)v\rangle\label{e.3.34}\\
&  =\langle\left[  dQ(z),dQ(w)\right]  v,u\rangle=\langle R\left(  z,w\right)
v,u\rangle=\langle R\left(  w,z\right)  u,v\rangle\nonumber
\end{align}
where we have used the anti-symmetry properties in 4a. and 4b. By Eq.
(\ref{e.3.34}) with $v=w=a,$%

\begin{align*}
\langle\operatorname*{Ric}u,z\rangle &  =\sum_{a\in S}\langle R(u,a)a,z\rangle
\\
&  =\sum_{a\in S}\left[  \langle dQ(a)a,dQ(u)z\rangle-\langle
dQ(u)a,dQ(a)z\rangle\right] \\
&  =\sum_{a\in S}\left[  \langle a,dQ(a)dQ(u)z\rangle-\langle
dQ(u)a,dQ(z)a\rangle\right] \\
&  =\sum_{a\in S}\left[  \langle a,dQ(dQ(u)z)a\rangle-\langle
dQ(z)dQ(u)a,a\rangle\right] \\
&  =\operatorname*{tr}(dQ(dQ(u)z)-dQ(z)dQ(u))
\end{align*}
which proves Eq. (\ref{e.3.32}). The assertion that $\operatorname{Ric}%
_{m}:T_{m}M\rightarrow T_{m}M$ is a symmetric operator follows easily from
this formula and item 3.
\end{proof}

\begin{notation}
\label{n.3.37}To each $v\in\mathbb{R}^{N},$ let $\partial_{v}$ denote the
vector field on $\mathbb{R}^{N}$ defined by
\[
\partial_{v}(\text{at }x)=v_{x}=\frac{d}{dt}|_{0}(x+tv).
\]
So if $F\in C^{\infty}(\mathbb{R}^{N}),$ then
\[
(\partial_{v}F)(x):=\frac{d}{dt}|_{0}F(x+tv)=F^{\prime}\left(  x\right)  v
\]
and%
\[
\left(  \partial_{v}\partial_{w}F\right)  \left(  x\right)  =F^{\prime\prime
}\left(  x\right)  \left(  v,w\right)  ,
\]
see Notation \ref{n.2.1}.
\end{notation}

Notice that if $w:\mathbb{R}^{N}\rightarrow\mathbb{R}^{N}$ is a function and
$v\in\mathbb{R}^{N},$ then%
\[
\left(  \partial_{v}\partial_{w}F\right)  \left(  x\right)  =\partial
_{v}\left[  F^{\prime}\left(  \cdot\right)  w\left(  \cdot\right)  \right]
\left(  x\right)  =F^{\prime}\left(  x\right)  \partial_{v}w\left(  x\right)
+F^{\prime\prime}\left(  x\right)  \left(  v,w\left(  x\right)  \right)  .
\]
The following variant of item 4. of Proposition \ref{p.3.36} will be useful in
proving the key Bochner-Weitenb\"{o}ck identity in Theorem \ref{t.3.49} below.

\begin{proposition}
\label{p.3.38}Suppose that $Z\in\Gamma\left(  TM\right)  ,$ $v,w\in T_{m}M$
and let $X,Y\in\Gamma\left(  TM\right)  $ such that $X\left(  m\right)  =v$
and $Y\left(  m\right)  =w.$ Then

\begin{enumerate}
\item $\nabla_{v\otimes w}^{2}Z$ defined by%
\begin{equation}
\nabla_{v\otimes w}^{2}Z:=\left(  \nabla_{X}\nabla_{Y}Z-\nabla_{\nabla_{X}%
Y}Z\right)  \left(  m\right)  \label{e.3.35}%
\end{equation}
is well defined, independent of the possible choices for $X$ and $Y.$

\item If $Z(m)=(m,z(m))$ with $z:\mathbb{R}^{N}\rightarrow\mathbb{R}^{N}$ a
smooth function such $z\left(  m\right)  \in\tau_{m}M$ for all $m\in M,$ then
\begin{equation}
\nabla_{v\otimes w}^{2}Z=dQ\left(  v\right)  dQ\left(  w\right)  z\left(
m\right)  +P\left(  m\right)  z^{\prime\prime}\left(  m\right)  \left(
v,w\right)  -P\left(  m\right)  z^{\prime}\left(  m\right)  \left[  dQ\left(
v\right)  w\right]  . \label{e.3.36}%
\end{equation}

\item The curvature tensor $R\left(  v,w\right)  $ may be computed as%
\begin{equation}
\nabla_{v\otimes w}^{2}Z-\nabla_{w\otimes v}^{2}Z=R\left(  v,w\right)
Z\left(  m\right)  . \label{e.3.37}%
\end{equation}

\item If $V$ is a smooth vector field along a path $\sigma\left(  s\right)  $
in $M,$ then the following product rule holds,
\begin{equation}
\frac{\nabla}{ds}\left(  \nabla_{V\left(  s\right)  }Z\right)  =\left(
\nabla_{\frac{\nabla}{ds}V\left(  s\right)  }Z\right)  +\nabla_{\sigma
^{\prime}\left(  s\right)  \otimes V\left(  s\right)  }^{2}Z. \label{e.3.38}%
\end{equation}

\end{enumerate}
\end{proposition}

\begin{proof}
We will prove items 1. and 2. by showing the right sides of Eq. (\ref{e.3.35})
and Eq. (\ref{e.3.36}) are equal. To do this write $X(m)=(m,x(m)),$
$Y(m)=(m,y(m)),$ and $Z(m)=(m,z(m))$ where $x,y,z:\mathbb{R}^{N}%
\rightarrow\mathbb{R}^{N}$ are smooth functions such that $x(m),$ $y(m),$ and
$z(m)$ are in $\tau_{m}M$ for all $m\in M.$ Then, suppressing $m$ from the
notation,%
\begin{align*}
\nabla_{X}\nabla_{Y}Z-\nabla_{\nabla_{X}Y}Z  &  =P\partial_{x}\left[
P\partial_{y}z\right]  -P\partial_{P\partial_{x}y}z\\
&  =P\left(  \partial_{x}P\right)  \partial_{y}z+P\partial_{x}\partial
_{y}z-P\partial_{P\partial_{x}y}z\\
&  =P\left(  \partial_{x}P\right)  \partial_{y}z+Pz^{\prime\prime}\left(
x,y\right)  +Pz^{\prime}\left[  \partial_{x}y-P\partial_{x}y\right] \\
&  =\left(  \partial_{x}P\right)  Q\partial_{y}z+Pz^{\prime\prime}\left(
x,y\right)  +Pz^{\prime}\left[  Q\partial_{x}y\right]  .
\end{align*}
Differentiating the identity, $Qy=0$ on $M$ shows $Q\partial_{x}y=-\left(
\partial_{x}Q\right)  y$ which combined with the previous equation gives%
\begin{align}
\nabla_{X}\nabla_{Y}Z-\nabla_{\nabla_{X}Y}Z  &  =\left(  \partial_{x}P\right)
Q\partial_{y}z+Pz^{\prime\prime}\left(  x,y\right)  -Pz^{\prime}\left[
\left(  \partial_{x}Q\right)  Y\right] \label{e.3.39}\\
&  =-\left(  \partial_{x}P\right)  \left(  \partial_{y}Q\right)
z+Pz^{\prime\prime}\left(  X,Y\right)  -Pz^{\prime}\left[  \left(
\partial_{x}Q\right)  Y\right]  .\nonumber
\end{align}
Evaluating this expression at $m$ proves the right side of Eq. (\ref{e.3.36}).

Equation (\ref{e.3.37}) now follows from Eqs. (\ref{e.3.36}) and
(\ref{e.3.28}), item 3. of Proposition \ref{p.3.36} and the fact the
$z^{\prime\prime}\left(  v,w\right)  =z^{\prime\prime}\left(  w,v\right)  $
because mixed partial derivatives commute.

We give two proofs of Eq. (\ref{e.3.38}). For the first proof, choose local
vector fields $\left\{  E_{i}\right\}  _{i=1}^{d}$ defined in a neighborhood
of $\sigma\left(  s\right)  $ such that $\left\{  E_{i}\left(  \sigma\left(
s\right)  \right)  \right\}  _{i=1}^{d}$ is a basis for $T_{\sigma\left(
s\right)  }M$ for each $s.$ We may then write $V\left(  s\right)  =\sum
_{i=1}^{d}V_{i}\left(  s\right)  E_{i}\left(  \sigma\left(  s\right)  \right)
$ and therefore,%
\begin{equation}
\frac{\nabla}{ds}V\left(  s\right)  =\sum_{i=1}^{d}\left\{  V_{i}^{\prime
}\left(  s\right)  E_{i}\left(  \sigma\left(  s\right)  \right)  +V_{i}\left(
s\right)  \nabla_{\sigma^{\prime}\left(  s\right)  }E_{i}\right\}
\label{e.3.40}%
\end{equation}
and
\begin{align*}
\frac{\nabla}{ds}\left(  \nabla_{V\left(  s\right)  }Z\right)   &
=\frac{\nabla}{ds}\left(  \sum_{i=1}^{d}V_{i}\left(  s\right)  \left(
\nabla_{E_{i}}Z\right)  \left(  \sigma\left(  s\right)  \right)  \right) \\
&  =\sum_{i=1}^{d}V_{i}^{\prime}\left(  s\right)  \left(  \nabla_{E_{i}%
}Z\right)  \left(  \sigma\left(  s\right)  \right)  +\sum_{i=1}^{d}%
V_{i}\left(  s\right)  \nabla_{\sigma^{\prime}\left(  s\right)  }\left(
\nabla_{E_{i}}Z\right)  .
\end{align*}
Using Eq. (\ref{e.3.35}),%
\[
\nabla_{\sigma^{\prime}\left(  s\right)  }\left(  \nabla_{E_{i}}Z\right)
=\nabla_{\sigma^{\prime}\left(  s\right)  \otimes E_{i}\left(  \sigma\left(
s\right)  \right)  }^{2}Z+\left(  \nabla_{\nabla_{\sigma^{\prime}\left(
s\right)  }E_{i}}Z\right)
\]
and using this in the previous equation along with Eq. (\ref{e.3.40}) shows%
\begin{align*}
\frac{\nabla}{ds}\left(  \nabla_{V\left(  s\right)  }Z\right)   &
=\nabla_{\sum_{i=1}^{d}\left\{  V_{i}^{\prime}\left(  s\right)  E_{i}\left(
\sigma\left(  s\right)  \right)  +V_{i}\left(  s\right)  \nabla_{\sigma
^{\prime}\left(  s\right)  }E_{i}\right\}  }Z+\sum_{i=1}^{d}V_{i}\left(
s\right)  \nabla_{\sigma^{\prime}\left(  s\right)  \otimes E_{i}\left(
\sigma\left(  s\right)  \right)  }^{2}Z\\
&  =\left(  \nabla_{\frac{\nabla}{ds}V\left(  s\right)  }Z\right)
+\nabla_{\sigma^{\prime}\left(  s\right)  \otimes V\left(  s\right)  }^{2}Z.
\end{align*}

For the second proof, write $V\left(  s\right)  =\left(  \sigma\left(
s\right)  ,v\left(  s\right)  \right)  =v\left(  s\right)  _{\sigma\left(
s\right)  }$ and $p\left(  s\right)  :=P\left(  \sigma\left(  s\right)
\right)  ,$ then
\begin{align*}
\frac{\nabla}{ds}\left(  \nabla_{V}Z\right)  -\left(  \nabla_{\frac{\nabla
}{ds}V}Z\right)   &  =p\frac{d}{ds}\left(  pz^{\prime}\left(  v\right)
\right)  -pz^{\prime}\left(  pv^{\prime}\right) \\
&  =p\left[  p^{\prime}z^{\prime}\left(  v\right)  +pz^{\prime\prime}\left(
\sigma^{\prime},v\right)  +pz^{\prime}\left(  v^{\prime}\right)  \right]
-pz^{\prime}\left(  pv^{\prime}\right) \\
&  =pp^{\prime}z^{\prime}\left(  v\right)  +pz^{\prime\prime}\left(
\sigma^{\prime},v\right)  +pz^{\prime}\left(  qv^{\prime}\right) \\
&  =p^{\prime}qz^{\prime}\left(  v\right)  +pz^{\prime\prime}\left(
\sigma^{\prime},v\right)  -pz^{\prime}\left(  q^{\prime}v\right) \\
&  =\nabla_{\sigma^{\prime}\left(  s\right)  \otimes V\left(  s\right)  }^{2}Z
\end{align*}
wherein the last equation we have made use of Eq. (\ref{e.3.39}).
\end{proof}

\subsection{Formulas for the Divergence and the Laplacian}

\begin{theorem}
\label{t.3.39} Let $Y$ be a vector field on $M,$ then
\begin{equation}
\operatorname*{div}Y=\operatorname*{tr}(\nabla Y). \label{e.3.41}%
\end{equation}
(Note: $(v_{m}\rightarrow\nabla_{v_{m}}Y)\in\operatorname*{End}(T_{m}M)$ for
each $m\in M,$ so it makes sense to take the trace.) Consequently, if $f$ is a
smooth function on $M,$ then
\begin{equation}
\Delta f=\text{$\operatorname*{tr}$}(\nabla\operatorname*{grad}f).
\label{e.3.42}%
\end{equation}

\end{theorem}

\begin{proof}
Let $x$ be a chart on $M$, $\partial_{i}:=\partial/\partial x^{i},$
$\nabla_{i}:=\nabla_{\partial_{i}}$, and $Y^{i}:=dx^{i}(Y).$ Then by the
product rule and the fact that $\nabla$ is Torsion free (item 2. of the
Proposition \ref{p.3.36}),
\[
\nabla_{i}Y=\sum_{j=1}^{d}\nabla_{i}(Y^{j}\partial_{j})=\sum_{j=1}%
^{d}(\partial_{i}Y^{j}\partial_{j}+Y^{j}\nabla_{i}\partial_{j}),
\]
and $\nabla_{i}\partial_{j}=\nabla_{j}\partial_{i}.$ Hence,
\begin{align*}
\text{$\operatorname*{tr}$}(\nabla Y)  &  =\sum_{i=1}^{d}dx^{i}(\nabla
_{i}Y)=\sum_{i=1}^{d}\partial_{i}Y^{i}+\sum_{i,j=1}^{d}dx^{i}(Y^{j}\nabla
_{i}\partial_{j})\\
&  =\sum_{i=1}^{d}\partial_{i}Y^{i}+\sum_{i,j=1}^{d}dx^{i}(Y^{j}\nabla
_{j}\partial_{i}).
\end{align*}
Therefore, according to Eq. \text{(\ref{e.3.20})}, to finish the proof it
suffices to show that
\[
\sum_{i=1}^{d}dx^{i}(\nabla_{j}\partial_{i})=\partial_{j}\log\sqrt{g}.
\]
From Lemma \ref{l.2.7},%
\[
\partial_{j}\log\sqrt{g}={\frac{1}{2}}\partial_{j}\log(\det g)={\frac{1}{2}%
}\text{$\operatorname*{tr}$}(g^{-1}\partial_{j}g)={\frac{1}{2}}\sum
_{k,l=1}^{d}g^{kl}\partial_{j}g_{kl},
\]
and using Eq. \text{(\ref{e.3.26}) we have}
\[
\partial_{j}g_{kl}=\partial_{j}\langle\partial_{k},\partial_{l}\rangle
=\langle\nabla_{j}\partial_{k},\partial_{l}\rangle+\langle\partial_{k}%
,\nabla_{j}\partial_{l}\rangle.
\]
Combining the last two equations along with the symmetry of $g^{kl}$ implies%
\[
\partial_{j}\log\sqrt{g}=\sum_{k,l=1}^{d}g^{kl}\langle\nabla_{j}\partial
_{k},\partial_{l}\rangle=\sum_{k=1}^{d}dx^{k}(\nabla_{j}\partial_{k}),
\]
where we have used
\[
\sum_{k=1}^{d}g^{kl}\langle\cdot,\partial_{l}\rangle=dx^{k}.
\]
This last equality is easily verified by applying both sides of this equation
to $\partial_{i}$ for $i=1,2,\ldots,n.$
\end{proof}

\begin{definition}
[One forms]\label{d.3.40}A \textbf{one form} $\omega$ on $M$ is a smooth
function $\omega:TM\rightarrow\mathbb{R}$ such that $\omega_{m}:=\omega
|_{T_{m}M}$ is linear for all $m\in M.$ Note: if $x$ is a chart of $M$ with
$m\in\mathcal{D}(x),$ then
\[
\omega_{m}=\sum_{i=1}^{d}\omega_{i}(m)dx^{i}|_{T_{m}M},
\]
where $\omega_{i}:=\omega(\partial/\partial x^{i}).$ The condition that
$\omega$ is smooth is equivalent to the condition that each of the functions
$\omega_{i}$ is smooth on $M.$ Let $\Omega^{1}(M)$ denote the smooth one-forms
on $M.$
\end{definition}

Given a one form, $\omega\in\Omega^{1}(M),$ there is a unique vector field $X$
on $M$ such that $\omega_{m}=\langle X(m),\cdot\rangle_{m}$ for all $m\in M.$
Using this observation, we may extend the definition of $\nabla$ to one forms
by requiring
\begin{equation}
\nabla_{v_{m}}\omega:=\langle\nabla_{v_{m}}X,\cdot\rangle\in T_{m}^{\ast
}M:=(T_{m}M)^{\ast}. \label{e.3.43}%
\end{equation}

\begin{lemma}
[Product Rule]\label{l.3.41} Keep the notation of the above paragraph. Let
$Y\in\Gamma(TM),$ then
\begin{equation}
v_{m}\left[  \omega(Y)\right]  =(\nabla_{v_{m}}\omega)(Y(m))+\omega
(\nabla_{v_{m}}Y). \label{e.3.44}%
\end{equation}
Moreover, if $\theta:M\rightarrow(\mathbb{R}^{N})^{\ast}$ is a smooth function
and
\[
\omega(v_{m}):=\theta(m)v
\]
for all $v_{m}\in TM,$ then
\begin{equation}
(\nabla_{v_{m}}\omega)(w_{m})=d\theta(v_{m})w-\theta(m)dQ(v_{m})w=(d(\theta
P)(v_{m}))w, \label{e.3.45}%
\end{equation}
where $(\theta P)(m):=\theta(m)P(m)\in(\mathbb{R}^{N})^{\ast}.$
\end{lemma}

\begin{proof}
Using the metric compatibility of $\nabla,$
\begin{align*}
v_{m}(\omega(Y))  &  =v_{m}(\langle X,Y\rangle)=\langle\nabla_{v_{m}%
}X,Y(m)\rangle+\langle X(m),\nabla_{v_{m}}Y\rangle\\
&  =(\nabla_{v_{m}}\omega)(Y(m))+\omega(\nabla_{v_{m}}Y).
\end{align*}
Writing $Y(m)=(m,y(m))=y(m)_{m}$ and using Eq. \text{(\ref{e.3.44})}, it
follows that
\begin{align*}
(\nabla_{v_{m}}\omega)(Y(m))  &  =v_{m}(\omega(Y))-\omega(\nabla_{v_{m}}Y)\\
&  =v_{m}(\theta(\cdot)y(\cdot))-\theta(m)(dy(v_{m})+dQ(v_{m})y(m))\\
&  =(d\theta(v_{m}))y(m)-\theta(m)(dQ(v_{m}))y(m).
\end{align*}
Choosing $Y$ such that $Y(m)=w_{m}$ proves the first equality in Eq.
\text{(\ref{e.3.45})}. The second equality in Eq. \text{(\ref{e.3.45})} is a
simple consequence of the formula
\[
d(\theta P)=d\theta(\cdot)P+\theta dP=d\theta(\cdot)P-\theta dQ.
\]

\end{proof}

Before continuing, let us record the following useful corollary of the
previous proof.

\begin{corollary}
\label{c.3.42}To every one -- form $\omega$ on $M,$ there exists $f_{i}%
,g_{i}\in C^{\infty}(M)$ for $i=1,2,\dots,N$ such that%
\begin{equation}
\omega=\sum_{i=1}^{N}f_{i}dg_{i}. \label{e.3.46}%
\end{equation}

\end{corollary}

\begin{proof}
Let $f_{i}\left(  m\right)  :=\theta(m)P(m)e_{i}$ and $g_{i}\left(  m\right)
=x^{i}\left(  m\right)  =\langle m,e_{i}\rangle_{\mathbb{R}^{N}}$ where
$\left\{  e_{i}\right\}  _{i=1}^{N}$ is the standard basis for $\mathbb{R}%
^{N}$ and $P\left(  m\right)  $ is orthogonal projection of $\mathbb{R}^{N}$
onto $\tau_{m}M$ for each $m\in M.$
\end{proof}

\begin{definition}
\label{d.3.43} For $f\in C^{\infty}(M)$ and $v_{m},$ $w_{m}$ in $T_{m}M$, let
\[
\nabla df(v_{m},w_{m}):=(\nabla_{v_{m}}df)(w_{m}),
\]
so that
\[
\nabla df:\cup_{m\in M}(T_{m}M\times T_{m}M)\rightarrow\mathbb{R}.
\]
We call $\nabla df$ the \textbf{Hessian} of $f.$
\end{definition}

\begin{lemma}
\label{l.3.44}Let $f\in C^{\infty}(M),$ $F\in C^{\infty}(\mathbb{R}^{N})$ such
that $f=F|_{M},$ $X,Y\in\Gamma(TM)$ and $v_{m},w_{m}\in T_{m}M.$ Then:

\begin{enumerate}
\item $\nabla df(X,Y)=XYf-df(\nabla_{X}Y).$

\item $\nabla df(v_{m},w_{m})=F^{\prime\prime}(m)(v,w)-F^{\prime}%
(m)dQ(v_{m})w.$

\item $\nabla df(v_{m},w_{m})=\nabla df(w_{m},v_{m})$ -- another manifestation
of zero torsion.
\end{enumerate}
\end{lemma}

\begin{proof}
Using the product rule (see Eq. \text{(\ref{e.3.44})}):
\[
XYf=X(df(Y))=(\nabla_{X}df)(Y)+df(\nabla_{X}Y),
\]
and hence
\[
\nabla df(X,Y)=(\nabla_{X}df)(Y)=XYf-df(\nabla_{X}Y).
\]
This proves item 1. From this last equation and Proposition \ref{p.3.36}
$(\nabla$ has zero torsion), it follows that
\[
\nabla df(X,Y)-\nabla df(Y,X)=[X,Y]f-df(\nabla_{X}Y-\nabla_{Y}X)=0.
\]
This proves the third item upon choosing $X$ and $Y$ such that $X(m)=v_{m}$
and $Y(m)=w_{m}.$ Item 2 follows easily from Lemma \ref{l.3.41} applied with
$\theta:=F^{\prime}.$
\end{proof}

\begin{definition}
\label{d.3.45}Given a point $m\in M,$ a \textbf{local orthonormal frame}
$\left\{  E_{i}\right\}  _{i=1}^{d}$ at $m$ is a collection of local vector
fields defined near $m$ such that $\left\{  E_{i}\left(  p\right)  \right\}
_{i=1}^{d}$ is an orthonormal basis for $T_{p}M$ for all $p$ near $m.$
\end{definition}

\begin{corollary}
\label{c.3.46} Suppose that $F\in C^{\infty}(\mathbb{R}^{N}),$ $f:=F|_{M},$
and $m\in M.$ Let $\{e_{i}\}_{i=1}^{d}$ be an orthonormal basis for $\tau
_{m}M$ and let $\{E_{i}\}_{i=1}^{d}$ be an orthonormal frame near $m\in M.$
Then
\begin{equation}
\Delta f(m)=\sum_{i=1}^{d}\nabla df(E_{i}(m),E_{i}(m)), \label{e.3.47}%
\end{equation}%
\begin{equation}
\Delta f(m)=\sum_{i=1}^{d}\{E_{i}E_{i}f)(m)-df(\nabla_{E_{i}(m)}E_{i})\},
\label{e.3.48}%
\end{equation}
and
\begin{equation}
\Delta f(m)=\sum_{i=1}^{d}F^{\prime\prime}(m)(e_{i},e_{i})-F^{\prime
}(m)(dQ(E_{i}\left(  m\right)  )e_{i}) \label{e.3.49}%
\end{equation}
where $E_{i}\left(  m\right)  :=\left(  m,e_{i}\right)  .$
\end{corollary}

\begin{proof}
By Theorem \ref{t.3.39}, $\Delta f=\sum_{i=1}^{d}\left\langle \nabla_{E_{i}%
}\operatorname*{grad}f,E_{i}\right\rangle $ and by Eq. \text{(\ref{e.3.43})},
$\nabla_{E_{i}}df=\left\langle \nabla_{E_{i}}\operatorname*{grad}%
f,\cdot\right\rangle .$ Therefore
\[
\Delta f=\sum_{i=1}^{d}(\nabla_{E_{i}}df)(E_{i})=\sum_{i=1}^{d}\nabla
df(E_{i},E_{i}),
\]
which proves Eq. \text{(\ref{e.3.47})}. Equations \text{(\ref{e.3.48})} and
\text{(\ref{e.3.49})} follows from Eq. \text{(\ref{e.3.47})} and Lemma
\ref{l.3.44}.
\end{proof}

\begin{notation}
\label{n.3.47}Let $\left\{  e_{i}\right\}  _{i=1}^{N}$ be the standard basis
on $\mathbb{R}^{N}$ and define $X_{i}\left(  m\right)  :=P\left(  m\right)
e_{i}$ for all $m\in M$ and $i=1,2,\dots,N.$
\end{notation}

In the next proposition we will express the gradient, divergence and the
Laplacian in terms of the vector fields, $\left\{  X_{i}\right\}  _{i=1}%
^{N.}.$ These formula will prove very useful when we start discussing Brownian
motion on $M.$

\begin{proposition}
\label{p.3.48}Let $f\in C^{\infty}\left(  M\right)  $ and $Y\in\Gamma\left(
TM\right)  $ then

\begin{enumerate}
\item $v_{m}=\sum_{i=1}^{N}\langle v_{m},X_{i}\left(  m\right)  \rangle
X_{i}\left(  m\right)  $ for all $v_{m}\in T_{m}M.$

\item $\vec{\nabla}f=\operatorname*{grad}f=\sum_{i=1}^{N}X_{i}f\cdot X_{i}$

\item $\vec{\nabla}\cdot Y=\operatorname*{div}(Y)=\sum_{i=1}^{N}\langle
\nabla_{X_{i}}Y,X_{i}\rangle$

\item $\sum_{i=1}^{N}\nabla_{X_{i}}X_{i}=0$

\item $\Delta f=\sum_{i=1}^{N}X_{i}^{2}f.$
\end{enumerate}
\end{proposition}

\begin{proof}
1. The main point is to show%
\begin{equation}
\sum_{i=1}^{N}X_{i}\left(  m\right)  \otimes X_{i}\left(  m\right)
=\sum_{i=1}^{d}u_{i}\otimes u_{i} \label{e.3.50}%
\end{equation}
where $\left\{  u_{i}\right\}  _{i=1}^{d}$ is an orthonormal basis for
$T_{m}M.$ But this is easily proved since
\[
\sum_{i=1}^{N}X_{i}\left(  m\right)  \otimes X_{i}\left(  m\right)
=\sum_{i=1}^{N}P\left(  m\right)  e_{i}\otimes P\left(  m\right)  e_{i}%
\]
and the latter expression is independent of the choice of orthonormal basis
$\left\{  e_{i}\right\}  _{i=1}^{N}$ for $\mathbb{R}^{N}.$ Hence if we choose
$\left\{  e_{i}\right\}  _{i=1}^{N}$ so that $e_{i}=u_{i}$ for $i=1,\dots,d,$
then
\[
\sum_{i=1}^{N}P\left(  m\right)  e_{i}\otimes P\left(  m\right)  e_{i}%
=\sum_{i=1}^{d}u_{i}\otimes u_{i}%
\]
as desired. Since $\sum_{i=1}^{N}\langle v_{m},X_{i}\left(  m\right)  \rangle
X_{i}\left(  m\right)  $ is quadratic in $X_{i},$ it now follows that
\[
\sum_{i=1}^{N}\langle v_{m},X_{i}\left(  m\right)  \rangle X_{i}\left(
m\right)  =\sum_{i=1}^{d}\langle v_{m},u_{i}\rangle u_{i}=v_{m}.
\]

2. This is an immediate consequence of item 1:%
\[
\operatorname*{grad}f\left(  m\right)  =\sum_{i=1}^{N}\langle
\operatorname*{grad}f\left(  m\right)  ,X_{i}\left(  m\right)  \rangle
X_{i}\left(  m\right)  =\sum_{i=1}^{N}X_{i}f\left(  m\right)  \cdot
X_{i}\left(  m\right)  .
\]

3. Again $\sum_{i=1}^{N}\langle\nabla_{X_{i}}Y,X_{i}\rangle\left(  m\right)  $
is quadratic in $X_{i}$ and so by Eq. (\ref{e.3.50}) and Theorem \ref{t.3.39},%
\[
\sum_{i=1}^{N}\langle\nabla_{X_{i}}Y,X_{i}\rangle\left(  m\right)  =\sum
_{i=1}^{d}\langle\nabla_{u_{i}}Y,u_{i}\rangle\left(  m\right)
=\operatorname*{div}(Y).
\]

4. By definition of $X_{i}$ and $\nabla$ and using Lemma \ref{l.3.30},%
\begin{equation}
\sum_{i=1}^{N}\left(  \nabla_{X_{i}}X_{i}\right)  \left(  m\right)
=\sum_{i=1}^{N}P\left(  m\right)  dP\left(  X_{i}\left(  m\right)  \right)
e_{i}=\sum_{i=1}^{N}dP\left(  P\left(  m\right)  e_{i}\right)  Q\left(
m\right)  e_{i}. \label{e.3.51}%
\end{equation}
The latter expression is independent of the choice of orthonormal basis
$\left\{  e_{i}\right\}  _{i=1}^{N}$ for $\mathbb{R}^{N}.$ So again we may
choose $\left\{  e_{i}\right\}  _{i=1}^{N}$ so that $e_{i}=u_{i}$ for
$i=1,\dots,d,$ in which case $P\left(  m\right)  e_{j}=0$ for $j>d$ and so
each summand in the right member of Eq. (\ref{e.3.51}) is zero.

5. To compute $\Delta f,$ use items 2.-- 4., the definition of $\vec{\nabla}f$
and the product rule to find%
\begin{align*}
\Delta f  &  =\vec{\nabla}\cdot(\vec{\nabla}f)=\sum_{i=1}^{N}\langle
\nabla_{X_{i}}\vec{\nabla}f,X_{i}\rangle\\
&  =\sum_{i=1}^{N}X_{i}\langle\vec{\nabla}f,X_{i}\rangle-\sum_{i=1}^{N}%
\langle\vec{\nabla}f,\nabla_{X_{i}}X_{i}\rangle=\sum_{i=1}^{N}X_{i}X_{i}f.
\end{align*}

\end{proof}

The following commutation formulas are at the heart of many of the results to
appear in the latter sections of these note.

\begin{theorem}
[The Bochner-Weitenb\"{o}ck Identity]\label{t.3.49}Let $f\in C^{\infty}\left(
M\right)  $ and $a,b,c\in T_{m}M,$ then%
\begin{equation}
\langle\nabla_{a\otimes b}^{2}\vec{\nabla}f,c\rangle=\langle\nabla_{a\otimes
c}^{2}\vec{\nabla}f,b\rangle\label{e.3.52}%
\end{equation}
and if $S\subset T_{m}M$ is an orthonormal basis, then%
\begin{equation}
\sum_{a\in S}\nabla_{a\otimes a}^{2}\vec{\nabla}f=\left(  \operatorname*{grad}%
\mathrm{\ }\Delta f\right)  \left(  m\right)  +\operatorname*{Ric}\vec{\nabla
}f\left(  m\right)  . \label{e.3.53}%
\end{equation}

\end{theorem}

This result is the first indication that the Ricci tensor is going to play an
important role in later developments. The proof will be given after the next
technical lemma which will be helpful in simplifying the proof of the theorem.

\begin{lemma}
\label{l.3.50}Given $m\in M$ and $v\in T_{m}M$ there exists $V\in\Gamma\left(
TM\right)  $ such that $V\left(  m\right)  =v$ and $\nabla_{w}V=0$ for all
$w\in T_{m}M.$ Moreover if $\left\{  e_{i}\right\}  _{i=1}^{d}$ is an
orthonormal basis for $T_{m}M,$ there exists a local orthonormal frame
$\left\{  E_{i}\right\}  _{i=1}^{d}$ near $m$ such that $\nabla_{w}E_{i}=0$
for all $w\in T_{m}M.$
\end{lemma}

\begin{proof}
In the proof to follow it is assume that $V,$ $Q$ and $P$ have all been
extended off $M$ to smooth function on the ambient space. If $V$ is to exist,
we must have%
\[
0=\nabla_{w}V=V^{\prime}\left(  m\right)  w+\partial_{w}Q\left(  m\right)  v,
\]
i.e.
\[
V^{\prime}\left(  m\right)  w=-\partial_{w}Q\left(  m\right)  v\text{ for all
}w\in T_{m}M.
\]
This helps to motivate defining $V$ by
\[
V\left(  x\right)  :=P\left(  x\right)  \left(  v-\left(  \partial
_{x-m}Q\right)  \left(  m\right)  v\right)  \in T_{x}M\text{ for all }x\in M.
\]
By construction, $V\left(  m\right)  =v$ and making use of the identities in
Lemma \ref{l.3.30},%
\begin{align*}
\nabla_{w}V  &  =\partial_{w}\left[  P\left(  x\right)  \left(  v-\left(
\partial_{x-m}Q\right)  \left(  m\right)  v\right)  \right]  |_{x=m}+\left(
\partial_{w}Q\right)  \left(  m\right)  v\\
&  =\left(  \partial_{w}P\right)  \left(  m\right)  v-P\left(  m\right)
\left(  \partial_{w}Q\right)  \left(  m\right)  v+\left(  \partial
_{w}Q\right)  \left(  m\right)  v\\
&  =\left(  \partial_{w}P\right)  \left(  m\right)  v+Q\left(  m\right)
\left(  \partial_{w}Q\right)  \left(  m\right)  v=\left(  \partial
_{w}P\right)  \left(  m\right)  v+\left(  \partial_{w}Q\right)  \left(
m\right)  v=0
\end{align*}
as desired.

For the second assertion, choose a local frame $\left\{  V_{i}\right\}
_{i=1}^{d}$ such that $V_{i}\left(  m\right)  =e_{i}$ and $\nabla_{w}V_{i}=0$
for all $i$ and $w\in T_{m}M.$ The desired frame $\left\{  E_{i}\right\}
_{i=1}^{d}$ is now constructed by performing Gram-Schmidt orthogonalization on
$\left\{  V_{i}\right\}  _{i=1}^{d}.$ The resulting orthonormal frame,
$\left\{  E_{i}\right\}  _{i=1}^{d},$ still satisfies $\nabla_{w}E_{i}=0$ for
all $w\in T_{m}M.$ For example, $E_{1}=\langle V_{1},V_{1}\rangle^{-1/2}V_{1}$
and since
\[
w\langle V_{1},V_{1}\rangle=2\langle\nabla_{w}V_{1},V_{1}\left(  m\right)
\rangle=0
\]
it follows that
\[
\nabla_{w}E_{1}=w\left(  \langle V_{1},V_{1}\rangle^{-1/2}\right)  \cdot
V_{1}\left(  m\right)  +\langle V_{1},V_{1}\rangle^{-1/2}\left(  m\right)
\nabla_{w}V_{1}\left(  m\right)  =0.
\]
The similar verifications that $\nabla_{w}E_{j}=0$ for $j=2,\dots,d$ will be
left to the reader.
\end{proof}

\begin{proof}
(\emph{Proof of Theorem \ref{t.3.49}.}) Let $a,b,c\in T_{m}M$ and suppose
$A,B,C\in\Gamma\left(  TM\right)  $ have been chosen as in Lemma \ref{l.3.50},
so that $A\left(  m\right)  =a,$ $B\left(  m\right)  =b$ and $C\left(
m\right)  =c$ with $\nabla_{w}A=\nabla_{w}B=\nabla_{w}C=0$ for all $w\in
T_{m}M.$ Then
\begin{align*}
ABCf  &  =AB\langle\vec{\nabla}f,C\rangle=A\langle\nabla_{B}\vec{\nabla
}f,C\rangle+A\langle\vec{\nabla}f,\nabla_{B}C\rangle\\
&  =\langle\nabla_{A}\nabla_{B}\vec{\nabla}f,C\rangle+\langle\nabla_{B}%
\vec{\nabla}f,\nabla_{A}C\rangle+A\langle\vec{\nabla}f,\nabla_{B}C\rangle
\end{align*}
which evaluated at $m$ gives%
\begin{align*}
\left(  ABCf\right)  \left(  m\right)   &  =\left(  \langle\nabla_{A}%
\nabla_{B}\vec{\nabla}f,C\rangle+A\langle\vec{\nabla}f,\nabla_{B}%
C\rangle\right)  \left(  m\right) \\
&  =\langle\nabla_{a\otimes b}^{2}\vec{\nabla}f,c\rangle+\left(  A\langle
\vec{\nabla}f,\nabla_{B}C\rangle\right)  \left(  m\right)
\end{align*}
wherein the last equality we have used $\left(  \nabla_{A}B\right)  \left(
m\right)  =0.$ Interchanging $B$ and $C$ in this equation and subtracting then
implies%
\begin{align*}
\left(  A\left[  B,C\right]  f\right)  \left(  m\right)   &  =\langle
\nabla_{a\otimes b}^{2}\vec{\nabla}f,c\rangle-\langle\nabla_{a\otimes c}%
^{2}\vec{\nabla}f,b\rangle+\left(  A\langle\vec{\nabla}f,\nabla_{B}%
C-\nabla_{C}B\rangle\right)  \left(  m\right) \\
&  =\langle\nabla_{a\otimes b}^{2}\vec{\nabla}f,c\rangle-\langle
\nabla_{a\otimes c}^{2}\vec{\nabla}f,b\rangle+\left(  A\langle\vec{\nabla
}f,[B,C]\rangle\right)  \left(  m\right) \\
&  =\langle\nabla_{a\otimes b}^{2}\vec{\nabla}f,c\rangle-\langle
\nabla_{a\otimes c}^{2}\vec{\nabla}f,b\rangle+\left(  A[B,C]f\right)  \left(
m\right)
\end{align*}
and this equation implies Eq. (\ref{e.3.52}).

Now suppose that $\left\{  E_{i}\right\}  _{i=1}^{d}\subset T_{m}M$ is an
orthonormal frame as in Lemma \ref{l.3.50} and $e_{i}=E_{i}\left(  m\right)
.$ Then, using Proposition \ref{p.3.38},%
\begin{equation}
\sum_{i=1}^{d}\langle\nabla_{e_{i}\otimes e_{i}}^{2}\vec{\nabla}%
f,c\rangle=\sum_{i=1}^{d}\langle\nabla_{e_{i}\otimes c}^{2}\vec{\nabla}%
f,e_{i}\rangle=\sum_{i=1}^{d}\langle\nabla_{c\otimes e_{i}}^{2}\vec{\nabla
}f+R\left(  e_{i},c\right)  \vec{\nabla}f\left(  m\right)  ,e_{i}\rangle.
\label{e.3.54}%
\end{equation}
Since%
\begin{align*}
\sum_{i=1}^{d}\langle\nabla_{c\otimes e_{i}}^{2}\vec{\nabla}f,e_{i}\rangle &
=\sum_{i=1}^{d}\left(  \langle\nabla_{C}\nabla_{E_{i}}\vec{\nabla}%
f,E_{i}\rangle\right)  \left(  m\right)  =\sum_{i=1}^{d}\left(  C\langle
\nabla_{E_{i}}\vec{\nabla}f,E_{i}\rangle\right)  \left(  m\right) \\
&  =\left(  C\Delta f\right)  \left(  m\right)  =\langle\left(  \vec{\nabla
}\Delta f\right)  \left(  m\right)  ,c\rangle
\end{align*}
and (using $R\left(  e_{i},c\right)  ^{\operatorname*{tr}}=R\left(
c,e_{i}\right)  )$%
\begin{align*}
\sum_{i=1}^{d}\langle R\left(  e_{i},c\right)  \vec{\nabla}f\left(  m\right)
,e_{i}\rangle &  =\sum_{i=1}^{d}\langle\vec{\nabla}f\left(  m\right)
,R\left(  c,e_{i}\right)  e_{i}\rangle\\
&  =\langle\vec{\nabla}f\left(  m\right)  ,\operatorname*{Ric}c\rangle
=\langle\operatorname*{Ric}\vec{\nabla}f\left(  m\right)  ,c\rangle,
\end{align*}
Eq. (\ref{e.3.54}) is implies%
\[
\sum_{i=1}^{d}\langle\nabla_{e_{i}\otimes e_{i}}^{2}\vec{\nabla}%
f,c\rangle=\langle\left(  \vec{\nabla}\Delta f\right)  \left(  m\right)
+\operatorname*{Ric}\vec{\nabla}f\left(  m\right)  ,c\rangle
\]
which proves Eq. (\ref{e.3.53}) since $c\in T_{m}M$ was arbitrary.
\end{proof}

\subsection{Parallel Translation\label{s.3.6}}

\begin{definition}
\label{d.3.51}Let $V$ be a smooth path in $TM.$ $V$ is said to
\textbf{parallel} or \textbf{covariantly constant} if $\nabla V(s)/ds\equiv0.$
\end{definition}

\begin{theorem}
\label{t.3.52} Let $\sigma$ be a smooth path in $M$ and $(v_{0})_{\sigma
(0)}\in T_{\sigma(0)}M.$ Then there exists a unique smooth vector field $V$
along $\sigma$ such that $V$ is parallel and $V(0)=(v_{0})_{\sigma(0)}.$
Moreover if $V\left(  s\right)  $ and $W\left(  s\right)  $ are parallel along
$\sigma,$ then $\langle V(s),W(s)\rangle=\langle V\left(  0\right)  ,W\left(
0\right)  \rangle$ for all $s.$
\end{theorem}

\begin{proof}
If $V$ and $W$ are parallel, then%
\[
\frac{d}{ds}\langle V(s),W(s)\rangle=\left\langle \frac{\nabla}{ds}%
V(s),W(s)\right\rangle +\left\langle V(s),\frac{\nabla}{ds}W(s)\right\rangle
=0
\]
which proves the last assertion of the theorem. If a parallel vector field
$V(s)=(\sigma(s),v(s))$ along $\sigma(s)$ is to exist, then
\begin{equation}
dv(s)/ds+dQ(\sigma^{\prime}(s))v(s)=0\quad\text{\textrm{\ and }}\quad
v(0)=v_{0}. \label{e.3.55}%
\end{equation}
By existence and uniqueness of solutions to ordinary differential equations,
there is exactly one solution to Eq. \text{(\ref{e.3.55})}. Hence, if $V$
exists it is unique.

Now let $v$ be the unique solution to Eq. \text{(\ref{e.3.55})} and set
$V(s):=(\sigma(s),v(s)).$ To finish the proof it suffices to show that
$v(s)\in\tau_{\sigma(s)}M.$ Equivalently, we must show that $w(s):=q(s)v(s)$
is identically zero, where $q(s):=Q(\sigma(s)).$ Letting $v^{\prime
}(s)=dv(s)/ds$ and $p(s)=P(\sigma(s)),$ then Eq. (\ref{e.3.55}) states
$v^{\prime}=-q^{\prime}v$ and from Lemma \ref{l.3.30} we have $pq^{\prime
}=q^{\prime}q.$ Thus the function $w$ satisfies%
\[
w^{\prime}=q^{\prime}v+qv^{\prime}=q^{\prime}v-qq^{\prime}v=pq^{\prime
}v=q^{\prime}qv=q^{\prime}w
\]
with $w(0)=0.$ But this linear ordinary differential equation has $w\equiv0$
as its unique solution.
\end{proof}

\begin{definition}
[Parallel Translation]\label{d.3.53}Given a smooth path $\sigma,$ let
$//_{s}(\sigma):T_{\sigma(0)}M\rightarrow T_{\sigma(s)}M$ be defined by
$//_{s}(\sigma)(v_{0})_{\sigma(0)}=V(s),$ where $V$ is the unique parallel
vector field along $\sigma$ such that $V(0)=(v_{0})_{\sigma(0)}.$ We call
$//_{s}(\sigma)$ \textbf{parallel translation} along $\sigma$ up to time $s.$
\end{definition}

\begin{remark}
\label{r.3.54}Notice that $//_{s}(\sigma)v_{\sigma(0)}=(u(s)v)_{\sigma(0)},$
where $s\rightarrow u(s)\in\operatorname*{Hom}(\tau_{\sigma(0)}M,\mathbb{R}%
^{N})$ is the unique solution to the differential equation
\begin{equation}
u^{\prime}(s)+dQ(\sigma^{\prime}(s))u(s)=0\quad\text{\textrm{\ with }}\quad
u(0)=P\left(  \sigma\left(  0\right)  \right)  . \label{e.3.56}%
\end{equation}
Because of Theorem \ref{t.3.52}, $u(s):\tau_{\sigma(0)}M\rightarrow
\mathbb{R}^{N}$ is an isometry for all $s$ and the range of $u(s)$ is
$\tau_{\sigma(s)}M.$ Moreover, if we let $\bar{u}\left(  s\right)  $ denote
the solution to
\begin{equation}
\bar{u}^{\prime}\left(  s\right)  -\bar{u}(s)dQ(\sigma^{\prime}(s))=0\text{
with }\bar{u}\left(  0\right)  =P\left(  \sigma\left(  0\right)  \right)  ,
\label{e.3.57}%
\end{equation}
then
\begin{align*}
\frac{d}{ds}\left[  \bar{u}\left(  s\right)  u\left(  s\right)  \right]   &
=\bar{u}^{\prime}\left(  s\right)  u\left(  s\right)  +\bar{u}\left(
s\right)  u^{\prime}\left(  s\right) \\
&  =\bar{u}(s)dQ(\sigma^{\prime}(s))u\left(  s\right)  -\bar{u}\left(
s\right)  dQ(\sigma^{\prime}(s))u(s)=0.
\end{align*}
Hence $\bar{u}\left(  s\right)  u\left(  s\right)  =P\left(  \sigma\left(
0\right)  \right)  $ for all $s$ and therefore $\bar{u}\left(  s\right)  $ is
the inverse to $u\left(  s\right)  $ thought of as an linear operator from
$\tau_{\sigma(0)}M$ to $\tau_{\sigma(s)}M.$ See also Lemma \ref{l.3.57} below.
\end{remark}

The following techniques for computing covariant derivatives will be useful in
the sequel.

\begin{lemma}
\label{l.3.55}Suppose $Y\in\Gamma\left(  TM\right)  ,$ $\sigma\left(
s\right)  $ is a path in $M,$ $W\left(  s\right)  =\left(  \sigma\left(
s\right)  ,w\left(  s\right)  \right)  $ is a vector field along $\sigma$ and
let $//_{s}=//_{s}\left(  \sigma\right)  $ be parallel translation along
$\sigma.$ Then

\begin{enumerate}
\item $\frac{\nabla}{ds}W\left(  s\right)  =//_{s}\frac{d}{ds}\left[
//_{s}^{-1}W\left(  s\right)  \right]  .$

\item For any $v\in T_{\sigma\left(  0\right)  }M,$%
\begin{equation}
\frac{\nabla}{ds}\nabla_{//_{s}v}Y=\nabla_{\sigma^{\prime}\left(  s\right)
\otimes//_{s}v}^{2}Y. \label{e.3.58}%
\end{equation}
where $\nabla_{\sigma^{\prime}\left(  s\right)  \otimes//_{s}v}^{2}Y$ was
defined in Proposition \ref{p.3.38}.
\end{enumerate}
\end{lemma}

\begin{proof}
Let $\bar{u}$ be as in Eq. (\ref{e.3.57}). From Eq. (\ref{e.3.25}),
\[
\frac{\nabla W(s)}{ds}=\left(  \frac{d}{ds}w(s)+dQ(\sigma^{\prime
}(s)))w(s)\right)  _{\sigma\left(  s\right)  }%
\]
while, using Remark \ref{r.3.54},%
\begin{align*}
\frac{d}{ds}\left[  //_{s}^{-1}W\left(  s\right)  \right]   &  =\left(
\frac{d}{ds}\left[  \bar{u}\left(  s\right)  w\left(  s\right)  \right]
\right)  _{\sigma\left(  s\right)  }\\
&  =\left(  \bar{u}^{\prime}\left(  s\right)  W\left(  s\right)  +\bar
{u}\left(  s\right)  w^{\prime}\left(  s\right)  \right)  _{\sigma\left(
s\right)  }\\
&  =\left(  \bar{u}\left(  s\right)  dQ\left(  \sigma^{\prime}(s)\right)
w\left(  s\right)  +\bar{u}\left(  s\right)  w^{\prime}\left(  s\right)
\right)  _{\sigma\left(  s\right)  }\\
&  =//_{s}^{-1}\frac{\nabla W(s)}{ds}.
\end{align*}
This proves the first item. We will give two proves of the second item, the
first proof being extrinsic while the second will be intrinsic. In each of
these proofs there will be an implied sum on repeated indices.

\textbf{First proof.} Let $\left\{  X_{i}\right\}  _{i=1}^{N}\subset
\Gamma\left(  TM\right)  $ be as in Notation \ref{n.3.47}, then by Proposition
\ref{p.3.48},%
\begin{equation}
//_{s}v=\langle//_{s}v,X_{i}\left(  \sigma\left(  s\right)  \right)  \rangle
X_{i}\left(  \sigma\left(  s\right)  \right)  =\langle v,//_{s}^{-1}%
X_{i}\left(  \sigma\left(  s\right)  \right)  \rangle X_{i}\left(
\sigma\left(  s\right)  \right)  \label{e.3.59}%
\end{equation}
and therefore,%
\begin{align}
\frac{\nabla}{ds}\nabla_{//_{s}v}Y  &  =\frac{\nabla}{ds}\left[  \langle
//_{s}v,X_{i}\left(  \sigma\left(  s\right)  \right)  \rangle\cdot\left(
\nabla_{X_{i}}Y\right)  \left(  \sigma\left(  s\right)  \right)  \right]
\nonumber\\
&  =\langle//_{s}v,X_{i}\left(  \sigma\left(  s\right)  \right)  \rangle
\cdot\nabla_{\sigma^{\prime}\left(  s\right)  }\left(  \nabla_{X_{i}}Y\right)
+\langle//_{s}v,\nabla_{\sigma^{\prime}\left(  s\right)  }X_{i}\rangle
\cdot\left(  \nabla_{X_{i}}Y\right)  \left(  \sigma\left(  s\right)  \right)
. \label{e.3.60}%
\end{align}
Now
\[
\nabla_{\sigma^{\prime}\left(  s\right)  }\left(  \nabla_{X_{i}}Y\right)
=\nabla_{\sigma^{\prime}\left(  s\right)  \otimes X_{i}}^{2}Y+\nabla
_{\sigma^{\prime}\left(  s\right)  X_{i}}Y
\]
and so again using Proposition \ref{p.3.48},
\begin{equation}
\langle//_{s}v,X_{i}\left(  \sigma\left(  s\right)  \right)  \rangle
\cdot\nabla_{\sigma^{\prime}\left(  s\right)  }\left(  \nabla_{X_{i}}Y\right)
=\nabla_{\sigma^{\prime}\left(  s\right)  \otimes//_{s}v}^{2}Y+\langle
//_{s}v,X_{i}\left(  \sigma\left(  s\right)  \right)  \rangle\cdot
\nabla_{\sigma^{\prime}\left(  s\right)  X_{i}}Y. \label{e.3.61}%
\end{equation}
Taking $\nabla/ds$ of Eq. (\ref{e.3.59}) shows%
\[
0=\langle//_{s}v,\nabla_{\sigma^{\prime}\left(  s\right)  }X_{i}\rangle
X_{i}\left(  \sigma\left(  s\right)  \right)  +\langle//_{s}v,X_{i}\left(
\sigma\left(  s\right)  \right)  \rangle\nabla_{\sigma^{\prime}\left(
s\right)  }X_{i}.
\]
and so%
\begin{equation}
\langle//_{s}v,X_{i}\left(  \sigma\left(  s\right)  \right)  \rangle
\cdot\nabla_{\sigma^{\prime}\left(  s\right)  X_{i}}Y=-\langle//_{s}%
v,\nabla_{\sigma^{\prime}\left(  s\right)  }X_{i}\rangle\cdot\left(
\nabla_{X_{i}}Y\right)  \left(  \sigma\right)  \left(  s\right)  .
\label{e.3.62}%
\end{equation}
Assembling Eqs. (\ref{e.3.59}), (\ref{e.3.61}) and (\ref{e.3.62}) proves Eq.
(\ref{e.3.58}).

\textbf{Second proof. } Let $\left\{  E_{i}\right\}  _{i=1}^{d}$ be an
orthonormal frame near $\sigma\left(  s\right)  ,$ then%
\begin{align}
\frac{\nabla}{ds}\nabla_{//_{s}v}Y  &  =\frac{\nabla}{ds}\left[  \langle
//_{s}v,E_{i}\left(  \sigma\left(  s\right)  \right)  \rangle\cdot\left(
\nabla_{E_{i}}Y\right)  \left(  \sigma\left(  s\right)  \right)  \right]
\nonumber\\
&  =\langle//_{s}v,\nabla_{\sigma^{\prime}\left(  s\right)  }E_{i}\rangle
\cdot\left(  \nabla_{E_{i}}Y\right)  \left(  \sigma\left(  s\right)  \right)
+\langle//_{s}v,E_{i}\left(  \sigma\left(  s\right)  \right)  \rangle
\cdot\nabla_{\sigma^{\prime}\left(  s\right)  }\nabla_{E_{i}}Y. \label{e.3.63}%
\end{align}
Working as in the first proof,%
\begin{align*}
\langle//_{s}v,E_{i}\left(  \sigma\left(  s\right)  \right)  \rangle
\cdot\nabla_{\sigma^{\prime}\left(  s\right)  }\nabla_{E_{i}}Y  &
=\langle//_{s}v,E_{i}\left(  \sigma\left(  s\right)  \right)  \rangle
\cdot\left(  \nabla_{\sigma^{\prime}\left(  s\right)  \otimes E_{i}}%
^{2}Y+\nabla_{\nabla_{\sigma^{\prime}\left(  s\right)  }E_{i}}Y\right) \\
&  =\nabla_{\sigma^{\prime}\left(  s\right)  \otimes//_{s}v}^{2}%
Y+\nabla_{\langle//_{s}v,E_{i}\left(  \sigma\left(  s\right)  \right)
\rangle_{\nabla_{\sigma^{\prime}\left(  s\right)  }E_{i}}}Y
\end{align*}
and using%
\[
0=\frac{\nabla}{ds}//_{s}v=\langle//_{s}v,\nabla_{\sigma^{\prime}\left(
s\right)  }E_{i}\rangle\cdot E_{i}\left(  \sigma\left(  s\right)  \right)
+\langle//_{s}v,E_{i}\left(  \sigma\left(  s\right)  \right)  \rangle
\cdot\nabla_{\sigma^{\prime}\left(  s\right)  }E_{i}%
\]
we learn%
\[
\langle//_{s}v,E_{i}\left(  \sigma\left(  s\right)  \right)  \rangle
\cdot\nabla_{\sigma^{\prime}\left(  s\right)  }\nabla_{E_{i}}Y=\nabla
_{\sigma^{\prime}\left(  s\right)  \otimes//_{s}v}^{2}Y-\langle//_{s}%
v,\nabla_{\sigma^{\prime}\left(  s\right)  }E_{i}\rangle\cdot\left(
\nabla_{E_{i}}Y\right)  \left(  \sigma\left(  s\right)  \right)  .
\]
This equation combined with Eq. (\ref{e.3.63}) again proves Eq. (\ref{e.3.58}).
\end{proof}

The remainder of this section discusses a covariant derivative on
$M\times\mathbb{R}^{N}$ which \textquotedblleft extends\textquotedblright%
\ $\nabla$ defined above. This will be needed in Section \ref{s.5}, where it
will be convenient to have a covariant derivative on the normal bundle:%
\[
N(M):=\cup_{m\in M}(\{m\}\times\tau_{m}M^{\perp})\subset M\times\mathbb{R}%
^{N}.
\]

Analogous to the definition of $\nabla$ on $TM,$ it is reasonable to extend
$\nabla$ to the normal bundle $N(M)$ by setting
\[
\frac{\nabla V(s)}{ds}=(\sigma(s),Q(\sigma(s))v^{\prime}(s))=(\sigma
(s),v^{\prime}(s)+dP(\sigma^{\prime}(s))v(s)),
\]
for all smooth paths $s\rightarrow V(s)=(\sigma(s),v(s))$ in $N(M).$ Then this
covariant derivative on the normal bundle satisfies analogous properties to
$\nabla$ on the tangent bundle $TM.$ The covariant derivatives on $TM$ and
$N\left(  M\right)  $ can be put together to make a covariant derivative on
$M\times\mathbb{R}^{N}.$ Explicitly, if $V(s)=(\sigma(s),v(s))$ is a smooth
path in $M\times\mathbb{R}^{N},$ let $p(s):=P(\sigma(s)),$ $q(s):=Q(\sigma
(s))$ and then define%
\[
\frac{\nabla V(s)}{ds}:=(\sigma(s),p(s)\frac{d}{ds}\{p(s)v(s)\}+q(s)\frac
{d}{ds}\{q(s)v(s)\}).
\]
Since
\begin{align*}
\frac{\nabla V(s)}{ds}  &  =(\sigma(s),\frac{d}{ds}\{p(s)v(s)\}+q^{\prime
}(s)p(s)v(s)\\
&  \qquad+\frac{d}{ds}\{q(s)v(s)\}+p^{\prime}(s)q(s)v(s))\\
&  =(\sigma(s),v^{\prime}(s)+q^{\prime}(s)p(s)v(s)+p^{\prime}(s)q(s)v(s))\\
&  =(\sigma(s),v^{\prime}(s)+dQ(\sigma^{\prime}(s))P(\sigma(s))v(s)+dP(\sigma
^{\prime}(s))Q(\sigma(s))v(s))
\end{align*}
we may write $\nabla V(s)/ds$ as%
\begin{equation}
\frac{\nabla V(s)}{ds}=(\sigma(s),v^{\prime}(s)+\Gamma(\sigma^{\prime
}(s))v(s)) \label{e.3.64}%
\end{equation}
where
\begin{equation}
\Gamma(w_{m})v:=dQ(w_{m})P(m)v+dP(w_{m})Q(m)v \label{e.3.65}%
\end{equation}
for all $w_{m}\in TM$ and $v\in\mathbb{R}^{N}.$

It should be clear from the above computation that the covariant derivative
defined in \text{(\ref{e.3.64})} agrees with those already defined on $TM$ and
$N(M).$ Many of the properties of the covariant derivative on $TM$ follow
quite naturally from this fact and Eq. \text{(\ref{e.3.64})}.

\begin{lemma}
\label{l.3.56}For each $w_{m}\in TM,$ $\Gamma(w_{m})$ is a skew symmetric
$N\times N$ -- matrix. Hence, if $u(s)$ is the solution to the differential
equation
\begin{equation}
u^{\prime}(s)+\Gamma(\sigma^{\prime}(s))u(s)=0\quad\text{\textrm{\ with }%
}\quad u(0)=I, \label{e.3.66}%
\end{equation}
then $u$ is an orthogonal matrix for all $s.$
\end{lemma}

\begin{proof}
Since $\Gamma=dQP+dPQ$ and $P$ and $Q$ are orthogonal projections and hence
symmetric, the adjoint $\Gamma^{\mathrm{tr}}$ of $\Gamma$ is given by
\[
\Gamma^{\mathrm{tr}}=PdQ+QdP=-dPQ-dQP=-\Gamma.
\]
where Lemma \ref{l.3.30} was used in the second equality. Hence $\Gamma$ is a
skew-symmetric valued one form. Now let $u$ denote the solution to
\text{(\ref{e.3.66})}\ and $A(s):=\Gamma(\sigma^{\prime}(s)).$ Then
\[
\frac{d}{ds}u^{\mathrm{tr}}u=(-Au)^{\mathrm{tr}}u+u^{\mathrm{tr}%
}(-Au)=u^{\mathrm{tr}}(A-A)u=0,
\]
which shows that $u^{\mathrm{tr}}(s)u(s)=u^{\mathrm{tr}}(0)u(0)=I.$
\end{proof}

\begin{lemma}
\label{l.3.57}Let $u$ be the solution to \text{(\ref{e.3.66})}. Then
\begin{equation}
u(s)(\tau_{\sigma(0)}M)=\tau_{\sigma(s)}M \label{e.3.67}%
\end{equation}
and
\begin{equation}
u(s)(\tau_{\sigma(0)}M)^{\perp}=\tau_{\sigma(s)}M^{\perp}. \label{e.3.68}%
\end{equation}
In particular, if $v\in\tau_{\sigma(0)}M$ $(v\in\tau_{\sigma(0)}M^{\perp})$
then $V(s):=(\sigma(s),u(s)v)$ is the parallel vector field along $\sigma$ in
$TM$ $(N(M))$ such that $V(0)=v_{\sigma(0)}.$
\end{lemma}

\begin{proof}
By the product rule,%
\begin{equation}
\frac{d}{ds}\{u^{\mathrm{tr}}P\left(  \sigma\right)  u\}=u^{\mathrm{tr}%
}\{\Gamma\left(  \sigma^{\prime}\right)  P\left(  \sigma\right)  +dP\left(
\sigma^{\prime}\right)  -P\left(  \sigma\right)  \Gamma\left(  \sigma^{\prime
}\right)  \}u. \label{e.3.69}%
\end{equation}
Moreover, making use of Lemma \ref{l.3.30},
\begin{align*}
\Gamma\left(  \sigma^{\prime}\right)  P\left(  \sigma\right)   &  -P\left(
\sigma\right)  \Gamma\left(  \sigma^{\prime}\right)  +dP\left(  \sigma
^{\prime}\right) \\
&  =dP\left(  \sigma^{\prime}\right)  +\left[  dQ(\sigma^{\prime}%
)P(\sigma)+dP(\sigma^{\prime})Q(\sigma)\right]  P\left(  \sigma\right) \\
&  \qquad-P\left(  \sigma\right)  \left[  dQ(\sigma^{\prime})P(\sigma
)+dP(\sigma^{\prime})Q(\sigma)\right] \\
&  =dP\left(  \sigma^{\prime}\right)  +dQ(\sigma^{\prime})P(\sigma
)-dP(\sigma^{\prime})Q(\sigma)\\
&  =dP\left(  \sigma^{\prime}\right)  +dQ(\sigma^{\prime})=0,
\end{align*}
which combined with Eq. (\ref{e.3.69}) shows $\frac{d}{ds}\{u^{\mathrm{tr}%
}P\left(  \sigma\right)  u\}=0.$ Therefore,
\[
u^{\mathrm{tr}}(s)P\left(  \sigma\left(  s\right)  \right)  u(s)=P(\sigma
\left(  0\right)  )
\]
for all $s.$ Combining this with Lemma \ref{l.3.56}, shows%
\[
P\left(  \sigma\left(  s\right)  \right)  u(s)=u(s)P(\sigma\left(  0\right)
).
\]
This last equation is equivalent to Eq. \text{(\ref{e.3.67})}. Eq.
\text{(\ref{e.3.68})} has completely analogous proof or can be seen easily
from the fact that $P+Q=I.$
\end{proof}

\subsection{More References}

I recommend \cite{GHL}\ and \cite{DoC} for more details on Riemannian
geometry. The references,
\cite{Ab2,BC,Davies90,DoC,GHL,Hi,Kl1,Kl2,Kl3,KN1,KN2,O'} and the complete five
volume set of Spivak's books on differential geometry starting with
\cite{Spivak1} are also very useful.

\section{Flows and Cartan's Development Map\label{s.4}}

The results of this section will serve as a warm-up for their stochastic
counter parts. These types of theorems will be crucial for the path space
analysis results to be developed in Sections \ref{s.7} and \ref{s.8} below.

\subsection{Time - Dependent Smooth Flows\label{s.4.1}}

\begin{figure}
[ptbh]
\begin{center}
\includegraphics[
height=1.702in,
width=4.1393in
]%
{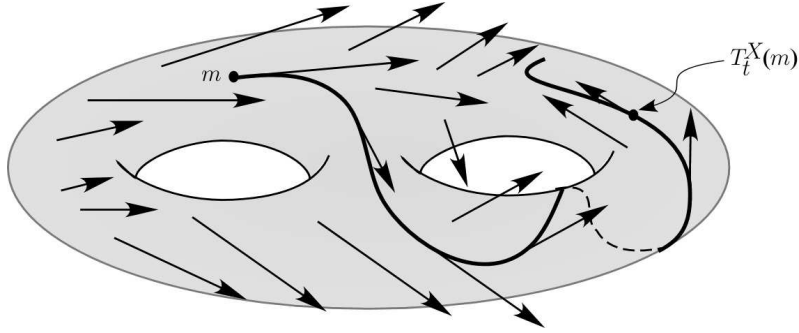}%
\caption{Going with the flow. Here we suppose that $X$ is a time independent
vector field which is indicated by the arrows in the picture and the curve is
the corresponding flow line starting at $m\in M.$}%
\label{fig.10}%
\end{center}
\end{figure}

\begin{notation}
\label{n.4.1}Given a smooth \textbf{time dependent vector} field,
$(t,m)\rightarrow X_{t}\left(  m\right)  \in T_{m}M$ on a manifold $M,$ let
$T_{t}^{X}(m)$ denote the solution to the ordinary differential equation,%
\[
\frac{d}{dt}T_{t}^{X}(m)=X_{t}\circ T_{t}^{X}(m)\text{ with }T_{0}^{X}(m)=m.
\]
If $X$ is \textbf{time independent }we will write $e^{tX}(m)$ for $T_{t}%
^{X}(m).$ We call $T^{X}$ the \textbf{flow} of $X.$ See Figure \ref{fig.10}.
\end{notation}

%

%TCIMACRO{\FRAME{ftbphFU}{4.1393in}{1.702in}{0pt}{\Qcb{Going with the flow.
%Here we suppose that $X$ is a time independent vector field which is indicated
%by the arrows in the picture and the curve is the corresponding flow line
%starting at $m\in M.$}}{\Qlb{fig.10}}{fdflow.eps}%
%width 4.1393in;  height 1.702in;  depth 0pt;  original-width 5.5786in;
%original-height 2.2755in;  cropleft "0";  croptop "1";  cropright "1";
%filename 'fdflow.eps';file-properties "XNPEU";}}}%
%BeginExpansion
%EndExpansion

\begin{theorem}
[Flow Theorem]\label{t.4.2}Suppose that $X_{t}$ is a smooth time dependent
vector field on $M.$ Then for each $m\in M,$ there exists a maximal open
interval $J_{m}\subset\mathbb{R}$ such that $0\in J_{m}$ and $t\rightarrow
T_{t}^{X}(m)$ exists for $t\in J_{m}.$ Moreover the set $\mathcal{D}\left(
X\right)  :=\cup_{m}\left(  J_{m}\times\left\{  m\right\}  \right)
\subset\mathbb{R}\times M$ is open and the map $(t,m)\in\mathcal{D}\left(
X\right)  \rightarrow T_{t}^{X}(m)\in M$ is a smooth map.
\end{theorem}

\begin{proof}
Let $Y_{t}$ be a smooth extension of $X_{t}$ to a vector field on $E$ where
$E$ is the Euclidean space in which $M$ is imbedded. The stated results with
$X$ replaced by $Y$ follows from the standard theory of ordinary differential
equations on Euclidean spaces. Let $T_{t}^{Y}$ denote the flow of $Y$ on $E.$
We will construct $T^{X}$ by setting $T_{t}^{X}\left(  m\right)  :=T_{t}%
^{Y}\left(  m\right)  $ for all $m\in M$ and $t\in J_{m}.$ In order for this
to work we must show that $T_{t}^{Y}\left(  m\right)  \in M$ whenever $m\in
M.$

To verify this last assertion, let $x$ be a chart on $M$ such that
$m\in\mathcal{D}\left(  x\right)  ,$ then $\sigma\left(  t\right)  $ solves
$\dot{\sigma}\left(  t\right)  =X_{t}\left(  \sigma\left(  t\right)  \right)
$ with $\sigma\left(  0\right)  =m$ iff
\[
\frac{d}{dt}\left[  x\circ\sigma\left(  t\right)  \right]  =dx\left(
\dot{\sigma}\left(  t\right)  \right)  =dx\left(  X_{t}\left(  \sigma\left(
t\right)  \right)  \right)  =dx\left(  X_{t}\circ x^{-1}\left(  x\circ
\sigma\left(  t\right)  \right)  \right)
\]
with $x\circ\sigma\left(  0\right)  =m.$ Since this is a differential equation
for $x\circ\sigma\left(  t\right)  \in\mathcal{R}\left(  z\right)  $ and
$\mathcal{R}\left(  z\right)  $ is an open subset $\mathbb{R}^{d},$ the
standard local existence theorem for ordinary differential equations implies
$x\circ\sigma\left(  t\right)  $ exists for small time. This then implies
$\sigma\left(  t\right)  \in M$ exists for small $t$ and satisfies
\[
\dot{\sigma}\left(  t\right)  =X_{t}\left(  \sigma\left(  t\right)  \right)
=Y_{t}\left(  \sigma\left(  t\right)  \right)  \text{ with }\sigma\left(
0\right)  =m.
\]
By uniqueness of solutions to ordinary differential equations, we must have
$T_{t}^{Y}\left(  m\right)  =\sigma\left(  t\right)  $ for small $t$ and in
particular $T_{t}^{Y}\left(  m\right)  \in M$ for small $t.$ Let
\[
\tau:=\sup\left\{  t\in J_{m}:T_{s}^{Y}\left(  m\right)  \in M\text{ for
}0\leq s\leq t\right\}
\]
and for sake of contradiction suppose that $[0,\tau]\subset J_{m}.$ Then by
continuity, $T_{\tau}^{Y}\left(  m\right)  \in M$ and by repeating the above
argument using a chart $x$ on $M$ centered at $T_{\tau}^{Y}\left(  m\right)
,$ we would find that $T_{t}^{Y}(m)\in M$ for $t$ in a neighborhood of $\tau.$
This contradicts the definition of $\tau$ and hence we may conclude that
$\tau$ is the right end point of $J_{m}.$ A similar argument works for $t\in
J_{m}$ with $t<0$ and hence $T_{t}^{Y}\left(  m\right)  \in M$ for all $t\in
J_{m}.$
\end{proof}

\begin{assumption}
[Completeness]\label{a.1}For simplicity in these notes it will always be
assumed that $X$ is \textbf{complete}, i.e. $J_{m}=\mathbb{R}$ for all $m\in
M$ and hence $\mathcal{D}\left(  X\right)  =\mathbb{R}\times M.$ This will be
the case if, for example, $M$ is compact or $M$ is imbedded in $\mathbb{R}%
^{N}$ and the vector field $X$ satisfies a Lipschitz condition. (Later we will
restrict to the compact case.)
\end{assumption}

\begin{notation}
\label{n.4.3}For $g,h\in\mathrm{Diff}(M)$ let $Ad_{g}h:=g\circ h\circ g^{-1}.$
We will also write $Ad_{g}$ for the linear transformation on $\Gamma\left(
TM\right)  $ defined by
\[
Ad_{g}Y=\frac{d}{ds}|_{0}Ad_{g}e^{sY}=\frac{d}{ds}|_{0}g\circ e^{sY}\circ
g^{-1}=g_{\ast}\left(  Y\circ g^{-1}\right)
\]
for all $Y\in\Gamma\left(  TM\right)  .$ (The vector space $\Gamma\left(
TM\right)  $ should be interpreted as the Lie algebra of the diffeomorphism
group, $\mathrm{Diff}(M).)$
\end{notation}

In order to verify $T_{t}^{X}$ is invertible, let $T_{t,s}^{X}$ denote the
solution to%
\[
\frac{d}{dt}T_{t,s}^{X}=X_{t}\circ T_{t,s}^{X}\text{ with }T_{s,s}^{X}=id.
\]

\begin{lemma}
\label{l.4.4}Suppose that $X_{t}$ is a complete time dependent vector field on
$M,$ then $T_{t}^{X}\in\mathrm{Diff}(M)$ for all $t$ and
\begin{equation}
\left(  T_{t}^{X}\right)  ^{-1}=T_{0,t}^{X}=T_{t}^{-Ad_{\left(  T^{X}\right)
^{-1}}X}, \label{e.4.1}%
\end{equation}
where
\[
\left(  Ad_{\left(  T^{X}\right)  ^{-1}}X\right)  _{t}:=Ad_{\left(  T_{t}%
^{X}\right)  ^{-1}}X_{t}.
\]

\end{lemma}

\begin{proof}
If $s,t,u\in\mathbb{R},$ then $S_{t}:=T_{t,s}^{X}\circ T_{s,u}^{X}$ solves%
\[
\dot{S}_{t}=X_{t}\circ S_{t}\text{ with }S_{s}=T_{s,u}^{X}%
\]
which is the same equation that $t\rightarrow T_{t,u}^{X}$ solves and
therefore $T_{t,u}^{X}=T_{t,s}^{X}\circ T_{s,u}^{X}.$ In particular,
$T_{0,t}^{X}$ is the inverse to $T_{t}^{X}.$ Moreover if we let $T_{t}%
:=T_{t}^{X}$ and $S_{t}:=T_{t}^{-1}$ then%
\[
0=\frac{d}{dt}id=\frac{d}{dt}\left[  T_{t}\circ S_{t}\right]  =X_{t}\circ
T_{t}\circ S_{t}+T_{t\ast}\dot{S}_{t}.
\]
So it follows that $S_{t}$ solves%
\[
\dot{S}_{t}=-T_{t\ast}^{-1}X_{t}\circ T_{t}\circ S_{t}=-\left(  Ad_{T_{t}%
^{-1}}X_{t}\right)  \circ S_{t}%
\]
which proves the second equality in Eq. (\ref{e.4.1}).
\end{proof}

\subsection{Differentials of $T_{t}^{X}$\label{s.4.2}}

In the later sections of this article, we will make heavy use of the
stochastic analogues of the following two differentiation theorems.

\begin{theorem}
[Differentiating $m\rightarrow T_{t}^{X}\left(  m\right)  $]\label{t.4.5}%
Suppose $\nabla$ is the Levi-Civita\footnote{Actually, for those in the know,
any torsion zero covariant derivative could be used here.} covariant
derivative on $TM$ and $T_{t}=T_{t}^{X}$ as above, then%
\begin{equation}
\frac{\nabla}{dt}T_{t\ast}v=\nabla_{T_{t\ast}v}X_{t}\text{ for all }v\in TM.
\label{e.4.2}%
\end{equation}
If we further let $m\in M,$ $//_{t}=//_{t}\left(  \tau\rightarrow T_{\tau
}\left(  m\right)  \right)  $ be parallel translation relative to $\nabla$
along the flow line $\tau\rightarrow T_{\tau}\left(  m\right)  $ and
$z_{t}:=//_{t}^{-1}T_{t\ast m},$ then%
\begin{equation}
\frac{d}{dt}z_{t}v=//_{t}^{-1}\nabla_{//_{t}z_{t}v}X_{t}\text{ for all }v\in
T_{m}M. \label{e.4.3}%
\end{equation}
(This is a linear differential equation for $z_{t}\in\operatorname*{End}%
\left(  T_{m}M\right)  .)$
\end{theorem}

\begin{proof}
Let $\sigma\left(  s\right)  $ be smooth path in $M$ such that $\sigma
^{\prime}\left(  0\right)  =v,$ then%
\begin{align*}
\frac{\nabla}{dt}T_{t\ast}v  &  =\frac{\nabla}{dt}\frac{d}{ds}|_{0}%
T_{t}\left(  \sigma\left(  s\right)  \right)  =\frac{\nabla}{ds}|_{0}\frac
{d}{dt}T_{t}\left(  \sigma\left(  s\right)  \right) \\
&  =\frac{\nabla}{ds}|_{0}X_{t}\left(  T_{t}\left(  \sigma\left(  s\right)
\right)  \right)  =\nabla_{T_{t\ast}v}X_{t}%
\end{align*}
wherein the second equality we have used $\nabla$ has zero torsion. Eq.
(\ref{e.4.3}) follows directly from Eq. (\ref{e.4.2}) using $\frac{\nabla}%
{dt}=//_{t}\frac{d}{dt}//_{t}^{-1},$ see Lemma \ref{l.3.55}.
\end{proof}

\begin{remark}
\label{r.4.6}As a warm up for writing the stochastic version of Eq.
(\ref{e.4.3}) in It\^{o} form let us pause to compute $\frac{\nabla}%
{dt}\left(  \nabla_{T_{t\ast}v}Y\right)  $ for $Y\in\Gamma\left(  TM\right)
.$ Using Eqs. (\ref{e.3.38}), (\ref{e.3.37}) and (\ref{e.3.35}) of Proposition
\ref{p.3.38},%
\begin{align}
\frac{\nabla}{dt}\nabla_{T_{t\ast}v}Y  &  =\nabla_{\dot{T}_{t}(m)\otimes
T_{t\ast}v}^{2}Y+\nabla_{\frac{\nabla}{dt}T_{t\ast}v}Y=\nabla_{X_{t}\left(
T_{t}\left(  m\right)  \right)  \otimes T_{t\ast}v}^{2}Y+\nabla_{\nabla
_{T_{t\ast}v}X_{t}}Y\nonumber\\
&  =\nabla_{T_{t\ast}v\otimes X_{t}\left(  T_{t}\left(  m\right)  \right)
}^{2}Y+R^{\nabla}\left(  X_{t}\left(  T_{t}\left(  m\right)  \right)
,T_{t\ast}v\right)  Y\left(  T_{t}\left(  m\right)  \right)  +\nabla
_{\nabla_{T_{t\ast}v}X_{t}}Y\nonumber\\
&  =R^{\nabla}\left(  X_{t}\left(  T_{t}\left(  m\right)  \right)  ,T_{t\ast
}v\right)  Y\left(  T_{t}\left(  m\right)  \right)  +\nabla_{T_{t\ast}%
v}\left(  \nabla_{X_{t}}Y\right)  . \label{e.4.4}%
\end{align}

\end{remark}

\begin{theorem}
[Differentiating $T_{t}^{X}$ in $X$]\label{t.4.7}Suppose $(t,m)\rightarrow
X_{t}\left(  m\right)  $ and $(t,m)\rightarrow Y_{t}\left(  m\right)  $ are
smooth time dependent vector fields on $M$ and let
\begin{equation}
\partial_{Y}T_{t}^{X}:=\frac{d}{ds}|_{0}T_{t}^{X+sY}. \label{e.4.5}%
\end{equation}
Then%
\begin{equation}
\partial_{Y}T_{t}^{X}=T_{t\ast}^{X}\int_{0}^{t}\left(  T_{\tau\ast}%
^{X}\right)  ^{-1}Y_{\tau}\circ T_{\tau}^{X}d\tau=T_{t\ast}^{X}\int_{0}%
^{t}Ad_{T_{\tau}^{X}}^{-1}Y_{\tau}d\tau. \label{e.4.6}%
\end{equation}
This formula may also be written as%
\begin{equation}
\partial_{Y}T_{t}^{X}=\left(  \int_{0}^{t}Ad_{T_{t,\tau}^{X}}Y_{\tau}%
d\tau\right)  \circ T_{t}^{X}=\left(  \int_{0}^{t}Ad_{T_{t}^{X}\circ\left(
T_{\tau}^{X}\right)  ^{-1}}Y_{\tau}d\tau\right)  \circ T_{t}^{X}.
\label{e.4.7}%
\end{equation}

\end{theorem}

\begin{proof}
To simplify notation, let $T_{t}:=T_{t}^{X}$ and define $V_{t}:=\left(
T_{t\ast}^{X}\right)  ^{-1}\partial_{Y}T_{t}^{X}.$ Then $V_{0}=0$ and
$\partial_{Y}T_{t}^{X}=T_{t\ast}^{X}V_{t}$ or equivalently, for all $f\in
C^{\infty}(M),$%
\[
\frac{d}{ds}|_{0}f\circ T_{t}^{X+sY}=\left(  T_{t\ast}^{X}V_{t}\right)
f=V_{t}\left(  f\circ T_{t}^{X}\right)  .
\]
Given $f\in C^{\infty}(M),$ on one hand we have%
\begin{align*}
\frac{d}{dt}\frac{d}{ds}|_{0}f\circ T_{t}^{X+sY}  &  =\frac{d}{dt}\left[
V_{t}(f\circ T_{t}^{X})\right]  =\dot{V}_{t}(f\circ T_{t}^{X})+V_{t}%
(X_{t}f\circ T_{t}^{X})\\
&  =\left(  T_{t\ast}^{X}\dot{V}_{t}\right)  f+V_{t}(X_{t}f\circ T_{t}^{X})
\end{align*}
while on the other hand%
\begin{align*}
\frac{d}{ds}|_{0}\frac{d}{dt}f\circ T_{t}^{X+sY}  &  =\frac{d}{ds}|_{0}\left[
\left(  \left(  X_{t}+sY_{t}\right)  f\right)  \circ T_{t}^{X+sY}\right]
=\left(  Y_{t}f\right)  \circ T_{t}^{X}+V_{t}\left(  X_{t}f\circ T_{t}%
^{X}\right) \\
&  =\left(  Y_{t}\circ T_{t}^{X}\right)  f+V_{t}\left(  X_{t}f\circ T_{t}%
^{X}\right)  .
\end{align*}
Since $\left[  \frac{d}{dt},\frac{d}{ds}|_{0}\right]  =0,$ the previous two
displayed equations imply $\left(  T_{t\ast}^{X}\dot{V}_{t}\right)  f=\left(
Y_{t}\circ T_{t}^{X}\right)  f$ and because this holds for all $f\in
C^{\infty}(M),$
\begin{equation}
T_{t\ast}^{X}\dot{V}_{t}=Y_{t}\circ T_{t}^{X}. \label{e.4.8}%
\end{equation}
Solving Eq. (\ref{e.4.8}) for $\dot{V}_{t}$ and then integrating on $t$ shows%
\[
V_{t}=\int_{0}^{t}\left(  T_{\tau\ast}^{X}\right)  ^{-1}Y_{\tau}\circ T_{\tau
}^{X}d\tau.
\]
which along with the relation, $\partial_{Y}T_{t}^{X}=T_{t\ast}^{X}V_{t},$
implies Eq. (\ref{e.4.6}).

We may now rewrite the formula in Eq. (\ref{e.4.6}) as%
\begin{align*}
\partial_{Y}T_{t}^{X}  &  =T_{t\ast}^{X}\left(  \int_{0}^{t}Ad_{T_{\tau}^{X}%
}^{-1}Y_{\tau}d\tau\right)  \circ\left(  T_{t}^{X}\right)  ^{-1}\circ
T_{t}^{X}=Ad_{T_{t}^{X}}\left(  \int_{0}^{t}Ad_{T_{\tau}^{X}}^{-1}Y_{\tau
}d\tau\right)  \circ T_{t}^{X}\\
&  =\left(  \int_{0}^{t}Ad_{T_{t}^{X}}Ad_{T_{\tau}^{X}}^{-1}Y_{\tau}%
d\tau\right)  \circ T_{t}^{X}=\left(  \int_{0}^{t}Ad_{T_{t}^{X}\circ\left(
T_{\tau}^{X}\right)  ^{-1}}Y_{\tau}d\tau\right)  \circ T_{t}^{X}\\
&  =\left(  \int_{0}^{t}Ad_{T_{t,\tau}^{X}}Y_{\tau}d\tau\right)  \circ
T_{t}^{X}%
\end{align*}
which gives Eq. (\ref{e.4.7}).
\end{proof}

\begin{example}
\label{ex.4.8}Suppose that $G$ is a Lie group, $\mathfrak{g}:=\mathrm{Lie}%
\left(  G\right)  ,$ $A_{t}$ and $B_{t}$ are two smooth $\mathfrak{g}$ --
valued functions and $g_{t}^{A}\in G$ solves the equation%
\[
\frac{d}{dt}g_{t}^{A}=\tilde{A}_{t}\left(  g_{t}^{A}\right)  \text{ with
}g_{0}^{A}=e\in G
\]
where $\tilde{A}_{t}\left(  x\right)  :=L_{x\ast}A_{t}$ is the \textbf{left
invariant }vector field on $G$ associated to $A_{t}\in\mathfrak{g}%
,\mathfrak{\ }$see Examples \ref{ex.2.34} and \ref{ex.3.27}. Then
\[
\partial_{B}g_{t}^{A}=R_{g_{t}^{A}\ast}\int_{0}^{t}Ad_{g_{\tau}^{A}}B_{\tau
}d\tau
\]
where
\[
Ad_{g}A=R_{g^{-1}\ast}L_{g\ast}A\text{ for all }g\in G\text{ and }%
A\in\mathfrak{g}.
\]

\end{example}

\begin{proof}
Let $T_{t}^{A}$ denote the flow of $A_{t}.$ Because $A_{t}$ is left invariant,%
\[
T_{t}^{A}\left(  x\right)  =xg_{t}^{A}=R_{g_{t}^{A}}x
\]
as the reader should verify. Thus%
\begin{align*}
\partial_{B}g_{t}^{A}  &  =\partial_{B}T_{t}^{A}\left(  e\right)
=R_{g_{t}^{A}\ast}\int_{0}^{t}\left(  R_{g_{\tau}^{A}\ast}\right)  ^{-1}%
\tilde{B}_{\tau}\circ R_{g_{\tau}^{A}}\left(  e\right)  d\tau\\
&  =R_{g_{t}^{A}\ast}\int_{0}^{t}\left(  R_{g_{\tau}^{A}\ast}\right)
^{-1}\tilde{B}_{\tau}\left(  g_{\tau}^{A}\right)  d\tau=R_{g_{t}^{A}\ast}%
\int_{0}^{t}\left(  R_{g_{\tau}^{A}\ast}\right)  ^{-1}L_{g_{\tau}^{A}\ast
}B_{\tau}d\tau\\
&  =R_{g_{t}^{A}\ast}\int_{0}^{t}Ad_{g_{\tau}^{A}}B_{\tau}d\tau.
\end{align*}

The next theorem expresses $\left[  X_{t},Y\right]  $ using the flow $T^{X}.$
The stochastic analog of this theorem is a key ingredient in the
\textquotedblleft Malliavin calculus,\textquotedblright\ see Proposition
\ref{p.8.14} below.
\end{proof}

\begin{theorem}
\label{t.4.9}If $X_{t}$ and $T_{t}^{X}$ are as above and $Y\in\Gamma\left(
TM\right)  ,$ then%
\begin{equation}
\frac{d}{dt}\left[  \left(  T_{t\ast}^{X}\right)  ^{-1}Y\circ T_{t}%
^{X}\right]  =\left(  T_{t\ast}^{X}\right)  ^{-1}\left[  X_{t},Y\right]  \circ
T_{t}^{X} \label{e.4.9}%
\end{equation}
or equivalently put%
\begin{equation}
\frac{d}{dt}Ad_{T_{t}^{X}}^{-1}=Ad_{T_{t}^{X}}^{-1}L_{X_{t}} \label{e.4.10}%
\end{equation}
where $L_{X}Y:=\left[  X,Y\right]  .$
\end{theorem}

\begin{proof}
Let $V_{t}:=\left(  T_{t\ast}^{X}\right)  ^{-1}Y\circ T_{t}^{X}$ which is
equivalent to $T_{t\ast}^{X}V_{t}=Y\circ T_{t}^{X},$ or more explicitly to%
\[
Yf\circ T_{t}^{X}=\left(  Y\circ T_{t}^{X}\right)  f=\left(  T_{t\ast}%
^{X}V_{t}\right)  f=V_{t}\left(  f\circ T_{t}^{X}\right)  \text{ for all }f\in
C^{\infty}(M).
\]
Differentiating this equation in $t$ then shows
\begin{align*}
\left(  X_{t}Yf\right)  \circ T_{t}^{X}  &  =\dot{V}_{t}\left(  f\circ
T_{t}^{X}\right)  +V_{t}\left(  X_{t}f\circ T_{t}^{X}\right) \\
&  =\left(  T_{t\ast}^{X}\dot{V}_{t}\right)  f+\left(  T_{t\ast}^{X}%
V_{t}\right)  X_{t}f\\
&  =\left(  T_{t\ast}^{X}\dot{V}_{t}\right)  f+\left(  Y\circ T_{t}%
^{X}\right)  X_{t}f\\
&  =\left(  T_{t\ast}^{X}\dot{V}_{t}\right)  f+\left(  YX_{t}f\right)  \circ
T_{t}^{X}.
\end{align*}
Therefore%
\[
\left(  T_{t\ast}^{X}\dot{V}_{t}\right)  f=\left(  \left[  X_{t},Y\right]
f\right)  \circ T_{t}^{X}%
\]
from which we conclude $T_{t\ast}^{X}\dot{V}_{t}=\left[  X_{t},Y\right]  \circ
T_{t}^{X}$ and therefore%
\[
\dot{V}_{t}=\left(  T_{t\ast}^{X}\right)  ^{-1}\left[  X_{t},Y\right]  \circ
T_{t}^{X}.
\]

\end{proof}

\subsection{Cartan's Development Map\label{s.4.3}}

For this section assume that $M$ is compact\footnote{It would actually be
sufficient to assume that $M$ is a \textquotedblleft
complete\textquotedblright\ Riemannian manifold for this section.} Riemannian
manifold and let $W^{\infty}(T_{0}M)$ be the collection of piecewise smooth
paths, $b:$ $[0,1]\rightarrow T_{o}M$ such that $b\left(  0\right)  =0_{o}\in
T_{o}M$ and let $W_{o}^{\infty}(M)$ be the collection of piecewise smooth
paths, $\sigma:$ $[0,1]\rightarrow M$ such that $\sigma\left(  0\right)  =o\in
M.$

\begin{theorem}
[Development Map]\label{t.4.10}To each $b\in W^{\infty}(T_{0}M)$ there is a
unique $\sigma\in W_{o}^{\infty}(M)$ such that%
\begin{equation}
\sigma^{\prime}(s):=(\sigma(s),d\sigma(s)/ds)=//_{s}(\sigma)b^{\prime}%
(s)\quad\text{\textrm{\ and }}\quad\sigma(0)=o, \label{e.4.11}%
\end{equation}
where $//_{s}(\sigma)$ denotes parallel translation along $\sigma.$
\end{theorem}

\begin{proof}
Suppose that $\sigma$ is a solution to Eq. \text{(\ref{e.4.11})} and
$//_{s}(\sigma)v_{o}=(o,u(s)v),$ where $u(s):\tau_{o}M\rightarrow
\mathbb{R}^{N}.$ Then $u$ satisfies the differential equation
\begin{equation}
u^{\prime}\left(  s\right)  +dQ(\sigma^{\prime}(s))u(s)=0\quad
\text{\textrm{\ with }}\quad u(0)=u_{0}, \label{e.4.12}%
\end{equation}
where $u_{0}v:=v$ for all $v\in\tau_{o}M,$ see Remark \ref{r.3.54}. Hence Eq.
\text{(\ref{e.4.11})} is equivalent to the following pair of coupled ordinary
differential equations:
\begin{equation}
\sigma^{\prime}\left(  s\right)  =u(s)b^{\prime}(s)\quad\text{\textrm{\ with
}}\quad\sigma(0)=o, \label{e.4.13}%
\end{equation}
and
\begin{equation}
u^{\prime}\left(  s\right)  +dQ((\sigma(s),u(s)b^{\prime}(s))u(s)=0\quad
\text{\textrm{\ with }}\quad u(0)=u_{0}. \label{e.4.14}%
\end{equation}
Therefore the uniqueness assertion follows from standard uniqueness theorems
for ordinary differential equations. The slickest prove of existence to Eq.
(\ref{e.4.11}) is to first introduce the orthogonal frame bundle, $O\left(
M\right)  ,$ on $M$ defined by $O\left(  M\right)  :=\cup_{m\in M}O_{m}(M)$
where $O_{m}(M)$ is the set of all isometries, $u:T_{o}M\rightarrow T_{m}M.$
It is then possible to show that $O\left(  M\right)  $ is an imbedded
submanifold in $\mathbb{R}^{N}\times\mathrm{\operatorname*{Hom}}\left(
\tau_{o}M,\mathbb{R}^{N}\right)  $ and that coupled pair of ordinary
differential equations (\ref{e.4.13}) and (\ref{e.4.14}) may be viewed as a
flow equation on $O(M).$ Hence the existence of solutions may be deduced from
the Theorem \ref{t.4.2}, see, for example, \cite{D5} for details of this
method. Here I will sketch a proof which does not require us to develop the
frame bundle formalism in detail.

Looking at the proof of Lemma \ref{l.2.30}, $Q$ has an extension to a
neighborhood in $\mathbb{R}^{N}$ of $m\in M$ in such a way that $Q(x)$ is
still an orthogonal projection onto $\mathrm{\mathrm{\operatorname*{Nul}}%
}(F^{\prime}(x)),$ where $F(x)=z_{>}(x)$ is as in Lemma \ref{l.2.30}. Hence
for small $s$, we may define $\sigma$ and $u$ to be the unique solutions to
Eq. \text{(\ref{e.4.13})} and Eq. \text{(\ref{e.4.14})} with values in
$\mathbb{R}^{N}$ and $\operatorname*{Hom}(\tau_{o}M,\mathbb{R}^{N})$
respectively. The key point now is to show that $\sigma(s)\in M$ and that the
range of $u(s)$ is $\tau_{\sigma(s)}M.$

Using the same proof as in Theorem \ref{t.3.52}, $w(s):=Q(\sigma(s))u(s)$
satisfies,%
\begin{align*}
w^{\prime}  &  =dQ\left(  \sigma^{\prime}\right)  u+Q\left(  \sigma\right)
u^{\prime}=dQ\left(  \sigma^{\prime}\right)  u-Q\left(  \sigma\right)
dQ(\sigma^{\prime})u\\
&  =P\left(  \sigma\right)  dQ\left(  \sigma^{\prime}\right)  u=dQ\left(
\sigma^{\prime}\right)  Q\left(  \sigma\right)  u=dQ\left(  \sigma^{\prime
}\right)  w,
\end{align*}
where Lemma \ref{l.3.30} was used in the last equality. Since $w(0)=0,$ it
follows by uniqueness of solutions to linear ordinary differential equations
that $w\equiv0$ and hence%
\[
\operatorname*{Ran}\left[  u(s)\right]  \subset\text{$\operatorname*{Nul}$%
}\,\left[  Q(\sigma(s))\right]  =\text{$\operatorname*{Nul}$}\,\left[
F^{\prime}(\sigma(s))\right]  .
\]
Consequently%
\[
dF(\sigma(s))/ds=F^{\prime}(\sigma(s))d\sigma(s)/ds=F^{\prime}(\sigma
(s))u(s)b^{\prime}(s)=0
\]
for small $s$ and since $F(\sigma(0))=F(o)=0,$ it follows that $F(\sigma
(s))=0,$ i.e. $\sigma(s)\in M.$ So we have shown that there is a solution
$(\sigma,u)$ to \text{(\ref{e.4.13})} and \text{(\ref{e.4.14})} for small $s$
such that $\sigma$ stays in $M$ and $u(s)$ is parallel translation along $s.$
By standard ordinary differential equation methods, there is a maximal
solution $(\sigma,u)$ with these properties. Notice that $(\sigma,u)$ is a
path in $M\times\text{\textrm{Iso}}(T_{o}M,\mathbb{R}^{N}),$ where
$\text{\textrm{Iso}}(T_{o}M,\mathbb{R}^{N})$ is the set of isometries from
$T_{o}M$ to $\mathbb{R}^{N}.$ Since $M\times\text{\textrm{Iso}}(T_{o}%
M,\mathbb{R}^{N})$ is a compact space, $(\sigma,u)$ can not explode. Therefore
$(\sigma,u)$ is defined on the same interval where $b$ is defined.
\end{proof}

\begin{figure}
[ptbh]
\begin{center}
\includegraphics[
height=2.6347in,
width=3.3438in
]%
{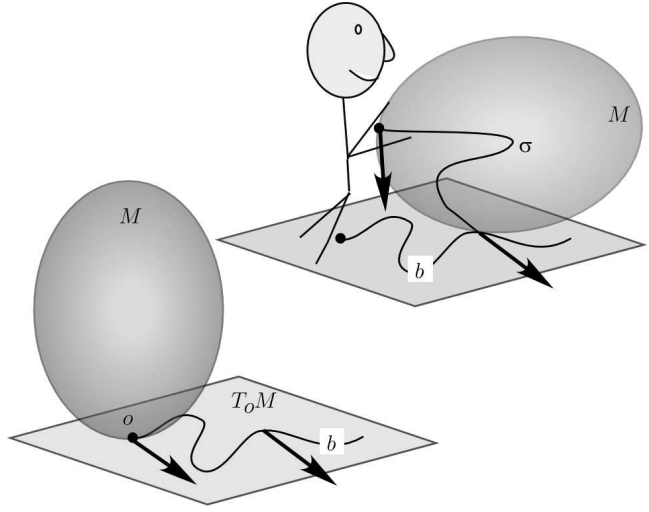}%
\caption{Monsieur Cartan is shown here rolling, without \textquotedblleft
slipping,\textquotedblright\ a manifold $M$ along a curve, $b,$ in $T_{o}M$ to
produce a curve, $\sigma,$ on $M.$}%
\label{fig.11}%
\end{center}
\end{figure}

The geometric interpretation of Cartan's map is to roll the manifold $M$ along
a freshly painted curve $b$ in $T_{o}M$ to produce a curve $\sigma$ on $M,$
see Figure \ref{fig.11}.%
%TCIMACRO{\FRAME{ftbphFU}{3.3438in}{2.6347in}{0pt}{\Qcb{Monsieur Cartan is
%shown here rolling, without \textquotedblleft slipping,\textquotedblright\ a
%manifold $M$ along a curve, $b,$ in $T_{o}M$ to produce a curve, $\sigma,$ on
%$M.$}}{\Qlb{fig.11}}{cartan.eps}{\special{ language "Scientific Word";
%type "GRAPHIC";  maintain-aspect-ratio TRUE;  display "USEDEF";
%valid_file "F";  width 3.3438in;  height 2.6347in;  depth 0pt;
%original-width 4.7458in;  original-height 3.7335in;  cropleft "0";
%filename 'cartan.eps';file-properties "XNPEU";}}}%
%BeginExpansion
%EndExpansion

\begin{notation}
\label{n.4.11}Let $\phi:W^{\infty}(T_{0}M)$ $\rightarrow W_{o}^{\infty}(M)$ be
the map $b\rightarrow\sigma,$ where $\sigma$ is the solution to
\text{(\ref{e.4.11})}. It is easy to construct the inverse map $\Psi
:=\phi^{-1}.$ Namely, $\Psi(\sigma)=b,$ where
\[
\Psi_{s}(\sigma)=b(s):=\int_{0}^{s}//_{r}(\sigma)^{-1}\sigma^{\prime}(r)dr.
\]
We now conclude this section by computing the differentials of $\Psi$ and
$\phi.$ For more details on computations of this nature the reader is referred
to \cite{Driver89,D5} and the references therein.
\end{notation}

\begin{figure}
[ptbh]
\begin{center}
\includegraphics[
height=1.6973in,
width=3.4805in
]%
{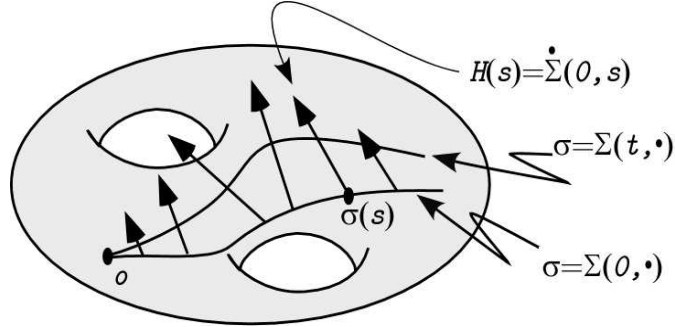}%
\caption{A variation of $\sigma$ giving rise to a vector field along $\sigma
.$}%
\label{fig.12}%
\end{center}
\end{figure}

\begin{theorem}
[Differential of $\Psi$]\label{t.4.12}Let $(t,s)\rightarrow\Sigma(t,s)$ be a
smooth map into $M$ such that $\Sigma(t,\cdot)\in W_{o}^{\infty}(M)$ for all
$t.$ Let
\[
H(s):=\dot{\Sigma}(0,s):=(\Sigma(0,s),d\Sigma(t,s)/dt|_{t=0}),
\]
so that $H$ is a vector-field along $\sigma:=\Sigma(0,\cdot).$ One should view
$H$ as an element of the \textquotedblleft tangent space\textquotedblright\ to
$W_{o}^{\infty}(M)$ at $\sigma,$ see Figure \ref{fig.12}. Let $u(s):=//_{s}%
(\sigma),$ $h(s):=//_{s}(\sigma)^{-1}H(s)$ $b:=\Psi_{s}(\sigma)$ and, for all
$a,c\in T_{o}M,$ let%
\begin{equation}
(R_{u}(a,c))(s):=u(s)^{-1}R(u(s)a,u(s)c)u(s). \label{e.4.15}%
\end{equation}
Then
\begin{equation}
d\Psi(H)=d\Psi(\Sigma(t,\cdot))/dt|_{t=0}=h+\int_{0}\left(  \int_{0}%
R_{u}(h,\delta b)\right)  \delta b, \label{e.4.16}%
\end{equation}
where $\delta b(s)$ is short hand notation for $b^{\prime}(s)ds,$ and
$\int_{0}f\delta b$ denotes the function $s\rightarrow\int_{0}^{s}%
f(r)b^{\prime}(r)dr$ when $f$ is a path of matrices.
\end{theorem}

%

%TCIMACRO{\FRAME{ftbphFU}{3.4805in}{1.6973in}{0pt}{\Qcb{A variation of $\sigma$
%giving rise to a vector field along $\sigma.$}}{\Qlb{fig.12}}{vfield.eps}%
%width 3.4805in;  height 1.6973in;  depth 0pt;  original-width 3.3616in;
%original-height 1.6249in;  cropleft "0";  croptop "1";  cropright "1";
%filename 'vfield.eps';file-properties "XNPEU";}}}%
%BeginExpansion
%EndExpansion

\begin{proof}
To simplify notation let \textquotedblleft\ $\overset{\cdot}{}$
\textquotedblright$=\frac{d}{dt}|_{0},$ \textquotedblleft\ $^{\prime}$
\textquotedblright$=\frac{d}{ds},$ $B(t,s):=\Psi(\Sigma(t,\cdot))(s)$,
$U(t,s):=//_{s}(\Sigma(t,\cdot))$, $u(s):=//_{s}(\sigma)=U(0,s)$ and
\[
\dot{b}(s):=(d\Psi(H))(s):=dB(t,s)/dt|_{t=0}.
\]
I will also suppress $(t,s)$ from the notation when possible. With this
notation
\begin{equation}
\Sigma^{\prime}=UB^{\prime},\quad\dot{\Sigma}=H=uh, \label{e.4.17}%
\end{equation}
and
\begin{equation}
\frac{\nabla U}{ds}=0. \label{e.4.18}%
\end{equation}
In Eq. (\ref{e.4.18}), $\frac{\nabla U}{ds}:T_{o}M\rightarrow T_{\Sigma}M$ is
defined by $\frac{\nabla U}{ds}=P\left(  \Sigma\right)  U^{\prime}$ or
equivalently by%
\[
\frac{\nabla U}{ds}a:=\frac{\nabla\left(  Ua\right)  }{ds}\text{ for all }a\in
T_{o}M.
\]
Taking $\nabla/dt$ of \text{(\ref{e.4.17})} at $t=0$ gives, with the aid of
Proposition \ref{p.3.32},
\[
\frac{\nabla U}{dt}|_{t=0}b^{\prime}+u\dot{b}^{\prime}=\nabla\Sigma^{\prime
}/dt|_{t=0}=\nabla\dot{\Sigma}/ds=uh^{\prime}.
\]
Therefore,
\begin{equation}
\dot{b}^{\prime}=h^{\prime}+Ab^{\prime}, \label{e.4.19}%
\end{equation}
where $A:=-U^{-1}\frac{\nabla U}{dt}|_{t=0},$ i.e.
\[
\frac{\nabla U}{dt}(0,\cdot)=-uA.
\]
Taking $\nabla/ds$ of this last equation and using $\nabla u/ds=0$ along with
Proposition \ref{p.3.32} gives%
\[
-uA^{\prime}=\left.  \frac{\nabla}{ds}\frac{\nabla}{dt}U\right\vert
_{t=0}=\left.  \left[  \frac{\nabla}{ds},\frac{\nabla}{dt}\right]
U\right\vert _{t=0}=R(\sigma^{\prime},H)u
\]
and hence $A^{\prime}=R_{u}(h,b^{\prime}).$ By integrating this identity using
$A(0)=0$ ($\nabla U(t,0)/dt=0$ since $U(t,0):=//_{0}(\Sigma(t,\cdot))=I$ is
independent of $t)$ shows%
\begin{equation}
A=\int_{0}R_{u}(h,\delta b) \label{e.4.20}%
\end{equation}
The theorem now follows by integrating \text{(\ref{e.4.19})} relative to $s$
making use of Eq. \text{(\ref{e.4.20})} and the fact that $\dot{b}(0)=0.$
\end{proof}

\begin{theorem}
[Differential of $\phi$]\label{t.4.13} Let $b,k\in W^{\infty}(T_{0}M)$ and
$(t,s)\rightarrow B(t,s)$ be a smooth map into $T_{o}M$ such that
$B(t,\cdot)\in W^{\infty}(T_{0}M)$ $,$ $B(0,s)=b(s),$ and $\dot{B}(0,s)=k(s).$
(For example take $B(t,s)=b(s)+tk(s).)$ Then
\[
\phi_{\ast}(k_{b}):=\frac{d}{dt}|_{0}\phi(B(t,\cdot))=//_{\cdot}(\sigma)h,
\]
where $\sigma:=\phi(b)$ and $h$ is the first component in the solution $(h,A)$
to the pair of coupled differential equations:
\begin{equation}
k^{\prime}=h^{\prime}+Ab^{\prime},\quad\text{\textrm{\ with }}\quad h(0)=0
\label{e.4.21}%
\end{equation}
and
\begin{equation}
A^{\prime}=R_{u}(h,b^{\prime})\quad\text{\textrm{\ with }}\quad A(0)=0.
\label{e.4.22}%
\end{equation}

\end{theorem}

\begin{proof}
This theorem has an analogous proof to that of Theorem \ref{t.4.12}. We can
also deduce the result from Theorem \ref{t.4.12}\ by defining $\Sigma$ by
$\Sigma(t,s):=\phi_{s}(B(t,\cdot)).$ We now assume the same notation used in
Theorem \ref{t.4.12}\ and its proof. Then $B(t,\cdot)=\Psi(\Sigma(t,\cdot))$
and hence by Theorem \ref{t.4.13}\
\[
k=\frac{d}{dt}|_{0}\Psi(\Sigma(t,\cdot))=d\Psi(H)=h+\int_{0}(\int_{0}%
R_{u}(h,\delta b))\delta b.
\]
Therefore, defining $A:=\int_{0}R_{u}(h,\delta b)$ and differentiating this
last equation relative to $s,$ it follows that $A$ solves \text{(\ref{e.4.22}%
)}\ and that $h$ solves \text{(\ref{e.4.21})}.
\end{proof}

The following theorem is a mild extension of Theorem \ref{t.4.12} to include
the possibility that $\Sigma(t,\cdot)\notin W_{o}^{\infty}(M)$ when $t\neq0,$
i.e. the base point may change.

\begin{theorem}
\label{t.4.14} Let $(t,s)\rightarrow\Sigma(t,s)$ be a smooth map into $M$ such
that $\sigma:=\Sigma(0,\cdot)\in W_{o}^{\infty}(M)$. Define $H(s):=d\Sigma
(t,s)/dt|_{t=0},$ $\sigma:=\Sigma(0,\cdot),$ and $h(s):=//_{s}(\sigma
)^{-1}H(s)$. (\textbf{Note:} $H(0)$ and $h(0)$ are no longer necessarily equal
to zero.) Let
\[
U(t,s):=//_{s}(\Sigma(t,\cdot))//_{t}(\Sigma(\cdot,0)):T_{o}M\rightarrow
T_{\Sigma(t,s)}M,
\]
so that $\nabla U(t,0)/dt=0$ and $\nabla U(t,s)/ds\equiv0.$ Set $B(t,s):=\int
_{0}^{s}U(t,r)^{-1}\Sigma^{\prime}(t,r)dr$, then
\begin{equation}
\dot{b}(s):=\frac{d}{dt}|_{0}B(t,s)=h_{s}+\int_{0}^{s}\left(  \int_{0}%
R_{u}(h,\delta b)\right)  \delta b, \label{e.4.23}%
\end{equation}
where as before $b:=\Psi(\sigma).$
\end{theorem}

\begin{proof}
The proof is almost identical to the proof of Theorem \ref{t.4.12} and hence
will be omitted.
\end{proof}

\section{Stochastic Calculus on Manifolds\label{s.5}}

In this section and the rest of the text the reader is assumed to be well
versed in stochastic calculus in the Euclidean context.

\begin{notation}
\label{n.5.1}In the sequel we will always assume there is any underlying
filtered probability space $(\Omega,\{\mathcal{F}_{s}\}_{s\geq0}%
,\mathcal{F},\mu)$ satisfying the \textquotedblleft usual
hypothesis.\textquotedblright\ Namely, $\mathcal{F}$ is $\mu$ -- complete,
$\mathcal{F}_{s}$ contains all of the null sets in $\mathcal{F},$ and
$\mathcal{F}_{s}$ is right continuous. As usual $\mathbb{E}$ will be used to
denote the expectation relative to the probability measure $\mu.$
\end{notation}

\begin{definition}
\label{d.5.2}For simplicity, we will call a function $\Sigma:\mathbb{R}%
_{+}\times\Omega\rightarrow V$ $(V$ a vector space) a \textbf{process} if
$\Sigma_{s}=\Sigma(s):=\Sigma(s,\cdot)$ is $\mathcal{F}_{s}$ -- measurable for
all $s\in\mathbb{R}_{+}:=[0,\infty),$ i.e. a process will mean an adapted
process unless otherwise stated. As above, we will always assume that $M$ is
an imbedded submanifold of $\mathbb{R}^{N}$ with the induced Riemannian
structure. An $M$ -- \textbf{valued semi-martingale }is a \textbf{continuous}
$\mathbb{R}^{N}$-valued semi-martingale $(\Sigma)$ such that $\Sigma
(s,\omega)\in M$ for all $(s,\omega)\in\mathbb{R}_{+}\times\Omega.$ It will be
convenient to let $\lambda$ be the distinguished process: $\lambda\left(
s\right)  =\lambda_{s}:=s.$
\end{definition}

Since $f\in C^{\infty}(M)$ is the restriction of a smooth function $F$ on
$\mathbb{R}^{N},$ it follows by It\^{o}'s lemma that $f\circ\Sigma
=F\circ\Sigma$ is a real-valued semi-martingale if $\Sigma$ is an $M$ --
valued semi-martingale. Conversely, if $\Sigma$ is an $M$ -- valued process
and $f\circ\Sigma$ is a real-valued semi-martingale for all $f\in C^{\infty
}(M)$ then $\Sigma$ is an $M$ -- valued semi-martingale. Indeed, let
$x=(x^{1},\ldots,x^{N})$ be the standard coordinates on $\mathbb{R}^{N},$ then
$\Sigma^{i}:=x^{i}\circ\Sigma$ is a real semi-martingale for each $i,$ which
implies that $\Sigma$ is a $\mathbb{R}^{N}$- valued semi-martingale.

\begin{notation}
[Fisk-Stratonovich Integral]\label{n.5.3}Suppose $V$ is a finite dimensional
vector space and%
\[
\pi=\{0=s_{0}<s_{1}<s_{2}<\cdots\}
\]
is a partition of $\mathbb{R}_{+}\ $with $\lim_{n\rightarrow\infty}%
s_{n}=\infty.$ To such a partition $\pi,$ let $|\pi|:=\sup_{i}|s_{i+1}-s_{i}|$
be the \textbf{mesh size }of $\pi$ and $s\wedge s_{i}:=\min\{s,s_{i}\}.$ To
each $\mathrm{\operatorname*{Hom}}\left(  \mathbb{R}^{N},V\right)  $ -- valued
semi-martingale $Z_{t}$ and each $M$ -- valued semi-martingale $\Sigma_{t},$
the \textbf{Fisk-Stratonovich integral }of $Z$ relative to $\Sigma$ is defined
by%
\begin{align*}
\int_{0}^{s}Z\delta\Sigma &  =\lim_{|\pi|\rightarrow0}\sum_{i=0}^{\infty
}{\frac{1}{2}}\left(  Z_{s\wedge s_{i}}+Z_{s\wedge s_{i+1}}\right)
(\Sigma_{s\wedge s_{i+1}}-\Sigma_{s\wedge s_{i}})\\
&  =\int_{0}^{s}Zd\Sigma+\frac{1}{2}\int_{0}^{s}dZd\Sigma\in V
\end{align*}
where
\[
\int_{0}^{s}Zd\Sigma=\lim_{|\pi|\rightarrow0}\sum_{i=0}^{\infty}Z_{s\wedge
s_{i}}(\Sigma_{s\wedge s_{i+1}}-\Sigma_{s\wedge s_{i}})\in V
\]
is the \textbf{It\^{o} integral }and
\[
\lbrack Z,\Sigma]_{s}=\int_{0}^{s}dZd\Sigma:=\lim_{|\pi|\rightarrow0}%
\sum_{i=0}^{\infty}\left(  Z_{s\wedge s_{i}}-Z_{s\wedge s_{i+1}}\right)
(\Sigma_{s\wedge s_{i+1}}-\Sigma_{s\wedge s_{i}})\in V
\]
is the \textbf{mutual variation} of $Z$ and $\Sigma.$ (All limits may be taken
in the sense of uniform convergence on compact subsets of $\mathbb{R}_{+}$ in probability.)
\end{notation}

\subsection{Stochastic Differential Equations on Manifolds\label{s.5.1}}

\begin{notation}
\label{n.5.4}Suppose that $\left\{  X_{i}\right\}  _{i=0}^{n}\subset
\Gamma\left(  TM\right)  $ are vector fields on $M.$ For $a\in\mathbb{R}^{n}$
let
\[
X_{a}\left(  m\right)  :=\mathbf{X}\left(  m\right)  a:=\sum_{i=1}^{n}%
a_{i}X_{i}\left(  m\right)
\]
With this notation, \textbf{$X$}$\left(  m\right)  :\mathbb{R}^{n}\rightarrow
T_{m}M$ is a linear map for each $m\in M.$
\end{notation}

\begin{definition}
\label{d.5.5}Given an $\mathbb{R}^{n}$ -- valued semi-martingale, $\beta_{s},$
we say an $M$ -- valued semi-martingale $\Sigma_{s}$ solves the
Fisk-Stratonovich stochastic differential equation%
\begin{equation}
\delta\Sigma_{s}=\mathbf{X}\left(  \Sigma_{s}\right)  \delta\beta_{s}%
+X_{0}\left(  \Sigma_{s}\right)  ds:=\sum_{i=1}^{n}X_{i}\left(  \Sigma
_{s}\right)  \delta\beta_{s}^{i}+X_{0}\left(  \Sigma_{s}\right)
ds\label{e.5.1}%
\end{equation}
if for all $f\in C^{\infty}(M),$
\[
\delta f\left(  \Sigma_{s}\right)  =\sum_{i=1}^{n}\left(  X_{i}f\right)
\left(  \Sigma_{s}\right)  \delta\beta_{s}^{i}+X_{0}f\left(  \Sigma
_{s}\right)  ds,
\]
i.e. if
\[
f\left(  \Sigma_{s}\right)  =f\left(  \Sigma_{0}\right)  +\sum_{i=1}^{n}%
\int_{0}^{s}\left(  X_{i}f\right)  \left(  \Sigma_{r}\right)  \delta\beta
_{r}^{i}+\int_{0}^{s}X_{0}f\left(  \Sigma_{r}\right)  dr.
\]

\end{definition}

\begin{lemma}
[It\^{o} Form of Eq. (\ref{e.5.1})]\label{l.5.6}Suppose that $\beta=B$ is an
$\mathbb{R}^{n}$ -- valued Brownian motion and let $L:=\frac{1}{2}\sum
_{i=1}^{n}X_{i}^{2}+X_{0}.$ Then an $M$ -- valued semi-martingale $\Sigma_{s}$
solves Eq. (\ref{e.5.1}) iff%
\begin{equation}
f\left(  \Sigma_{s}\right)  =f\left(  \Sigma_{0}\right)  +\sum_{i=1}^{n}%
\int_{0}^{s}\left(  X_{i}f\right)  \left(  \Sigma_{r}\right)  dB_{r}^{i}%
+\int_{0}^{s}Lf\left(  \Sigma_{r}\right)  dr \label{e.5.2}%
\end{equation}
for all $f\in C^{\infty}(M).$
\end{lemma}

\begin{proof}
Suppose that $\Sigma_{s}$ solves Eq. (\ref{e.5.1}), then
\begin{align*}
d\left[  \left(  X_{i}f\right)  \left(  \Sigma_{r}\right)  \right]   &
=\sum_{j=1}^{n}\left(  X_{j}X_{i}f\right)  \left(  \Sigma_{r}\right)  \delta
B_{s}^{j}+X_{0}X_{i}f\left(  \Sigma_{s}\right)  ds\\
&  =\sum_{j=1}^{n}\left(  X_{j}X_{i}f\right)  \left(  \Sigma_{r}\right)
dB_{s}^{j}+d\left(  BV\right)
\end{align*}
where $BV$ denotes a process of bounded variation. Hence%
\begin{align*}
\int_{0}^{s}\left(  X_{i}f\right)  \left(  \Sigma_{r}\right)  \delta
B_{r}^{i}  &  =\sum_{i=1}^{n}\int_{0}^{s}\left(  X_{i}f\right)  \left(
\Sigma_{r}\right)  dB_{r}^{i}+\frac{1}{2}\int_{0}^{s}d\left[  \left(
X_{i}f\right)  \left(  \Sigma_{r}\right)  \right]  dB_{r}^{i}\\
&  =\sum_{i=1}^{n}\int_{0}^{s}\left(  X_{i}f\right)  \left(  \Sigma
_{r}\right)  dB_{r}^{i}+\frac{1}{2}\sum_{i,j=1}^{n}\int_{0}^{s}\left(
X_{j}X_{i}f\right)  \left(  \Sigma_{r}\right)  dB_{s}^{j}dB_{r}^{i}\\
&  =\sum_{i=1}^{n}\int_{0}^{s}\left(  X_{i}f\right)  \left(  \Sigma
_{r}\right)  dB_{r}^{i}+\frac{1}{2}\int_{0}^{s}\sum_{i=1}^{n}X_{i}^{2}f\left(
\Sigma_{r}\right)  dr.
\end{align*}

Similarly if Eq. (\ref{e.5.2}) holds for all $f\in C^{\infty}\left(  M\right)
$ we have%
\[
d\left[  \left(  X_{i}f\right)  \left(  \Sigma_{r}\right)  \right]  =\left(
X_{j}X_{i}f\right)  \left(  \Sigma_{r}\right)  dB_{s}^{j}+LX_{i}f\left(
\Sigma_{s}\right)  ds
\]
and so as above%
\[
\int_{0}^{s}\left(  X_{i}f\right)  \left(  \Sigma_{r}\right)  \delta B_{r}%
^{i}=\sum_{i=1}^{n}\int_{0}^{s}\left(  X_{i}f\right)  \left(  \Sigma
_{r}\right)  dB_{r}^{i}+\frac{1}{2}\int_{0}^{s}\sum_{i=1}^{n}X_{i}^{2}f\left(
\Sigma_{r}\right)  dr.
\]
Solving for $\int_{0}^{s}\left(  X_{i}f\right)  \left(  \Sigma_{r}\right)
dB_{r}^{i}$ and putting the result into Eq. (\ref{e.5.2}) shows%
\begin{align*}
f\left(  \Sigma_{s}\right)   &  =f\left(  \Sigma_{0}\right)  +\sum_{i=1}%
^{n}\int_{0}^{s}\left(  X_{i}f\right)  \left(  \Sigma_{r}\right)  \delta
B_{r}^{i}-\frac{1}{2}\int_{0}^{s}\sum_{i=1}^{n}X_{i}^{2}f\left(  \Sigma
_{r}\right)  dr+\int_{0}^{s}Lf\left(  \Sigma_{r}\right)  dr\\
&  =f\left(  \Sigma_{0}\right)  +\sum_{i=1}^{n}\int_{0}^{s}\left(
X_{i}f\right)  \left(  \Sigma_{r}\right)  \delta B_{r}^{i}+\int_{0}^{s}%
X_{0}f\left(  \Sigma_{r}\right)  dr.
\end{align*}

\end{proof}

To avoid technical problems with possible explosions of stochastic
differential equations in the sequel, we make the following assumption.

\begin{assumption}
{\label{ass.2}}Unless otherwise stated, in the remainder of these notes, $M$
will be a compact manifold imbedded in $E:=\mathbb{R}^{N}.$
\end{assumption}

To shortcut the development of a number of issues here it is useful to recall
the following Wong and Zakai type approximation theorem for solutions to
Fisk-Stratonovich stochastic differential equations.

\begin{notation}
\label{n.5.7}Let $\{B_{s}\}_{s\in\lbrack0,T]}$ be a standard $\mathbb{R}^{n}%
$---valued Brownian motion. Given a partition
\[
\pi=\{0=s_{0}<s_{1}<s_{2}<...<s_{k}=T\}
\]
of $[0,T],$ let
\[
\left\vert \pi\right\vert =\max\left\{  s_{i}-s_{i-1}:i=1,2,\dots,k\right\}
\]
and%
\[
B_{\pi}(s)=B(s_{i-1})+(s-s_{i-1})\frac{\Delta_{i}B}{\Delta_{i}s}%
\text{\textrm{\ if }}s\in(s_{i-1},s_{i}],
\]
where $\Delta_{i}B:=B(s_{i})-B(s_{i-1})$ and $\Delta_{i}s:=s_{i}-s_{i-1}.$
Notice that $B_{\pi}\left(  s\right)  $ is a continuous piecewise linear path
in $\mathbb{R}^{n}.$
\end{notation}

\begin{theorem}
[Wong-Zakai type approximation theorem]\label{t.5.8}Let $a\in\mathbb{R}^{N},$%
\[
f:\mathbb{R}^{n}\times\mathbb{R}^{N}\rightarrow\mathrm{\operatorname*{Hom}%
}(\mathbb{R}^{n},\mathbb{R}^{N})\text{ and }f_{0}:\mathbb{R}^{n}%
\times\mathbb{R}^{N}\rightarrow\mathbb{R}^{N}%
\]
be twice differentiable functions with bounded continuous derivatives. Let
$\pi$ and $B_{\pi}$ be as in Notation \ref{n.5.7} and $\xi_{\pi}(s)$ denote
the solution to the ordinary differential equation:
\begin{equation}
\xi_{\pi}^{\prime}(s)=f(B_{\pi}(s),\xi_{\pi}(s))B_{\pi}^{\prime}%
(s)+f_{0}(B_{\pi}(s),\xi_{\pi}(s)),\qquad\xi_{\pi}(0)=a \label{e.5.3}%
\end{equation}
and $\xi$ denote the solution to the Fisk-Stratonovich stochastic differential
equation,
\begin{equation}
d\xi_{s}=f(B_{s},\xi_{s})\delta B_{s}+f_{0}(B_{s},\xi_{s})ds,\qquad\xi_{0}=a.
\label{e.5.4}%
\end{equation}
Then, for any $\gamma\in(0,\frac{1}{2})$ and $p\in\lbrack1,\infty)$, there is
a constant $C(p,\gamma)<\infty$ such that
\begin{equation}
\lim_{|\pi|\rightarrow0}\mathbb{E}\left[  \sup_{s\leq T}|\xi_{\pi}(s)-\xi
_{s}|^{p}\right]  \leq C(p,\gamma)|\pi|^{\gamma p}. \label{e.5.5}%
\end{equation}

\end{theorem}

This theorem is a special case of Theorem 5.7.3 and Example 5.7.4 in Kunita
\cite{Kunita90}. Theorems of this type have a long history starting with Wong
and Zakai \cite{Wong:Zakai65b,Wong:Zakai67}. The reader may also find this and
related results in the following \emph{partial}\textbf{\ }list of references:
\cite{Amit91,Bally89,Bally89b,Bismut81,Blum84,Doss79,Elworthy78,Guo-Shu82,Ikeda81,Jorgensen75,Kaneko:Nakao,Kurtz:Protter:1991a,Kurtz:Protter:1991b,Lyons96,Malliavin78b,Malliavin97,McShane72,McShane74,Moulinier88,Nakao:Yamato,Pinsky78,Stroock:Taniguchi1994,StVar69b,Stroock:Varadhan72,Sussmann91,Tw93}%
. Also see \cite{Driver98a,Driver04} and the references therein for more of
the geometry associated to the Wong and Zakai approximation scheme.
%BRUCE: Comment out: The theorem as stated here may be found in <cite>DHu</cite>.

\begin{remark}
[Transfer Principle]\label{r.5.9}Theorem \ref{t.5.8} is a manifestation of the
\textbf{transfer principle} (coined by Malliavin) which loosely states: to get
a correct stochastic formula one should take the corresponding deterministic
smooth formula and replace all derivatives by Fisk-Stratonovich differentials.
We will see examples of this principle over and over again in the sequel.
\end{remark}

\begin{theorem}
\label{t.5.10}Given a point $m\in M$ there exits a unique $M$ -- valued semi
martingale $\Sigma$ which solves Eq. (\ref{e.5.1}) with the initial condition,
$\Sigma_{0}=m.$ We will write $T_{s}\left(  m\right)  $ for $\Sigma_{s}$ if we
wish to emphasize the dependence of the solution on the initial starting point
$m\in M.$
\end{theorem}

\begin{proof}
\textbf{Existence.} If for the moment we assumed that the Brownian motion
$B_{s}$ were differentiable in $s,$ Eq. (\ref{e.5.1}) could be written as%
\[
\Sigma_{s}^{\prime}=X_{s}\left(  \Sigma_{s}\right)  \text{ with }\Sigma_{0}=m
\]
where
\[
X_{s}\left(  m\right)  :=\sum_{i=1}^{n}X_{i}\left(  m\right)  \left(
B^{i}\right)  ^{\prime}\left(  s\right)  +X_{0}\left(  m\right)
\]
and the existence of $\Sigma_{s}$ could be deduced from Theorem \ref{t.4.2}.
We will make this rigorous with an application of Theorem \ref{t.5.8}.

Let $\left\{  Y_{i}\right\}  _{i=0}^{n}$ be smooth vector fields on $E$ with
compact support such that $Y_{i}=X_{i}$ on $M$ for each $i$ and let $B_{\pi
}\left(  s\right)  $ be as in Notation \ref{n.5.7} and define%
\begin{align*}
X_{s}^{\pi}\left(  m\right)   &  :=\sum_{i=1}^{n}X_{i}\left(  m\right)
\left(  B_{\pi}^{i}\right)  ^{\prime}\left(  s\right)  +X_{0}\left(  m\right)
\text{ and}\\
Y_{s}^{\pi}\left(  m\right)   &  :=\sum_{i=1}^{n}Y_{i}\left(  m\right)
\left(  B_{\pi}^{i}\right)  ^{\prime}\left(  s\right)  +Y_{0}\left(  m\right)
.
\end{align*}
Then by Theorem \ref{t.4.2} we may use $X^{\pi}$ and $Y^{\pi}$ to generate
(random) flows $T^{\pi}:=T^{X^{\pi}}$ on $M$ and $\tilde{T}^{\pi}:=T^{Y^{\pi}%
}$ on $E$ respectively. Moreover, as in the proof of Theorem \ref{t.4.2} we
know $T_{s}^{\pi}(m)=$ $\tilde{T}_{s}^{\pi}(m)$ for all $m\in M.$ An
application of Theorem \ref{t.5.8} now shows that $\Sigma_{s}:=\tilde{T}%
_{s}\left(  m\right)  :=\lim_{\left\vert \pi\right\vert \rightarrow0}\tilde
{T}_{s}^{\pi}(m)=\lim_{\left\vert \pi\right\vert \rightarrow0}T_{s}^{\pi
}(m)\in M$ exists\footnote{Here we have used the fact that $M$ is a closed
subset of $\mathbb{R}^{N}.$} and satisfies the Fisk-Stratonovich differential
equation \textbf{on }$E,$%
\begin{equation}
d\Sigma_{s}=\sum_{i=1}^{n}Y_{i}\left(  \Sigma_{s}\right)  \delta B_{s}%
^{i}+Y_{0}\left(  \Sigma_{s}\right)  ds\text{ with }\Sigma_{0}=m.
\label{e.5.6}%
\end{equation}
Given $f\in C^{\infty}(M),$ let $F\in C^{\infty}(E)$ be chosen so that
$f=F|_{M}.$ Then Eq. (\ref{e.5.6}) implies%
\begin{equation}
d\left[  F\left(  \Sigma_{s}\right)  \right]  =\sum_{i=1}^{n}Y_{i}F\left(
\Sigma_{s}\right)  \delta B_{s}^{i}+Y_{0}F\left(  \Sigma_{s}\right)  ds.
\label{e.5.7}%
\end{equation}
Since we have already seen $\Sigma_{s}\in M$ and by construction $Y_{i}=X_{i}$
on $M,$ we have $F\left(  \Sigma_{s}\right)  =f\left(  \Sigma_{s}\right)  $
and $Y_{i}F\left(  \Sigma_{s}\right)  =X_{i}f\left(  \Sigma_{s}\right)  .$
Therefore Eq. (\ref{e.5.7}) implies%
\[
d\left[  f\left(  \Sigma_{s}\right)  \right]  =\sum_{i=1}^{n}X_{i}f\left(
\Sigma_{s}\right)  \delta B_{s}^{i}+Y_{0}F\left(  \Sigma_{s}\right)  ds,
\]
i.e. $\Sigma_{s}$ solves Eq. (\ref{e.5.1}) as desired.

\textbf{Uniqueness. }If $\Sigma$ is a solution to Eq. (\ref{e.5.1}), then for
$F\in C^{\infty}(E),$ we have%
\begin{align*}
dF\left(  \Sigma_{s}\right)   &  =\sum_{i=1}^{n}X_{i}F\left(  \Sigma
_{s}\right)  \delta B_{s}^{i}+X_{0}F\left(  \Sigma_{s}\right)  ds\\
&  =\sum_{i=1}^{n}Y_{i}F\left(  \Sigma_{s}\right)  \delta B_{s}^{i}%
+Y_{0}F\left(  \Sigma_{s}\right)  ds
\end{align*}
which shows, by taking $F$ to be the standard linear coordinates on $E,$
$\Sigma_{s}$ also solves Eq. (\ref{e.5.6}). But this is a stochastic
differential equation on a Euclidean space $E$ with smooth compactly supported
coefficients and therefore has a \textbf{unique }solution.
\end{proof}

\subsection{Line Integrals}

For $a,b\in\mathbb{R}^{N},$ let $\langle a,b\rangle_{\mathbb{R}^{N}}%
:=\sum_{i=1}^{N}a_{i}b_{i}$ denote the standard inner product on
$\mathbb{R}^{N}.$ Also let $\mathfrak{g}l(N)=\mathfrak{g}l(N,\mathbb{R})$ be
the set of $N\times N$ real matrices. (It is not necessary to assume $M$ is
compact for most of the results in this section.)

\begin{theorem}
\label{t.5.11}As above, for $m\in M,$ let $P\left(  m\right)  $ and $Q\left(
m\right)  $ denote orthogonal projection or $\mathbb{R}^{N}$ onto $\tau_{m}M$
and $\tau_{m}M^{\perp}$ respectively. Then for any $M$ -- valued
semi-martingale $\Sigma,$
\[
0=Q(\Sigma)\delta\Sigma\text{ and }d\Sigma=P\left(  \Sigma\right)
\delta\Sigma,
\]
i.e.
\[
\Sigma_{s}-\Sigma_{0}=\int_{0}^{s}P(\Sigma_{r})\delta\Sigma_{r}.
\]

\end{theorem}

\begin{proof}
We will first assume that $M$ is the level set of a function $F$ as in Theorem
\ref{t.2.5}. Then we may assume that
\[
Q(x)=\phi(x)F^{\prime}(x)^{\ast}(F^{\prime}(x)F^{\prime}(x)^{\ast}%
)^{-1}F^{\prime}(x),
\]
where $\phi$ is smooth function on $\mathbb{R}^{N}$ such that $\phi:=1$ in a
neighborhood of $M$ and the support of $\phi$ is contained in the set:
$\{x\in\mathbb{R}^{N}|F^{\prime}(x)\text{\textrm{\ is surjective}}\}.$ By
It\^{o}'s lemma
\[
0=d0=d(F(\Sigma))=F^{\prime}(\Sigma)\delta\Sigma.
\]
The lemma follows in this special case by multiplying the above equation
through by $\phi(\Sigma)F^{\prime}(\Sigma)^{\ast}(F^{\prime}(\Sigma)F^{\prime
}(\Sigma)^{\ast})^{-1},$ see the proof of Lemma \ref{l.2.30}.

For the general case, choose two open covers $\{V_{i}\}$ and $\{U_{i}\}$ of
$M$ such that each $\bar{V}_{i}$ is compactly contained in $U_{i},$ there is a
smooth function $F_{i}\in C_{c}^{\infty}(U_{i}\rightarrow\mathbb{R}^{N-d})$
such that $V_{i}\cap M=V_{i}\cap\{F_{i}^{-1}(\{0\})\}$ and $F_{i}$ has a
surjective differential on $V_{i}\cap M.$ Choose $\phi_{i}\in C_{c}^{\infty
}(\mathbb{R}^{N})$ such that the support of $\phi_{i}$ is contained in $V_{i}$
and $\sum\phi_{i}=1$ on $M,$ with the sum being locally finite. (For the
existence of such covers and functions, see the discussion of partitions of
unity in any reasonable book about manifolds.) Notice that $\phi_{i}\cdot
F_{i}\equiv0$ and that $F_{i}\cdot\phi_{i}^{\prime}\equiv0$ on $M$ so that
\begin{align*}
0  &  =d\{\phi_{i}(\Sigma)F_{i}(\Sigma)\}=(\phi_{i}^{\prime}(\Sigma
)\delta\Sigma)F_{i}(\Sigma)+\phi_{i}(\Sigma)F_{i}^{\prime}(\Sigma)\delta
\Sigma\\
&  =\phi_{i}(\Sigma)F_{i}^{\prime}(\Sigma)\delta\Sigma.
\end{align*}
Multiplying this equation by $\Psi_{i}(\Sigma)F_{i}^{\prime}(\Sigma)^{\ast
}(F_{i}^{\prime}(\Sigma)F_{i}^{\prime}(\Sigma)^{\ast})^{-1},$ where each
$\Psi_{i}$ is a smooth function on $\mathbb{R}^{N}$ such that $\Psi_{i}%
\equiv1$ on the support of $\phi_{i}$ and the support of $\Psi_{i}$ is
contained in the set where $F_{i}^{\prime}$ is surjective, we learn that
\begin{equation}
0=\phi_{i}(\Sigma)F_{i}^{\prime}(\Sigma)^{\ast}(F_{i}^{\prime}(\Sigma
)F_{i}^{\prime}(\Sigma)^{\ast})^{-1}F_{i}^{\prime}(\Sigma)\delta\Sigma
=\phi_{i}(\Sigma)Q(\Sigma)\delta\Sigma\label{e.5.8}%
\end{equation}
for all $i.$ By a stopping time argument we may assume that $\Sigma$ never
leaves a compact set, and therefore we may choose a finite subset $I$ of the
indices $\{i\}$ such that $\sum_{i\in I}\phi_{i}(\Sigma)Q(\Sigma)=Q(\Sigma).$
Hence summing over $i\in I$ in equation \text{(\ref{e.5.8})}\ shows that
$0=Q(\Sigma)\delta\Sigma.$ Since $Q+P=I,$ it follows that
\[
d\Sigma=I\delta\Sigma=\left[  Q(\Sigma)+P(\Sigma)\right]  \delta
\Sigma=P\left(  \Sigma\right)  \delta\Sigma.
\]

\end{proof}

The following notation will be needed to define line integrals along a
semi-martingale $\Sigma.$

\begin{notation}
\label{n.5.12}Let $P\left(  m\right)  $ be orthogonal projection of
$\mathbb{R}^{N}$ onto $\tau_{m}M$ as above.

\begin{enumerate}
\item Given a one-form $\alpha$ on $M$ let $\tilde{\alpha}:M\rightarrow
(\mathbb{R}^{N})^{\ast}$ be defined by
\begin{equation}
\tilde{\alpha}(m)v:=\alpha((P(m)v)_{m}) \label{e.5.9}%
\end{equation}
for all $m\in M$ and $v\in\mathbb{R}^{N}.$

\item Let $\Gamma(T^{\ast}M\otimes T^{\ast}M)$ denote the set of functions
$\rho:\cup_{m\in M}T_{m}M\otimes T_{m}M\rightarrow\mathbb{R}$ such that
$\rho_{m}:=\rho|_{T_{m}M\otimes T_{m}M}$ is linear, and $m\rightarrow
\rho(X(m)\otimes Y(m))$ is a smooth function on $M$ for all smooth
vector-fields $X,Y\in\Gamma(TM).$ (Riemannian metrics and Hessians of smooth
functions are examples of elements of $\Gamma(T^{\ast}M\otimes T^{\ast}M).)$

\item For $\rho\in\Gamma(T^{\ast}M\otimes T^{\ast}M),$ let $\tilde{\rho
}:M\rightarrow(\mathbb{R}^{N}\otimes\mathbb{R}^{N})^{\ast}$ be defined by
\begin{equation}
\tilde{\rho}(m)(v\otimes w):=\rho((P(m)v)_{m}\otimes(P(m)w)_{m}).
\label{e.5.10}%
\end{equation}

\end{enumerate}
\end{notation}

\begin{definition}
\label{d.5.13}Let $\alpha$ be a one form on $M,$ $\rho\in\Gamma(T^{\ast
}M\otimes T^{\ast}M),$ and $\Sigma$ be an $M$ -- valued semi-martingale. Then
the \textbf{Fisk-Stratonovich} integral of $\alpha$ along $\Sigma$ is:
\begin{equation}
\int_{0}^{\cdot}\alpha(\delta\Sigma):=\int_{0}^{\cdot}\tilde{\alpha}%
(\Sigma)\delta\Sigma, \label{e.5.11}%
\end{equation}
and the \textbf{It\^{o}} integral is given by:
\begin{equation}
\int_{0}^{\cdot}\alpha({\bar{d}}\Sigma):=\int_{0}^{\cdot}\tilde{\alpha}%
(\Sigma)d\Sigma, \label{e.5.12}%
\end{equation}
where the stochastic integrals on the right hand sides of Eqs.
\text{(\ref{e.5.11})}\ and \text{(\ref{e.5.12})}\ are Fisk-Stratonovich and
It\^{o} integrals respectively. Formally, ${\bar{d}}\Sigma:=P(\Sigma)d\Sigma.$
We also define \textbf{quadratic } \textbf{integral}:
\begin{equation}
\int_{0}^{\cdot}\rho(d\Sigma\otimes d\Sigma):=\int_{0}^{\cdot}\tilde{\rho
}(\Sigma)(d\Sigma\otimes d\Sigma):=\sum_{i,j=1}^{N}\int_{0}^{\cdot}\tilde
{\rho}(\Sigma)(e_{i}\otimes e_{j})d[\Sigma^{i},\Sigma^{j}], \label{e.5.13}%
\end{equation}
where $\{e_{i}\}_{i=1}^{N}$ is an orthonormal basis for $\mathbb{R}^{N},$
$\Sigma^{i}:=\langle e_{i},\Sigma\rangle,$ and $d[\Sigma^{i},\Sigma^{j}]$ is
the differential of the mutual quadratic variation of $\Sigma^{i}$ and
$\Sigma^{j}.$
\end{definition}

So as not to confuse $[\Sigma^{i},\Sigma^{j}]$ with a commutator or a Lie
bracket, in the sequel we will write $d\Sigma^{i}d\Sigma^{j}$ for
$d[\Sigma^{i},\Sigma^{j}].$

\begin{remark}
\label{r.5.14} The above definitions may be generalized as follows. Suppose
that $\alpha$ is now a $T^{\ast}M$ -- valued semi-martingale and $\Sigma$ is
the $M$ valued semi-martingale such that $\alpha_{s}\in T_{\Sigma_{s}}^{\ast
}M$ for all $s.$ Then we may define
\[
\tilde{\alpha}_{s}v:=\alpha_{s}((P(\Sigma_{s})v)_{\Sigma_{s}}),
\]%
\begin{equation}
\int_{0}^{\cdot}\alpha(\delta\Sigma):=\int_{0}^{\cdot}\tilde{\alpha}%
\delta\Sigma, \label{e.5.14}%
\end{equation}
and
\begin{equation}
\int_{0}^{\cdot}\alpha({\bar{d}}\Sigma):=\int_{0}^{\cdot}\tilde{\alpha}%
d\Sigma. \label{e.5.15}%
\end{equation}
Similarly, if $\rho$ is a process in $T^{\ast}M\otimes T^{\ast}M$ such that
$\rho_{s}\in T_{\Sigma_{s}}^{\ast}M\otimes T_{\Sigma_{s}}^{\ast}M$, let
\begin{equation}
\int_{0}^{\cdot}\rho(d\Sigma\otimes d\Sigma)=\int_{0}^{\cdot}\tilde{\rho
}(d\Sigma\otimes d\Sigma), \label{e.5.16}%
\end{equation}
where
\[
\tilde{\rho}_{s}(v\otimes w):=\rho_{s}((P(\Sigma_{s})v)_{\Sigma_{s}}%
\otimes(P(\Sigma_{s})v)_{\Sigma_{s}})
\]
and
\begin{equation}
d\Sigma\otimes d\Sigma=\sum_{i,j=1}^{N}e_{i}\otimes e_{j}d\Sigma^{i}%
d\Sigma^{j} \label{e.5.17}%
\end{equation}
as in Eq. \text{(\ref{e.5.13})}.
\end{remark}

\begin{lemma}
\label{l.5.15}Suppose that $\alpha=fdg$ for some functions $f,g\in C^{\infty
}(M),$ then
\[
\int_{0}^{\cdot}\alpha(\delta\Sigma)=\int_{0}^{\cdot}f(\Sigma)\delta\lbrack
g(\Sigma)].
\]
Since, by Corollary \ref{c.3.42}, any one form $\alpha$ on $M$ may be written
as $\alpha=\sum_{i=1}^{N}f_{i}dg_{i}$ with $f_{i},g_{i}\in C^{\infty}(M),$ it
follows that the Fisk-Stratonovich integral is intrinsically defined
independent of how $M$ is imbedded into a Euclidean space.
\end{lemma}

\begin{proof}
Let $G$ be a smooth function on $\mathbb{R}^{N}$ such that $g=G|_{M}.$ Then
$\tilde{\alpha}(m)=f(m)G^{\prime}(m)P(m),$ so that
\begin{align*}
\int_{0}^{\cdot}\alpha(\delta\Sigma)  &  =\int_{0}^{\cdot}f(\Sigma)G^{\prime
}(\Sigma)P(\Sigma)\delta\Sigma\\
&  =\int_{0}^{\cdot}f(\Sigma)G^{\prime}(\Sigma)\delta\Sigma\qquad
\text{\textrm{(by Theorem \ref{t.5.11})}}\\
&  =\int_{0}^{\cdot}f(\Sigma)\delta\lbrack G(\Sigma)]\qquad\text{\textrm{(by
It\^{o}'s Lemma)}}\\
&  =\int_{0}^{\cdot}f(\Sigma)\delta\lbrack g(\Sigma)].\qquad(g(\Sigma
)=G(\Sigma))
\end{align*}

\end{proof}

\begin{lemma}
\label{l.5.16} Suppose that $\rho=fdh\otimes dg$, where $f,g,h\in C^{\infty
}(M),$ then
\[
\int_{0}^{\cdot}\rho(d\Sigma\otimes d\Sigma)=\int_{0}^{\cdot}f(\Sigma
)d[h(\Sigma),g(\Sigma)]=:\int_{0}^{\cdot}f(\Sigma)d\left[  h(\Sigma)\right]
d\left[  g(\Sigma)\right]  .
\]
Since, by an argument similar to that in Corollary \ref{c.3.42}, any $\rho
\in\Gamma(T^{\ast}M\otimes T^{\ast}M)$ may be written as a finite linear
combination $\rho=\sum_{i}f_{i}dh_{i}\otimes dg_{i}$ with $f_{i},h_{i}%
,g_{i}\in C^{\infty}(M),$ it follows that the quadratic integral is
intrinsically defined independent of the imbedding.
\end{lemma}

\begin{proof}
By Theorem \ref{t.5.11}, $\delta\Sigma=P(\Sigma)\delta\Sigma,$ so that
\begin{align*}
\Sigma_{s}^{i}  &  =\Sigma_{0}^{i}+\int_{0}^{\cdot}(e_{i},P(\Sigma
)d\Sigma)+B.V.\\
&  =\Sigma_{0}^{i}+\sum_{k}\int_{0}^{\cdot}(e_{i},P(\Sigma)e_{k})d\Sigma
^{k}+B.V.,
\end{align*}
where $B.V.$ denotes a process of bounded variation. Therefore
\begin{equation}
d[\Sigma^{i},\Sigma^{j}]=\sum_{k,l}(e_{i},P(\Sigma)e_{k})(e_{i},P(\Sigma
)e_{l})d\Sigma^{k}d\Sigma^{l}. \label{e.5.18}%
\end{equation}
Now let $H$ and $G$ be in $C^{\infty}(\mathbb{R}^{N})$ such that $h=H|_{M}$
and $g=G|_{M}.$ By It\^{o}'s lemma and Eq. (\ref{e.5.18}),%
\begin{align*}
d[h(\Sigma),g(\Sigma)]  &  =\sum_{i,j}(H^{\prime}(\Sigma)e_{i})(G^{\prime
}(\Sigma)e_{j})d[\Sigma^{i},\Sigma^{j}]\\
&  =\sum_{i,j,k,l}(H^{\prime}(\Sigma)e_{i})(G^{\prime}(\Sigma)e_{j}%
)(e_{i},P(\Sigma)e_{k})(e_{i},P(\Sigma)e_{l})d\Sigma^{k}d\Sigma^{l}\\
&  =\sum_{k,l}(H^{\prime}(\Sigma)P(\Sigma)e_{k})(G^{\prime}(\Sigma
)P(\Sigma)e_{l})d\Sigma^{k}d\Sigma^{l}.
\end{align*}
Since
\[
\tilde{\rho}(m)=f(m)\cdot(H^{\prime}(m)P(m))\otimes(G^{\prime}(m)P(m)),
\]
it follows from Eq. \text{(\ref{e.5.13})}\ and the two above displayed
equations that
\begin{align*}
\int_{0}^{\cdot}f(\Sigma)d[h(\Sigma),g(\Sigma)]  &  :=\int_{0}^{\cdot}%
\sum_{k,l}f(\Sigma)(H^{\prime}(\Sigma)P(\Sigma)e_{k})(G^{\prime}%
(\Sigma)P(\Sigma)e_{l})d\Sigma^{k}d\Sigma^{l}\\
&  =\int_{0}^{\cdot}\tilde{\rho}(\Sigma)(d\Sigma\otimes d\Sigma)=:\int
_{0}^{\cdot}\rho(d\Sigma\otimes d\Sigma).
\end{align*}

\end{proof}

\begin{theorem}
\label{t.5.17}Let $\alpha$ be a one form on $M$, and $\Sigma$ be a $M$ --
valued semi-martingale. Then
\begin{equation}
\int_{0}^{\cdot}\alpha(\delta\Sigma)=\int_{0}^{\cdot}\alpha({\bar{d}}%
\Sigma)+{\frac{1}{2}}\int_{0}^{\cdot}\nabla\alpha(d\Sigma\otimes d\Sigma),
\label{e.5.19}%
\end{equation}
where $\nabla\alpha(v_{m}\otimes w_{m}):=(\nabla_{v_{m}}\alpha)(w_{m})$ and
$\nabla\alpha$ is defined in Definition \ref{d.3.40}, also see Lemma
\ref{l.3.41}. (This shows that the It\^{o} integral depends not only on the
manifold structure of $M$ but on the geometry of $M$ as reflected in the
Levi-Civita covariant derivative $\nabla.)$
\end{theorem}

\begin{proof}
Let $\tilde{\alpha}$ be as in Eq. \text{(\ref{e.5.9})}. For the purposes of
the proof, suppose that $\tilde{\alpha}:M\rightarrow(\mathbb{R}^{N})^{\ast}$
has been extended to a smooth function from $\mathbb{R}^{N}\rightarrow
(\mathbb{R}^{N})^{\ast}.$ We still denote this extension by $\tilde{\alpha}.$
Then using Eq. \text{(\ref{e.5.18})},
\begin{align*}
\int_{0}^{\cdot}  &  \alpha(\delta\Sigma):=\int_{0}^{\cdot}\tilde{\alpha
}(\Sigma)\delta\Sigma\\
&  =\int_{0}^{\cdot}\tilde{\alpha}(\Sigma)d\Sigma+{\frac{1}{2}}\int_{0}%
^{\cdot}\tilde{\alpha}^{\prime}(\Sigma)(d\Sigma)d\Sigma\\
&  =\int_{0}^{\cdot}\alpha({\bar{d}}\Sigma)+{\frac{1}{2}}\sum_{i,j,k,l}%
\int_{0}^{\cdot}\tilde{\alpha}^{\prime}(\Sigma)(e_{i})e_{j}(e_{i}%
,P(\Sigma)e_{k})(e_{i},P(\Sigma)e_{l})d\Sigma^{k}d\Sigma^{l}\\
&  =\int_{0}^{\cdot}\alpha({\bar{d}}\Sigma)+{\frac{1}{2}}\sum_{k,l}\int
_{0}^{\cdot}\tilde{\alpha}^{\prime}(\Sigma)(P(\Sigma)e_{k})P(\Sigma
)e_{l}d\Sigma^{k}d\Sigma^{l}\\
&  =\int_{0}^{\cdot}\alpha({\bar{d}}\Sigma)+{\frac{1}{2}}\sum_{k,l}\int
_{0}^{\cdot}d\tilde{\alpha}((P(\Sigma)e_{k})_{\Sigma})P(\Sigma)e_{l}%
d\Sigma^{k}d\Sigma^{l}.
\end{align*}
But by Eq. \text{(\ref{e.3.45})}, we know for all $v_{m},w_{m}\in TM$ that
\[
\nabla\alpha(v_{m}\otimes w_{m})=d\tilde{\alpha}(v_{m})w-\tilde{\alpha
}(m)dQ(v_{m})w.
\]
Since $\tilde{\alpha}(m)=\tilde{\alpha}(m)P(m)$ and $PdQ=dQQ$ (Lemma
\ref{l.3.30}), we find
\[
\tilde{\alpha}(m)dQ(v_{m})w=\tilde{\alpha}(m)dQ(v_{m})Q(m)w=0~\forall
~v_{m},w_{m}\in TM.
\]
Hence combining the three above displayed equations shows that
\begin{align*}
\int_{0}^{\cdot}\alpha(\delta\Sigma)  &  =\int_{0}^{\cdot}\alpha({\bar{d}%
}\Sigma)+{\frac{1}{2}}\sum_{k,l}\int_{0}^{\cdot}\nabla\alpha((P(\Sigma
)e_{k})_{\Sigma}\otimes(P(\Sigma)e_{l})_{\Sigma})d\Sigma^{k}d\Sigma^{l}\\
&  =\int_{0}^{\cdot}\alpha({\bar{d}}\Sigma)+{\frac{1}{2}}\sum_{k,l}\int
_{0}^{\cdot}\nabla\alpha(d\Sigma\otimes d\Sigma).
\end{align*}

\end{proof}

\begin{corollary}
[It\^{o}'s Lemma for Manifolds]\label{c.5.18}If $u\in C^{\infty}\left(
\left(  0,T\right)  \times M\right)  $ and $\Sigma$ is an $M-$valued
semi-martingale, then%
\begin{align}
d\left[  u\left(  s,\Sigma_{s}\right)  \right]   &  =\left(  \partial
_{s}u\right)  \left(  s,\Sigma_{s}\right)  ds\nonumber\\
&  +d_{M}\left[  u\left(  s,\cdot\right)  \right]  ({\bar{d}}\Sigma
_{s})+{\frac{1}{2}}\left(  \nabla d_{M}u\left(  s,\cdot\right)  \right)
(d\Sigma_{s}\otimes d\Sigma_{s}), \label{e.5.20}%
\end{align}
where, as in Notation \ref{n.2.20}, $d_{M}u\left(  s,\cdot\right)  $ is being
used to denote the differential of the map: $m\in M\rightarrow u\left(
s,m\right)  .$
\end{corollary}

\begin{proof}
Let $U\in C^{\infty}(\left(  0,T\right)  \times\mathbb{R}^{N})$ such that
$u\left(  s,\cdot\right)  =U\left(  s,\cdot\right)  |_{M}.$ Then by It\^{o}'s
lemma and Theorem \ref{t.5.11},
\begin{align*}
d\left[  u\left(  s,\Sigma_{s}\right)  \right]   &  =d\left[  U\left(
s,\Sigma_{s}\right)  \right]  =\left(  \partial_{s}U\right)  \left(
s,\Sigma_{s}\right)  ds+D_{\Sigma}U(s,\Sigma_{s})\delta\Sigma_{s}\\
&  =\left(  \partial_{s}U\right)  \left(  s,\Sigma_{s}\right)  ds+D_{\Sigma
}U(s,\Sigma_{s})P(\Sigma_{s})\delta\Sigma_{s}\\
&  =\left(  \partial_{s}u\right)  \left(  s,\Sigma_{s}\right)  ds+d_{M}\left[
u\left(  s,\cdot\right)  \right]  (\delta\Sigma_{s})\\
&  =\left(  \partial_{s}u\right)  \left(  s,\Sigma_{s}\right)  ds+d_{M}\left[
u\left(  s,\cdot\right)  \right]  ({\bar{d}}\Sigma_{s})\\
&  \qquad+{\frac{1}{2}}\left(  \nabla d_{M}u\left(  s,\cdot\right)  \right)
(d\Sigma_{s}\otimes d\Sigma_{s}),
\end{align*}
wherein the last equality is a consequence of Theorem \ref{t.5.17}.
\end{proof}

\subsection{$M$ -- valued Martingales and Brownian Motions\label{s.5.3}}

\begin{definition}
\label{d.5.19}An $M$ -- valued semi-martingale $\Sigma$ is said to be a
(local) \textbf{martingale }(more precisely a $\nabla$-martingale) if
\begin{equation}
\int_{0}^{\cdot}df({\bar{d}}\Sigma)=f(\Sigma)-f(\Sigma_{0})-{\frac{1}{2}}%
\int_{0}^{\cdot}\nabla df(d\Sigma\otimes d\Sigma) \label{e.5.21}%
\end{equation}
is a (local) martingale for all $f\in C^{\infty}(M).$ (See Theorem
\ref{t.5.17} for the truth of the equality in Eq. (\ref{e.5.21}).) The process
$\Sigma$ is said to be a \textbf{Brownian motion} if
\begin{equation}
f(\Sigma)-f(\Sigma_{0})-{\frac{1}{2}}\int_{0}^{\cdot}\Delta f(\Sigma
)d\lambda\label{e.5.22}%
\end{equation}
is a local martingale for all $f\in C^{\infty}(M),$ where $\lambda(s):=s$ and
$\int_{0}^{\cdot}\Delta f(\Sigma)d\lambda$ denotes the process $s\rightarrow
\int_{0}^{s}\Delta f(\Sigma)d\lambda.$
\end{definition}

\begin{figure}
[ptbh]
\begin{center}
\includegraphics[
height=2.3162in,
width=4.412in
]%
{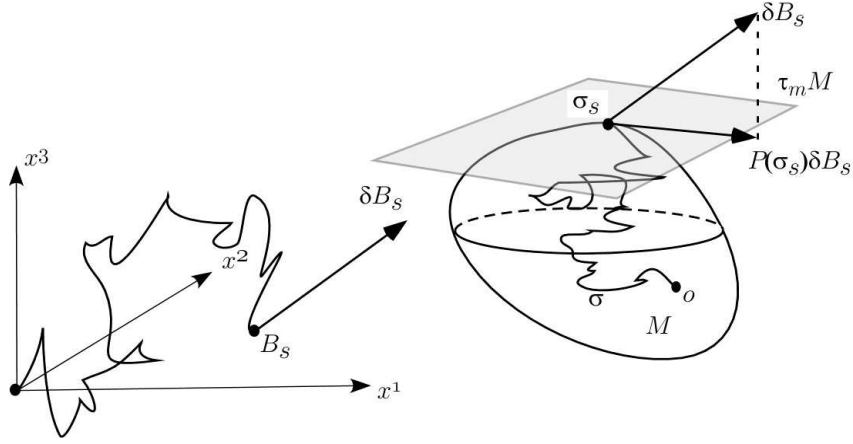}%
\caption{Projection construction of Brownian motion on $M.$}%
\label{f.13}%
\end{center}
\end{figure}

\begin{theorem}
[Projection Construction of Brownian Motion]\label{t.5.20}Suppose that
$B=\left(  B^{1},B^{2},\dots,B^{N}\right)  $ is an $N$ -- dimensional Brownian
motion. The there is a unique $M$ -- valued semi-martingale $\Sigma$ which
solves the Fisk-Stratonovich stochastic differential equation,%
\begin{equation}
\delta\Sigma=P(\Sigma)\delta B\quad\text{\textrm{\ with }}\quad\Sigma_{0}=o\in
M, \label{e.5.23}%
\end{equation}
see Figure \ref{f.13}. Moreover, $\Sigma$ is an $M$ -- valued Brownian motion.
\end{theorem}

%

%TCIMACRO{\FRAME{ftbphFU}{4.412in}{2.3162in}{0pt}{\Qcb{Projection construction
%of Brownian motion on $M.$}}{\Qlb{f.13}}{stroock.eps}%
%width 4.412in;  height 2.3162in;  depth 0pt;  original-width 5.3744in;
%original-height 2.8084in;  cropleft "0";  croptop "1";  cropright "1";
%filename 'stroock.eps';file-properties "XNPEU";}}}%
%BeginExpansion
%EndExpansion

\begin{proof}
Let $\left\{  e_{i}\right\}  _{i=1}^{N}$ be the standard basis for
$\mathbb{R}^{N}$ and $X_{i}\left(  m\right)  :=P\left(  m\right)  e_{i}\in
T_{m}M$ for each $i=1,2,\dots,N$ and $m\in M.$ Then Eq. (\ref{e.5.23}) is
equivalent to the Stochastic differential equation.,
\[
\delta\Sigma=\sum_{i=1}^{N}X_{i}(\Sigma)\delta B^{i}\quad\text{\textrm{\ with
}}\quad\Sigma_{0}=o\in M
\]
which has a unique solution by Theorem \ref{t.5.10}. Using Lemma \ref{l.5.6},
this equation may be rewritten in It\^{o} form as
\[
d\left[  f\left(  \Sigma\right)  \right]  =\sum_{i=1}^{N}X_{i}f(\Sigma
)dB^{i}+\frac{1}{2}\sum_{i=1}^{N}X_{i}^{2}f\left(  \Sigma\right)  ds\text{ for
all }f\in C^{\infty}(M).
\]
This completes the proof since $\sum_{i=1}^{N}X_{i}^{2}=\Delta$ by Proposition
\ref{p.3.48}.
\end{proof}

\begin{lemma}
[L\'{e}vy's Criteria]\label{l.5.21}For each $m\in M,$ let $\mathcal{I}%
(m):=\sum_{i=1}^{d}E_{i}\otimes E_{i},$ where $\{E_{i}\}_{i=1}^{d}$ is an
orthonormal basis for $T_{m}M.$ An $M$ -- valued semi-martingale, $\Sigma,$ is
a Brownian motion iff $\Sigma$ is a martingale and%
\begin{equation}
d\Sigma\otimes d\Sigma=\mathcal{I}(\Sigma)d\lambda. \label{e.5.24}%
\end{equation}
More precisely, this last condition is to be interpreted as:
\begin{equation}
\int_{0}^{\cdot}\rho(d\Sigma\otimes d\Sigma)=\int_{0}^{\cdot}\rho
(\mathcal{I}(\Sigma))d\lambda~\forall~\rho\in\Gamma(T^{\ast}M\otimes T^{\ast
}M). \label{e.5.25}%
\end{equation}

\end{lemma}

\begin{proof}
$(\Rightarrow)$ Suppose that $\Sigma$ is a Brownian motion on $M$ (so Eq.
\text{(\ref{e.5.22}) holds) }and $f,g\in C^{\infty}(M).$ Then on one hand
\begin{align*}
d(f(\Sigma)g(\Sigma))  &  =d\left[  f(\Sigma)\right]  \cdot g(\Sigma
)+f(\Sigma)d\left[  g(\Sigma)\right]  +d[f(\Sigma),g(\Sigma)]\\
&  \cong{\frac{1}{2}}\{\Delta f(\Sigma)g(\Sigma)+f(\Sigma)\Delta
g(\Sigma)\}d\lambda+d[f(\Sigma),g(\Sigma)],
\end{align*}
where $``\cong$\textquotedblright\ denotes equality up to the differential of
a martingale. On the other hand,
\begin{align*}
d(f(\Sigma)g(\Sigma))  &  \cong{\frac{1}{2}}\Delta(fg)(\Sigma)d\lambda\\
&  ={\frac{1}{2}}\{\Delta f(\Sigma)g(\Sigma)+f(\Sigma)\Delta g(\Sigma
)+2\langle\operatorname*{grad}f,\text{$\operatorname*{grad}$}g\rangle
(\Sigma)\}d\lambda.
\end{align*}
Comparing the above two equations implies that
\[
d[f(\Sigma),g(\Sigma)]=\langle\operatorname*{grad}%
f,\text{$\operatorname*{grad}$}g\rangle(\Sigma)d\lambda=df\otimes
dg(\mathcal{I}(\Sigma)\rangle d\lambda.
\]
Therefore by Lemma \ref{l.5.16}, if $\rho=h\cdot df\otimes dg$ then
\begin{align*}
\int_{0}^{\cdot}\rho(d\Sigma\otimes d\Sigma)  &  =\int_{0}^{\cdot}%
h(\Sigma)d[f(\Sigma),g(\Sigma)]\\
&  =\int_{0}^{\cdot}h(\Sigma)(df\otimes dg)(\mathcal{I}(\Sigma))d\lambda
=\int_{0}^{\cdot}\rho(\mathcal{I}(\Sigma))d\lambda.
\end{align*}
Since the general element $\rho$ of $\Gamma(T^{\ast}M\otimes T^{\ast}M)$ is a
finite linear combination of expressions of the form $hdf\otimes dg,$ it
follows that Eq. \text{(\ref{e.5.24})}\ holds. Moreover, Eq. (\ref{e.5.24})
implies%
\begin{equation}
\left(  \nabla df\right)  (d\Sigma\otimes d\Sigma)=\left(  \nabla df\right)
(\mathcal{I}(\Sigma))d\lambda=\Delta f(\Sigma)d\lambda\label{e.5.26}%
\end{equation}
and therefore,
\begin{align}
f(\Sigma)  &  -f(\Sigma_{0})-{\frac{1}{2}}\int_{0}^{\cdot}\nabla
df(d\Sigma\otimes d\Sigma)\nonumber\\
&  =f(\Sigma)-f(\Sigma_{0})-{\frac{1}{2}}\int_{0}^{\cdot}\Delta f(\Sigma
)d\lambda\label{e.5.27}%
\end{align}
is a martingale and so by definition $\Sigma$ is a martingale.

Conversely assume $\Sigma$ is a martingale and Eq. (\ref{e.5.24}) holds. Then
Eq. (\ref{e.5.26}) and Eq. (\ref{e.5.27}) hold and they imply $\Sigma$ is a
Brownian motion, see Definition \ref{d.5.19}.
\end{proof}

\begin{definition}
[$\delta^{\nabla}V:=P\delta V$]\label{d.5.22}Suppose $\alpha$ is a one form on
$M$ and $V$ is a $TM-$valued semi-martingale, i.e. $V_{s}=(\Sigma_{s},v_{s}),$
where $\Sigma$ is an $M$ -- valued semi-martingale and $v$ is a $\mathbb{R}%
^{N}$-valued semi-martingale such that $v_{s}\in\tau_{\Sigma_{s}}M$ for all
$s.$ Then we define:
\begin{equation}
\int_{0}^{\cdot}\alpha(\delta^{\nabla}V):=\int_{0}^{\cdot}\tilde{\alpha
}(\Sigma)\delta v=\int_{0}^{\cdot}\alpha(\Sigma)\left(  P\left(
\Sigma\right)  \delta v\right)  . \label{e.5.28}%
\end{equation}

\end{definition}

\begin{remark}
\label{r.5.23}Suppose that $\alpha(v_{m})=\theta(m)v,$ where $\theta
:M\rightarrow(\mathbb{R}^{N})^{\ast}$ is a smooth function. Then
\[
\int_{0}^{\cdot}\alpha(\delta^{\nabla}V):=\int_{0}^{\cdot}\theta
(\Sigma)P(\Sigma)\delta v=\int_{0}^{\cdot}\theta(\Sigma)\{\delta
v+dQ(\delta\Sigma)v\},
\]
where we have used the identity:
\begin{equation}
\delta^{\nabla}V=P(\Sigma)\delta v=\delta v+dQ(\delta\Sigma)v. \label{e.5.29}%
\end{equation}
This last identity follows by taking the differential of the identity,
$v=P(\Sigma)v,$ as in the proof of Proposition \ref{p.3.32}.
\end{remark}

\begin{proposition}
[Product Rule]\label{p.5.24}Keeping the notation of above, we have
\begin{equation}
\delta(\alpha(V))=\nabla\alpha(\delta\Sigma\otimes V)+\alpha(\delta^{\nabla
}V), \label{e.5.30}%
\end{equation}
where $\nabla\alpha(\delta\Sigma\otimes V):=\gamma(\delta\Sigma)$ and $\gamma$
is the $T^{\ast}M$ -- valued semi-martingale defined by
\[
\gamma_{s}\left(  w\right)  :=\nabla\alpha(w\otimes V_{s})=\left(  \nabla
_{w}\alpha\right)  (V_{s})\text{ for any }w\in T_{\Sigma_{s}}M.
\]

\end{proposition}

\begin{proof}
Let $\theta:\mathbb{R}^{N}\rightarrow(\mathbb{R}^{N})^{\ast}$ be a smooth map
such that $\tilde{\alpha}(m)=\theta(m)|_{\tau_{m}M}$ for all $m\in M.$ By
Lemma \ref{l.5.15},\ $\delta(\theta(\Sigma)P(\Sigma))=d(\theta P)(\delta
\Sigma)$ and hence by Lemma \ref{l.3.41},\ $\delta(\theta(\Sigma
)P(\Sigma))v=\nabla\alpha(\delta\Sigma\otimes V),$ where $\nabla\alpha
(v_{m}\otimes w_{m}):=(\nabla_{v_{m}}\alpha)(w_{m})$ for all $v_{m},w_{m}\in
TM.$ Therefore:
\begin{align*}
\delta(\alpha(V))  &  =\delta(\theta(\Sigma)v)=\delta(\theta(\Sigma
)P(\Sigma)v)=(d(\theta P)(\delta\Sigma))v+\theta(\Sigma)P(\Sigma)\delta v\\
&  =(d(\theta P)(\delta\Sigma))v+\tilde{\alpha}(\Sigma)\delta v=\nabla
\alpha(\delta\Sigma\otimes V)+\alpha(\delta^{\nabla}V).
\end{align*}

\end{proof}

\subsection{Stochastic Parallel Translation and Development Maps\label{s.5.4}}

\begin{definition}
\label{d.5.25}A $TM$ -- valued semi-martingale $V$ is said to be parallel if
$\delta^{\nabla}V\equiv0,$ i.e. $\int_{0}^{\cdot}\alpha(\delta^{\nabla
}V)\equiv0$ for all one forms $\alpha$ on $M.$
\end{definition}

\begin{proposition}
\label{p.5.26}A $TM$ -- valued semi-martingale $V=(\Sigma,v)$ is parallel iff
\begin{equation}
\int_{0}^{\cdot}P(\Sigma)\delta v=\int_{0}^{\cdot}\{\delta v+dQ(\delta
\Sigma)v\}\equiv0. \label{e.5.31}%
\end{equation}

\end{proposition}

\begin{proof}
Let $x=(x^{1},\ldots,x^{N})$ denote the standard coordinates on $\mathbb{R}%
^{N}.$ If $V$ is parallel then,%
\[
0\equiv\int_{0}^{\cdot}dx^{i}(\delta^{\nabla}V)=\int_{0}^{\cdot}\langle
e_{i},P(\Sigma)\delta v\rangle
\]
for each $i\ $which implies Eq. \text{(\ref{e.5.31})}. The converse follows
from Remark \ref{r.5.23}.
\end{proof}

In the following theorem, $V_{0}$ is said to be a measurable vector-field on
$M$ if $V_{0}(m)=(m,v(m))$ with $v:M\rightarrow\mathbb{R}^{N}$ being a
measurable function such that $v(m)\in\tau_{m}M$ for all $m\in M.$

\begin{theorem}
[Stochastic Parallel Translation on $M\times\mathbb{R}^{N}$]\label{t.5.27}Let
$\Sigma$ be an $M$ -- valued semi-martingale, and $V_{0}(m)=(m,v(m))$ be a
measurable vector-field on $M,$ then there is a unique parallel $TM$-valued
semi-martingale $V$ such that $V_{0}=V_{0}(\Sigma_{0})$ and $V_{s}\in
T_{\Sigma_{s}}M$ for all $s.$ Moreover, if $u$ denotes the solution to the
stochastic differential equation:
\begin{equation}
\delta u+\Gamma(\delta\Sigma)u=0\quad\text{\textrm{\ with }}\quad u_{0}=I\in
O(N), \label{e.5.32}%
\end{equation}
(where $O\left(  N\right)  $ is as in Example \ref{ex.2.6} and $\Gamma$ is as
in Eq. (\ref{e.3.65})) then $V_{s}=(\Sigma_{s},u_{s}v(\Sigma_{0}).$ The
process $u$ defined in \text{(\ref{e.5.32})}\ is orthogonal for all $s$ and
satisfies $P(\Sigma_{s})u_{s}=u_{s}P(\Sigma_{0}).$ Moreover if $\Sigma
_{0}=o\in M$ a.e. and $v\in\tau_{o}M$ and $w\perp\tau_{o}M,$ then $u_{s}v$ and
$u_{s}w$ satisfy
\begin{equation}
\delta\left[  u_{s}v\right]  +dQ\left(  \delta\Sigma\right)  u_{s}v=P\left(
\Sigma\right)  \delta\left[  u_{s}v\right]  =0 \label{e.5.33}%
\end{equation}
and%
\begin{equation}
\delta\left[  u_{s}w\right]  +dP\left(  \delta\Sigma\right)  u_{s}v=Q\left(
\Sigma\right)  \delta\left[  u_{s}v\right]  =0. \label{e.5.34}%
\end{equation}

\end{theorem}

\begin{proof}
The assertions prior to Eq. (\ref{e.5.33}) are the stochastic analogs of
Lemmas \ref{l.3.56} and \ref{l.3.57}. The proof may be given by replacing
$\frac{d}{ds}$ everywhere in the proofs of Lemmas \ref{l.3.56} and
\ref{l.3.57} by $\delta_{s}$ to get a proof in this stochastic setting. Eqs.
(\ref{e.5.33}) and (\ref{e.5.34}) are now easily verified, for example using
and $P\left(  \Sigma\right)  uv=uv,$ we have%
\[
\delta\left[  uv\right]  =\delta\left[  P\left(  \Sigma\right)  uv\right]
=P\left(  \delta\Sigma\right)  uv+P\left(  \Sigma\right)  \delta\left[
uv\right]
\]
which proves the first equality in Eq. (\ref{e.5.33}). For the second equality
in Eq. (\ref{e.5.33}),%
\begin{align*}
P\left(  \Sigma\right)  \delta\left[  uv\right]   &  =-P\left(  \Sigma\right)
\Gamma\left(  \delta\Sigma\right)  \left[  uv\right] \\
&  =-P\left(  \Sigma\right)  \left[  dQ(\delta\Sigma)P(\Sigma)+dP(\delta
\Sigma)Q(\Sigma)\right]  \left[  uv\right] \\
&  =-dQ(\delta\Sigma)Q\left(  \Sigma\right)  P(\Sigma)\delta\left[  uv\right]
=0
\end{align*}
where Lemma \ref{l.3.30} was used in the third equality. The proof of Eq.
(\ref{e.5.34}) is completely analogous. The skeptical reader is referred to
Section 3 of Driver \cite{D5}\ for more details.
\end{proof}

\begin{definition}
[Stochastic Parallel Translation]\label{d.5.28} Given $v\in\mathbb{R}^{N}$ and
an $M$ -- valued semi-martingale $\Sigma,$ let $//_{s}(\Sigma)v_{\Sigma_{0}%
}=(\Sigma_{s},u_{s}v),$ where $u$ solves \text{(\ref{e.5.32})}. (Note:
$V_{s}=//_{s}(\Sigma)V_{0}.$)
\end{definition}

In the remainder of these notes, I will often abuse notation and write $u_{s}$
instead of $//_{s}:=//_{s}(\Sigma)$ and $v_{s}$ rather than $V_{s}=(\Sigma
_{s},v_{s}).$ For example, the reader should sometimes interpret $u_{s}v$ as
$//_{s}(\Sigma)v_{\Sigma_{0}}$ depending on the context. Essentially, we will
be identifying $\tau_{m}M$ with $T_{m}M$ when no particular confusion will arise.

\textbf{Convention. }Let us now fix a\textbf{ base point} $o\in M$ and unless
otherwise noted, we will assume that all $M$ -- valued semi-martingales,
$\Sigma,$ start of $o\in M,$ i.e. $\Sigma_{0}=o$ a.e.

To each $M$ -- valued semi-martingale, $\Sigma,$ let $\Psi(\Sigma):=b$ where
\[
b:=\int_{0}^{\cdot}//^{-1}\delta\Sigma=\int_{0}^{\cdot}u^{-1}\delta\Sigma
=\int_{0}^{\cdot}u^{\mathrm{tr}}\delta\Sigma.
\]
Then $b=\Psi(\Sigma)$ is a $T_{o}M$ -- valued semi-martingale such that
$b_{0}=0_{o}\in T_{o}M.$ The converse holds as well.

\begin{theorem}
[Stochastic Development Map]\label{t.5.29}Suppose that $o\in M$ is given and
$b$ is a $T_{o}M$ -- valued semi-martingale. Then there exists a unique $M$ --
valued semi-martingale $\Sigma$ such that
\begin{equation}
\delta\Sigma_{s}=//_{s}\delta b_{s}=u_{s}\delta b_{s}\quad\text{\textrm{\ with
}}\quad\Sigma_{0}=o \label{e.5.35}%
\end{equation}
where $u$ solves \text{(\ref{e.5.32})}. $.$
\end{theorem}

\begin{proof}
This theorem is a stochastic analog of Theorem \ref{t.4.10} and the reader is
again referred to Figure \ref{fig.11}. To prove the existence and uniqueness,
we may follow the method in the proof of Theorem \ref{t.4.10}. Namely, the
pair $\left(  \Sigma,u\right)  \in M\times O\left(  N\right)  $ solves an
Stochastic differential equation. of the form
\begin{align*}
\delta\Sigma &  =u\delta b\quad\text{\textrm{with}}\quad\Sigma_{0}=o\\
\delta u  &  =-\Gamma\left(  \delta\Sigma\right)  u=-\Gamma\left(  u\delta
b\right)  u\quad\text{\textrm{with}}\quad u_{0}=I\in O(N)
\end{align*}
which after a little effort can be expressed in a form for which Theorem
\ref{t.5.10} may be applied. The details will be left to the reader, or see
(for example) Section 3 of Driver \cite{D5}.
\end{proof}

\begin{notation}
\label{n.5.30}As in the smooth case, define $\Sigma=\phi(b),$ so that%
\[
\Psi\left(  \Sigma\right)  :=\phi^{-1}\left(  b\right)  =\int_{0}^{\cdot
}//_{r}\left(  \Sigma\right)  ^{-1}\delta\Sigma_{r}.
\]

\end{notation}

In what follows, we will assume that $b_{s},$ $u_{s}$ $($or equivalently
$//_{s}(\Sigma)),$ and $\Sigma_{s}$ are related by Equations
\text{(\ref{e.5.35})} and \text{(\ref{e.5.32})}$,$ i.e. $\Sigma=\phi\left(
b\right)  $ and $u=//=//\left(  \Sigma\right)  .$ Recall that ${\bar{d}}%
\Sigma=P\left(  \Sigma\right)  d\Sigma$ is the It\^{o} differential of
$\Sigma,$ see Definition \ref{d.5.13}.

\begin{proposition}
\label{p.5.31}Let $\Sigma=\phi\left(  b\right)  ,$ then%
\begin{equation}
{\bar{d}}\Sigma=P(\Sigma)d\Sigma=udb. \label{e.5.36}%
\end{equation}
Also
\begin{equation}
d\Sigma\otimes d\Sigma=udb\otimes udb:=\sum_{i,j=1}^{d}ue_{i}\otimes
ue_{j}db^{i}db^{j}, \label{e.5.37}%
\end{equation}
where $\{e_{i}\}_{i=1}^{d}$ is an orthonormal basis for $T_{o}M$ and
$b=\sum_{i=1}^{d}b^{i}e_{i}.$ More precisely
\[
\int_{0}^{\cdot}\rho(d\Sigma\otimes d\Sigma)=\int_{0}^{\cdot}\sum_{i,j=1}%
^{d}\rho(ue_{i}\otimes ue_{j})db^{i}db^{j},
\]
for all $\rho\in\Gamma(T^{\ast}M\otimes T^{\ast}M).$
\end{proposition}

\begin{proof}
Consider the identity:
\begin{align*}
d\Sigma &  =u\delta b=udb+{\frac{1}{2}}dudb\\
&  =udb-{\frac{1}{2}}\Gamma(\delta\Sigma)udb=udb-{\frac{1}{2}}\Gamma(udb)udb
\end{align*}
where $\Gamma$ is as defined in Eq. (\ref{e.3.65}). Hence
\[
{\bar{d}}\Sigma=P(\Sigma)d\Sigma=udb-{\frac{1}{2}}\sum_{i,j=1}^{d}%
P(\Sigma)\Gamma((ue_{i})_{\Sigma})ue_{j}db^{i}db^{j}.
\]
The proof of Eq. \text{(\ref{e.5.36})}\ is finished upon observing,
\[
P\Gamma P=P\{dQP+dPQ\}P=PdQP=PQdQ=0.
\]
The proof of Eq. \text{(\ref{e.5.37})}\ is easy and will be left for the reader.
\end{proof}

\begin{fact}
\label{fact.5.32}If $(M,g)$ is a complete Riemannian manifold and the Ricci
curvature tensor is bounded from below\footnote{These assumptions are always
satisfied when $M$ is compact.}, then $\Delta=\Delta_{g}$ acting on
$C_{c}^{\infty}(M)$ is essentially self-adjoint, i.e. the closure $\bar
{\Delta}$ of $\Delta$ is an unbounded self-adjoint operator on $L^{2}(M,dV).$
(Here $dV=\sqrt{g}dx^{1}\dots dx^{n}$ is being used to denote the Riemann
volume measure on $M.)$ Moreover, the semi-group $e^{t\bar{\Delta}/2}$ has a
smooth integral kernel, $p_{t}(x,y),$ such that
\begin{align*}
p_{t}(x,y)  &  \geq0\text{ for all }x,y\in M\\
\int_{M}p_{t}(x,y)dV(y)  &  =1\text{ for all }x\in M\text{ and}%
\end{align*}%
\[
\left(  e^{t\bar{\Delta}/2}f\right)  (x)=\int_{M}p_{t}(x,y)f(y)dV(y)\text{ for
all }f\in L^{2}(M).
\]
If $f\in C_{c}^{\infty}(M),$ the function $u\left(  t,x\right)  :=e^{t\bar
{\Delta}/2}f\left(  x\right)  $ is smooth for $t>0$ and $x\in M$ and
$Le^{t\bar{\Delta}/2}f\left(  x\right)  $ is continuous for $t\geq0$ and $x\in
M$ for any smooth linear differential operator $L$ on $C^{\infty}\left(
M\right)  .$ For these results, see for example Strichartz \cite{Strichartz83}%
, Dodziuk \cite{Dodziuk83} and Davies \cite{Davies90}.
\end{fact}

\begin{theorem}
[Stochastic Rolling Constructions]\label{t.5.33}Assume $M$ is compact and let
$\Sigma,$ $u_{s}=//_{s},$ and $b$ be as in Theorem \ref{t.5.29}, then:

\begin{enumerate}
\item \ $\Sigma$ is a martingale iff $b$ is a $T_{o}M$ -- valued martingale.

\item \ $\Sigma$ is a Brownian motion iff $b$ is a $T_{o}M$ -- valued Brownian motion.
\end{enumerate}

Furthermore if $\Sigma$ is a Brownian motion, $T\in\left(  0,\infty\right)  $
and $f\in C^{\infty}(M),$ then
\[
M_{s}:=\left(  e^{(T-s)\bar{\Delta}/2}f\right)  \left(  \Sigma_{s}\right)
\]
is a martingale for $s\in\lbrack0,T]$ and%
\begin{equation}
dM_{s}=\left(  de^{(T-s)\bar{\Delta}/2}f\right)  (u_{s}db_{s})_{\Sigma_{s}%
}=\left(  de^{(T-s)\bar{\Delta}/2}f\right)  (//_{s}db_{s}). \label{e.5.38}%
\end{equation}

\end{theorem}

\begin{proof}
Keep the same notation as in Proposition \ref{p.5.31} and let $f\in C^{\infty
}(M).$ By Proposition \ref{p.5.31}, if $b$ is a martingale, then $\int
_{0}^{\cdot}df({\bar{d}}\Sigma)=\int_{0}^{\cdot}df(udb)$ is also a martingale
and hence $\Sigma$ is a martingale. Combining this with Corollary \ref{c.5.18}
and Proposition \ref{p.5.31},
\begin{align*}
d[f(\Sigma)]  &  =df({\bar{d}}\Sigma)+{\frac{1}{2}}\nabla df(d\Sigma\otimes
d\Sigma)\\
&  =df(udb)+{\frac{1}{2}}\nabla df(udb\otimes udb).
\end{align*}
Since $u$ is an isometry and $b$ is a Brownian motion, $udb\otimes
udb=\mathcal{I}(\Sigma)d\lambda.$ Hence
\[
d[f(\Sigma)]=df(udb)+{\frac{1}{2}}\Delta f(\Sigma)d\lambda
\]
from which it follows that $\Sigma$ is a Brownian motion.

Conversely, if $\Sigma$ is a $M$ -- valued martingale, then
\begin{equation}
N:=\sum_{i=1}^{N}\int_{0}^{\cdot}dx^{i}({\bar{d}}\Sigma)e_{i}=\sum_{i=1}%
^{N}\int_{0}^{\cdot}\langle e_{i},udb\rangle e_{i}=\int_{0}^{\cdot}udb
\label{e.5.39}%
\end{equation}
is a martingale, where $x=(x^{1},\ldots,x^{N})$ are standard coordinates on
$\mathbb{R}^{N}$ and $\{e_{i}\}_{i=1}^{N}$ is the standard basis for
$\mathbb{R}^{N}.$ From Eq. (\ref{e.5.39}), it follows that $b=\int_{0}^{\cdot
}u^{-1}dN$ is also a martingale.

Now suppose that $\Sigma$ is an $M$ -- valued Brownian motion, then we have
already proved that $b$ is a martingale. To finish the proof it suffices by
L\'{e}vy's criteria (Lemma \ref{l.5.21}) to show that $db\otimes
db=\mathcal{I}(o)d\lambda.$ But $\Sigma=N+(\text{\textrm{bounded variation}})$
and hence
\begin{align*}
db\otimes db  &  =u^{-1}d\Sigma\otimes u^{-1}d\Sigma=u^{-1}dN\otimes
u^{-1}dN\\
&  =(u^{_{-1}}\otimes u^{-1})(d\Sigma\otimes d\Sigma)\\
&  =(u^{_{-1}}\otimes u^{-1})\mathcal{I}(\Sigma)d\lambda=\mathcal{I}%
(o)d\lambda,
\end{align*}
wherein Eq. (\ref{e.5.24}) was used in the fourth equality and the
orthogonality of $u$ was used in the last equality.\textrm{ }To prove Eq.
(\ref{e.5.38}), let $M_{s}=u\left(  s,\Sigma_{s}\right)  $ where $u\left(
s,x\right)  :=\left(  e^{(T-s)\bar{\Delta}/2}f\right)  \left(  x\right)  $
which satisfies
\[
\partial_{s}u\left(  s,x\right)  +\frac{1}{2}\Delta u\left(  s,x\right)
=0\text{ with }u\left(  T,x\right)  =f\left(  x\right)
\]
By It\^{o}'s Lemma (see Corollary \ref{c.5.18}) along with Lemma \ref{l.5.21}
and Proposition \ref{p.5.31},%
\begin{align*}
dM_{s}  &  =\partial_{s}u\left(  s,\Sigma_{s}\right)  ds+d_{M}\left[  u\left(
s,\cdot\right)  \right]  ({\bar{d}}\Sigma_{s})+{\frac{1}{2}}\nabla
d_{M}\left[  u\left(  s,\cdot\right)  \right]  (d\Sigma_{s}\otimes d\Sigma
_{s})\\
&  =\partial_{s}u\left(  s,\Sigma_{s}\right)  ds+\frac{1}{2}\Delta u\left(
s,\Sigma_{s}\right)  ds+\left(  d_{M}e^{(T-s)\bar{\Delta}/2}f\right)  \left(
(u_{s}db_{s})_{\Sigma_{s}}\right) \\
&  =\left(  d_{M}e^{(T-s)\bar{\Delta}/2}f\right)  \left(  (u_{s}%
db_{s})_{\Sigma_{s}}\right)  .
\end{align*}

\end{proof}

The rolling construction of Brownian motion seems to have first been
discovered by Eells and Elworthy \cite{EE1971} who used ideas of Gangolli
\cite{Gangolli1964}. The relationship of the stochastic development map to
stochastic differential equations on the orthogonal frame bundle $O(M)$ of $M$
is pointed out in Elworthy \cite{Elworthy74,Elworthy75,Elworthy78}. The frame
bundle point of view has also been extensively developed by Malliavin, see for
example \cite{Malliavin78,Malliavin78b,Malliavin79}. For a more detailed
history of the stochastic development map, see pp. 156--157 in Elworthy
\cite{Elworthy78}. The reader may also wish to consult
\cite{Emery,Ikeda81,Kunita90,Malliavin97,Stroock2001,Hsu2003c}.

\begin{corollary}
\label{c.5.34}If $\Sigma$ is a Brownian motion on $M,$%
\[
\pi=\left\{  0=s_{0}<s_{1}<\dots<s_{n}=T\right\}
\]
is a partition of $[0,T]$ and $f\in C^{\infty}\left(  M^{n}\right)  ,$ then%
\begin{equation}
\mathbb{E}f\left(  \Sigma_{s_{1}},\dots,\Sigma_{s_{n}}\right)  =\int_{M^{n}%
}f\left(  x_{1},x_{2},\dots,x_{n}\right)  \prod_{i=1}^{n}p_{\Delta_{i}%
s}\left(  x_{i-1},x_{i}\right)  d\lambda\left(  x_{i}\right)  \label{e.5.40}%
\end{equation}
where $\Delta_{i}s:=s_{i}-s_{i-1},$ $x_{0}:=o$ and $\lambda:=\lambda_{M}.$ In
particular $\Sigma$ is a Markov process relative to the filtration, $\left\{
\mathcal{F}_{s}\right\}  $ where $\mathcal{F}_{s}$ is the $\sigma$ -- algebra
generated by $\left\{  \Sigma_{\tau}:\tau\leq s\right\}  .$
\end{corollary}

\begin{proof}
By standard measure theoretic arguments, it suffices to prove Eq.
(\ref{e.5.40}) when $f$ is a product function of the form $f\left(
x_{1},x_{2},\dots,x_{n}\right)  =\prod_{i=1}^{n}f_{i}\left(  x_{i}\right)  $
with $f_{i}\in C^{\infty}(M).$ By Theorem \ref{t.5.33}, $M_{s}:=e^{\left(
T-s\right)  \bar{\Delta}/2}f_{n}\left(  \Sigma_{s}\right)  $ is a martingale
for $s\leq T$ and therefore%
\begin{align}
\mathbb{E}\left[  f\left(  \Sigma_{s_{1}},\dots,\Sigma_{s_{n}}\right)
\right]   &  =\mathbb{E}\left[  \prod_{i=1}^{n-1}f_{i}\left(  \Sigma_{s_{i}%
}\right)  \cdot M_{T}\right]  =\mathbb{E}\left[  \prod_{i=1}^{n-1}f_{i}\left(
\Sigma_{s_{i}}\right)  \cdot M_{s_{n-1}}\right] \nonumber\\
&  =\mathbb{E}\left[  \prod_{i=1}^{n-1}f_{i}\left(  \Sigma_{s_{i}}\right)
\cdot\left(  P_{\Delta_{n}s}f_{n}\right)  \left(  \Sigma_{s_{n-1}}\right)
\right]  . \label{e.5.41}%
\end{align}
In particular if $n=1,$ it follows that
\[
\mathbb{E}\left[  f_{1}\left(  \Sigma_{T}\right)  \right]  =\mathbb{E}\left[
\left(  e^{T\bar{\Delta}/2}f_{1}\right)  \left(  \Sigma_{0}\right)  \right]
=\int_{M}p_{T}(o,x_{1})f_{1}\left(  x_{1}\right)  d\lambda\left(
x_{1}\right)  .
\]
Now assume we have proved Eq. (\ref{e.5.40}) with $n$ replaced by $n-1$ and to
simplify notation let $g\left(  x_{1},x_{2},\dots,x_{n-1}\right)
:=\prod_{i=1}^{n-1}f_{i}\left(  x_{i}\right)  .$ It would then follow from Eq.
(\ref{e.5.41}) that
\begin{align*}
\mathbb{E}  &  \left[  f\left(  \Sigma_{s_{1}},\dots,\Sigma_{s_{n}}\right)
\right] \\
&  =\int_{M^{n-1}}g\left(  x_{1},x_{2},\dots,x_{n-1}\right)  \left(
e^{\frac{s_{n}-s_{n-1}}{2}\bar{\Delta}}f_{n}\right)  \left(  x_{n-1}\right)
\prod_{i=1}^{n-1}p_{\Delta_{i}s}\left(  x_{i-1},x_{i}\right)  d\lambda\left(
x_{i}\right) \\
&  =\int_{M^{n-1}}g\left(  x_{1},x_{2},\dots,x_{n-1}\right)  \left[  \int
_{M}f_{n}\left(  x_{n}\right)  p_{\Delta_{n}s}\left(  x_{n-1},x_{n}\right)
d\lambda\left(  x_{n}\right)  \right]  \times\\
&  \qquad\qquad\qquad\qquad\qquad\qquad\qquad\qquad\qquad\times\prod
_{i=1}^{n-1}p_{\Delta_{i}s}\left(  x_{i-1},x_{i}\right)  d\lambda\left(
x_{i}\right) \\
&  =\int_{M^{n}}f\left(  x_{1},x_{2},\dots,x_{n}\right)  \prod_{i=1}%
^{n}p_{\Delta_{i}s}\left(  x_{i-1},x_{i}\right)  d\lambda\left(  x_{i}\right)
.
\end{align*}
This completes the induction step and hence also the proof of the theorem.
\end{proof}

\subsection{More Constructions of Semi-Martingales and Brownian
Motions\label{s.5.5}}

Let $\Gamma$ be the one form on $M$ with values in the skew symmetric $N\times
N$ matrices defined by $\Gamma=dQP+dPQ$ as in Eq. \text{(\ref{e.3.65})}. Given
an $M-$valued semi-martingale $\Sigma,$ let $u$ denote parallel translation
along $\Sigma$ as defined in Eq. \text{(\ref{e.5.32})}\ of Theorem
\ref{t.5.27}.

\begin{lemma}
[Orthogonality Lemma]\label{l.5.35}Suppose that $B$ is an $\mathbb{R}^{N}$ --
valued semi-martingale and $\Sigma$ is the solution to
\begin{equation}
\delta\Sigma=P(\Sigma)\delta B\quad\text{\textrm{\ with }}\quad\Sigma_{0}=o\in
M. \label{e.5.42}%
\end{equation}
Let $\{e_{i}\}_{i=1}^{N}$ be any orthonormal basis for $\mathbb{R}^{N}$ and
define $B^{i}:=\langle e_{i},B\rangle$ then%
\[
P(\Sigma)dB\otimes Q(\Sigma)dB:=\sum_{i,j=1}^{N}P(\Sigma)e_{i}\otimes
Q(\Sigma)e_{j}\left(  dB^{i}dB^{j}\right)  =0.
\]

\end{lemma}

\begin{proof}
Suppose $\left\{  v_{i}\right\}  _{i=1}^{N}$ is another orthonormal basis for
$\mathbb{R}^{N}.$ Using the bilinearity of the joint quadratic variation,%
\begin{align*}
\lbrack\langle e_{i},B\rangle,\langle e_{j},B\rangle]  &  =\sum_{k,l}[\langle
e_{i},v_{k}\rangle\langle v_{k},B\rangle,\langle e_{j},v_{l}\rangle\langle
v_{l},B\rangle]\\
&  =\sum_{k,l}\langle e_{i},v_{k}\rangle\langle e_{j},v_{l}\rangle
\lbrack\langle v_{k},B\rangle,\langle v_{l},B\rangle].
\end{align*}
Therefore,%
\begin{align*}
\sum_{i,j=1}^{N}  &  P(\Sigma)e_{i}\otimes Q(\Sigma)e_{j}\cdot d\left[
B^{i},B^{j}\right] \\
&  =\sum_{i,j,k,l=1}^{N}\left[  P(\Sigma)e_{i}\otimes Q(\Sigma)e_{j}\right]
\langle e_{i},v_{k}\rangle\langle e_{j},v_{l}\rangle d[\langle v_{k}%
,B\rangle,\langle v_{l},B\rangle]\\
&  =\sum_{k,l=1}^{N}\left[  P(\Sigma)v_{k}\otimes Q(\Sigma)v_{l}\right]
d[\langle v_{k},B\rangle,\langle v_{l},B\rangle]
\end{align*}
which shows $P(\Sigma)dB\otimes Q(\Sigma)dB$ is well defined.

Now define%
\[
\tilde{B}:=\int_{0}^{\cdot}u^{-1}dB\text{ and }\tilde{B}^{i}:=\langle
e_{i},\tilde{B}\rangle=\int_{0}^{\cdot}\langle ue_{i},dB\rangle
\]
where $u$ is parallel translation along $\Sigma$ in $M\times\mathbb{R}^{N}$ as
defined in Eq. \text{(\ref{e.5.32}). }Then%
\begin{align*}
P(\Sigma)dB\otimes Q(\Sigma)dB  &  =\sum_{i,j,k,l=1}^{N}P(\Sigma)ue_{k}\otimes
Q(\Sigma)ue_{l}\langle e_{i},ue_{k}\rangle\langle e_{j},ue_{l}\rangle\left(
dB^{i}dB^{j}\right) \\
&  =\sum_{k,l=1}^{N}P(\Sigma)ue_{k}\otimes Q(\Sigma)ue_{l}\left(  d\tilde
{B}^{k}d\tilde{B}^{l}\right) \\
&  =\sum_{k,l=1}^{N}uP(o)e_{k}\otimes uQ(o)e_{l}\left(  d\tilde{B}^{k}%
d\tilde{B}^{l}\right)
\end{align*}
wherein we have used $P(\Sigma)u=uP(o)$ and $Q(\Sigma)u=uQ(o),$ see Theorem
\ref{t.5.27}. This last expression is easily seen to be zero by choosing
$\{e_{i}\}$ such that $P(o)e_{i}=e_{i}$ for $i=1,2,\ldots,d$ and $Q\left(
o\right)  e_{j}=e_{j}$ for $j=d+1,\dots,N.$
\end{proof}

The next proposition is a stochastic analogue of Lemma \ref{l.3.55} and the
proof is very similar to that of Lemma \ref{l.3.55}.

\begin{proposition}
\label{p.5.36}Suppose that $V$ is a $TM$ -- valued semi-martingale,
$\Sigma=\pi\left(  V\right)  $ so that $\Sigma$ is an $M$ -- valued
semi-martingale and $V_{s}\in T_{\Sigma_{s}}M$ for all $s\geq0.$ Then
\begin{equation}
//_{s}\delta_{s}\left[  //_{s}^{-1}V_{s}\right]  =\delta_{s}^{\nabla}%
V_{s}=:P\left(  \Sigma_{s}\right)  \delta V_{s} \label{e.5.43}%
\end{equation}
where $//_{s}$ is stochastic parallel translation along $\Sigma.$ If $Y_{s}%
\in\Gamma\left(  TM\right)  $ is a time dependent vector field, then%
\begin{equation}
\delta_{s}\left[  //_{s}^{-1}Y_{s}\left(  \Sigma_{s}\right)  \right]
=//_{s}^{-1}\left(  \frac{d}{ds}Y_{s}\right)  \left(  \Sigma_{s}\right)
ds+//_{s}^{-1}\nabla_{\delta\Sigma_{s}}Y_{s} \label{e.5.44}%
\end{equation}
and for $w\in T_{o}M,$%
\begin{align}
//_{s}^{-1}\delta_{s}^{\nabla}\left[  \nabla_{//_{s}w}Y_{s}\right]   &
=\delta_{s}\left[  //_{s}^{-1}\nabla_{//_{s}w}Y_{s}\right] \nonumber\\
&  =//_{s}^{-1}\nabla_{\delta\Sigma_{s}\otimes//_{s}w}^{2}Y_{s}+//_{s}%
^{-1}\left[  \nabla_{//_{s}w}\left(  \frac{d}{ds}Y_{s}\right)  \right]  .
\label{e.5.45}%
\end{align}
Furthermore if $\Sigma_{s}$ is a Brownian motion, then%
\begin{align}
d\left[  //_{s}^{-1}Y_{s}\left(  \Sigma_{s}\right)  \right]  =  &  //_{s}%
^{-1}\nabla_{//_{s}db_{s}}Y_{s}+//_{s}^{-1}\left(  \frac{d}{ds}Y_{s}\right)
\left(  \Sigma_{s}\right)  ds\nonumber\\
&  \qquad+\frac{1}{2}\sum_{i=1}^{d}//_{s}^{-1}\nabla_{//_{s}e_{i}\otimes
//_{s}e_{i}}^{2}Y_{s}ds \label{e.5.46}%
\end{align}
where $\left\{  e_{i}\right\}  _{i=1}^{d}$ is an orthonormal basis for
$T_{o}M.$
\end{proposition}

\begin{proof}
We will use the convention of summing on repeated indices and write $u_{s}$
for $//_{s},$ i.e. stochastic parallel translation along $\Sigma$ on $TM.$
Recall that $u_{s}$ solves%
\[
\delta u_{s}+dQ\left(  \delta\Sigma_{s}\right)  u_{s}=0\text{ with }%
u_{0}=I_{T_{o}M}.
\]
Define $\bar{u}_{s}$ as the solution to:%
\[
\delta\bar{u}_{s}=\bar{u}_{s}dQ\left(  \delta\Sigma_{s}\right)  \text{ with
}\bar{u}_{0}=I_{T_{o}M}.
\]
Then%
\[
\delta\left(  \bar{u}_{s}u_{s}\right)  =-\bar{u}_{s}dQ\left(  \delta\Sigma
_{s}\right)  u_{s}+\bar{u}_{s}dQ\left(  \delta\Sigma_{s}\right)  u_{s}=0
\]
from which it follows that $\bar{u}_{s}u_{s}=I$ for all $s$ and hence $\bar
{u}_{s}=u_{s}^{-1}.$
%(BRUCE: what about domains and ranges. Should this be done in Lemma <ref>l.z3.50</ref>.)
This proves Eq. (\ref{e.5.43}) since%
\begin{align*}
u_{s}\delta_{s}\left[  u_{s}^{-1}V_{s}\right]   &  =u_{s}\left[  u_{s}%
^{-1}dQ\left(  \delta\Sigma_{s}\right)  V_{s}+u_{s}^{-1}\delta V_{s}\right] \\
&  =dQ\left(  \delta\Sigma_{s}\right)  V_{s}+\delta V_{s}=\delta^{\nabla}%
V_{s},
\end{align*}
where the last equality comes from Eq. (\ref{e.5.29}).

Applying Eq. (\ref{e.5.43}) to $V_{s}:=Y_{s}\left(  \Sigma_{s}\right)  $
gives
\begin{align*}
\delta_{s}\left[  //_{s}^{-1}Y_{s}\left(  \Sigma_{s}\right)  \right]   &
=//_{s}^{-1}P\left(  \Sigma_{s}\right)  \delta_{s}\left[  Y_{s}\left(
\Sigma_{s}\right)  \right] \\
&  =//_{s}^{-1}P\left(  \Sigma_{s}\right)  \left(  \frac{d}{ds}Y_{s}\right)
\left(  \Sigma_{s}\right)  ds+//_{s}^{-1}P\left(  \Sigma_{s}\right)
Y_{s}^{\prime}\left(  \Sigma_{s}\right)  \delta_{s}\Sigma_{s}\\
&  =//_{s}^{-1}\left(  \frac{d}{ds}Y_{s}\right)  \left(  \Sigma_{s}\right)
ds+//_{s}^{-1}\nabla_{\delta_{s}\Sigma_{s}}Y_{s},
\end{align*}
which proves Eq. (\ref{e.5.44}).

To prove Eq. (\ref{e.5.45}), let $X_{i}\left(  m\right)  =P\left(  m\right)
e_{i}$ for $i=1,2,\dots,N.$ By Proposition \ref{p.3.48},
\begin{align}
\nabla_{//_{s}w}Y_{s}  &  =\langle//_{s}w,X_{i}\left(  \Sigma_{s}\right)
\rangle\left(  \nabla_{X_{i}}Y_{s}\right)  \left(  \Sigma_{s}\right)
\label{e.5.47}\\
&  =\langle w,//_{s}^{-1}X_{i}\left(  \Sigma_{s}\right)  \rangle\left(
\nabla_{X_{i}}Y_{s}\right)  \left(  \Sigma_{s}\right) \nonumber
\end{align}
and%
\[
//_{s}w=\langle//_{s}w,X_{i}\left(  \Sigma_{s}\right)  \rangle X_{i}\left(
\Sigma_{s}\right)  =\langle w,//_{s}^{-1}X_{i}\left(  \Sigma_{s}\right)
\rangle X_{i}\left(  \Sigma_{s}\right)
\]
or equivalently,%
\begin{equation}
w=\langle w,//_{s}^{-1}X_{i}\left(  \Sigma_{s}\right)  \rangle//_{s}^{-1}%
X_{i}\left(  \Sigma_{s}\right)  . \label{e.5.48}%
\end{equation}
Taking the covariant differential of Eq. (\ref{e.5.47}), making use of Eq.
(\ref{e.5.44}), gives%
\begin{align}
\delta_{s}^{\nabla}  &  \left[  \nabla_{//_{s}w}Y_{s}\right] \nonumber\\
&  =\langle//_{s}w,\nabla_{\delta_{s}\Sigma_{s}}X_{i}\rangle\left(
\nabla_{X_{i}}Y_{s}\right)  \left(  \Sigma_{s}\right)  +\langle//_{s}%
w,X_{i}\left(  \Sigma_{s}\right)  \rangle\nabla_{\delta_{s}\Sigma_{s}}%
\nabla_{X_{i}}Y_{s}\nonumber\\
&  \qquad+\langle//_{s}w,X_{i}\left(  \Sigma_{s}\right)  \rangle\left(
\nabla_{X_{i}}\left(  \frac{d}{ds}Y_{s}\right)  \right)  \left(  \Sigma
_{s}\right) \nonumber\\
&  =\langle//_{s}w,\nabla_{\delta_{s}\Sigma_{s}}X_{i}\rangle\left(
\nabla_{X_{i}}Y_{s}\right)  \left(  \Sigma_{s}\right)  +\langle//_{s}%
w,X_{i}\left(  \Sigma_{s}\right)  \rangle\nabla_{\delta_{s}\Sigma_{s}\otimes
X_{i}}^{2}Y_{s}\nonumber\\
&  \qquad+\langle//_{s}w,X_{i}\left(  \Sigma_{s}\right)  \rangle\nabla
_{\nabla_{\delta_{s}\Sigma_{s}}X_{i}}Y_{s}+\left(  \nabla_{//_{s}w}\left(
\frac{d}{ds}Y_{s}\right)  \right)  \left(  \Sigma_{s}\right) \nonumber\\
&  =\left(  \nabla_{\langle//_{s}w,\nabla_{\delta_{s}\Sigma_{s}}X_{i}\rangle
X_{i}\left(  \Sigma_{s}\right)  +\langle//_{s}w,X_{i}\left(  \Sigma
_{s}\right)  \rangle\nabla_{\delta_{s}\Sigma_{s}}X_{i}}Y_{s}\right)  \left(
\Sigma_{s}\right) \nonumber\\
&  \qquad+\nabla_{\delta_{s}\Sigma_{s}\otimes//_{s}w}^{2}Y_{s}+\left(
\nabla_{//_{s}w}\left(  \frac{d}{ds}Y_{s}\right)  \right)  \left(  \Sigma
_{s}\right)  , \label{e.5.49}%
\end{align}
Taking the differential of Eq. (\ref{e.5.48}) implies%
\[
0=\delta v=\langle v,//_{s}^{-1}\nabla_{\delta_{s}\Sigma_{s}}X_{i}%
\rangle//_{s}^{-1}X_{i}\left(  \Sigma_{s}\right)  +\langle v,//_{s}^{-1}%
X_{i}\left(  \Sigma_{s}\right)  \rangle//_{s}^{-1}\nabla_{\delta_{s}\Sigma
_{s}}X_{i}%
\]
which upon multiplying by $//_{s}$ shows%
\[
\langle//_{s}v,\nabla_{\delta_{s}\Sigma_{s}}X_{i}\rangle X_{i}\left(
\Sigma_{s}\right)  +\langle//_{s}v,X_{i}\left(  \Sigma_{s}\right)
\rangle\nabla_{\delta_{s}\Sigma_{s}}X_{i}=0.
\]
Using this identity in Eq. (\ref{e.5.49}) completes the proof of Eq.
(\ref{e.5.45}).

Now suppose that $\Sigma_{s}$ is a Brownian motion and $b_{s}=\Psi_{s}\left(
\Sigma\right)  $ is the anti-developed $T_{o}M$ -- valued Brownian motion
associated to $\Sigma.$ Then by Eq. (\ref{e.5.44}),%
\begin{align*}
d\left[  //_{s}^{-1}Y_{s}\left(  \Sigma_{s}\right)  \right]   &  =//_{s}%
^{-1}\left(  \frac{d}{ds}Y_{s}\right)  \left(  \Sigma_{s}\right)
ds+//_{s}^{-1}\nabla_{//_{s}\delta b_{s}}Y_{s}\\
&  =//_{s}^{-1}\left(  \frac{d}{ds}Y_{s}\right)  \left(  \Sigma_{s}\right)
ds+\left(  //_{s}^{-1}\nabla_{//_{s}e_{i}}Y_{s}\right)  \delta b_{s}^{i}.
\end{align*}
Using Eq. (\ref{e.5.45}),%
\begin{align*}
\left(  //_{s}^{-1}\nabla_{//_{s}e_{i}}Y_{s}\right)  \delta b_{s}^{i}  &
=\left(  //_{s}^{-1}\nabla_{//_{s}e_{i}}Y_{s}\right)  db_{s}^{i}+\frac{1}%
{2}d\left(  //_{s}^{-1}\nabla_{//_{s}e_{i}}Y_{s}\right)  db_{s}^{i}\\
&  =//_{s}^{-1}\nabla_{//_{s}db_{s}}Y_{s}+\frac{1}{2}//_{s}^{-1}\nabla
_{\delta\Sigma_{s}\otimes//_{s}e_{i}}^{2}Y_{s}db_{s}^{i}\\
&  =//_{s}^{-1}\nabla_{//_{s}db_{s}}Y_{s}+\frac{1}{2}//_{s}^{-1}\nabla
_{//_{s}e_{j}\otimes//_{s}e_{i}}^{2}Y_{s}db_{s}^{i}db_{s}^{j}\\
&  =//_{s}^{-1}\nabla_{//_{s}db_{s}}Y_{s}+\frac{1}{2}//_{s}^{-1}\nabla
_{//_{s}e_{i}\otimes//_{s}e_{i}}^{2}Y_{s}ds.
\end{align*}
Combining the last two equations proves Eq. (\ref{e.5.46}).
\end{proof}

\begin{theorem}
\label{t.5.37}Let $\Sigma_{s}$ denote the solution to Eq. (\ref{e.5.1}) with
$\Sigma_{0}=o\in M$ and $b_{s}=\Psi_{s}\left(  \Sigma\right)  \in T_{o}M.$
Then%
\begin{align}
b_{s}=  &  \int_{0}^{s}//_{r}^{-1}\left(  \Sigma\right)  \left[
\mathbf{X}\left(  \Sigma_{r}\right)  \delta B_{r}+X_{0}\left(  \Sigma
_{r}\right)  dr\right] \nonumber\\
=  &  \int_{0}^{s}//_{r}^{-1}\left(  \Sigma\right)  \mathbf{X}\left(
\Sigma_{r}\right)  dB_{r}\nonumber\\
&  \quad+\int_{0}^{s}//_{r}^{-1}\left[  \frac{1}{2}\sum_{i,j=1}^{n}\left(
\nabla_{X_{i}}X_{j}\right)  \left(  \Sigma_{r}\right)  dB_{r}^{i}dB_{r}%
^{j}+X_{0}\left(  \Sigma_{r}\right)  dr\right]  . \label{e.5.50}%
\end{align}
Hence if $B$ is a Brownian motion, then%
\begin{align}
b_{s}=  &  \int_{0}^{s}//_{r}^{-1}\left(  \Sigma\right)  \mathbf{X}\left(
\Sigma_{r}\right)  dB_{r}\nonumber\\
&  \quad+\int_{0}^{s}//_{r}^{-1}\left[  \frac{1}{2}\sum_{i=1}^{n}\left(
\nabla_{X_{i}}X_{i}\right)  \left(  \Sigma_{r}\right)  +X_{0}\left(
\Sigma_{r}\right)  \right]  dr. \label{e.5.51}%
\end{align}

\end{theorem}

\begin{proof}
By the definition of $b,$%
\begin{align*}
db_{s}  &  =//_{s}^{-1}\left(  \Sigma\right)  \left[  \mathbf{X}\left(
\Sigma_{s}\right)  \delta B_{s}+X_{0}\left(  \Sigma_{s}\right)  ds\right] \\
&  =//_{s}^{-1}\left(  \Sigma\right)  \left[  \mathbf{X}\left(  \Sigma
_{s}\right)  dB_{s}+X_{0}\left(  \Sigma_{s}\right)  ds\right]  +\frac{1}%
{2}d\left[  //_{s}^{-1}\left(  \Sigma\right)  \mathbf{X}\left(  \Sigma
_{s}\right)  \right]  dB_{s}\\
&  =//_{s}^{-1}\left(  \Sigma\right)  \left[  \mathbf{X}\left(  \Sigma
_{s}\right)  dB_{s}+X_{0}\left(  \Sigma_{s}\right)  ds\right]  +\frac{1}%
{2}\left[  //_{s}^{-1}\left(  \Sigma\right)  \nabla_{\mathbf{X}\left(
\Sigma_{s}\right)  dB_{s}}\mathbf{X}\right]  dB_{s}\\
&  =//_{s}^{-1}\left(  \Sigma\right)  \left[  \mathbf{X}\left(  \Sigma
_{s}\right)  dB_{s}+ds\right]  +\frac{1}{2}//_{s}^{-1}\left(  \Sigma\right)
\sum_{i,j=1}^{n}\left(  \nabla_{X_{i}}X_{j}\right)  \left(  \Sigma_{s}\right)
dB_{s}^{i}dB_{s}^{j}%
\end{align*}
which combined with the identity,
\begin{align*}
d\left[  //_{s}^{-1}\left(  \Sigma\right)  \mathbf{X}\left(  \Sigma
_{s}\right)  \right]  dB_{s}  &  =\left[  //_{s}^{-1}\left(  \Sigma\right)
\nabla_{d\Sigma_{s}}\mathbf{X}\right]  dB_{s}=\left[  //_{s}^{-1}\left(
\Sigma\right)  \nabla_{\mathbf{X}\left(  \Sigma_{s}\right)  dB_{s}}%
\mathbf{X}\right]  dB_{s}\\
&  =\sum_{i,j=1}^{n}\left(  \nabla_{X_{i}}X_{j}\right)  \left(  \Sigma
_{s}\right)  dB_{s}^{i}dB_{s}^{j}%
\end{align*}
proves Eq. (\ref{e.5.50}).
\end{proof}

\begin{corollary}
\label{c.5.38}Suppose $B_{s}$ is an $\mathbb{R}^{n}$ -- valued Brownian
motion, $\Sigma_{s}$ is the solution to Eq. (\ref{e.5.1}) with $\beta=B$ and
$\frac{1}{2}\sum_{k=1}^{n}\left(  \nabla_{X_{k}}X_{k}\right)  +X_{0}=0,$ then
$\Sigma$ is an $M$ -- valued martingale with quadratic variation,%
\begin{equation}
d\Sigma_{s}\otimes d\Sigma_{s}=\sum_{k=1}^{n}X_{k}\left(  \Sigma_{s}\right)
\otimes X_{k}\left(  \Sigma_{s}\right)  ds. \label{e.5.52}%
\end{equation}

\end{corollary}

\begin{proof}
By Eq. (\ref{e.5.51}) and Theorem \ref{t.5.33}, $\Sigma$ is a martingale and
from Eq. (\ref{e.5.1}),
\[
d\Sigma^{i}d\Sigma^{j}=\sum_{k,l=1}^{n}X_{k}^{i}\left(  \Sigma\right)
X_{l}^{j}\left(  \Sigma\right)  dB^{k}dB^{l}=\sum_{k=1}^{n}X_{k}^{i}\left(
\Sigma\right)  X_{k}^{j}\left(  \Sigma\right)  ds
\]
where $\left\{  e_{i}\right\}  _{i=1}^{N}$ is the standard basis for
$\mathbb{R}^{N},$ $\Sigma^{i}:=\langle\Sigma,e_{i}\rangle$ and $X_{k}%
^{i}\left(  \Sigma\right)  =\langle X_{k}\left(  \Sigma\right)  ,e_{i}%
\rangle.$ Using this identity in Eq. (\ref{e.5.17}), shows%
\[
d\Sigma_{s}\otimes d\Sigma_{s}=\sum_{i,j=1}^{N}\sum_{k=1}^{n}e_{i}\otimes
e_{j}X_{k}^{i}\left(  \Sigma\right)  X_{k}^{j}\left(  \Sigma\right)
ds=\sum_{k=1}^{n}X_{k}\left(  \Sigma_{s}\right)  \otimes X_{k}\left(
\Sigma_{s}\right)  ds.
\]

\end{proof}

\begin{corollary}
\label{c.5.39}Suppose now that $B_{s}$ is an $\mathbb{R}^{N}$ -- valued
semi-martingale and $\Sigma_{s}$ is the solution to Eq. (\ref{e.5.42}) in
Lemma \ref{l.5.35}. If $B$ is a martingale, then $\Sigma$ is a martingale and
if $B$ is a Brownian motion, then $\Sigma$ is a Brownian motion.
\end{corollary}

\begin{proof}
Solving Eq. (\ref{e.5.42}) is the same as solving Eq. (\ref{e.5.1}) with
$n=N,$ $\beta=B,$ $X_{0}\equiv0$ and $X_{i}\left(  m\right)  =P\left(
m\right)  e_{i}$ for all $i=1,2,\dots,N.$ Since
\[
\nabla_{X_{i}}X_{j}=PdP\left(  X_{i}\right)  e_{j}=dP\left(  X_{i}\right)
Qe_{j}=dP\left(  Pe_{i}\right)  Qe_{j},
\]
it follows from orthogonality Lemma \ref{l.5.35} that%
\[
\sum_{i,j=1}^{n}\left(  \nabla_{X_{i}}X_{j}\right)  \left(  \Sigma_{r}\right)
dB_{r}^{i}dB_{r}^{j}=0.
\]
Therefore from Eq. (\ref{e.5.50}), $b_{s}:=\int_{0}^{s}//_{r}^{-1}\delta
\Sigma_{r}$ is a $T_{o}M$ -- martingale which is equivalent to $\Sigma_{s}$
being a $M$ -- valued martingale. Finally if $B$ is a Brownian motion, then
from Eq. (\ref{e.5.52}), $\Sigma$ has quadratic variation given by%
\begin{equation}
d\Sigma_{s}\otimes d\Sigma_{s}=\sum_{i=1}^{N}P\left(  \Sigma_{s}\right)
e_{i}\otimes P\left(  \Sigma_{s}\right)  e_{i}ds \label{e.5.53}%
\end{equation}
Since $\sum_{i=1}^{N}P(m)e_{i}\otimes P(m)e_{i}$ is independent of the choice
of orthonormal basis for $\mathbb{R}^{N},$ we may choose $\{e_{i}\}$ such that
$\{e_{i}\}_{i=1}^{d}$ is an orthonormal basis for $\tau_{m}M$ to learn%
\[
\sum_{i=1}^{N}P(m)e_{i}\otimes P(m)e_{i}=\mathcal{I}(m).
\]
Using this in Eq. (\ref{e.5.53}) we learn that $d\Sigma_{s}\otimes d\Sigma
_{s}=\mathcal{I}\left(  \Sigma_{s}\right)  ds$ and hence $\Sigma$ is a
Brownian motion on $M$ by the L\'{e}vy criteria, see Lemma \ref{l.5.21}.
\end{proof}

\begin{theorem}
\label{t.5.40}Let $B$ be any $\mathbb{R}^{N}$-valued semi-martingale, $\Sigma$
be the solution to Eq. (\ref{e.5.42}),%
\begin{equation}
b:=\int_{0}^{\cdot}u^{-1}\delta\Sigma=\int_{0}^{\cdot}u^{-1}P(\Sigma)\delta B
\label{e.5.54}%
\end{equation}
be the anti-development of $\Sigma$ and
\begin{equation}
\beta:=\int_{0}^{\cdot}u^{-1}Q(\Sigma)dB=Q(o)\int_{0}^{\cdot}u^{-1}dB
\label{e.5.55}%
\end{equation}
be the \textquotedblleft normal\textquotedblright\ process. Then
\begin{equation}
b=\int_{0}^{\cdot}u^{-1}P(\Sigma)dB=P(o)\int_{0}^{\cdot}u^{-1}dB,
\label{e.5.56}%
\end{equation}
i.e. the Fisk-Stratonovich integral may be replaced by the It\^{o} integral.
Moreover if $B$ is a standard $\mathbb{R}^{N}$ -- valued Brownian motion then
$(b,\beta)$ is also a standard $\mathbb{R}^{N}$ -- valued Brownian and the
processes, $b_{s},$ $\Sigma_{s}$ and $//_{s}$ are all independent of $\beta.$
\end{theorem}

\begin{proof}
Let $p=P(\Sigma)$ and $u$ be parallel translation on $M\times\mathbb{R}^{N}$
(see Eq. (\ref{e.5.32})), then
\begin{align*}
d(u^{-1}P(\Sigma))\cdot dB  &  =u^{-1}\left[  \Gamma(\delta\Sigma
)P(\Sigma)dB+dP(\delta\Sigma)dB\right] \\
&  =u^{-1}\left[  \left(  dQ(\delta\Sigma)P\left(  \Sigma\right)
+dP(\delta\Sigma)Q\left(  \Sigma\right)  \right)  P(\Sigma)dB+dP(\delta
\Sigma)dB\right] \\
&  =u^{-1}\left[  dQ(\delta\Sigma)P(\Sigma)dB-dQ(\delta\Sigma)dB\right] \\
&  =-u^{-1}dQ(\delta\Sigma)Q(\Sigma)dB=-u^{-1}dQ(P\left(  \Sigma\right)
dB)Q(\Sigma)dB=0
\end{align*}
where we have again used $P\left(  \Sigma\right)  dB\otimes Q\left(
\Sigma\right)  dB=0.$ This proves \text{(\ref{e.5.56})}.

Now suppose that $B$ is a Brownian motion. Since $(b,\beta)=\int_{0}^{\cdot
}u^{-1}dB$ and $u$ is an orthogonal process, it easily follow's using
L\'{e}vy's criteria that $(b,\beta)$ is a standard Brownian motion and in
particular, $\beta$ is independent of $b.$ Since $(\Sigma,u)$ satisfies the
coupled pair of stochastic differential equations
\begin{align*}
d\Sigma &  =u\delta b\text{ and }du+\Gamma(u\delta b)u=0\text{ with}\\
\Sigma_{0}  &  =o\text{ and }u_{0}=I\in\operatorname*{End}(\mathbb{R}^{N}),
\end{align*}
it follows that $(\Sigma,u)$ is a functional of $b$ and hence the process
$(\Sigma,u)$ are independent of $\beta.$
\end{proof}

\subsection{The Differential in the Starting Point of a Stochastic
Flow\label{d.5.6}}

In this section let $B_{s}$ be an $\mathbb{R}^{n}$ -- valued Brownian motion
and for each $m\in M$ let $T_{s}\left(  m\right)  =\Sigma_{s}$ where
$\Sigma_{s}$ is the solution to Eq. (\ref{e.5.1}) with $\Sigma_{0}=m.$ It is
well known, see Kunita \cite{Kunita90} that there is a version of
$T_{s}\left(  m\right)  $ which is continuous in $s$ and smooth in $m,$
moreover the differential of $T_{s}\left(  m\right)  $ relative to $m$ solves
the stochastic differential equation found by differentiating Eq.
\text{(}\ref{e.5.1}\text{)}. Let
\begin{equation}
Z_{s}:=T_{s\ast o}\text{ and }z_{s}:=//_{s}^{-1}Z_{s}\in\operatorname{End}%
\left(  T_{o}M\right)  \label{e.5.57}%
\end{equation}
where $//_{s}$ is stochastic parallel translation along $\Sigma_{s}%
:=T_{s}\left(  o\right)  .$

\begin{theorem}
\label{t.5.41}For all $v\in T_{o}M$%
\begin{equation}
\delta_{s}^{\nabla}Z_{s}v=\left(  \nabla_{Z_{s}v}\mathbf{X}\right)  \delta
B_{s}+\left(  \nabla_{Z_{s}v}X_{0}\right)  ds\text{ with }Z_{0}v=v.
\label{e.5.58}%
\end{equation}
Alternatively $z_{s}$ satisfies%
\begin{equation}
dz_{s}v=//_{s}^{-1}\left(  \nabla_{//_{s}z_{s}v}\mathbf{X}\right)  \delta
B_{s}+//_{s}^{-1}\left(  \nabla_{//_{s}z_{s}v}X_{0}\right)  ds. \label{e.5.59}%
\end{equation}

\end{theorem}

\begin{proof}
Equations (\ref{e.5.58}) and (\ref{e.5.59}) are the formal analogues Eqs.
(\ref{e.4.2}) and (\ref{e.4.3}) respectively. Because of Proposition
\ref{p.5.36}, Eq. (\ref{e.5.58}) is equivalent to Eq. (\ref{e.5.59}). To prove
Eq. (\ref{e.5.58}), differentiate Eq. (\ref{e.5.1}) in $m$ in the direction
$v\in T_{o}M$ to find%
\[
\delta_{s}Z_{s}v=DX_{i}\left(  \Sigma_{s}\right)  Z_{s}v\circ\delta B_{s}%
^{i}+DX_{0}\left(  \Sigma_{s}\right)  Z_{s}vds\text{ with }Z_{0}v=v.
\]
Multiplying this equation through by $P\left(  \Sigma_{s}\right)  $ on the
left then gives Eq. (\ref{e.5.58}).
\end{proof}

\begin{notation}
\label{n.5.42}The pull back, $\operatorname{Ric}_{//_{s}},$ of the Ricci
tensor by parallel translation is defined by
\begin{equation}
\operatorname{Ric}_{//_{s}}:=//_{s}^{-1}\operatorname*{Ric}\nolimits_{\Sigma
_{s}}//_{s}. \label{e.5.60}%
\end{equation}

\end{notation}

\begin{theorem}
[It\^{o} form of Eq. (\ref{e.5.59})]\label{t.5.43}The It\^{o} form of Eq.
(\ref{e.5.59}) is%
\begin{equation}
dz_{s}v=//_{s}^{-1}\left(  \nabla_{//_{s}z_{s}v}\mathbf{X}\right)
dB_{s}+\alpha_{s}ds \label{e.5.61}%
\end{equation}
where%
\begin{equation}
\alpha_{s}:=//_{s}^{-1}\left[  \nabla_{//_{s}z_{s}v}\left(  \sum_{i=1}%
^{n}\nabla_{X_{i}}X_{i}+X_{0}\right)  -\frac{1}{2}\sum_{i=1}^{n}R^{\nabla
}\left(  //_{s}z_{s}v,X_{i}\left(  \Sigma_{s}\right)  \right)  X_{i}\left(
\Sigma_{s}\right)  \right]  ds. \label{e.5.62}%
\end{equation}
If we further assume that $n=N$ and $X_{i}\left(  m\right)  =P\left(
m\right)  e_{i}$ (so that Eq. (\ref{e.5.1}) is equivalent to Eq.
(\ref{e.5.42}) if $X_{0}\equiv0)$, then $\alpha_{s}=-\frac{1}{2}%
\operatorname*{Ric}\nolimits_{//_{s}}z_{s}vds,$ i.e. Eq. (\ref{e.5.59}) is
equivalent to%
\begin{equation}
dz_{s}v=//_{s}^{-1}P\left(  \Sigma_{s}\right)  dP\left(  //_{s}z_{s}v\right)
dB_{s}+\left[  //_{s}^{-1}\nabla_{//_{s}z_{s}v}X_{0}-\frac{1}{2}%
\operatorname{Ric}_{//_{s}}z_{s}v\right]  ds. \label{e.5.63}%
\end{equation}

\end{theorem}

\begin{proof}
In this proof there will always be an implied sum on repeated indices. Using
Proposition \ref{p.5.36},%
\begin{align}
d\left[  //_{s}^{-1}\left(  \nabla_{//_{s}z_{s}v}\mathbf{X}\right)  \right]
dB_{s}  &  =//_{s}^{-1}\left[  \nabla_{\mathbf{X}\left(  \Sigma_{s}\right)
dB_{s}\otimes//_{s}z_{s}v}^{2}\mathbf{X}+\nabla_{//_{s}dz_{s}v}\mathbf{X}%
\right]  dB_{s}\nonumber\\
&  =//_{s}^{-1}\left[  \nabla_{\mathbf{X}\left(  \Sigma_{s}\right)
dB_{s}\otimes//_{s}z_{s}v}^{2}\mathbf{X}+\nabla_{\left(  \nabla_{//_{s}z_{s}%
v}\mathbf{X}\right)  dB_{s}}\mathbf{X}\right]  dB_{s}\nonumber\\
&  =//_{s}^{-1}\left[  \nabla_{X_{i}\left(  \Sigma_{s}\right)  \otimes
//_{s}z_{s}v}^{2}X_{i}+\nabla_{\left(  \nabla_{//_{s}z_{s}v}X_{i}\right)
}X_{i}\right]  ds. \label{e.5.64}%
\end{align}
Now by Proposition \ref{p.3.38},
\begin{align*}
\nabla_{X_{i}\left(  \Sigma_{s}\right)  \otimes//_{s}z_{s}v}^{2}X_{i}  &
=\nabla_{//_{s}z_{s}v\otimes X_{i}\left(  \Sigma_{s}\right)  }^{2}%
X_{i}ds+R^{\nabla}\left(  X_{i}\left(  \Sigma_{s}\right)  ,//_{s}%
z_{s}v\right)  X_{i}\left(  \Sigma_{s}\right) \\
&  =\nabla_{//_{s}z_{s}v\otimes X_{i}\left(  \Sigma_{s}\right)  }^{2}%
X_{i}ds-R^{\nabla}\left(  //_{s}z_{s}v,X_{i}\left(  \Sigma_{s}\right)
\right)  X_{i}\left(  \Sigma_{s}\right) \\
&  =\left[  \nabla_{//_{s}z_{s}v}\nabla_{X_{i}}X_{i}-\nabla_{\nabla
_{//_{s}z_{s}v}X_{i}}X_{i}\right] \\
&  \qquad-R^{\nabla}\left(  //_{s}z_{s}v,X_{i}\left(  \Sigma_{s}\right)
\right)  X_{i}\left(  \Sigma_{s}\right)
\end{align*}
which combined with Eq. (\ref{e.5.64}) implies%
\begin{equation}
d\left[  //_{s}^{-1}\left(  \nabla_{//_{s}z_{s}v}\mathbf{X}\right)  \right]
dB_{s}=//_{s}^{-1}\left[  \nabla_{//_{s}z_{s}v}\nabla_{X_{i}}X_{i}-R^{\nabla
}\left(  //_{s}z_{s}v,X_{i}\left(  \Sigma_{s}\right)  \right)  X_{i}\left(
\Sigma_{s}\right)  \right]  ds. \label{e.5.65}%
\end{equation}
Eq. (\ref{e.5.61}) is now a follows directly from this equation and Eq.
(\ref{e.5.59}).

If we further assume $n=N,$ $X_{i}\left(  m\right)  =P\left(  m\right)  e_{i}$
and $X_{0}\left(  m\right)  =0,$ then
\begin{equation}
\left(  \nabla_{//_{s}z_{s}v}\mathbf{X}\right)  dB_{s}=//_{s}^{-1}P\left(
\Sigma_{s}\right)  dP\left(  //_{s}z_{s}v\right)  dB_{s}. \label{e.5.66}%
\end{equation}
Moreover, from the definition of the Ricci tensor in Eq. (\ref{e.3.31}) and
making use of Eq. (\ref{e.3.50}) in the proof of Proposition \ref{p.3.48} we
have%
\begin{equation}
R^{\nabla}\left(  //_{s}z_{s}v,X_{i}\left(  \Sigma_{s}\right)  \right)
X_{i}\left(  \Sigma_{s}\right)  =\operatorname{Ric}_{//_{s}}//_{s}z_{s}v.
\label{e.5.67}%
\end{equation}
Combining Eqs. (\ref{e.5.66}) and (\ref{e.5.67}) along with $\nabla_{X_{i}%
}X_{i}=0$ (from Proposition \ref{p.3.48}) with Eqs. (\ref{e.5.61}) and
(\ref{e.5.62}) implies Eq. (\ref{e.5.63}).
\end{proof}

In the next result, we will filter out the \textquotedblleft redundant
noise\textquotedblright\ in Eq. (\ref{e.5.63}). This is useful for deducing
intrinsic formula from their extrinsic cousins, see, for example, Corollary
\ref{c.6.4} and Theorem \ref{t.7.39} below.

\begin{theorem}
[Filtering out the Redundant Noise]\label{t.5.44}Keep the same setup in
Theorem \ref{t.5.43} with $n=N$ and $X_{i}\left(  m\right)  =P\left(
m\right)  e_{i}.$ Further let $\mathcal{M}$ be the $\sigma$ -- algebra
generated by the solution $\Sigma=\left\{  \Sigma_{s}:s\geq0\right\}  .$ Then
there is a version, $\bar{z}_{s},$ of $\mathbb{E}\left[  z_{s}|\mathcal{M}%
\right]  $ such that $s\rightarrow\bar{z}_{s}$ is continuous and $\bar{z}$
satisfies,%
\begin{equation}
\bar{z}_{s}v=v+\int_{0}^{s}\left[  //_{r}^{-1}\left(  \nabla_{//_{r}\bar
{z}_{r}v}X_{0}\right)  -\frac{1}{2}\operatorname{Ric}_{//_{r}}\bar{z}%
_{r}v\right]  dr. \label{e.5.68}%
\end{equation}
In particular if $X_{0}=0,$ then
\begin{equation}
\frac{d}{ds}\bar{z}_{s}=-\frac{1}{2}\operatorname{Ric}_{//_{s}}\bar{z}%
_{s}\text{ with }\bar{z}_{0}=id, \label{e.5.69}%
\end{equation}

\end{theorem}

\begin{proof}
In this proof, we let $b_{s}$ be the martingale part of the anti-development
map, $\Psi_{s}\left(  \Sigma\right)  ,$ i.e.%
\[
b_{s}:=\int_{0}^{s}//_{r}^{-1}P\left(  \Sigma_{r}\right)  \delta B_{r}%
=\int_{0}^{s}//_{r}^{-1}P\left(  \Sigma_{r}\right)  dB_{r}.
\]
Since $(\Sigma_{s},u_{s})$ solves the stochastic differential equation,%
\begin{align*}
\delta\Sigma_{s}  &  =u_{s}\delta b_{s}+X_{0}\left(  \Sigma_{s}\right)
ds\text{ with }\Sigma_{0}=o\\
\delta u  &  =-\Gamma\left(  \delta\Sigma\right)  u=-\Gamma\left(  u\delta
b\right)  u\text{ with }u_{0}=I\in O(N)
\end{align*}
it follows that $\left(  \Sigma,u\right)  $ may be expressed as a function of
the Brownian motion, $b.$ Therefore by the martingale representation property,
see Corollary \ref{c.7.20} below, any measurable function, $f\left(
\Sigma\right)  ,$ of $\Sigma$ may be expressed as
\[
f\left(  \Sigma\right)  =f_{0}+\int_{0}^{1}\langle a_{r},db_{r}\rangle
=f_{0}+\int_{0}^{1}\langle a_{r},//_{r}^{-1}\left[  P\left(  \Sigma
_{r}\right)  dB_{r}\right]  \rangle.
\]
Hence, using $PdP=dPQ,$ the previous equation and the isometry property of the
It\^{o} integral,%
\begin{align*}
\mathbb{E}  &  \left\{  \int_{0}^{s}\left[  P\left(  \Sigma_{r}\right)
dP\left(  //_{r}z_{r}v\right)  dB_{r}\right]  f\left(  \Sigma\right)  \right\}
\\
&  =\mathbb{E}\left\{  \int_{0}^{s}\left[  dP\left(  //_{r}z_{r}v\right)
Q\left(  \Sigma_{r}\right)  dB_{r}\right]  \int_{0}^{1}\left\langle P\left(
\Sigma_{r}\right)  //_{r}a_{r},dB_{r}\right\rangle \right\} \\
&  =\mathbb{E}\left\{  \int_{0}^{s}\left[  dP\left(  //_{r}z_{r}v\right)
Q\left(  \Sigma_{r}\right)  P\left(  \Sigma_{r}\right)  //_{r}a_{r}\right]
dr\right\}  =0.
\end{align*}
This shows that
\[
\mathbb{E}\left[  \int_{0}^{s}P\left(  \Sigma_{r}\right)  dP\left(
//_{r}z_{r}v\right)  dB_{r}|\mathcal{M}\right]  =0
\]
and hence taking the conditional expectation, $\mathbb{E}\left[
\cdot|\mathcal{M}\right]  ,$ of the integrated version of Eq. (\ref{e.5.63})
implies Eq. (\ref{e.5.68}). In performing this operation we have used the fact
that $(\Sigma,//)$ is $\mathcal{M}$ -- measurable and that $z_{s}$ appears
linearly in Eq. (\ref{e.5.63}). I have also glossed over the technicality of
passing the conditional expectation past the integrals involving a $ds$ term.
For this detail and a much more general presentation of these ideas the reader
is referred to Elworthy, Li and Le Jan \cite{El-LJ-Li}.
\end{proof}

\subsection{More References\label{s.5.7}}

For more details on the sorts of results in this section, the books by
Elworthy \cite{El}, Emery \cite{Emery}, and Ikeda and Watanabe \cite{IW2},
Malliavin \cite{Malliavin97}, Stroock \cite{Stroock2001}, and Hsu
\cite{Hsu2003c} are highly recommended. The following articles and books are
also relevant,
\cite{BeDa,Bismut81,Bismut84a,Darling82,Ee2,EE1971,El2,Kendall87a,Malliavin78b,Me2,No2,Sc1,Sc2,Sc3,Wa1}%
.

\section{Heat Kernel Derivative Formula\label{s.6}}

In this short section we will illustrate how to derive Bismut type formulas
for derivatives of heat kernels. For more details and much more general
formula see, for example, Driver and Thalmaier \cite{Driver01b}, Elworthy, Le
Jan and Li \cite{El-LJ-Li}, Stroock and Turetsky \cite{StT1,StT2} and Hsu
\cite{Hsu99e} and the references therein. Throughout this section $\Sigma_{s}$
will be an $M$ -- valued semi-martingale, $//_{s}$ will be stochastic parallel
translation along $\Sigma$ and
\[
b_{s}=\Psi_{s}\left(  \Sigma\right)  :=\int_{0}^{s}//_{r}^{-1}\delta\Sigma
_{r}.
\]
Furthermore, let $Q_{s}$ denote the unique solution to the differential
equation:
\begin{equation}
\frac{dQ_{s}}{ds}=-\frac{1}{2}Q_{s}\text{$\operatorname*{Ric}$}_{//_{s}%
}\text{\textrm{\ with }}Q_{0}=I. \label{e.6.1}%
\end{equation}
See Eq. (\ref{e.5.60}) for the definition of $\operatorname*{Ric}$$_{//_{s}}.$

\begin{lemma}
\label{l.6.1}Let $f:M\rightarrow\mathbb{R}$ be a smooth function, $t>0$ and
for $s\in\lbrack0,t]$ let
\begin{equation}
F(s,m):=(e^{(t-s)\bar{\Delta}/2}f)(m). \label{e.6.2}%
\end{equation}
If $\Sigma_{s}$ is an $M$ -- valued Brownian motion, then the process
$s\in\lbrack0,t]\rightarrow Q_{s}//_{s}^{-1}\vec{\nabla}F(s,\Sigma_{s})$ is a
martingale and%
\begin{equation}
d\left[  Q_{s}//_{s}^{-1}\vec{\nabla}F(s,\Sigma_{s})\right]  =Q_{s}//_{s}%
^{-1}\nabla_{//_{s}db_{s}}\vec{\nabla}F(s,\cdot). \label{e.6.3}%
\end{equation}

\end{lemma}

\begin{proof}
Let $W_{s}:=//_{s}^{-1}\vec{\nabla}F(s,\Sigma_{s}).$ Then by Proposition
\ref{p.5.36} and Theorem \ref{t.3.49},%
\begin{align*}
dW_{s}=  &  \left[  //_{s}^{-1}\vec{\nabla}\partial_{s}F(s,\Sigma_{s}%
)+\frac{1}{2}//_{s}^{-1}\nabla_{//_{s}e_{i}\otimes//_{s}e_{i}}^{2}\vec{\nabla
}F(s,\cdot)\right]  ds\\
&  \qquad+//_{s}^{-1}\nabla_{//_{s}e_{i}}\vec{\nabla}F(s,\cdot)db_{s}^{i}\\
=  &  \frac{1}{2}//_{s}^{-1}\left[  \nabla_{//_{s}e_{i}\otimes//_{s}e_{i}}%
^{2}\vec{\nabla}F(s,\cdot)-\left(  \vec{\nabla}\Delta F(s,\cdot)\right)
\left(  \Sigma_{s}\right)  \right]  ds\\
&  \qquad+//_{s}^{-1}\nabla_{//_{s}e_{i}}\vec{\nabla}F(s,\cdot)db_{s}^{i}\\
=  &  \frac{1}{2}//_{s}^{-1}\operatorname*{Ric}\vec{\nabla}F(s,\Sigma
_{s})ds+//_{s}^{-1}\nabla_{//_{s}e_{i}}\vec{\nabla}F(s,\cdot)db_{s}^{i}\\
=  &  \frac{1}{2}\operatorname{Ric}_{//_{s}}W_{s}ds+//_{s}^{-1}\nabla
_{//_{s}e_{i}}\vec{\nabla}F(s,\cdot)db_{s}^{i}%
\end{align*}
where $\left\{  e_{i}\right\}  _{i=1}^{d}$ is an orthonormal basis for
$T_{o}M$ and there is an implied sum on repeated indices. Hence if $Q$ solves
Eq. \text{(\ref{e.6.1})}, then
\begin{align*}
d\left[  Q_{s}W_{s}\right]   &  =-\frac{1}{2}Q_{s}\operatorname*{Ric}%
\nolimits_{//_{s}}W_{s}ds+Q_{s}\left[  \frac{1}{2}\operatorname*{Ric}%
\nolimits_{//_{s}}W_{s}ds+//_{s}^{-1}\nabla_{//_{s}e_{i}}\vec{\nabla}%
F(s,\cdot)db_{s}^{i}\right] \\
&  =Q_{s}//_{s}^{-1}\nabla_{//_{s}e_{i}}\vec{\nabla}F(s,\cdot)db_{s}^{i}%
\end{align*}
which proves Eq. (\ref{e.6.3}) and shows that $Q_{s}W_{s}$ is a martingale as desired.
\end{proof}

\begin{theorem}
[Bismut]\label{t.6.2}Let $f:M\rightarrow\mathbb{R}$ be a smooth function and
$\Sigma$ be an $M$ -- valued Brownian motion with $\Sigma_{0}=o,$ then for
$0<t_{0}\leq t<\infty,$%
\begin{equation}
\vec{\nabla}(e^{t\Delta/2}f)(o)=\frac{1}{t_{0}}E\left[  \left(  \int
_{0}^{t_{0}}Q_{r}db_{r}\right)  f(\Sigma_{t})\right]  . \label{e.6.4}%
\end{equation}

\end{theorem}

\begin{proof}
The proof given here is modelled on Remark 6 on p. 84 in Bismut
\cite{Bismut84a}\ and the proof of Theorem 2.1 in Elworthy and Li
\cite{Elworthy94a}. Also see Norris \cite{No1, No2, No3}. For $(s,m)\in
\lbrack0,t]\times M$ let $F$ be defined as in Eq. (\ref{e.6.2}).$\ $ We wish
to compute the differential of $k_{s}:=\left(  \int_{0}^{s}Q_{r}db_{r}\right)
F(s,\Sigma_{s}).$ By Eq. (\ref{e.5.38}), $d\left[  F(s,\Sigma_{s})\right]
=\langle\vec{\nabla}(F(s,\cdot))(\Sigma_{s}),//_{s}db_{s}\rangle$ and
therefore:
\begin{align*}
dk_{s}  &  =F(s,\Sigma_{s})Q_{s}db_{s}+\left(  \int_{0}^{s}Q_{r}db_{r}\right)
\langle\vec{\nabla}(F(s,\cdot))(\Sigma_{s}),//_{s}db_{s}\rangle\\
&  \quad+\sum_{i=1}^{d}\langle\vec{\nabla}(F(s,\cdot))(\Sigma_{s}),//_{s}%
e_{i}\rangle Q_{s}e_{i}\,ds.
\end{align*}
From this we conclude that
\begin{align*}
\mathbb{E}\left[  k_{t_{0}}\right]   &  =\mathbb{E}\left[  k_{0}\right]
+\mathbb{E}\int_{0}^{t_{0}}\sum_{i=1}^{d}\langle//_{s}^{-1}\vec{\nabla
}(F(s,\cdot))(\Sigma_{s}),e_{i}\rangle Q_{s}e_{i}\,ds\\
&  =\int_{0}^{t_{0}}\mathbb{E}\left[  Q_{s}//_{s}^{-1}\vec{\nabla}%
(F(s,\cdot))(\Sigma_{s})\right]  ds\\
&  =\int_{0}^{t_{0}}\mathbb{E}\left[  Q_{0}//_{0}^{-1}\vec{\nabla}%
(F(0,\cdot))(\Sigma_{0})\right]  ds=t_{0}\vec{\nabla}(e^{t\Delta/2}f)(o)
\end{align*}
wherein the the third equality we have used (by Lemma \ref{l.6.1}) that
$s\rightarrow Q_{s}//_{s}^{-1}\vec{\nabla}(F(s,\cdot))(\Sigma_{s})$ is a
martingale. Hence
\[
\vec{\nabla}(e^{t\Delta/2}f)(o)=\frac{1}{t_{0}}\mathbb{E}\left[  \left(
\int_{0}^{t_{0}}Q_{s}db_{s}\right)  (e^{(t-t_{0})\Delta/2}f)(\Sigma_{t_{0}%
})\right]
\]
from which Eq. (\ref{e.6.4}) follows using either the Markov property of
$\Sigma_{s}$ or the fact that $s\rightarrow\left(  e^{(t-s)\Delta/2}f\right)
(\Sigma_{s})$ is a martingale.
\end{proof}

The following theorem is an non-intrinsic form of Theorem \ref{t.6.2}. In this
theorem we will be using the notation introduced before Theorem \ref{t.5.41}.
Namely, let $\left\{  X_{i}\right\}  _{i=0}^{n}\subset\Gamma\left(  TM\right)
$ be as in Notation \ref{n.5.4}, $B_{s}$ be an $\mathbb{R}^{n}$ -- valued
Brownian motion, and $T_{s}\left(  m\right)  =\Sigma_{s}$ where $\Sigma_{s}$
is the solution to Eq. (\ref{e.5.1}) with $\Sigma_{s}=m\in M$ and $\beta=B.$

\begin{theorem}
[Elworthy - Li]\label{t.6.3}Assume that $\mathbf{X}\left(  m\right)
:\mathbb{R}^{n}\rightarrow T_{m}M$ (recall $\mathbf{X}\left(  m\right)
a:=\sum_{i=1}^{n}X_{i}\left(  m\right)  a_{i})$ is surjective for all $m\in M$
and let%
\begin{equation}
\mathbf{X}\left(  m\right)  ^{\#}=\left[  \mathbf{X}\left(  m\right)
|_{\mathrm{Nul}(\mathbf{X}\left(  m\right)  )^{\perp}}\right]  ^{-1}%
:T_{m}M\rightarrow\mathbb{R}^{n}, \label{e.6.5}%
\end{equation}
where the orthogonal complement is taken relative to the standard inner
product on $\mathbb{R}^{n}.$ (See Lemma \ref{l.7.38} below for more on
$\mathbf{X}\left(  m\right)  ^{\#}.)$ Then for all $v\in T_{o}M,$
$0<t_{o}<t<\infty$ and $f\in C\left(  M\right)  $ we have%
\begin{equation}
v\left(  e^{tL/2}f\right)  =\frac{1}{t_{0}}\mathbb{E}\left[  f\left(
\Sigma_{t}\right)  \int_{0}^{t_{0}}\langle\mathbf{X}\left(  \Sigma_{s}\right)
^{\#}Z_{s}v,dB_{s}\rangle\right]  \label{e.6.6}%
\end{equation}
where $Z_{s}=T_{s\ast o}$ as in Eq. (\ref{e.5.57}).
\end{theorem}

\begin{proof}
Let $L=\sum_{i=1}^{n}X_{i}^{2}+2X_{0}$ be the generator of the diffusion,
$\left\{  T_{s}\left(  m\right)  \right\}  _{s\geq0}.$ Since $\mathbf{X}%
\left(  m\right)  :\mathbb{R}^{n}\rightarrow T_{m}M$ is surjective for all
$m\in M,$ $L$ is an elliptic operator on $C^{\infty}\left(  M\right)  .$ So,
using results similar to those in Fact \ref{fact.5.32}, it makes sense to
define $F_{s}\left(  m\right)  :=\left(  e^{\left(  t-s\right)  L/2}f\right)
\left(  m\right)  $ and $N_{s}^{m}=F_{s}\left(  T_{s}\left(  m\right)
\right)  .$ Then%
\[
\partial_{s}F_{s}+\frac{1}{2}LF_{s}=0\text{ with }F_{t}=f
\]
and by It\^{o}'s lemma,%
\begin{equation}
dN_{s}^{m}=d\left[  F_{s}\left(  T_{s}\left(  m\right)  \right)  \right]
=\sum_{i=1}^{n}\left(  X_{i}F_{s}\right)  (T_{s}\left(  m\right)  )dB_{s}^{i}.
\label{e.6.7}%
\end{equation}
This shows $N_{s}^{m}$ is a martingale for all $m\in M$ and, upon integrating
Eq. (\ref{e.6.7}) on $s,$ that%
\[
f\left(  T_{t}\left(  m\right)  \right)  =e^{tL/2}f(m)+\sum_{i=1}^{n}\int
_{0}^{t}\left(  X_{i}F_{s}\right)  (T_{s}\left(  m\right)  )dB_{s}^{i}.
\]
Hence if $a_{s}\in\mathbb{R}^{n}$ is a predictable process such that
$\mathbb{E}\int_{0}^{t}\left\vert a_{s}\right\vert ^{2}ds<\infty,$ then by the
It\^{o} isometry property,%
\begin{align}
\mathbb{E}\left[  f\left(  T_{t}\left(  m\right)  \right)  \int_{0}^{t}\langle
a,dB\rangle\right]   &  =\int_{0}^{t}\mathbb{E}\left[  \left(  X_{i}%
F_{s}\right)  (T_{s}\left(  m\right)  )a_{i}(s)\right]  ds\nonumber\\
&  =\int_{0}^{t}\mathbb{E}\left[  \left(  d_{M}F_{s}\right)  \left(
\mathbf{X}(T_{s}\left(  m\right)  )a_{s}\right)  \right]  ds. \label{e.6.8}%
\end{align}

Suppose that $\ell_{s}\in\mathbb{R}$ is a continuous piecewise differentiable
function and let $a_{s}:=\ell_{s}^{\prime}\mathbf{X}\left(  \Sigma_{s}\right)
^{\#}Z_{s}v.$ Then form Eq. (\ref{e.6.8}) we have%
\begin{equation}
\mathbb{E}\left[  f\left(  \Sigma_{t}\right)  \int_{0}^{t}\langle\ell
_{s}^{\prime}\mathbf{X}\left(  \Sigma_{s}\right)  ^{\#}Z_{s}v,dB_{s}%
\rangle\right]  =\int_{0}^{t}\ell_{s}^{\prime}\mathbb{E}\left[  \left(
d_{M}F_{s}\right)  \left(  Z_{s}v\right)  \right]  ds. \label{e.6.9}%
\end{equation}
Since $N_{s}^{m}=F_{s}\left(  T_{s}\left(  m\right)  \right)  $ is a
martingale for all $m,$ we may deduce that%
\begin{equation}
v\left(  m\rightarrow N_{s}^{m}\right)  =d_{M}F_{s}\left(  T_{s\ast
o}v\right)  =d_{M}F_{s}\left(  Z_{s}v\right)  \label{e.6.10}%
\end{equation}
is a martingale as well for any $v\in T_{o}M.$ In particular, $s\in
\lbrack0,t]\rightarrow\mathbb{E}\left[  \left(  d_{M}F_{s}\right)  \left(
Z_{s}v\right)  \right]  $ is constant and evaluating this expression at $s=0$
and $s=t$ implies%
\begin{equation}
\mathbb{E}\left[  \left(  d_{M}F_{s}\right)  \left(  Z_{s}v\right)  \right]
=v\left(  e^{tL/2}f\right)  =\mathbb{E}\left[  \left(  d_{M}f\right)  \left(
Z_{t}v\right)  \right]  . \label{e.6.11}%
\end{equation}
Using Eq. (\ref{e.6.11}) in Eq. (\ref{e.6.9}) then shows%
\[
\mathbb{E}\left[  f\left(  \Sigma_{t}\right)  \int_{0}^{t}\langle\ell
_{s}^{\prime}\mathbf{X}\left(  \Sigma_{s}\right)  ^{\#}Z_{s}v,dB_{s}%
\rangle\right]  =\left(  \ell_{t}-\ell_{0}\right)  v\left(  e^{tL/2}f\right)
\]
which, by taking $\ell_{s}=s\wedge t_{0},$ implies Eq. (\ref{e.6.6}).
\end{proof}

\begin{corollary}
\label{c.6.4}Theorem \ref{t.6.3} may be used to deduce Theorem \ref{t.6.2}.
\end{corollary}

\begin{proof}
Apply Theorem \ref{t.6.3} with $n=N,$ $X_{0}\equiv0$ and $X_{i}\left(
m\right)  =P\left(  m\right)  e_{i}$ for $i=1,\dots,N$ to learn%
\begin{equation}
v\left(  e^{t\Delta/2}f\right)  =\frac{1}{t_{0}}\mathbb{E}\left[  f\left(
\Sigma_{t}\right)  \int_{0}^{t_{0}}\langle Z_{s}v,dB_{s}\rangle\right]
=\frac{1}{t_{0}}\mathbb{E}\left[  f\left(  \Sigma_{t}\right)  \int_{0}^{t_{0}%
}\langle//_{s}z_{s}v,dB_{s}\rangle\right]  \label{e.6.12}%
\end{equation}
where we have used $L=\Delta$ (see Proposition \ref{p.3.48}) and
$\mathbf{X}\left(  m\right)  ^{\#}=P\left(  m\right)  $ in this setting. By
Theorem \ref{t.5.40},%
\begin{align*}
\int_{0}^{t_{0}}\langle//_{s}z_{s}v,dB_{s}\rangle &  =\int_{0}^{t_{0}}%
\langle//_{s}z_{s}v,P\left(  \Sigma_{s}\right)  dB_{s}\rangle\\
&  =\int_{0}^{t_{0}}\langle z_{s}v,//_{s}^{-1}P\left(  \Sigma_{s}\right)
dB_{s}\rangle=\int_{0}^{t_{0}}\langle z_{s}v,db_{s}\rangle
\end{align*}
and therefore Eq. (\ref{e.6.12}) may be written as
\[
v\left(  e^{t\Delta/2}f\right)  =\frac{1}{t_{0}}\mathbb{E}\left[  f\left(
\Sigma_{t}\right)  \int_{0}^{t_{0}}\langle z_{s}v,db_{s}\rangle\right]  .
\]
Using Theorem \ref{t.5.44}, this may also be expressed as%
\begin{equation}
v\left(  e^{t\Delta/2}f\right)  =\frac{1}{t_{0}}\mathbb{E}\left[  f\left(
\Sigma_{t}\right)  \int_{0}^{t_{0}}\langle\bar{z}_{s}v,db_{s}\rangle\right]
=\frac{1}{t_{0}}\mathbb{E}\left[  f\left(  \Sigma_{t}\right)  \int_{0}^{t_{0}%
}\langle v,\bar{z}_{s}^{\operatorname*{tr}}db_{s}\rangle\right]
\label{e.6.13}%
\end{equation}
where $\bar{z}_{s}$ solves Eq. (\ref{e.5.69}). By taking transposes of Eq.
(\ref{e.5.69}) it follows that $\bar{z}_{s}^{\operatorname*{tr}}$ satisfies
Eq. (\ref{e.6.1}) and hence $\bar{z}_{s}^{\operatorname*{tr}}=Q_{s}.$ Since
$v\in T_{o}M$ was arbitrary, Equation (\ref{e.6.4}) is now an easy consequence
of Eq. (\ref{e.6.13}) and the definition of $\vec{\nabla}(e^{t\Delta/2}f)(o).$
\end{proof}

\section{Calculus on $W(M)$\label{s.7}}

In this section, $\left(  M,o\right)  $ is assumed to be either a compact
Riemannian manifold equipped with a fixed point $o\in M$ or $M=\mathbb{R}^{d}$
with $o=0.$

\begin{notation}
\label{n.7.1}We will be interested in the following path spaces:%
\begin{align*}
W(T_{o}M)  &  :=\{\omega\in C([0,1]\rightarrow T_{o}M)|\omega(0)=0_{o}\in
T_{o}M\},\\
H\left(  T_{o}M\right)   &  :=\{h\in W(T_{o}M):h(0)=0,\text{ \&~}\langle
h,h\rangle_{H}:=\int_{0}^{1}|h^{\prime}(s)|_{T_{o}M}^{2}ds<\infty\}
\end{align*}
and%
\[
W(M):=\left\{  \sigma\in C([0,1]\rightarrow M):\sigma\left(  0\right)  =0\in
M\right\}  .
\]
(By convention $\langle h,h\rangle_{H}=\infty$ if $h\in W(T_{o}M)$ is not
absolutely continuous.) We refer to $W(T_{o}M)$ as \textbf{Wiener space},
$W\left(  M\right)  $ as \textbf{curved Wiener space} and $H\left(
T_{o}M\right)  $ or $H\left(  \mathbb{R}^{d}\right)  $ as the
\textbf{Cameron-Martin Hilbert space}.
\end{notation}

\begin{definition}
\label{d.7.2}Let $\mu$ and $\mu_{W\left(  M\right)  }$ denote the Wiener
measures on $W\left(  T_{o}M\right)  $ and $W\left(  M\right)  $ respectively,
i.e. $\mu=\mathrm{Law}\left(  b\right)  $ and $\mu_{W\left(  M\right)
}=\mathrm{Law}\left(  \Sigma\right)  $ where $b$ and $\Sigma$ are Brownian
motions on $T_{o}M$ and $M$ starting at $0\in T_{o}M$ and $o\in M$ respectively.
\end{definition}

\begin{notation}
\label{n.7.3}The probability space in this section will often be $\left(
W\left(  M\right)  ,\mathcal{F},\mu_{W\left(  M\right)  }\right)  ,$ where
$\mathcal{F}$ is the completion of the $\sigma$ -- algebra generated by the
projection maps, $\Sigma_{s}:W\left(  M\right)  \rightarrow M$ defined by
$\Sigma_{s}\left(  \sigma\right)  =\sigma_{s}$ for $s\in\lbrack0,1].$ We make
this into a filtered probability space by taking $\mathcal{F}_{s}$ to be the
$\sigma$ -- algebra generated by $\left\{  \Sigma_{r}:r\leq s\right\}  $ and
the null sets in $\mathcal{F}_{s}.$ Also let $//_{s}$ be stochastic parallel
translation along $\Sigma.$
\end{notation}

\begin{definition}
\label{d.7.4}A function $F:W(M)\rightarrow\mathbb{R}$ is called a $C^{k}%
$\textbf{ -- cylinder function} if there exists a partition
\begin{equation}
\pi:=\{0=s_{0}<s_{1}<s_{2}\cdots<s_{n}=1\} \label{e.7.1}%
\end{equation}
of $[0,1]$ and $f\in C^{k}(M^{n})$ such that
\begin{equation}
F(\sigma)=f(\sigma_{s_{1}},\ldots,\sigma_{s_{n}})\text{ for all }\sigma\in
W\left(  M\right)  . \label{e.7.2}%
\end{equation}
If $M=\mathbb{R}^{d},$ we further require that $f$ and all of its derivatives
up to order $k$ have at most polynomial growth at infinity. The collection of
$C^{k}$ -- cylinder functions will be denoted by $\mathcal{F}C^{k}\left(
W\left(  M\right)  \right)  .$
\end{definition}

\begin{definition}
\label{d.7.5}The \textbf{continuous tangent space} to $W(M)$ at $\sigma\in
W(M)$ is the set $CT_{\sigma}W(M)$ of continuous vector-fields along $\sigma$
which are zero at $s=0:$
\begin{equation}
CT_{\sigma}W(M)=\{X\in C([0,1],TM)|X_{s}\in T_{\sigma_{s}}M\text{\textrm{\ }%
$\forall$ $s\in\lbrack0,1]$\textrm{ and }$X(0)=0$}\}. \label{e.7.3}%
\end{equation}

\end{definition}

To motivate the above definition, consider a differentiable path in $\gamma\in
W(M)$ going through $\sigma$ at $t=0.$ Writing $\gamma\left(  t\right)
\left(  s\right)  $ as $\gamma\left(  t,s\right)  ,$ the derivative
$X_{s}:=\frac{d}{dt}|_{0}\gamma(t,s)\in T_{\sigma\left(  s\right)  }M$ of such
a path should, by definition, be a tangent vector to $W(M)$ at $\sigma.$

We now wish to define a \textquotedblleft Riemannian metric\textquotedblright%
\ on $W(M).$ It turns out that the continuous tangent space $CT_{\sigma}W(M)$
is too large for our purposes, see for example the Cameron-Martin Theorem
\ref{t.7.13} below. To remedy this we will introduce a Riemannian structure on
a an a.e. defined \textquotedblleft sub-bundle\textquotedblright\ of
$CTW\left(  M\right)  .$

\begin{definition}
\label{d.7.6}A \textbf{Cameron-Martin process}, $h,$ is a $T_{o}M$ -- valued
process on $W\left(  M\right)  $ such that $s\rightarrow h(s)$ is in $H,$
$\mu_{W\left(  M\right)  }$ -- a.e. Contrary to our earlier assumptions, we do
\textbf{not} assume that $h$ is adapted unless explicitly stated.
\end{definition}

\begin{definition}
\label{d.7.7}Suppose that $X$ is a $TM$ -- valued process on $\left(  W\left(
M\right)  ,\mu_{W\left(  M\right)  }\right)  $ such that the process
$\pi\left(  X_{s}\right)  =\Sigma_{s}\in M.$ We will say $X$ is a
\textbf{Cameron-Martin vector-field} if
\begin{equation}
h_{s}:=//_{s}^{_{-1}}X_{s} \label{e.7.4}%
\end{equation}
is a Cameron-Martin valued process and
\begin{equation}
\langle X,X\rangle_{\mathcal{X}}:=\mathbb{E}[\langle h,h\rangle_{H}]<\infty.
\label{e.7.5}%
\end{equation}
A Cameron-Martin vector field $X$ is said to be adapted if $h:=//^{-1}X$ is
adapted. The set of Cameron-Martin vector-fields will be denoted by
$\mathcal{X}$ and those which are adapted will be denoted by $\mathcal{X}%
_{a}.$
\end{definition}

\begin{remark}
\label{r.7.8}Notice that $\mathcal{X}$ is a Hilbert space with the inner
product determined by $\langle\cdot,\cdot\rangle_{\mathcal{X}}$ in
\text{(\ref{e.7.5})}. Furthermore, $\mathcal{X}_{a}$ is a Hilbert-subspace of
$\mathcal{\ X}.$
\end{remark}

\begin{notation}
\label{n.7.9}Given a Cameron-Martin process $h$, let $X^{h}:=//h.$ In this way
we may identify Cameron-Martin processes with Cameron-Martin vector fields.
\end{notation}

We define a \textquotedblleft metric\textquotedblright, $G,$\footnote{The
function $G$ is to be loosely interpreted as a Riemannian metric on $W(M).$}
on$\mathcal{\ X}$ by
\begin{equation}
G(X^{h},X^{h})=\langle h,h\rangle_{H}. \label{e.7.6}%
\end{equation}
With this notation we have $\langle X,X\rangle_{\mathcal{X}}=\mathbb{E}\left[
G(X,X)\right]  .$

\begin{remark}
\label{r.7.10}Notice, if $\sigma$ is a smooth path then the expression in
\text{(\ref{e.7.6})} could be written as
\[
G(X,X)=\int_{0}^{1}g\left(  \frac{\nabla}{ds}X(s),\frac{\nabla}{ds}%
X(s)\right)  ds,
\]
where ${\frac{\nabla}{ds}}$ denotes the covariant derivative along the path
$\sigma$ which is induced from the covariant derivative $\nabla.$ This is a
typical metric used by differential geometers on path and loop spaces.
\end{remark}

\begin{notation}
\label{n.7.11}Given a Cameron-Martin vector field $X$ on $\left(  W\left(
M\right)  ,\mu_{W\left(  M\right)  }\right)  $ and a cylinder function
$F\in\mathcal{F}C^{1}\left(  W\left(  M\right)  \right)  $ as in Eq.
(\ref{e.7.2}), let $XF$ denote the random variable
\begin{equation}
XF\left(  \sigma\right)  :=\sum_{i=1}^{n}(\text{$\operatorname*{grad}$}%
_{i}F(\sigma),X_{s_{i}}\left(  \sigma\right)  ), \label{e.7.7}%
\end{equation}
where%
\begin{equation}
\text{$\operatorname*{grad}$}_{i}F(\sigma):=\left(
\text{$\operatorname*{grad}$}_{i}f\right)  (\sigma_{s_{1}},\ldots
,\sigma_{s_{n}}) \label{e.7.8}%
\end{equation}
and $\left(  \operatorname*{grad}_{i}f\right)  $ denotes the gradient of $f$
relative to the $i^{\text{th }}$ variable.
\end{notation}

\begin{notation}
\label{n.7.12}The \textbf{gradient}, $DF,$ of a smooth cylinder functions,
$F,$ on $W(M)$ is the unique Cameron-Martin process such that $G\left(
DF,X\right)  =XF$ for all $X\in\mathcal{X}.$ The explicit formula for $D,$ as
the reader should verify, is
\begin{equation}
\left(  DF\right)  _{s}=//_{s}\left(  \sum_{i=1}^{n}s\wedge s_{i}//_{s_{i}}%
{}^{-1}\text{$\operatorname*{grad}$}_{i}F(\sigma)\right)  . \label{e.7.9}%
\end{equation}
The formula in Eq. (\ref{e.7.9}) defines a densely defined operator,
$D:L^{2}\left(  \mu\right)  \rightarrow\mathcal{X}$ with $\mathcal{D}\left(
D\right)  =\mathcal{F}C^{1}\left(  W\left(  M\right)  \right)  $ as its domain.
\end{notation}

\subsection{Classical Wiener Space Calculus\label{s.7.1}}

In this subsection (which is a warm up for the sequel) we will specialize to
the case where $M=\mathbb{R}^{d},$ $o=0\in\mathbb{R}^{d}.$ To simplify
notation let $W:=W(\mathbb{R}^{d}),$ $H:=H\left(  \mathbb{R}^{d}\right)  ,$
$\mu=\mu_{W\left(  \mathbb{R}^{d}\right)  },$ $b_{s}\left(  \omega\right)
=\omega_{s}$ for all $s\in\lbrack0,1]$ and $\omega\in W.$ Recall that
$\left\{  \mathcal{F}_{s}:s\in\lbrack0,1]\right\}  $ is the filtration on $W$
as explained in Notation \ref{n.7.3} where we are now writing $b$ for
$\Sigma.$ Cameron and Martin \cite{CM0,CM1,CM2,CM3} and Cameron \cite{CM2}
began the study of calculus on this classical Wiener space. They proved the
following two results, see Theorem 2, p. 387 of \cite{CM1} and Theorem II, p.
919 of \cite{CM2} respectively. (There have been many extensions of these
results partly initiated by Gross' work in \cite{Gr2,Gr3}.)

\begin{theorem}
[Cameron \& Martin 1944]\label{t.7.13}Let $(W,\mathcal{F},\mu)$ be the
classical Wiener space described above and for $h\in W,$ define $T_{h}%
:W\rightarrow W$ by $T_{h}(\omega)=\omega+h$ for all $\omega\in W.$ If $h$ is
$C^{1},$ then $\mu T_{h}^{-1}$ is absolutely continuous relative to $\mu.$
\end{theorem}

This theorem was extended by Maruyama \cite{Mar} and Girsanov \cite{Gir} to
allow the same conclusion for $h\in H$ and more general Cameron-Martin
processes. Moreover it is now well known $\mu T_{h}^{-1}\ll\mu$ iff $h\in H.$
From the Cameron and Martin theorem one may prove Cameron's integration by
parts formula.

\begin{theorem}
[Cameron 1951]\label{t.7.14}Let $h\in H$ and $F,G\in L^{\infty-}(\mu
):=\cap_{1\leq p<\infty}L^{p}\left(  \mu\right)  $ such that $\partial
_{h}F:=\frac{d}{d\varepsilon}F\circ T_{\varepsilon h}|_{\varepsilon=0}$ and
$\partial_{h}G:=\frac{d}{d\varepsilon}G\circ T_{\varepsilon h}|_{\varepsilon
=0}$ where the derivatives are supposed to exist\footnote{The notion of
derivative stated here is weaker than the notion given in \cite{CM2}.
Nevertheless Cameron's proof covers this case without any essential change.}
in $L^{p}(\mu)$ for all $1\leq p<\infty.$ Then
\[
\int_{W}\partial_{h}F\cdot G\,d\mu=\int_{W}F\partial_{h}^{\ast}G\,d\mu,
\]
where $\partial_{h}^{\ast}G=-\partial_{h}G+z_{h}G$ and $z_{h}:=\int_{0}%
^{1}\langle h^{\prime}\left(  s\right)  ,db_{s}\rangle_{\mathbb{R}^{d}}.$
\end{theorem}

In this flat setting parallel translation is trivial, i.e. $//_{s}=id$ for all
$s.$ Hence the gradient operator $D$ in Eq. (\ref{e.7.9}) reduces to the
equation,%
\[
\left(  DF\right)  _{s}\left(  \omega\right)  =\left(  \sum_{i=1}^{n}s\wedge
s_{i}\text{$\operatorname*{grad}$}_{i}F(\omega_{s})\right)  .
\]
Similarly the association of a Cameron-Martin vector field $X$ on
$W(\mathbb{R}^{d})$ with a Cameron-Martin valued process $h$ in Eq.
(\ref{e.7.4}) is simply that $X=h.$

We will now recall that adapted Cameron-Martin vector fields, $X=h,$ are in
the domain of $D^{\ast}.$ From this fact it will easily follow that $D^{\ast}$
is densely defined.

\begin{theorem}
\label{t.7.15}Let $h$ be an adapted Cameron-Martin process (vector field) on
$W.$ Then $h\in\mathcal{D}(D^{\ast})$ and
\[
D^{\ast}h=\int_{0}^{1}\langle h^{\prime},db\rangle.
\]

\end{theorem}

\begin{proof}
We start by proving the theorem under the additional assumption that
\begin{equation}
\sup_{s\in\lbrack0,1]}\left\vert h_{s}^{\prime}\right\vert \leq
C,\label{e.7.10}%
\end{equation}
where $C$ is a non-random constant. For each $t\in\mathbb{R}$ let
$b(t,s)=b_{s}(t)=b_{s}+th_{s}.$ By Girsanov's theorem, $s\rightarrow b_{s}(t)$
(for fixed $t)$ is a Brownian motion relative to $Z_{t}\cdot\mu,$ where
\[
Z_{t}:=\exp\left(  -\int_{0}^{1}t\langle h_{s}^{\prime},db_{s}\rangle
-{\frac{1}{2}}t^{2}\int_{0}^{1}\langle h_{s}^{\prime},h_{s}^{\prime}\rangle
ds\right)  .
\]
Hence if $F$ is a smooth cylinder function on $W,$%
\[
\mathbb{E}\left[  F\left(  b(t,\cdot)\right)  \cdot Z_{t}\right]
=\mathbb{E}\left[  F(b)\right]  .
\]
Differentiating this equation in $t$ at $t=0,$ using%
\[
\langle DF,h\rangle_{H}=\frac{d}{dt}|_{0}F\left(  b\left(  t,\cdot\right)
\right)  \text{ and }\frac{d}{dt}|_{0}Z_{t}=-\int_{0}^{1}\langle h^{\prime
},db\rangle,
\]
shows
\[
\mathbb{E}\left[  \langle DF,h\rangle_{H}\right]  -\mathbb{E}\left[  F\int
_{0}^{1}\langle h^{\prime},db\rangle\right]  =0.
\]
From this equation it follows that $h\in\mathcal{D}(D^{\ast})$ and $D^{\ast
}h=\int_{0}^{1}\langle h^{\prime},db\rangle.$ So it now only remains to remove
the restriction placed on $h$ in Eq. (\ref{e.7.10}).

Let $h$ be a general adapted Cameron-Martin vector-field and for each
$n\in\mathbb{N},$ let
\begin{equation}
h_{n}(s\mathbb{)}:=\int_{0}^{s}h^{\prime}(r\mathbb{)}\cdot1_{|h^{\prime
}(r\mathbb{)}|\leq n}dr. \label{e.7.11}%
\end{equation}
(Notice that $h_{n}$ is still adapted.) By the special case above we know that
$h_{n}\in\mathcal{D}(D^{\ast})$ and $D^{\ast}h_{n}=\int_{0}^{1}\langle
h_{n}^{\prime},db\rangle.$ Therefore,
\[
\mathbb{E}\left\vert D^{\ast}(h_{m}-h_{n})\right\vert ^{2}=\mathbb{E}\int
_{0}^{1}|h_{m}^{\prime}-h_{n}^{\prime}|^{2}ds\rightarrow0\text{\textrm{\ as }%
}m,n\rightarrow\infty
\]
from which it follows that $D^{\ast}h_{n}$ is convergent. Because $D^{\ast}$
is a closed operator, $h\in\mathcal{D}(D^{\ast})$ and
\[
D^{\ast}h=\lim_{n\rightarrow\infty}D^{\ast}h_{n}=\lim_{n\rightarrow\infty}%
\int_{0}^{1}\langle h_{n}^{\prime},db\rangle=\int_{0}^{1}\langle h^{\prime
},db\rangle.
\]

\end{proof}

\begin{corollary}
\label{c.7.16}The operator $D^{\ast}$ is densely defined and hence $D$ is
closable. (Let $\bar{D}$ denote the closure of $D.)$
\end{corollary}

\begin{proof}
Let $h\in H$ and $F$ and $K$ be smooth cylinder functions. Then, by the
product rule,%
\begin{align*}
\langle DF,Kh\rangle_{\mathcal{X}}  &  =\mathbb{E}[\langle KDF,h\rangle
_{H}]=\mathbb{E}[\langle D\left(  KF\right)  -FDK,h\rangle_{H}]\\
&  =\mathbb{E}[F\cdot KD^{\ast}h-F\langle DK,h\rangle_{H}].
\end{align*}
Therefore $Kh\in\mathcal{D}(D^{\ast})$ ($\mathcal{D}(D^{\ast})$ is the domain
of $D^{\ast})$ and
\[
D^{\ast}(Kh)=KD^{\ast}h-\langle DK,h\rangle_{H}.
\]
Since the subspace,
\[
\{Kh|h\in H\text{\textrm{\ and }}K\text{\textrm{\ is a smooth cylinder
function}}\},
\]
is a dense subspace of $\mathcal{X},$ $D^{\ast}$ is densely defined.
\end{proof}

\subsubsection{Martingale Representation Property and the Clark-Ocone Formula}

\begin{lemma}
\label{l.7.17}Let $F(b)=f(b_{s_{1}},\ldots,b_{s_{n}})$ be the smooth cylinder
function on $W$ as in Definition \ref{d.7.4}, then
\begin{equation}
F=\mathbb{E}F+\int_{0}^{1}\langle a_{s},db_{s}\rangle, \label{e.7.12}%
\end{equation}
where $a_{s}$ is a bounded, piecewise-continuous (in $s)$ and predictable
process. Furthermore, the jumps points of $a_{s}$ are contained in the set
$\{s_{1},\ldots,s_{n}\}$ and $a_{s}\equiv0$ is $s\geq s_{n}.$
\end{lemma}

\begin{proof}
The proof will be by induction on $n.$ First assume that $n=1,$ so that
$F(b)=f(b_{t})$ for some $0<t\leq1.$ Let $H(s,m):=(e^{(t-s)\Delta/2}f)(m)$ for
$0\leq s\leq t$ and $m\in\mathbb{R}^{d}.$ Then, by It\^{o}'s formula (or see
Eq. (\ref{e.5.38})),%
\[
dH(s,b_{s})=\langle\operatorname*{grad}H(s,b_{s}),db_{s}\rangle
\]
which upon integrating on $s\in\lbrack0,t]$ gives%
\[
F(b)=(e^{t\Delta/2}f)(o)+\int_{0}^{t}\langle\text{$\operatorname*{grad}$%
}H(s,b_{s}),db_{s}\rangle=\mathbb{E}F+\int_{0}^{1}\langle a_{s},db_{s}%
\rangle,
\]
where $a_{s}=1_{s\leq t}//_{s}^{-1}\operatorname*{grad}H(s,b_{s}).$ This
proves the $n=1$ case. To finish the proof it suffices to show that we may
reduce the assertion of the lemma at the level $n$ to the assertion at the
level $n-1.$

Let $F(b)=f(b_{s_{1}},\ldots,b_{s_{n}}),$
\begin{align*}
(\Delta_{n}f)(x_{1},x_{2},\ldots,x_{n}) &  =(\Delta g)(x_{n})\text{ and}\\
(\operatorname{grad}_{n}f)(x_{1},x_{2},\ldots,x_{n}) &  =\vec{\nabla}g\left(
x_{n}\right)
\end{align*}
where $g(x):=f(x_{1},x_{2},\ldots,x_{n-1},x).$ (So $\Delta_{n}f$ and
$\operatorname*{grad}_{n}f$ is the Laplacian and the gradient of $f$ in the
$n^{\text{th}}$ -- variable.) It\^{o}'s lemma applied to the process,
\[
s\in\lbrack s_{n-1},s_{n}]\rightarrow H(s,b):=(e^{(s_{n}-s)\Delta_{n}%
/2}f)(b_{s_{1}},\ldots,b_{s_{n-1}},b_{s})
\]
gives%
\[
dH(s,b)=\langle\text{$\operatorname*{grad}$}_{n}e^{(s_{n}-s)\Delta_{n}%
/2}f)(b_{s_{1}},\ldots,b_{s_{n-1}},b_{s},db_{s}\rangle
\]
and hence%
\begin{align}
F(b) &  =(e^{(s_{n}-s_{n-1})\Delta_{n}/2}f)(b_{s_{1}},\ldots,b_{s_{n-1}%
},b_{s_{n-1}})\nonumber\\
&  \qquad+\int_{s_{n-1}}^{s_{n}}\langle\text{$\operatorname*{grad}$}%
_{n}e^{(s_{n}-s)\Delta_{n}/2}f)(b_{s_{1}},\ldots,b_{s_{n-1}},b_{s}%
,db_{s}\rangle\nonumber\\
&  =(e^{(s_{n}-s_{n-1})\Delta_{n}/2}f)(b_{s_{1}},\ldots,b_{s_{n-1}}%
,b_{s_{n-1}})+\int_{s_{n-1}}^{s_{n}}\langle\alpha_{s},db_{s}\rangle
,\label{e.7.13}%
\end{align}
where $\alpha_{s}:=(\operatorname*{grad}_{n}e^{(s_{n}-s)\Delta_{n}%
/2}f)(b_{s_{1}},\ldots,b_{s_{n-1}},b_{s})$ for $s\in(s_{n-1},s_{n}).$ By
induction we know that the smooth cylinder function
\[
(e^{(s_{n}-s_{n-1})\Delta_{n}/2}f)(b_{s_{1}},\ldots,b_{s_{n-1}},b_{s_{n-1}})
\]
may be written as a constant plus $\int_{0}^{1}\langle a_{s},db_{s}\rangle,$
where $a_{s}$ is bounded and piecewise continuous and $a_{s}\equiv0$ if $s\geq
s_{n-1}.$ Hence it follows by replacing $a_{s}$ by $a_{s}+1_{(s_{n-1},s_{n}%
)s}\alpha_{s}$ that
\[
F(b)=C+\int_{0}^{s_{n}}\langle a_{s},db_{s}\rangle
\]
for some constant $C.$ Taking expectations of both sides of this equation then
shows $C=\mathbb{E}\left[  F(b)\right]  .$
\end{proof}

\begin{remark}
\label{r.7.18}By being more careful in the proof of the Lemma \ref{l.7.17} (as
is done in more generality later in Theorem \ref{t.7.47}) it is possible to
show $a_{s}$ in Eq. (\ref{e.7.12}) may be written as%
\begin{equation}
a_{s}=\mathbb{E}\left[  \left.  \sum_{i=1}^{n}1_{s\leq s_{i}}%
\text{$\operatorname*{grad}$}_{i}f\left(  b_{s_{1}},\ldots,b_{s_{n}}\right)
\right\vert \mathcal{F}_{s}\right]  . \label{e.7.14}%
\end{equation}
This will also be explained, by indirect means, in Theorem \ref{t.7.21} below.
\end{remark}

\begin{corollary}
\label{c.7.19}Let $F$ be a smooth cylinder function on $W,$ then there is a
predictable, piecewise continuously differentiable Cameron-Martin process $h$
such that $F=\mathbb{E}F+D^{\ast}h.$
\end{corollary}

\begin{proof}
Let $h_{s}:=\int_{0}^{s}a_{r}dr$ where $a$ is the process as in Lemma
\ref{l.7.17}.
\end{proof}

\begin{corollary}
[Martingale Representation Property]\label{c.7.20}Let $F\in L^{2}(\mu),$ then
there is a predictable process, $a_{s},$ such that $\mathbb{E}\int_{0}%
^{1}|a_{s}|^{2}ds<\infty,$ and
\begin{equation}
F=\mathbb{E}F+\int_{0}^{1}\langle a,db\rangle. \label{e.7.15}%
\end{equation}

\end{corollary}

\begin{proof}
Choose a sequence of smooth cylinder functions $\{F_{n}\}$ such that
$F_{n}\rightarrow F$ as $n\rightarrow\infty.$ By replacing $F$ by
$F-\mathbb{E}F$ and $F_{n}$ by $F_{n}-\mathbb{E}F_{n},$ we may assume that
$\mathbb{E}F=0$ and $\mathbb{E}F_{n}=0.$ Let $a^{n}$ be predictable processes
such that $F_{n}=\int_{0}^{1}\langle a^{n},db\rangle$ for all $n.$ Notice
that
\[
\mathbb{E}\int_{0}^{1}|a_{s}^{n}-a_{s}^{m}|^{2}ds=\mathbb{E}(F_{n}-F_{m}%
)^{2}\rightarrow0\text{\textrm{\ as }}m,n\rightarrow\infty.
\]
Hence, if $a:=L^{2}(ds\times d\mu)-\lim_{n\rightarrow\infty}a^{n},$ then
\[
F_{n}=\int_{0}^{1}a^{n}\cdot db\rightarrow\int_{0}^{1}\langle a,db\rangle
\text{\textrm{\ as }}n\rightarrow\infty.
\]
This show that $F=\int_{0}^{1}\langle a,db\rangle.$
\end{proof}

\begin{theorem}
[Clark -- Ocone Formula]\label{t.7.21}Suppose that $F\in\mathcal{D}\left(
\bar{D}\right)  ,$ then\footnote{Here we are abusing notation and writing
$\mathbb{E}\left[  \frac{d}{ds}\left.  \bar{D}F_{s}\left(  b\right)
\right\vert \mathcal{F}_{s}\right]  $ for the \textquotedblleft
predictable\textquotedblright\ projection of the process $s\rightarrow\frac
{d}{ds}\bar{D}F_{s}\left(  b\right)  .$ Since we will only really use Eq.
(\ref{e.7.17}) in these notes, this technicality need not concern us here.}%
\begin{equation}
F=\mathbb{E}F+\int_{0}^{1}\left\langle \mathbb{E}\left[  \frac{d}{ds}\left.
\left(  \bar{D}F\right)  _{s}\left(  b\right)  \right\vert \mathcal{F}%
_{s}\right]  ,db_{s}\right\rangle . \label{e.7.16}%
\end{equation}
In particular if $F=f\left(  b_{s_{1}},\ldots,b_{s_{n}}\right)  $ is a smooth
cylinder function on $W\left(  M\right)  $ then%
\begin{equation}
F=\mathbb{E}F+\int_{0}^{1}\left\langle \mathbb{E}\left[  \left.  \sum
_{i=1}^{n}1_{s\leq s_{i}}\text{$\operatorname*{grad}$}_{i}f\left(  b_{s_{1}%
},\ldots,b_{s_{n}}\right)  \right\vert \mathcal{F}_{s}\right]  ,db_{s}%
\right\rangle . \label{e.7.17}%
\end{equation}

\end{theorem}

\begin{proof}
Let $h$ be a predictable Cameron-Martin valued process such that
$\mathbb{E}\int_{0}^{1}\left\vert h_{s}^{\prime}\right\vert ^{2}ds<\infty.$
Then using Theorem \ref{t.7.15} and the It\^{o} isometry property,%
\begin{align}
\mathbb{E}\langle\bar{D}F,h\rangle_{H}  &  =\mathbb{E}\left[  FD^{\ast
}h\right]  =\mathbb{E}\left[  F\int_{0}^{1}\langle h_{s}^{\prime}%
,db_{s}\rangle\right] \nonumber\\
&  =\mathbb{E}\left[  \left(  \mathbb{E}F+\int_{0}^{1}\langle a,db\rangle
\right)  \int_{0}^{1}\langle h_{s}^{\prime},db_{s}\rangle\right]
=\mathbb{E}\left[  \int_{0}^{1}\langle a_{s},h_{s}^{\prime}\rangle ds\right]
\label{e.7.18}%
\end{align}
where $a$ is the predictable process in Corollary \ref{c.7.20}. Since $h$ is
predictable,%
\begin{align}
\mathbb{E}\langle\bar{D}F,h\rangle_{H}  &  =\mathbb{E}\left[  \int_{0}%
^{1}\left\langle \frac{d}{ds}\left(  \bar{D}F\right)  _{s},h_{s}^{\prime
}\right\rangle ds\right] \nonumber\\
&  =\mathbb{E}\left[  \int_{0}^{1}\left\langle \mathbb{E}\left[  \left.
\frac{d}{ds}\left(  \bar{D}F\right)  _{s}\right\vert \mathcal{F}_{s}\right]
,h_{s}^{\prime}\right\rangle ds\right]  . \label{e.7.19}%
\end{align}
Since $h$ is an arbitrary predictable Cameron-Martin valued process, comparing
Eqs. (\ref{e.7.18}) and (\ref{e.7.18}) shows%
\[
a_{s}=\mathbb{E}\left[  \left.  \frac{d}{ds}\left(  \bar{D}F\right)
_{s}\right\vert \mathcal{F}_{s}\right]
\]
which combined with Eq. (\ref{e.7.12}) completes the proof.
\end{proof}

\begin{remark}
\label{r.7.22}As mentioned in Remark \ref{r.7.18} it is possible to prove Eq.
(\ref{e.7.17}) by an inductive procedure. On the other hand if we were to know
that Eq. (\ref{e.7.17}) was valid for all $F\in\mathcal{F}C^{1}\left(
W\right)  ,$ then for $h\in\mathcal{X}_{a},$%
\begin{align*}
\mathbb{E}\left[  F\int_{0}^{1}\langle h_{s}^{\prime},db_{s}\rangle\right]
&  =\mathbb{E}\left[  \left(  \mathbb{E}F+\int_{0}^{1}\left\langle
\mathbb{E}\left[  \frac{d}{ds}\left.  DF_{s}\right\vert \mathcal{F}%
_{s}\right]  ,db_{s}\right\rangle \right)  \int_{0}^{1}\langle h_{s}^{\prime
},db_{s}\rangle\right] \\
&  =\mathbb{E}\left[  \int_{0}^{1}\left\langle \mathbb{E}\left[  \frac{d}%
{ds}\left.  DF_{s}\right\vert \mathcal{F}_{s}\right]  ,h_{s}^{\prime
}\right\rangle ds\right] \\
&  =\mathbb{E}\left[  \int_{0}^{1}\left\langle \frac{d}{ds}DF_{s}%
,h_{s}^{\prime}\right\rangle ds\right]  =\langle DF,h\rangle_{\mathcal{X}}.
\end{align*}
This identity shows $h\in\mathcal{D}\left(  D^{\ast}\right)  $ and that
$D^{\ast}h=\int_{0}^{1}\langle h_{s}^{\prime},db_{s}\rangle,$ i.e. we have
recovered Theorem \ref{t.7.15}. In this way we see that the Clark-Ocone
formula may be used to recover integration by parts on Wiener space.
\end{remark}

Let $\mathcal{L}$ be the infinite dimensional Ornstein-Uhlenbeck operator
defined as the self-adjoint operator on $L^{2}(\mu)$ given by $\mathcal{L}%
=D^{\ast}\bar{D}.$ The following spectral gap inequality for $\mathcal{L}$ has
been known since the early days of quantum mechanics. This is because
$\mathcal{L}$ is unitarily equivalent to a \textquotedblleft harmonic
oscillator Hamiltonian\textquotedblright\ for which the full spectrum may be
found, see for example \cite{Shigekawa86}. However, these explicit
computations will not in general be available when we consider analogous
spectral gap inequalities when $\mathbb{R}^{d}$ is replaced by a general
compact Riemannian manifold $M.$

\begin{theorem}
[Ornstein Uhlenbeck Spectral Gap Inequality]\label{t.7.23}The null-space of
$\mathcal{L}$ consists of the constant functions on $W$ and $\mathcal{L}$ has
a spectral gap of size $1,$ i.e.%
\begin{equation}
\langle\mathcal{L}F,F\rangle_{L^{2}\left(  \mu\right)  }\geq\langle
F,F\rangle_{L^{2}\left(  \mu\right)  } \label{e.7.20}%
\end{equation}
for all $F\in\mathcal{D}(\mathcal{L})$ such that $F\in\operatorname*{Nul}%
(\mathcal{L})^{\perp}=\left\{  1\right\}  ^{\perp}.$
\end{theorem}

\begin{proof}
Let $F\in\mathcal{D}(\bar{D}),$ then by the Clark-Ocone formula in Eq.
(\ref{e.7.16}), the isometry property of the It\^{o} integral and the
contractive properties of conditional expectation,%
\begin{align*}
\mathbb{E}(F-\mathbb{E}F)^{2}  &  =\mathbb{E}\left[  \int_{0}^{1}\left\langle
\mathbb{E}\left[  \frac{d}{ds}\left.  \bar{D}F_{s}\left(  b\right)
\right\vert \mathcal{F}_{s}\right]  ,db_{s}\right\rangle \right]  ^{2}\\
&  =\mathbb{E}\left[  \int_{0}^{1}\left\vert \mathbb{E}\left[  \frac{d}%
{ds}\left.  \bar{D}F_{s}\left(  b\right)  \right\vert \mathcal{F}_{s}\right]
\right\vert ^{2}ds\right] \\
&  \leq\mathbb{E}\left[  \int_{0}^{1}\left(  \mathbb{E}\left[  \left\vert
\frac{d}{ds}\bar{D}F_{s}\left(  b\right)  \right\vert |\mathcal{F}_{s}\right]
\right)  ^{2}ds\right] \\
&  \leq\mathbb{E}\left[  \int_{0}^{1}\mathbb{E}\left[  \left\vert \frac{d}%
{ds}\bar{D}F_{s}\left(  b\right)  \right\vert ^{2}|\mathcal{F}_{s}\right]
ds\right] \\
&  =\mathbb{E}\left[  \int_{0}^{1}\left\vert \frac{d}{ds}\bar{D}F_{s}\left(
b\right)  \right\vert ^{2}ds\right]  =\langle\bar{D}F,\bar{D}F\rangle
_{\mathcal{X}}.
\end{align*}
In particular if $F\in\mathcal{D}(\mathcal{L}),$ then $\langle\bar{D}F,\bar
{D}F\rangle_{\mathcal{X}}=\mathbb{E}[\mathcal{L}F\cdot F],$ and hence
\begin{equation}
\langle\mathcal{L}F,F\rangle_{L^{2}\left(  \mu\right)  }\geq\langle
F-\mathbb{E}F,F-\mathbb{E}F\rangle_{L^{2}\left(  \mu\right)  }. \label{e.7.21}%
\end{equation}
Therefore, if $F\in\operatorname*{Nul}(\mathcal{L}),$ it follows that
$F=\mathbb{E}F,$ i.e. $F$ is a constant. Moreover if $F\perp1$ (i.e.
$\mathbb{E}F=0)$ then Eq. (\ref{e.7.20}) becomes Eq. (\ref{e.7.21}).
\end{proof}

It turns out that using a method which is attributed to Maurey and Neveu in
\cite{CHL97}, it is possible to use the Clark-Ocone formula as the starting
point for a proof of Gross' logarithmic Sobolev inequality which by general
theory is known to be stronger than the spectral gap inequality in Theorem
\ref{t.7.23}.

\begin{theorem}
[Gross' Logarithmic Sobolev Inequality for $W\left(  \mathbb{R}^{d}\right)  $%
]\label{t.7.24}For all $F\in\mathcal{D}\left(  \bar{D}\right)  ,$%
\begin{equation}
\mathbb{E}\left[  F^{2}\log F^{2}\right]  \leq2\mathbb{E}\left[  \langle
DF,DF\rangle_{H}\right]  +\mathbb{E}F^{2}\cdot\log\mathbb{E}F^{2}.
\label{e.7.22}%
\end{equation}

\end{theorem}

\begin{proof}
Let $F\in\mathcal{F}C^{1}\left(  W\right)  ,$ $\varepsilon>0,$ $H_{\varepsilon
}:=F^{2}+\varepsilon\in\mathcal{D}\left(  \bar{D}\right)  $ and $a_{s}%
=\mathbb{E}\left[  \frac{d}{ds}\left(  DH_{\varepsilon}\right)  _{s}%
|\mathcal{F}_{s}\right]  .$ By Theorem \ref{t.7.21},%
\[
H_{\varepsilon}=\mathbb{E}H_{\varepsilon}+\int_{0}^{1}\langle a,db\rangle
\]
and hence%
\[
M_{s}:=\mathbb{E}\left[  H_{\varepsilon}|\mathcal{F}_{s}\right]
=\mathbb{E}\left[  F^{2}+\varepsilon|\mathcal{F}_{s}\right]  \geq\varepsilon
\]
is a positive martingale which may be written as%
\[
M_{s}:=M_{0}+\int_{0}^{s}\langle a,db\rangle
\]
where $M_{0}=\mathbb{E}H_{\varepsilon}.$

Let $\phi\left(  x\right)  =x\ln x$ so that $\phi^{\prime}\left(  x\right)
=\ln x+1$ and $\phi^{\prime\prime}\left(  x\right)  =x^{-1}.$ Then by
It\^{o}'s formula,
\begin{align*}
d\left[  \phi\left(  M_{s}\right)  \right]   &  =\phi\left(  M_{0}\right)
+\phi^{\prime}\left(  M_{s}\right)  dM_{s}+\frac{1}{2}\phi^{\prime\prime
}\left(  M_{s}\right)  \left\vert a_{s}\right\vert ^{2}ds\\
&  =\phi\left(  M_{0}\right)  +\phi^{\prime}\left(  M_{s}\right)  dM_{s}%
+\frac{1}{2}\frac{1}{M_{s}}\left\vert a_{s}\right\vert ^{2}ds.
\end{align*}
Integrating this equation on $s$ and then taking expectations shows%
\begin{equation}
\mathbb{E}\left[  \phi\left(  M_{1}\right)  \right]  =\phi\left(
\mathbb{E}M_{1}\right)  +\frac{1}{2}\mathbb{E}\left[  \int_{0}^{1}\frac
{1}{M_{s}}\left\vert a_{s}\right\vert ^{2}ds\right]  . \label{e.7.23}%
\end{equation}
Since $\bar{D}H_{\varepsilon}=2F\bar{D}F,$ Eq. (\ref{e.7.23}) is equivalent to%
\[
\mathbb{E}\left[  \phi\left(  H_{\varepsilon}\right)  \right]  =\phi\left(
\mathbb{E}H_{\varepsilon}\right)  +\frac{1}{2}\mathbb{E}\left[  \int_{0}%
^{1}\frac{1}{\mathbb{E}\left[  H_{\varepsilon}|\mathcal{F}_{s}\right]
}\left\vert \mathbb{E}\left[  2F\left(  \bar{D}F\right)  _{s}^{\prime
}|\mathcal{F}_{s}\right]  \right\vert ^{2}ds\right]  .
\]
Using the Cauchy-Schwarz inequality and the contractive properties of
conditional expectations,
\begin{align*}
\left\vert \mathbb{E}\left[  2F\frac{d}{ds}\left(  \bar{D}F\right)
_{s}|\mathcal{F}_{s}\right]  \right\vert ^{2}  &  \leq4\left(  \mathbb{E}%
\left[  F\left\vert \frac{d}{ds}\left(  \bar{D}F\right)  _{s}\right\vert
|\mathcal{F}_{s}\right]  \right)  ^{2}\\
&  \leq4\mathbb{E}\left[  F^{2}|\mathcal{F}_{s}\right]  \cdot\mathbb{E}\left[
\left\vert \frac{d}{ds}\left(  \bar{D}F\right)  _{s}\right\vert ^{2}%
|\mathcal{F}_{s}\right]  .
\end{align*}
Combining the last two equations, using
\begin{equation}
\frac{\mathbb{E}\left[  F^{2}|\mathcal{F}_{s}\right]  }{\mathbb{E}\left[
H_{\varepsilon}|\mathcal{F}_{s}\right]  }=\frac{\mathbb{E}\left[
F^{2}|\mathcal{F}_{s}\right]  }{\mathbb{E}\left[  F^{2}|\mathcal{F}%
_{s}\right]  +\varepsilon}\leq1 \label{e.7.24}%
\end{equation}
gives,%
\begin{align*}
\mathbb{E}\left[  \phi\left(  H_{\varepsilon}\right)  \right]   &  \leq
\phi\left(  \mathbb{E}H_{\varepsilon}\right)  +2\mathbb{E}\int_{0}%
^{1}\mathbb{E}\left[  \left\vert \frac{d}{ds}\left(  \bar{D}F\right)
_{s}\right\vert ^{2}|\mathcal{F}_{s}\right]  ds\\
&  =\phi\left(  \mathbb{E}H_{\varepsilon}\right)  +2\mathbb{E}\int_{0}%
^{1}\left\vert \frac{d}{ds}\left(  \bar{D}F\right)  _{s}\right\vert ^{2}ds.
\end{align*}
We may now let $\varepsilon\downarrow0$ in this inequality to find Eq.
(\ref{e.7.22}) is valid for $F\in\mathcal{F}C^{1}\left(  W\right)  .$ Since
$\mathcal{F}C^{1}\left(  W\right)  $ is a core for $\bar{D},$ standard
limiting arguments show that Eq. (\ref{e.7.22}) is valid in general.

The main objective for the rest of this section is to generalize the previous
theorems to the setting of general compact Riemannian manifolds. Before doing
this we need to record the stochastic analogues of the differentiation formula
in Theorems \ref{t.4.7}, \ref{t.4.12}, and \ref{t.4.13}.
\end{proof}

\subsection{Differentials of Stochastic Flows and Developments\label{s.7.2}}

\begin{notation}
\label{n.7.25}Let $T_{s}^{\beta}\left(  m\right)  =\Sigma_{s}$ where
$\Sigma_{s}$ is the solution to Eq. (\ref{e.5.1}) with $\Sigma_{0}=m$ and
$\beta_{s}$ is an $\mathbb{R}^{n}$ -- valued semi-martingale, i.e.%
\[
\delta\Sigma_{s}=\sum_{i=1}^{n}X_{i}\left(  \Sigma_{s}\right)  \delta\beta
_{s}^{i}+X_{0}\left(  \Sigma_{s}\right)  ds\text{ with }\Sigma_{0}=m.
\]

\end{notation}

\begin{theorem}
[Differentiating $\Sigma$ in $B$]\label{t.7.26}Let $B_{s}$ be an
$\mathbb{R}^{n}$ -- valued Brownian motion and $h$ be an adapted
Cameron-Martin process, $h_{s}\in\mathbb{R}^{n}$ with $\left\vert
h_{s}^{\prime}\right\vert $ bounded.
%(BRUCE: this should be relaxed.)
Then there is a version of $T_{s}^{B+th}\left(  m\right)  $ which is
continuous in $s$ and differentiable in $\left(  t,m\right)  .$ Moreover if we
define $\partial_{h}T_{s}^{B}\left(  o\right)  :=\frac{d}{ds}|_{0}T_{s}%
^{B+sh}\left(  o\right)  ,$ then
\begin{equation}
\partial_{h}T_{s}^{B}\left(  o\right)  =Z_{s}\int_{0}^{s}Z_{r}^{-1}%
X_{h_{r}^{\prime}}\left(  \Sigma_{r}\right)  dr=//_{s}z_{s}\int_{0}^{s}%
z_{r}^{-1}//_{r}^{-1}X_{h_{r}^{\prime}}\left(  \Sigma_{r}\right)
dr\label{e.7.25}%
\end{equation}
where $Z_{s}:=\left(  T_{s}^{B}\right)  _{\ast o},$ $//_{s}$ is stochastic
parallel translation along $\Sigma,$ and $z_{s}:=//_{s}^{-1}Z_{s}.$ (See
Theorem \ref{t.5.41} for more on the processes $Z$ and $z.)$ Recall from
Notation \ref{n.5.4} that
\[
X_{a}\left(  m\right)  :=\sum_{i=1}^{n}a_{i}X_{i}\left(  m\right)
=\mathbf{X}\left(  m\right)  a.
\]

\end{theorem}

\begin{proof}
This is a stochastic analogue of Theorem \ref{t.4.7}. Formally, if $B_{s}$
were piecewise differentiable it would follow from Theorem \ref{t.4.7} with
$s=t,$%
\[
X_{s}\left(  m\right)  =\mathbf{X}\left(  m\right)  B_{s}^{\prime}%
+X_{0}\left(  m\right)  \text{ and }Y_{s}\left(  m\right)  =\mathbf{X}\left(
m\right)  h_{s}^{\prime}.
\]
(Notice that $\frac{d}{dt}|_{0}\left[  \mathbf{X}\left(  m\right)  \left(
B_{s}^{\prime}+th_{s}^{\prime}\right)  +X_{0}\left(  m\right)  \right]
=Y_{s}.)$ For a rigorous proof of this theorem in the flat case, which is
essentially applicable here because of $M$ is an imbedded submanifold, see
Bell \cite{Bell5a} or Nualart \cite{Nu1} for example. For this theorem in this
geometric context see Bismut \cite{Bismut81} or Driver \cite{D5} for example.
\end{proof}

\begin{notation}
\label{n.7.27}Let $b$ be an $T_{o}M\cong\mathbb{R}^{d}$ -- valued Brownian
motion. A $T_{o}M$ -- valued semi-martingale $Y$ is called an \textbf{adapted
vector field }or \textbf{tangent process }to $b$ if $Y$ can be written as
\begin{equation}
Y_{s}=\int_{0}^{s}q_{r}db_{r}+\int_{0}^{s}\alpha_{r}dr\label{e.7.26}%
\end{equation}
where $q_{r}$ is an $so(d)$ -- valued adapted process and $\alpha_{s}$ is a
$T_{o}M$ such that
\[
\int_{0}^{1}|\alpha_{s}|^{2}ds<\infty\text{ a.e.}%
\]

\end{notation}

A key point of a tangent process $Y$ as above is that it gives rise to natural
perturbations of the underlying Brownian motion $b.$ Namely, following Bismut
(also see Fang and Malliavin \cite{FM2}), for $t\in\mathbb{R}$ let $b_{s}^{t}$
be the process given by:
\begin{equation}
b_{s}^{t}:=\int_{0}^{s}e^{tq_{r}}b_{r}+t\int_{0}^{s}\alpha_{r}dr.
\label{e.7.27}%
\end{equation}
Then (under some integrability restrictions on $\alpha)$ by L\'{e}vy's
criteria and Girsanov's theorem, the law of $b^{t}$ is absolutely continuous
relative to the law of $b.$ Moreover $b^{0}=b$ and, with some additional
integrability assumptions on $q_{r},$ $\frac{d}{dt}|_{0}b^{t}=Y.$

Let $b$ be an $T_{o}M\cong\mathbb{R}^{d}$ -- valued Brownian motion,
$\Sigma:=\phi\left(  b\right)  $ be the stochastic development map as in
Notation \ref{n.5.30} and suppose that $X^{h}=//h$ is a Cameron-Martin vector
field on $W\left(  M\right)  .$ Using Theorem \ref{t.4.12} as motivation (see
Eq. (\ref{e.4.16})), the pull back of $X$ under the stochastic development map
should be the process $Y$ defined by%
\begin{equation}
Y_{s}=h_{s}+\int_{0}^{s}\left(  \int_{0}^{r}R_{//_{\rho}}(h_{\rho},\delta
b_{\rho})\right)  \delta b_{r} \label{e.7.28}%
\end{equation}
where
\begin{equation}
R_{//_{s}}(h_{s},\delta b_{s})=//_{s}^{-1}R(//_{s}h_{s},//_{s}\delta
b_{s})//_{s} \label{e.7.29}%
\end{equation}
like in Eq. (\ref{e.4.15}). Since%
\begin{align*}
\left(  \int_{0}^{r}R_{//_{\rho}}(h_{\rho},\delta b_{\rho})\right)  \delta
b_{r}  &  =\left(  \int_{0}^{r}R_{//_{\rho}}(h_{\rho},\delta b_{\rho})\right)
db_{r}+\frac{1}{2}R_{//_{\rho}}(h_{\rho},db_{\rho})db_{\rho}\\
&  =\left(  \int_{0}^{r}R_{//_{\rho}}(h_{\rho},\delta b_{\rho})\right)
db_{r}+\frac{1}{2}\sum_{i=1}^{d}R_{//_{\rho}}(h_{\rho},e_{i})e_{i}d\rho
\end{align*}
where $\left\{  e_{i}\right\}  _{i=1}^{d}$ is an orthonormal basis for
$T_{o}M,$ Eq. (\ref{e.7.28}) may be written in It\^{o} form as%
\begin{equation}
Y_{\cdot}=\int_{0}^{\cdot}C_{s}db_{s}+\int_{0}^{\cdot}r_{s}ds, \label{e.7.30}%
\end{equation}
where
\begin{align}
C_{s}  &  :=\int_{0}^{s}R_{//_{\sigma}}(h_{\sigma},\delta b_{\sigma}),\text{
}\ r_{s}=h_{s}^{\prime}+{\frac{1}{2}}\operatorname{Ric}_{//_{s}}h_{s}\text{
and}\label{e.7.31}\\
\operatorname{Ric}_{//_{s}}a  &  :=//_{s}^{-1}\operatorname*{Ric}//_{s}a\text{
}\forall\ a\in T_{o}M. \label{e.7.32}%
\end{align}
By the symmetry property in item 4b of Proposition \ref{p.3.36}, the matrix
$C_{s}$ is skew symmetric and therefore $Y$ is a tangent process. Here is a
theorem which relates $Y$ in Eq. (\ref{e.7.30}) to $X^{h}=//h.$

\begin{theorem}
[Differential of the development map]\label{t.7.28}Assume $M$ is compact
manifold, $o\in M$ is fixed, $b$ is $T_{o}M\cong\mathbb{R}^{d}$ -- valued
Brownian motion, $\Sigma:=\phi\left(  b\right)  ,$ $h$ is a Cameron-Martin
process with $\left\vert h_{s}^{\prime}\right\vert \leq K<\infty$ ($K$ is a
non-random constant) and $Y$ is as in Eq. (\ref{e.7.30}). As in Eq.
(\ref{e.7.27}) let%
\begin{equation}
b_{s}^{t}:=\int_{0}^{s}e^{tC_{r}}db_{r}+t\int_{0}^{s}r_{r}dr. \label{e.7.33}%
\end{equation}
Then there exists a version of $\phi_{s}\left(  b^{t}\right)  $ which is
continuous in $(s,t),$ differentiable in $t$ and $\frac{d}{dt}|_{0}\phi\left(
b^{t}\right)  =X^{h}.$
\end{theorem}

\begin{proof}
For the proof of this theorem and its generalization to more general $h,$ the
reader is referred to Section 3.1 of \cite{D12} and to \cite{D5}. Let me just
point out here that formally the proof is very analogous to the deterministic
version in Theorems \ref{t.4.12} and \ref{t.4.13}.
\end{proof}

\subsection{Quasi -- Invariance Flow Theorem for $W\left(  M\right)
$\label{s.7.3}}

In this section, we will discuss the $W\left(  M\right)  $ analogues of
Theorems \ref{t.7.13} and \ref{t.7.14}.

\begin{figure}
[ptbh]
\begin{center}
\includegraphics[
height=2.594in,
width=3.9223in
]%
{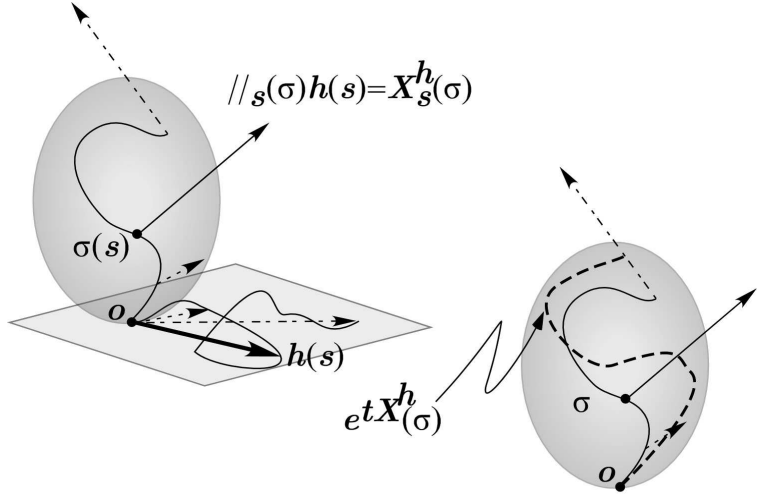}%
\caption{Constructing a vector field, $X^{h},$ on $W\left(  M\right)  $ from a
vector field $h$ on $W\left(  T_{o}M\right)  $. The dotted path indicates the
flow of $\sigma$ under this vector field.}%
\label{fig.14}%
\end{center}
\end{figure}

\begin{theorem}
[Cameron-Martin Theorem for $M$]\label{t.7.29}Let $h\in H(T_{o}M)$ and $X^{h}$
be the $\mu_{W(M)}$ -- a.e. well defined vector field on $W(M)$ given by%
\begin{equation}
X_{s}^{h}(\sigma)=//_{s}(\sigma)h_{s}\text{ for }s\in\lbrack0,1],
\label{e.7.34}%
\end{equation}
where $//_{s}\left(  \sigma\right)  $ is stochastic parallel translation along
$\sigma\in W\left(  M\right)  .$ Then $X^{h}$ admits a flow $e^{tX^{h}}$ on
$W(M)$ (see Figure \ref{fig.14}) and this flow leaves the Wiener measure,
$\mu_{W(M)},$ quasi-invariant.
\end{theorem}

%

%TCIMACRO{\FRAME{ftbphFU}{3.9223in}{2.594in}{0pt}{\Qcb{Constructing a vector
%field, $X^{h},$ on $W\left(  M\right)  $ from a vector field $h$ on $W\left(
%T_{o}M\right)  $. The dotted path indicates the flow of $\sigma$ under this
%vector field.}}{\Qlb{fig.14}}{flow.eps}{\special{ language "Scientific Word";
%type "GRAPHIC";  maintain-aspect-ratio TRUE;  display "USEDEF";
%valid_file "F";  width 3.9223in;  height 2.594in;  depth 0pt;
%original-width 7.3475in;  original-height 4.8483in;  cropleft "0";
%filename 'flow.eps';file-properties "XNPEU";}}}%
%BeginExpansion
%EndExpansion

This theorem first appeared in Driver \cite{D5} for $h\in H\left(
T_{o}M\right)  \cap C^{1}([0,1],T_{o}M)$ and was soon extended to all $h\in
H\left(  T_{o}M\right)  $ by E. Hsu \cite{Hsu95c,Hsu95d}. Other proofs may
also be found in \cite{ES1,Lyons96b,No3}. The proof of this theorem is rather
involved and will not be given here. A sketch of the argument and more
information on the technicalities involved may be found in \cite{D9}.

\begin{example}
\label{ex.7.30}When $M=\mathbb{R}^{d},$ $//_{s}(\sigma)v_{o}=v_{\sigma_{s}}$
for all $v\in\mathbb{R}^{d}$ and $\sigma\in W(\mathbb{R}^{d}).$ Thus
$X_{s}^{h}(\sigma)=(h_{s})_{\sigma_{s}}$ and $e^{tX^{h}}(\sigma)=\sigma+th$
and so Theorem \ref{t.7.29} becomes the classical Cameron-Martin Theorem
\ref{t.7.13}.
\end{example}

\begin{corollary}
[Integration by Parts for $\mu_{W(M)}$]\label{c.7.31}For $h\in H(T_{o}M)$ and
$F\in\mathcal{F}C^{1}(W(M))$ as in Eq. (\ref{e.7.2}), let%
\[
(X^{h}F)(\sigma)=\frac{d}{dt}|_{0}F(e^{tX^{h}}(\sigma))=G\left(
DF,X^{h}\right)
\]
as in Notation \ref{n.7.11}. Then%
\[
\int_{\mathrm{W}(M)}X^{h}F\,d\mu_{W\left(  M\right)  }=\int_{\mathrm{W}%
(M)}F\,z^{h}\,d\mu_{W\left(  M\right)  }%
\]
where
\[
z^{h}:=\int_{0}^{1}\langle h_{s}^{\prime}+\frac{1}{2}\operatorname*{Ric}%
\nolimits_{//_{s}}h_{s}^{\prime},db_{s}\rangle,
\]%
\[
b_{s}\left(  \sigma\right)  :=\Psi_{s}\left(  \sigma\right)  =\int_{0}%
^{s}//_{r}^{-1}\delta\sigma_{r}%
\]
and $\operatorname{Ric}_{//_{s}}\in\operatorname{End}(T_{o}M)$ is as in Eq.
(\ref{e.5.60}).
\end{corollary}

\begin{proof}
A special case of this Corollary \ref{c.7.31} with $F(\sigma)=f(\sigma_{s})$
for some $f\in C^{\infty}(M)$ first appeared in Bismut \cite{Bismut84a}. The
result stated here was proved in \cite{D5} as an infinitesimal form of the
flow Theorem \ref{t.7.29}. Other proofs of this corollary may be found in
\cite{Aida97,Airault96,Driver97b,Elworthy94a,Elworthy96a,Elworthy96c,ES1,FM2,Hsu95c,Hsu95d,Leandre93b,Leandre95e,Lyons96b,No3}%
. This corollary is a special case of Theorem \ref{t.7.32} below.
\end{proof}

\subsection{Divergence and Integration by Parts.}

In the next theorem, it will be shown that adapted Cameron-Martin vector
fields, $X,$ are in the domain of $D^{\ast}$ and consequently $D^{\ast}$ is
densely defined. For the purposes of this subsection, we assume that $b$ is a
$T_{o}M$ -- valued Brownian motion, $\Sigma=\phi\left(  b\right)  $ is the
evolved Brownian motion on $M$ and $//_{s}$ is stochastic parallel translation
along $\Sigma.$

\begin{theorem}
\label{t.7.32}Let $X\in\mathcal{X}_{a}$ be an adapted Cameron-Martin vector
field on $W(M)$ and $h:=//^{-1}X.$ Then $X\in\mathcal{D}(D^{\ast})$ and
\begin{equation}
X^{\ast}1=D^{\ast}X=\int_{0}^{1}\langle B(h),db\rangle=\int_{0}^{1}\langle
h_{s}^{\prime}+{\frac{1}{2}}\operatorname*{Ric}\nolimits_{//_{s}}h_{s}%
,db_{s}\rangle, \label{e.7.35}%
\end{equation}
where $B$ is the random linear operator mapping $H$ to $L^{2}(ds,T_{o}M)$
given by
\begin{equation}
\left[  B(h)\right]  _{s}:=h_{s}^{\prime}+{\frac{1}{2}}\operatorname*{Ric}%
\nolimits_{//_{s}}h_{s}. \label{e.7.36}%
\end{equation}

\end{theorem}

\begin{remark}
\label{r.7.33}There is a non-random constant $C<\infty$ depending only on the
bound on the Ricci tensor such that $\left\Vert B\right\Vert _{H\rightarrow
L^{2}(ds,T_{o}M)}\leq C.$
\end{remark}

\begin{proof}
I will give a sketch of the proof here, the interested reader may find
complete details of this proof in \cite{D12}. Moreover, we will give two more
proofs of this theorem, see Theorem \ref{t.7.40} and Corollary \ref{c.7.50} below.

We start by proving the theorem under the additional assumption that
$h:=//^{-1}X$ satisfies $\sup_{s\in\lbrack0,1]}\left\vert h_{s}^{\prime
}\right\vert \leq K,$ where $K$ is a non-random constant.

Let $b_{s}^{t}$ be defined as in Eq. (\ref{e.7.33}). (Notice that $b^{t}$ is
\textbf{not }the flow of the vector-field $Y$ in Eq. (\ref{e.7.30}) but does
have the property that $\frac{d}{dt}|_{0}b_{s}^{t}=Y_{s}.)$ Since $C_{s}$ is
skew-symmetric, $e^{tC_{s}}$ is orthogonal and so by Levy's criteria,
$s\rightarrow\int_{0}^{s}e^{tC_{r}}db_{r}$ is a Brownian motion. Combining
this with Girsanov's theorem, $s\rightarrow b_{s}^{t}$ (for fixed $t)$ is a
Brownian motion relative to the measure $Z_{t}\cdot\mu,$ where
\begin{equation}
Z_{t}:=\exp\left(  -\int_{0}^{1}t\langle r,e^{tC}db\rangle-{\frac{1}{2}}%
t^{2}\int_{0}^{1}\langle r,r\rangle ds\right)  .\label{e.7.37}%
\end{equation}
For $t\in\mathbb{R},$ let $\Sigma(t,\cdot):=\phi(b^{t})$ where $\phi$ is the
stochastic development map as in Theorem \ref{t.5.29}. Then by Theorem
\ref{t.7.28}, $X^{h}=\frac{d}{dt}|_{0}\Sigma(t,\cdot)$ and in particular if
$F$ is a smooth cylinder function then $X^{h}F=\frac{d}{dt}|_{0}%
F(\Sigma(t,\cdot)).$ So differentiating the identity,
\[
\mathbb{E}\left[  F(\Sigma(t,\cdot)Z_{t}\right]  =\mathbb{E}\left[  F\left(
\Sigma\right)  \right]  ,
\]
at $t=0$ gives:
\[
\mathbb{E}\left[  XF\right]  -\mathbb{E}\left[  F\int_{0}^{1}\langle
r,db\rangle\right]  =0.
\]
This last equation may be written alternatively as
\[
\langle DF,X\rangle_{\mathcal{X}}=\mathbb{E}\left[  G\left(  DF,X\right)
\right]  =\mathbb{E}\left[  F\cdot\int_{0}^{1}\langle B(h),db\rangle\right]  .
\]
Hence it follows that $X\in\mathcal{D}(D^{\ast})$ and
\[
D^{\ast}X=\int_{0}^{1}\langle B(h),db\rangle.
\]
This proves the theorem in the special case that $h^{\prime}$ is uniformly bounded.

Let $X$ be a general adapted Cameron-Martin vector-field and $h:=//^{-1}X.$
For each $n\in\mathbb{N},$ let $h_{n}\left(  s\right)  :=\int_{0}^{s}%
h^{\prime}(r\mathbb{)}\cdot1_{|h^{\prime}(r\mathbb{)}|\leq n}dr$ be as in Eq.
(\ref{e.7.11}). Set $X^{n}:=//h_{n},$ then by the special case above we know
that $X^{n}\in\mathcal{D}(D^{\ast})$ and $D^{\ast}X^{n}=\int_{0}^{1}\langle
B(h_{n}),db\rangle.$ It is easy to check that
\[
\langle X-X^{n},X-X^{n}\rangle_{\mathcal{X}}=\mathbb{E}\langle h-h_{n}%
,h-h_{n}\rangle_{H}\rightarrow0\text{ as }n\rightarrow\infty.
\]
Furthermore,
\[
\mathbb{E}\left\vert D^{\ast}(X^{m}-X^{n})\right\vert ^{2}=\mathbb{E}\int
_{0}^{1}|B(h_{m}-h_{n})|^{2}ds\leq C\mathbb{E}\langle h_{m}-h_{n},h_{m}%
-h_{n}\rangle_{H},
\]
from which it follows that $D^{\ast}X^{m}$ is convergent. Because $D^{\ast}$
is a closed operator, it follows that $X\in\mathcal{D}(D^{\ast})$ and
\[
D^{\ast}X=\lim_{n\rightarrow\infty}D^{\ast}X^{n}=\lim_{n\rightarrow\infty}%
\int_{0}^{1}\langle B(h_{n}),db\rangle=\int_{0}^{1}\langle B(h),db\rangle.
\]

\end{proof}

\begin{corollary}
\label{c.7.34}The operator $D^{\ast}:\mathcal{X}\rightarrow L^{2}\left(
W\left(  M\right)  ,\mu_{W\left(  M\right)  }\right)  $ is densely defined. In
particular $D$ is closable. (Let $\bar{D}$ denote the closure of $D.)$
\end{corollary}

\begin{proof}
Let $h\in H,$ $X^{h}:=//h,$ and $F$ and $K$ be smooth cylinder functions.
Then, by the product rule,%
\begin{align*}
\langle DF,KX^{h}\rangle_{\mathcal{X}}  &  =\mathbb{E}[G\left(  KDF,X^{h}%
\right)  ]=\mathbb{E}[G\left(  D\left(  KF\right)  -FDK,X^{h}\right)  ]\\
&  =\mathbb{E}[F\cdot KD^{\ast}X^{h}-FG\left(  DK,X^{h}\right)  ].
\end{align*}
Therefore $KX^{h}\in\mathcal{D}(D^{\ast})$ ($\mathcal{D}(D^{\ast})$ is the
domain of $D^{\ast})$ and
\[
D^{\ast}(KX^{h})=KD^{\ast}X^{h}-G(DK,X^{h}).
\]
Since
\[
\mathrm{span}\{KX^{h}|h\in H\text{\textrm{\ and }}K\in\mathcal{F}C^{\infty
}\}\subset\mathcal{D}(D^{\ast})
\]
is is a dense subspace of $\mathcal{X},$ $D^{\ast}$ is densely defined.
\end{proof}

\begin{corollary}
\label{c.7.35}Let $h$ be an adapted Cameron-Martin valued process and $Q_{s}$
be defined as in Eq. (\ref{e.6.1}). Then
\begin{equation}
\left(  X^{Q^{\mathrm{\mathrm{\mathrm{tr}}}}h}\right)  ^{\ast}1=\int_{0}%
^{1}\langle Q^{\mathrm{\mathrm{\mathrm{tr}}}}h^{\prime},db\rangle.
\label{e.7.38}%
\end{equation}

\end{corollary}

\begin{proof}
Taking the transpose of Eq. (\ref{e.6.1}) shows $Q^{\mathrm{tr}}$ solves,%
\begin{equation}
\frac{d}{ds}Q^{\mathrm{tr}}+\frac{1}{2}\operatorname{Ric}_{//}Q^{\mathrm{tr}%
}=0\text{ with }Q_{0}^{\mathrm{tr}}=Id. \label{e.7.39}%
\end{equation}
Therefore, from Eq. (\ref{e.7.35}),
\begin{align*}
\left(  X^{Q^{\mathrm{\mathrm{\mathrm{tr}}}}h}\right)  ^{\ast}1  &  =\int
_{0}^{1}\langle\left(  Q^{\mathrm{\mathrm{\mathrm{tr}}}}h\right)  ^{\prime
}+{\frac{1}{2}}\operatorname{Ric}_{//}Q^{\mathrm{\mathrm{\mathrm{tr}}}%
}h,db\rangle\\
&  =\int_{0}^{1}\left\langle \left[  \frac{d}{ds}+\frac{1}{2}%
\operatorname{Ric}_{//}\right]  \left(  Q^{\mathrm{\mathrm{\mathrm{tr}}}%
}h\right)  ,db\right\rangle \\
&  =\int_{0}^{1}\langle Q^{\mathrm{\mathrm{\mathrm{tr}}}}h^{\prime},db\rangle.
\end{align*}

\end{proof}

Theorem \ref{t.7.32}\ may be extended to allow for vector-fields on the paths
of $M$ which are not based. This theorem and it Corollary \ref{c.7.37} will
not be used in the sequel and may safely be skipped.

\begin{theorem}
\label{t.7.36}Let $h$ be an adapted $T_{o}M$ -- valued process such that
$h(0)$ is non-random and $h-h(0)$ is a Cameron-Martin process, $X:=X^{h}%
:=//h,$ $\mathbb{E}_{x}$ denote the path space expectation for a Brownian
motion starting at $x\in M,$ $F:C([0,1]\rightarrow M)\rightarrow\mathbb{R}$ be
a cylinder function as in Definition \ref{d.7.4} and $X^{h}F$ be defined as in
Eq. \text{(\ref{e.7.7})}. Then (writing $\langle df,v\rangle$ for $df\left(
v\right)  )$
\begin{equation}
\mathbb{E}_{o}[X^{h}F]=\mathbb{E}_{o}[FD^{\ast}X^{h}]+\langle d(\mathbb{E}%
_{(\cdot)}F),h(0)_{o}\rangle, \label{e.7.40}%
\end{equation}
where
\[
D^{\ast}X^{h}:=\int_{0}^{1}\langle h_{s}^{\prime}+{\frac{1}{2}}%
\operatorname{Ric}_{//_{s}}h_{s},db_{s}\rangle:=\int_{0}^{1}\langle
B(h),db\rangle,
\]
as in Eq. \text{(\ref{e.7.35})}\ and $B(h)$ is defined in Eq.
\text{(\ref{e.7.36})}.
\end{theorem}

\begin{proof}
Start by choosing a smooth path $\alpha$ in $M$ such that $\dot{\alpha
}(0)=h(0)_{o}.$ Let%
\begin{align*}
C  &  :=\int R_{//}(h,\delta b),\\
r  &  =h^{\prime}+{\frac{1}{2}}\operatorname{Ric}_{//}(h),\\
b_{s}^{t}  &  =\int_{0}^{s}e^{tC}db+t\int_{0}^{s}rd\lambda\text{ and}\\
Z_{t}  &  =\exp-\left\{  \int_{0}^{1}t\langle r,e^{tC}db\rangle+{\frac{1}{2}%
}t^{2}\int_{0}^{1}\langle r,r\rangle ds\right\}
\end{align*}
be defined by the same formulas as in the proof of Theorem \ref{t.7.32}. Let
$u_{0}(t)$ denote parallel translation along $\alpha,$ that is
\[
du_{0}(t)/dt+\Gamma(\dot{\alpha}(t))u_{0}(t)=0\quad\text{\textrm{\ with }%
}\quad u_{0}(0)=id.
\]
For $t\in\mathbb{R},$ define $\Sigma(t,\cdot)$ by
\[
\Sigma(t,\delta s)=u(t,s)\delta b_{s}^{t}\quad\text{\textrm{\ with }}%
\quad\Sigma(t,0)=\alpha(t)
\]
and
\[
u(t,\delta s)+\Gamma(u(t,s)\delta_{s}b_{s}^{t})u(t,s)=0\quad
\text{\textrm{\ with }}\quad u(t,0)=u_{o}(t).
\]
Appealing to a stochastic version of Theorem \ref{t.4.14} (after choosing a
good version of $\Sigma$) it is possible to show that $\dot{\Sigma}%
(0,\cdot)=X,$ so the $XF=\frac{d}{dt}|_{0}F\left[  \Sigma(t,\cdot)\right]  .$
As in the proof of Theorem \ref{t.7.32}, $b^{t}$ is a Brownian motion relative
to the expectation $\mathbb{E}_{t}$ defined by $\mathbb{E}_{t}(F):=\mathbb{E}%
\left[  Z_{t}F\right]  $. From this it is easy to see that $\Sigma(t,\cdot)$
is a Brownian motion on $M$ starting at $\alpha(t)$ relative to the
expectation $\mathbb{E}_{t}.$ Therefore, for all $t,$
\[
\mathbb{E}\left[  F\left(  \Sigma(t,\cdot)\right)  Z_{t}\right]
=\mathbb{E}_{\alpha(t)}F
\]
and differentiating this last expression at $t=0$ gives:
\[
\mathbb{E}\left[  XF(\Sigma)\right]  -\mathbb{E}\left[  F\int_{0}^{1}\langle
r,db\rangle\right]  =\langle d\mathbb{E}_{(\cdot)}F,h(0)_{o}\rangle.
\]
The rest of the proof is identical to the previous proof.
\end{proof}

As a corollary to Theorem \ref{t.7.36} we get Elton Hsu's derivative formula
which played a key role in the original proof of his logarithmic Sobolev
inequality on $W(M),$ see Theorem \ref{t.7.52} below and \cite{Hsu3}.

\begin{corollary}
[Hsu's Derivative Formula]\label{c.7.37} Let $v_{o}\in T_{o}M$. Define $h$ to
be the adapted $T_{o}M$ -- valued process solving the differential equation:
\begin{equation}
h_{s}^{\prime}+{\frac{1}{2}}\operatorname{Ric}_{//_{s}}h_{s}=0\quad
\text{\textrm{\ with }}\quad h_{0}=v_{o}. \label{e.7.41}%
\end{equation}
Then
\begin{equation}
\langle d(\mathbb{E}_{(\cdot)}F),v_{o}\rangle=\mathbb{E}_{o}[X^{h}F].
\label{e.7.42}%
\end{equation}

\end{corollary}

\begin{proof}
Apply Theorem \ref{t.7.36} to $X^{h}$ with $h$ defined by \text{(\ref{e.7.41}%
)}. Notice that $h$ has been constructed so that $B(h)\equiv0,$ i.e. $D^{\ast
}X^{h}=0.$
\end{proof}

The idea for the proof used here is similar Hsu's proof, the only question is
how one describes the perturbed process $\Sigma(t,\cdot)$ in the proof of
Theorem \ref{t.7.36} above. It is also possible to give a much more elementary
proof of Eq. (\ref{e.7.42}) based on the ideas in Section \ref{s.6}, see for
example \cite{Driver01b}.

\subsection{Elworthy-Li Integration by Parts Formula}

In this subsection, let $\left\{  X_{i}\right\}  _{i=0}^{n}\subset
\Gamma\left(  TM\right)  ,$ $B$ be a $\mathbb{R}^{n}$ -- valued Brownian
motion and $T_{s}^{B}\left(  m\right)  $ denote the solution to Eq.
(\ref{e.5.1}) with $\beta=B$ be as in Notation \ref{n.7.25}. We will further
assume that $\mathbf{X}\left(  m\right)  :\mathbb{R}^{n}\rightarrow T_{m}M$
(as in Notation \ref{n.5.4}) is surjective for all $m\in M$ and let
$\mathbf{X}\left(  m\right)  ^{\#}=\left[  \mathbf{X}\left(  m\right)
|_{\mathrm{Nul}(\mathbf{X}\left(  m\right)  )^{\perp}}\right]  ^{-1}$ as in
Eq. (\ref{e.6.5}). The following Lemma is an elementary exercise in linear algebra.

\begin{lemma}
\label{l.7.38}For $m\in M$ and $v,w\in T_{m}M$ let%
\[
\langle v,w\rangle_{m}:=\langle\mathbf{X}\left(  m\right)  ^{\#}%
v,\mathbf{X}\left(  m\right)  ^{\#}w\rangle_{\mathbb{R}^{n}}.
\]
Then

\begin{enumerate}
\item $m\rightarrow\langle\cdot,\cdot\rangle_{m}$ is a smooth Riemannian
metric on $M.$

\item $\mathbf{X}\left(  m\right)  ^{\mathrm{tr}}=\mathbf{X}\left(  m\right)
^{\#}$ and in particular $\mathbf{X}\left(  m\right)  \mathbf{X}\left(
m\right)  ^{\mathrm{tr}}=id_{T_{m}M}$ for all $m\in M.$

\item Every $v\in T_{m}M$ may be expanded as%
\begin{equation}
v=\sum_{j=1}^{n}\langle v,X_{j}\left(  m\right)  \rangle X_{j}\left(
m\right)  =\sum_{j=1}^{n}\langle v,\mathbf{X}\left(  m\right)  e_{j}%
\rangle\mathbf{X}\left(  m\right)  e_{j} \label{e.7.43}%
\end{equation}
where $\left\{  e_{j}\right\}  _{j=1}^{n}$ is the standard basis for
$\mathbb{R}^{n}.$
\end{enumerate}
\end{lemma}

The proof of this lemma is left to the reader with the comment that Eq.
(\ref{e.7.43}) is proved in the same manner as item (1) in Proposition
\ref{p.3.48}.

\begin{theorem}
[Elworthy - Li]\label{t.7.39}Suppose $k_{s}$ is a $T_{o}M$ valued
Cameron-Martin process such that $\mathbb{E}\int_{0}^{1}\left\vert
k_{s}^{\prime}\right\vert ^{2}ds<\infty$ and $F:W\left(  M\right)
\rightarrow\mathbb{R}$ is a bounded $C^{1}$ -- function with bounded
derivative on $W,$ for example $F$ could be a cylinder function. Then%
\begin{align}
\mathbb{E}\left[  \left(  d_{W\left(  M\right)  }F\right)  \left(  Z_{\cdot
}k_{\cdot}\right)  \right]   &  =\mathbb{E}\left[  F\left(  \Sigma\right)
\int_{0}^{T}\langle Z_{s}k_{s}^{\prime},\mathbf{X}\left(  \Sigma_{s}\right)
dB_{s}\rangle\right] \nonumber\\
&  =\mathbb{E}\left[  F\left(  \Sigma\right)  \int_{0}^{T}\langle
\mathbf{X}\left(  \Sigma_{s}\right)  ^{\mathrm{\mathrm{tr}}}Z_{s}k_{s}%
^{\prime},dB_{s}\rangle\right]  \label{e.7.44}%
\end{align}
where and $Z_{s}=\left(  T_{s}^{B}\right)  _{\ast o}$ is the differential of
$m\rightarrow T_{s}^{B}\left(  m\right)  $ at $o.$
\end{theorem}

\begin{proof}
Notice that $Z_{s}k_{s}\in T_{\Sigma_{s}}M$ for all $s$ as it should be. By a
reduction argument used in the proof of Theorem \ref{t.7.32}, it suffices to
consider the case where $\left\vert k_{s}^{\prime}\right\vert \leq K$ where
$K$ is a non-random constant. Let $h_{s}$ be the $T_{o}M$ -- valued
Cameron-Martin process defined by%
\[
h_{s}:=\int_{0}^{s}\mathbf{X}\left(  \Sigma_{r}\right)  ^{\mathrm{tr}}%
Z_{r}k_{r}^{\prime}dr.
\]
Then by Lemma \ref{l.7.38} and Theorem \ref{t.7.26},
\begin{align*}
\partial_{h}T_{s}^{B}\left(  o\right)   &  =Z_{s}\int_{0}^{s}Z_{r}%
^{-1}\mathbf{X}\left(  \Sigma_{r}\right)  h_{r}^{\prime}dr\\
&  =Z_{s}\int_{0}^{s}Z_{r}^{-1}\mathbf{X}\left(  \Sigma_{r}\right)
\mathbf{X}\left(  \Sigma_{r}\right)  ^{\mathrm{tr}}Z_{r}k_{r}^{\prime}%
dr=Z_{s}k_{s}.
\end{align*}
In particular this implies
\[
\partial_{h}F\left(  T_{\left(  \cdot\right)  }^{B}\left(  o\right)  \right)
=\langle dF\left(  \Sigma\right)  ,\partial_{h}T_{s}^{B}\left(  o\right)
\rangle=\langle d_{W\left(  M\right)  }F\left(  \Sigma\right)  ,Zk\rangle
\]
and therefore by integration by parts on the flat Wiener space (Theorem
\ref{t.7.32}) with $M=\mathbb{R}^{n})$ implies%
\begin{align*}
\mathbb{E}\left[  \left(  d_{W\left(  M\right)  }F\right)  \left(
\Sigma\right)  \left(  Z_{\cdot}k_{\cdot}\right)  \right]   &  =\mathbb{E}%
\left[  \partial_{h}\left[  F\left(  \Sigma\right)  \right]  \right]
=\mathbb{E}\left[  F\left(  \Sigma\right)  \int_{0}^{T}\langle h_{s}^{\prime
},dB_{s}\rangle\right] \\
&  =\mathbb{E}\left[  F\left(  \Sigma\right)  \int_{0}^{T}\langle
\mathbf{X}\left(  \Sigma_{s}\right)  ^{\mathrm{tr}}Z_{s}k_{s}^{\prime}%
,dB_{s}\rangle\right]  .
\end{align*}

%What about more general perturbations, k, not just Cameron-Martin valued ones.

\end{proof}

By factoring out the redundant noise in Theorem \ref{t.7.39}, we get yet
another proof of Corollary \ref{c.7.35} which also easily gives another proof
of Theorem \ref{t.7.32}.

\begin{theorem}
[Factoring out the redundant noise]\label{t.7.40}Assume $\mathbf{X}\left(
m\right)  =P\left(  m\right)  $ and $X_{0}=0,$ $k_{s}$ is a Cameron-Martin
valued process adapted to the filtration, $\mathcal{F}_{s}^{\Sigma}%
:=\sigma\left(  \Sigma_{r}:r\leq s\right)  ,$ then
\[
\mathbb{E}\left[  \left(  d_{W\left(  M\right)  }F\right)  \left(
//Q_{t}^{\mathrm{tr}}k\right)  \right]  =\mathbb{E}\left[  F\left(
\Sigma\right)  \int_{0}^{T}\langle//_{s}Q_{s}^{\mathrm{tr}}k_{s}^{\prime
},db_{s}\rangle\right]
\]
where $Q_{s}$ solves Eq. (\ref{e.6.1}).
\end{theorem}

\begin{proof}
Using%
\begin{align*}
\mathbb{E}\left[  \left(  d_{W\left(  M\right)  }F\right)  \left(
//zk\right)  \right]   &  =\mathbb{E}\left[  F\left(  \Sigma\right)  \int
_{0}^{T}\langle//_{s}z_{s}k_{s}^{\prime},P\left(  \Sigma_{s}\right)
dB_{s}\rangle\right] \\
&  =\mathbb{E}\left[  F\left(  \Sigma\right)  \int_{0}^{T}\langle//_{s}%
z_{s}k_{s}^{\prime},db_{s}\rangle\right]
\end{align*}
along with Theorem \ref{t.5.44} we find%
\[
\mathbb{E}\left[  \left(  d_{W\left(  M\right)  }F\right)  \left(  //\bar
{z}k\right)  \right]  =\mathbb{E}\left[  F\left(  \Sigma\right)  \int_{0}%
^{T}\langle//_{s}\bar{z}_{s}k_{s}^{\prime},db_{s}\rangle\right]  .
\]
As observed in the proof of Corollary \ref{c.6.4}, $\bar{z}_{t}=Q_{t}%
^{\mathrm{tr}}$ which completes the proof.
\end{proof}

The reader interested in seeing more of these type of arguments is referred to
Elworthy, Le Jan and Li \cite{El-LJ-Li} where these ideas are covered in much
greater detail and in full generality.

\subsection{Fang's Spectral Gap Theorem and Proof\label{s.7.6}}

%BRUCE: Should this section be simplified by placing after the Clark Ocone formulas?

As in the flat case we let $\mathcal{L}=D^{\ast}\bar{D}$ -- an unbounded
operator on $L^{2}\left(  W\left(  M\right)  ,\mu_{W\left(  M\right)
}\right)  $ which is a \textquotedblleft curved\textquotedblright\ analogue of
the Ornstein-Uhlenbeck operator used in Theorem \ref{t.7.23}. It has been
shown in Driver and R\"{o}ckner \cite{DR} that this operator generates a
diffusion on $W(M).$ This last result also holds for pinned paths on $M$ and
free loops on $\mathbb{R}^{N},$ see \cite{ALR}.

In this section, we will give a proof of S. Fang's \cite{Fang1} spectral gap
inequality for $\mathcal{L}.$ Hsu's stronger logarithmic Sobolev inequality
will be covered later in Theorem \ref{t.7.52} below.

\begin{theorem}
[Fang]\label{t.7.41}Let $\bar{D}$ be the closure of $D$ and $\mathcal{L}$ be
the self-adjoint operator on $L^{2}\left(  \mu_{W\left(  M\right)  }\right)  $
defined by $\mathcal{L}=D^{\ast}\bar{D}.$ (Note, if $M=\mathbb{R}^{d}$ then
$\mathcal{L}$ would be an infinite dimensional Ornstein-Uhlenbeck operator.)
Then the null-space of $\mathcal{L}$ consists of the constant functions on
$W(M)$ and $\mathcal{L}$ has a spectral gap, i.e. there is a constant $c>0$
such that $\langle\mathcal{L}F,F\rangle_{L^{2}\left(  \mu_{W\left(  M\right)
}\right)  }\geq c\langle F,F\rangle_{L^{2}\left(  \mu_{W\left(  M\right)
}\right)  }$ for all $F\in\mathcal{D}(\mathcal{L})$ which are perpendicular to
the constant functions.
\end{theorem}

This theorem is the $W\left(  M\right)  $ analogue of Theorem \ref{t.7.23}.
The proof of this theorem will be given at the end of this subsection. We
first will need to represent $F$ in terms of $DF.$ (Also see Section
\ref{s.7.7} below.)

\begin{lemma}
\label{l.7.42}For each $F\in L^{2}\left(  W\left(  M\right)  ,\mu_{W\left(
M\right)  }\right)  ,$ there is a unique adapted Cameron-Martin vector field
$X$ on $W(M)$ such that
\[
F=\mathbb{E}F+D^{\ast}X.
\]

\end{lemma}

\begin{proof}
By the martingale representation theorem (see Corollary \ref{c.7.20}), there
is a predictable $T_{o}M$--valued process, $a,$ (which is not in general
continuous) such that
\[
\mathbb{E}\int_{0}^{1}|a_{s}|^{2}ds<\infty,
\]
and
\begin{equation}
F=\mathbb{E}F+\int_{0}^{1}\langle a_{s},db_{s}\rangle. \label{e.7.45}%
\end{equation}
Define $h:=B^{-1}(a),$ i.e. let $h$ be the solution to the differential
equation:
\begin{equation}
h_{s}^{\prime}+{\frac{1}{2}}\operatorname{Ric}_{//_{s}}h_{s}=a_{s}\ \text{with
}h_{0}=0. \label{e.7.46}%
\end{equation}
\medskip\noindent\textbf{Claim:} $B_{\sigma}^{-1}$ is a bounded linear map
from $L^{2}(ds,T_{o}M)\rightarrow H$ for each $\sigma\in W\left(  M\right)  ,$
and furthermore the norm of $B_{\sigma}^{-1}$ is bounded independent of
$\sigma\in W\left(  M\right)  $.

To prove the claim, use Duhamel's principle to write the solution to
\text{(\ref{e.7.46})} as:
\begin{equation}
h_{s}=\int_{0}^{s}Q_{s}^{\mathrm{tr}}\left(  Q_{\tau}^{\mathrm{tr}}\right)
^{-1}a_{\tau}d\tau. \label{e.7.47}%
\end{equation}
Since, $W_{s}:=Q_{s}^{\mathrm{tr}}\left(  Q_{\tau}^{\mathrm{tr}}\right)
^{-1}$ solves the differential equation
\[
W_{s}^{\prime}+A_{s}W_{s}=0\text{ with }W_{\tau}=I
\]
it is easy to show from the boundedness of $\operatorname*{Ric}_{//_{s}}$ and
an application of Gronwall's inequality that%
\[
\left\vert Q_{s}^{\mathrm{tr}}\left(  Q_{\tau}^{\mathrm{tr}}\right)
^{-1}\right\vert =\left\vert W_{s}\right\vert \leq C,
\]
where $C$ is a non-random constant independent of $s$ and $\tau.$ Therefore,
\begin{align*}
\langle h,h\rangle_{H}  &  =\int_{0}^{1}|a_{s}-{\frac{1}{2}}\operatorname{Ric}%
_{//_{s}}h_{s}|^{2}ds\\
&  \leq2\int_{0}^{1}|a_{s}|^{2}ds+2\int_{0}^{1}|{\frac{1}{2}}%
\operatorname{Ric}_{//_{s}}h_{s}|^{2}ds\\
&  \leq2(1+C^{2}K^{2})\int_{0}^{1}|a_{s}|^{2}ds,
\end{align*}
where $K$ is a bound on the process ${\frac{1}{2}}\operatorname*{Ric}%
\nolimits_{//_{s}}.$ This proves the claim.

Because of the claim, $h:=B^{-1}(a)$ satisfies $\mathbb{E}\left[  \langle
h,h\rangle_{H}\right]  <\infty$ and because of Eq. \text{(\ref{e.7.47}), }$h$
is adapted. Hence, $X:=//h$ is an adapted Cameron-Martin vector field and
\[
D^{\ast}X=\int_{0}^{1}\langle B(h),db\rangle=\int_{0}^{1}\langle a,db\rangle.
\]
The existence part of the theorem now follows from this identity and Eq.
\text{(\ref{e.7.45})}.

The uniqueness assertion follows from the energy identity:
\[
\mathbb{E}\left[  D^{\ast}X\right]  ^{2}=\mathbb{E}\int_{0}^{1}|B(h)_{s}%
|^{2}ds\geq C\mathbb{E}\left[  \langle h,h\rangle_{H}\right]  .
\]
Indeed if $D^{\ast}X=0,$ then $h=0$ and hence $X=//h=0.$
\end{proof}

The next goal is to find an expression for the vector-field $X$ in the above
lemma in terms of the function $F$ itself. This will be the content of Theorem
\ref{t.7.45} below.

\begin{notation}
\label{n.7.43}Let $L_{a}^{2}(\mu_{W\left(  M\right)  }:L^{2}(ds,T_{o}M))$
denote the $T_{o}M$ -- valued predictable processes, $v_{s}$ on $W\left(
M\right)  $ such that $\mathbb{E}\int_{0}^{1}\left\vert v_{s}\right\vert
^{2}ds<\infty.$ Define the bounded linear operator $\bar{B}:\mathcal{X}%
_{a}\rightarrow L_{a}^{2}(\mu_{W\left(  M\right)  }:L^{2}(ds,T_{o}M))$ by
\[
\bar{B}(X)=B(//^{-1}X)=\frac{d}{ds}\left[  //_{s}^{-1}X_{s}\right]  +{\frac
{1}{2}}//_{s}^{-1}\operatorname*{Ric}X_{s}.
\]
Also let $\mathcal{\ Q}:\mathcal{X}\rightarrow\mathcal{X}$ denote the
orthogonal projection of $\mathcal{\ X}$ onto$\mathcal{\ X}_{a}.$
\end{notation}

\begin{remark}
\label{r.7.44}Notice that $D^{\ast}X=\int_{0}^{1}\langle\bar{B}(X),db\rangle$
for all $X\in\mathcal{X}_{a}.$ We have seen that $\bar{B}$ has a bounded
inverse, in fact $\bar{B}^{-1}(a)=//B^{-1}(a).$
\end{remark}

\begin{theorem}
\label{t.7.45} As above let $\bar{D}$ denote the closure of $D.$ Also let
$T:\mathcal{\ X}\rightarrow\mathcal{X}_{a}$ be the bounded linear operator
defined by
\[
T(X)=(\bar{B}^{\ast}\bar{B})^{-1}\mathcal{Q}X
\]
for all $X\in\mathcal{X}.$ Then for all $F\in\mathcal{D}(\bar{D}),$
\begin{equation}
F=\mathbb{E}F+D^{\ast}T\bar{D}F. \label{e.7.48}%
\end{equation}

\end{theorem}

It is worth pointing out that $\bar{B}^{\ast}$ is not $//B^{\ast}$ but is
instead given by $\mathcal{Q}//B^{\ast}.$ This is because $//B^{\ast}$ does
not take adapted processes to adapted processes. This is the reason it is
necessary to introduce the orthogonal projection, $\mathcal{Q}.$

\begin{proof}
Let $Y\in\mathcal{X}_{a}$ be given and $X\in\mathcal{X}_{a}$ be chosen so that
$F=\mathbb{E}F+D^{\ast}X.$ Then
\begin{align*}
\langle Y,\mathcal{\ Q}\bar{D}F\rangle_{\mathcal{X}}  &  =\langle Y,\bar
{D}F\rangle_{\mathcal{X}}=\mathbb{E}\left[  D^{\ast}Y\cdot F\right] \\
&  =\mathbb{E}\left[  D^{\ast}Y\cdot D^{\ast}X\right]  =\mathbb{E}\left[
\langle\bar{B}(Y),\bar{B}(X)\rangle_{L^{2}(ds)}\right] \\
&  =\langle Y,\bar{B}^{\ast}\bar{B}(X)\rangle_{\mathcal{X}},
\end{align*}
where in going from the first to the second line we have used $\mathbb{E}%
\left[  D^{\ast}Y\right]  =0.$ From the above displayed equation it follows
that $\mathcal{\ Q}\bar{D}F=\bar{B}^{\ast}\bar{B}(X)$ and hence $X=(\bar
{B}^{\ast}\bar{B})^{-1}\mathcal{Q}\bar{D}F=T(\bar{D}F).$
\end{proof}

\subsubsection{Proof of Theorem \ref{t.7.41}}

Let $F\in\mathcal{D}(\bar{D}).$ By Theorem \ref{t.7.45},%
\[
\mathbb{E}\left[  F-\mathbb{E}F\right]  ^{2}=\mathbb{E}\left[  D^{\ast}%
T\bar{D}F\right]  ^{2}=\mathbb{E}|\bar{B}(T\bar{D}F)|_{L^{2}(ds,T_{o}M)}%
^{2}\leq C\langle\bar{D}F,\bar{D}F\rangle_{\mathcal{X}}%
\]
where $C$ is the operator norm of $\bar{B}T.$ In particular if $F\in
\mathcal{D}(\mathcal{L}),$ then $\langle\bar{D}F,\bar{D}F\rangle_{\mathcal{X}%
}=\mathbb{E}[\mathcal{L}F\cdot F],$ and hence
\[
\langle\mathcal{L}F,F\rangle_{L^{2}\left(  \mu_{W\left(  M\right)  }\right)
}\geq C^{-1}\langle F-\mathbb{E}F,F-\mathbb{E}F\rangle_{L^{2}\left(
\mu_{W\left(  M\right)  }\right)  }.
\]
Therefore, if $F\in\operatorname*{Nul}(\mathcal{L}),$ it follows that
$F=\mathbb{E}F,$ i.e. $F$ is a constant. Moreover if $F\perp1$ (i.e.
$\mathbb{E}F=0)$ then
\[
\langle\mathcal{L}F,F\rangle_{L^{2}\left(  \mu_{W\left(  M\right)  }\right)
}\geq C^{-1}\langle F,F\rangle_{L^{2}\left(  \mu_{W\left(  M\right)  }\right)
},
\]
proving Theorem \ref{t.7.41}\ with $c=C^{-1}.$

\subsection{$W\left(  M\right)  $ -- Martingale Representation Theorem}

\label{s.7.7}

In this subsection, $\Sigma$ is a Brownian motion on $M$ starting at $o\in M,$
$//_{s}$ is stochastic parallel translation along $\Sigma$ and%
\[
b_{s}=\left[  \Psi(\Sigma)\right]  _{s}=\int_{0}^{s}//_{r}^{-1}\delta
\Sigma_{r}%
\]
is the undeveloped $T_{o}M$ -- valued Brownian motion associated to $\Sigma$
as described before Theorem \ref{t.5.29}.

\begin{lemma}
\label{l.7.46}If $f\in C^{\infty}\left(  M^{n+1}\right)  $ and $i\leq n,$ then%
\begin{align}
&  \mathbb{E}\left[  \left.  //_{s_{i}}^{-1}\text{$\operatorname*{grad}$}%
_{i}f\left(  \Sigma_{s_{1}},\ldots,\Sigma_{s_{n}},\Sigma_{s_{n+1}}\right)
\right\vert \mathcal{F}_{s_{n}}\right] \nonumber\\
&  =//_{s_{i}}^{-1}\text{$\operatorname*{grad}$}_{i}(e^{(s_{n+1}-s_{n}%
)\bar{\Delta}_{n+1}/2}f)\left(  \Sigma_{s_{1}},\ldots,\Sigma_{s_{n}}%
,\Sigma_{s_{n}}\right)  . \label{e.7.49}%
\end{align}

\end{lemma}

\begin{proof}
Let us begin with the special case where $f=g\otimes h$ for some $g\in
C^{\infty}\left(  M^{n}\right)  $ and $h\in C^{\infty}(M)$ where $g\otimes
h\left(  x_{1},\dots,x_{n+1}\right)  :=g\left(  x_{1},\dots,x_{n}\right)
h\left(  x_{n+1}\right)  .$ In this case%
\[
//_{s_{i}}^{-1}\text{$\operatorname*{grad}$}_{i}f\left(  \Sigma_{s_{1}}%
,\ldots,\Sigma_{s_{n}},\Sigma_{s_{n+1}}\right)  =//_{s_{i}}^{-1}%
\text{$\operatorname*{grad}$}_{i}g\left(  \Sigma_{s_{1}},\ldots,\Sigma_{s_{n}%
}\right)  \cdot h\left(  \Sigma_{s_{n+1}}\right)
\]
where $//_{s_{i}}^{-1}\operatorname*{grad}_{i}g\left(  \Sigma_{s_{1}}%
,\ldots,\Sigma_{s_{n}}\right)  $ is $\mathcal{F}_{s_{n}}$ -- measurable. Hence
by the Markov property we have%
\begin{align*}
&  \mathbb{E}\left[  \left.  //_{s_{i}}^{-1}\text{$\operatorname*{grad}$}%
_{i}f\left(  \Sigma_{s_{1}},\ldots,\Sigma_{s_{n}},\Sigma_{s_{n+1}}\right)
\right\vert \mathcal{F}_{s_{n}}\right] \\
&  =//_{s_{i}}^{-1}\text{$\operatorname*{grad}$}_{i}g\left(  \Sigma_{s_{1}%
},\ldots,\Sigma_{s_{n}}\right)  \mathbb{E}\left[  \left.  h\left(
\Sigma_{s_{n+1}}\right)  \right\vert \mathcal{F}_{s_{n}}\right] \\
&  =//_{s_{i}}^{-1}\text{$\operatorname*{grad}$}_{i}g\left(  \Sigma_{s_{1}%
},\ldots,\Sigma_{s_{n}}\right)  (e^{(s_{n+1}-s_{n})\bar{\Delta}/2}h)\left(
\Sigma_{s_{n}}\right) \\
&  =//_{s_{i}}^{-1}\text{$\operatorname*{grad}$}_{i}(e^{(s_{n+1}-s_{n}%
)\bar{\Delta}_{n+1}/2}f)\left(  \Sigma_{s_{1}},\ldots,\Sigma_{s_{n}}%
,\Sigma_{s_{n}}\right)  .
\end{align*}
\textbf{Alternatively}, as we have already seen, $M_{s}:=(e^{(s_{n+1}%
-s)\bar{\Delta}/2}h)\left(  \Sigma_{s}\right)  $ is a martingale for $s\leq
s_{n+1}$, and therefore,%
\[
\mathbb{E}\left[  \left.  h\left(  \Sigma_{s_{n+1}}\right)  \right\vert
\mathcal{F}_{s_{n}}\right]  =\mathbb{E}\left[  \left.  M_{s_{n+1}}\right\vert
\mathcal{F}_{s_{n}}\right]  =M_{s_{n}}=(e^{(s_{n+1}-s_{n})\bar{\Delta}%
/2}h)\left(  \Sigma_{s_{n}}\right)  .
\]
Since Eq. (\ref{e.7.49}) is linear in $f,$ this proves Eq. (\ref{e.7.49}) when
$f$ is a linear combination of functions of the form $g\otimes h$ as above.

Using a partition unity argument along with the standard convolution
approximation methods; to any $f\in C^{\infty}\left(  M^{n+1}\right)  $ there
exists a sequence of $f_{k}\in C^{\infty}\left(  M^{n+1}\right)  $ with each
$f_{k}$ being a linear combination of functions of the form $g\otimes h$ such
that $f_{k}$ along with all of its derivatives converges uniformly to $f.$
Passing to the limit in Eq. (\ref{e.7.49}) with $f$ being replaced by $f_{k},$
shows that Eq. (\ref{e.7.49}) holds for all $f\in C^{\infty}\left(
M^{n+1}\right)  .$
\end{proof}

Recall that $Q_{s}$ is the $\operatorname*{End}\left(  T_{o}M\right)  $ --
valued process determined in Eq. (\ref{e.6.1}) and since
\[
\frac{d}{ds}Q_{s}^{-1}=-Q_{s}^{-1}\left[  \frac{d}{ds}Q_{s}\right]  Q_{s}%
^{-1},
\]
$Q_{s}^{-1}$ solves the equation,%
\begin{equation}
\frac{d}{ds}Q_{s}^{-1}=\frac{1}{2}\operatorname{Ric}_{//_{s}}Q_{s}^{-1}\text{
with }Q_{0}^{-1}=I. \label{e.7.50}%
\end{equation}

\begin{theorem}
[Representation Formula]\label{t.7.47}Suppose that $F$ is a smooth cylinder
function of the form $F\left(  \sigma\right)  =f\left(  \sigma_{s_{1}}%
,\ldots,\sigma_{s_{n}}\right)  ,$ then
\begin{equation}
F(\Sigma)=\mathbb{E}F+\int_{0}^{1}\left\langle a_{s},db_{s}\right\rangle
\label{e.7.51}%
\end{equation}
where $a_{s}$ is a bounded predictable process, $a_{s}$ is zero if $s\geq
s_{n}$ and $s\rightarrow a_{s}$ is continuous off the partition set,
$\{s_{1},\ldots,s_{n}\}.$ Moreover $a_{s}$ may be expressed as
\begin{equation}
a_{s}:=Q_{s}^{-1}\mathbb{E}\left[  \left.  \sum_{i=1}^{n}1_{s\leq s_{i}%
}Q_{s_{i}}//_{s_{i}}^{-1}\text{$\operatorname*{grad}$}_{i}f\left(
\Sigma_{s_{1}},\ldots,\Sigma_{s_{n}}\right)  \right\vert \mathcal{F}%
_{s}\right]  . \label{e.7.52}%
\end{equation}

\end{theorem}

\begin{proof}
The proof will be by induction on $n.$ For $n=1$ suppose $F\left(
\Sigma\right)  =f\left(  \Sigma_{t}\right)  $ for some $t\in(0,1].$
Integrating Eq. (\ref{e.5.38}) from $[0,t]$ with $g=f$ implies%
\begin{equation}
F\left(  \Sigma\right)  =f\left(  \Sigma_{t}\right)  =e^{t\bar{\Delta}%
/2}f\left(  o\right)  +\int_{0}^{t}\langle//_{s}^{-1}\operatorname*{grad}%
\mathrm{\ }e^{(t-s)\bar{\Delta}/2}f\left(  \Sigma_{s}\right)  ,db_{s}\rangle.
\label{e.7.53}%
\end{equation}
Since $e^{t\bar{\Delta}/2}f\left(  o\right)  =\mathbb{E}F,$ Eq. (\ref{e.7.53})
shows Eq. (\ref{e.7.51}) holds with
\[
a_{s}=1_{0\leq s\leq t}//_{s}^{-1}\operatorname*{grad}\mathrm{\ }%
e^{(t-s)\bar{\Delta}/2}f\left(  \Sigma_{s}\right)  .
\]
By Lemma \ref{l.6.1}, $Q_{s}//_{s}^{-1}\operatorname*{grad}\mathrm{\ }%
e^{(t-s)\bar{\Delta}/2}f\left(  \Sigma_{s}\right)  $ is a martingale, and
hence%
\[
Q_{s}//_{s}^{-1}\operatorname*{grad}\mathrm{\ }e^{(t-s)\bar{\Delta}/2}f\left(
\Sigma_{s}\right)  =\mathbb{E}\left[  \left.  Q_{t}//_{t}^{-1}%
\operatorname*{grad}\mathrm{\ }f\left(  \Sigma_{t}\right)  \right\vert
\mathcal{F}_{s}\right]
\]
from which it follows that
\[
a_{s}=1_{0\leq s\leq t}//_{s}^{-1}\operatorname*{grad}\mathrm{\ }%
e^{(t-s)\bar{\Delta}/2}f\left(  \Sigma_{s}\right)  =1_{0\leq s\leq t}%
Q_{s}^{-1}\mathbb{E}\left[  \left.  Q_{t}//_{t}^{-1}\operatorname*{grad}%
\mathrm{\ }f\left(  \Sigma_{t}\right)  \right\vert \mathcal{F}_{s}\right]  .
\]
This shows that Eq. (\ref{e.7.52}) is valid for $n=1.$

To carry out the inductive step, suppose the result holds for level $n$ and
now suppose that%
\[
F\left(  \Sigma\right)  =f\left(  \Sigma_{s_{1}},\ldots,\Sigma_{s_{n+1}%
}\right)
\]
with $0<s_{1}<s_{2}\cdots<s_{n+1}\leq1.$ Let
\[
(\Delta_{n+1}f)(x_{1},x_{2},\ldots,x_{n+1})=(\Delta g)(x_{n+1})
\]
where $g(x):=f(x_{1},x_{2},\ldots,x_{n},x).$ Similarly, let
$\operatorname*{grad}_{n+1}$ denote the gradient acting on the $\left(
n+1\right)  ^{\text{th}}$ -- variable of a function $f\in C^{\infty}%
(M^{n+1}).$ Set
\[
H(s,\Sigma):=(e^{(s_{n+1}-s)\bar{\Delta}_{n+1}/2}f)(\Sigma_{s_{1}}%
,\ldots,\Sigma_{s_{n}},\Sigma_{s})
\]
for $s_{n}\leq s\leq s_{n+1}.$ By It\^{o}'s Lemma, (see Corollary \ref{c.5.18}
and also Eq. (\ref{e.5.38}),%
\[
d\left[  H(s,\Sigma_{s})\right]  =\langle\text{$\operatorname*{grad}$}%
_{n+1}e^{(s_{n+1}-s)\bar{\Delta}_{n+1}/2}f)(\Sigma_{s_{1}},\ldots
,\Sigma_{s_{n}},\Sigma_{s},//_{s}db_{s}\rangle
\]
for $s_{n}\leq s\leq s_{n+1}.$ Integrating this last expression from $s_{n}$
to $s_{n+1}$ yields:
\begin{align}
F(\Sigma)  &  =(e^{(s_{n+1}-s_{n})\bar{\Delta}_{n+1}/2}f)(\Sigma_{s_{1}%
},\ldots,\Sigma_{s_{n}},\Sigma_{s_{n}})\nonumber\\
&  \qquad+\int_{s_{n}}^{s_{n+1}}\langle//_{s}^{-1}\text{$\operatorname*{grad}%
$}_{n+1}e^{(s_{n+1}-s)\bar{\Delta}_{n+1}/2}f)\left(  \Sigma_{s_{1}}%
,\ldots,\Sigma_{s_{n}},\Sigma_{s}\right)  ,db_{s}\rangle\label{e.7.54}\\
&  =(e^{(s_{n+1}-s_{n})\bar{\Delta}_{n+1}/2}f)(\Sigma_{s_{1}},\ldots
,\Sigma_{s_{n}},\Sigma_{s_{n}})+\int_{s_{n}}^{s_{n+1}}\langle\alpha_{s}%
,db_{s}\rangle, \label{e.7.55}%
\end{align}
where $\alpha_{s}:=//_{s}^{-1}(\operatorname*{grad}_{n+1}e^{(s_{n+1}%
-s)\bar{\Delta}_{n+1}/2}f)(\Sigma_{s_{1}},\ldots,\Sigma_{s_{n}},\Sigma_{s}).$
By the induction hypothesis, the smooth cylinder function,
\[
(e^{(s_{n+1}-s_{n})\bar{\Delta}_{n+1}/2}f)(\Sigma_{s_{1}},\ldots,\Sigma
_{s_{n}},\Sigma_{s_{n}}),
\]
may be written as a constant plus $\int_{0}^{1}\langle\tilde{a}_{s}%
,db_{s}\rangle,$ where $\tilde{a}_{s}$ is bounded and piecewise continuous and
$\tilde{a}_{s}\equiv0$ if $s\geq s_{n}.$ Thus if we let $a_{s}:=\tilde{a}%
_{s}+1_{s_{n}<s\leq s_{n+1}}\alpha_{s},$ we have shown%
\[
F(\Sigma)=C+\int_{0}^{s_{n+1}}\langle a_{s},db_{s}\rangle
\]
for some constant $C.$ Taking expectations of both sides of this equation then
shows $C=\mathbb{E}\left[  F(\Sigma)\right]  $ and the proof of Eq.
(\ref{e.7.51}) is complete. So to finish the proof it only remains to verify
Eq. (\ref{e.7.52}).

Again by Lemma \ref{l.6.1},
\[
s\rightarrow M_{s}:=Q_{s}//_{s}^{-1}(\operatorname{grad}_{n+1}e^{(s_{n+1}%
-s)\bar{\Delta}_{n+1}/2}f)(\Sigma_{s_{1}},\ldots,\Sigma_{s_{n}},\Sigma_{s})
\]
is a martingale for $s\in\lbrack s_{n},s_{n+1}]$ and therefore,%
\begin{align}
M_{s}  &  =Q_{s}//_{s}^{-1}(\operatorname{grad}_{n+1}e^{(s_{n+1}-s)\bar
{\Delta}_{n+1}/2}f)(\Sigma_{s_{1}},\ldots,\Sigma_{s_{n}},\Sigma_{s}%
)\nonumber\\
&  =\mathbb{E}\left[  \left.  M_{s_{n+1}}\right\vert \mathcal{F}_{s}\right]
=\mathbb{E}\left[  \left.  Q_{s_{n+1}}//_{s_{n+1}}^{-1}\left(
\text{$\operatorname*{grad}$}_{n+1}f\right)  (\Sigma_{s_{1}},\ldots
,\Sigma_{s_{n}},\Sigma_{s_{n+1}})\right\vert \mathcal{F}_{s}\right]  ,
\label{e.7.56}%
\end{align}
i.e.%
\begin{align}
//_{s}^{-1}  &  (\operatorname{grad}_{n+1}e^{(s_{n+1}-s)\bar{\Delta}_{n+1}%
/2}f)(\Sigma_{s_{1}},\ldots,\Sigma_{s_{n}},\Sigma_{s})\nonumber\\
&  =Q_{s}^{-1}\mathbb{E}\left[  \left.  Q_{s_{n+1}}//_{s_{n+1}}^{-1}\left(
\text{$\operatorname*{grad}$}_{n+1}f\right)  (\Sigma_{s_{1}},\ldots
,\Sigma_{s_{n}},\Sigma_{s_{n+1}})\right\vert \mathcal{F}_{s}\right]  .
\label{e.7.57}%
\end{align}
Using this identity, Eq. (\ref{e.7.54}) may be written as%
\begin{align}
F(\Sigma)=  &  g(\Sigma_{s_{1}},\ldots,\Sigma_{s_{n}})\nonumber\\
&  +\int_{s_{n}}^{s_{n+1}}\left\langle Q_{s}^{-1}\mathbb{E}\left[  \left.
Q_{s_{n+1}}//_{s_{n+1}}^{-1}\left(  \text{$\operatorname*{grad}$}%
_{n+1}f\right)  (\Sigma_{s_{1}},\ldots,\Sigma_{s_{n}},\Sigma_{s_{n+1}%
})\right\vert \mathcal{F}_{s}\right]  ,db_{s}\right\rangle . \label{e.7.58}%
\end{align}

where%
\[
g\left(  x_{1},\dots,x_{n}\right)  :=(e^{(s_{n+1}-s_{n})\bar{\Delta}_{n+1}%
/2}f)\left(  x_{1},\dots,x_{n},x_{n}\right)  .
\]
By the induction hypothesis,%
\begin{align}
g  &  (\Sigma_{s_{1}},\ldots,\Sigma_{s_{n}})\nonumber\\
&  =C+\int_{0}^{1}\left\langle Q_{s}^{-1}\mathbb{E}\left[  \left.  \sum
_{i=1}^{n}1_{s\leq s_{i}}Q_{s_{i}}//_{s_{i}}^{-1}\text{$\operatorname*{grad}$%
}_{i}g\left(  \Sigma_{s_{1}},\ldots,\Sigma_{s_{n}}\right)  \right\vert
\mathcal{F}_{s}\right]  ,db_{s}\right\rangle \label{e.7.59}%
\end{align}
where $C=\mathbb{E}\left[  F(\Sigma)\right]  $ as we have already seen or
alternatively, by the Markov property,%
\begin{align}
C  &  :=\mathbb{E}(e^{(s_{n+1}-s_{n})\bar{\Delta}_{n+1}/2}f)(\Sigma_{s_{1}%
},\ldots,\Sigma_{s_{n}},\Sigma_{s_{n}})\nonumber\\
&  =\mathbb{E}f(\Sigma_{s_{1}},\ldots,\Sigma_{s_{n}},\Sigma_{s_{n+1}%
})=\mathbb{E}\left[  F(\Sigma)\right]  . \label{e.7.60}%
\end{align}
By Lemma \ref{l.7.46}, for $s\leq s_{n}$ and $i<n$%
\begin{align}
\mathbb{E}  &  \left[  \left.  Q_{s_{i}}//_{s_{i}}^{-1}%
\text{$\operatorname*{grad}$}_{i}g\left(  \Sigma_{s_{1}},\ldots,\Sigma_{s_{n}%
}\right)  \right\vert \mathcal{F}_{s}\right] \nonumber\\
&  =\mathbb{E}\left[  \left.  Q_{s_{i}}\mathbb{E}\left[  \left.  //_{s_{i}%
}^{-1}\text{$\operatorname*{grad}$}_{i}(e^{(s_{n+1}-s_{n})\bar{\Delta}%
_{n+1}/2}f)\left(  \Sigma_{s_{1}},\ldots,\Sigma_{s_{n}},\Sigma_{s_{n}}\right)
\right\vert \mathcal{F}_{s_{n}}\right]  \right\vert \mathcal{F}_{s}\right]
\nonumber\\
&  =\mathbb{E}\left[  \left.  Q_{s_{i}}//_{s_{i}}^{-1}%
\text{$\operatorname*{grad}$}_{i}f\left(  \Sigma_{s_{1}},\ldots,\Sigma_{s_{n}%
},\Sigma_{s_{n+1}}\right)  \right\vert \mathcal{F}_{s}\right]  .
\label{e.7.61}%
\end{align}
While for $s\leq s_{n}$ and $i=n,$ we have:%
\begin{align*}
\text{$\operatorname*{grad}$}_{n}g\left(  \Sigma_{s_{1}},\ldots,\Sigma_{s_{n}%
}\right)   &  =\text{$\operatorname*{grad}$}_{n}(e^{(s_{n+1}-s_{n})\bar
{\Delta}_{n+1}/2}f)\left(  \Sigma_{s_{1}},\ldots,\Sigma_{s_{n}},\Sigma_{s_{n}%
}\right) \\
&  +\text{$\operatorname*{grad}$}_{n+1}(e^{(s_{n+1}-s_{n})\bar{\Delta}%
_{n+1}/2}f)\left(  \Sigma_{s_{1}},\ldots,\Sigma_{s_{n}},\Sigma_{s_{n}}\right)
,
\end{align*}%
\begin{align*}
\mathbb{E}  &  \left[  \left.  Q_{s_{n}}//_{s_{n}}^{-1}%
\text{$\operatorname*{grad}$}_{n}(e^{(s_{n+1}-s_{n})\bar{\Delta}_{n+1}%
/2}f)\left(  \Sigma_{s_{1}},\ldots,\Sigma_{s_{n}},\Sigma_{s_{n}}\right)
\right\vert \mathcal{F}_{s}\right] \\
&  =\mathbb{E}\left[  \left.  Q_{s_{n}}\mathbb{E}\left[  \left.  //_{s_{n}%
}^{-1}\text{$\operatorname*{grad}$}_{n}(e^{(s_{n+1}-s_{n})\bar{\Delta}%
_{n+1}/2}f)\left(  \Sigma_{s_{1}},\ldots,\Sigma_{s_{n}},\Sigma_{s_{n}}\right)
\right\vert \mathcal{F}_{s_{n}}\right]  \right\vert \mathcal{F}_{s}\right] \\
&  =\mathbb{E}\left[  \left.  Q_{s_{n}}//_{s_{n}}^{-1}%
\text{$\operatorname*{grad}$}_{n}f\left(  \Sigma_{s_{1}},\ldots,\Sigma_{s_{n}%
},\Sigma_{s_{n+1}}\right)  \right\vert \mathcal{F}_{s}\right]
\end{align*}
by Lemma \ref{l.7.46} and%
\begin{align*}
\mathbb{E}  &  \left[  \left.  \mathbb{E}\left[  \left.  Q_{s_{n}}//_{s_{n}%
}^{-1}\text{$\operatorname*{grad}$}_{n+1}(e^{(s_{n+1}-s_{n})\bar{\Delta}%
_{n+1}/2}f)\left(  \Sigma_{s_{1}},\ldots,\Sigma_{s_{n}},\Sigma_{s_{n}}\right)
\right\vert \mathcal{F}_{s_{n}}\right]  \right\vert \mathcal{F}_{s}\right] \\
&  =\mathbb{E}\left[  \left.  Q_{s_{n+1}}//_{s_{n+1}}^{-1}\left(
\text{$\operatorname*{grad}$}_{n+1}f\right)  (\Sigma_{s_{1}},\ldots
,\Sigma_{s_{n}},\Sigma_{s_{n+1}})\right\vert \mathcal{F}_{s}\right]
\end{align*}
from Eq. (\ref{e.7.57}) with $s=s_{n}.$ Combining the previous three displayed
equations shows,%
\begin{align}
\mathbb{E}  &  \left[  \left.  Q_{s_{n}}//_{s_{n}}^{-1}%
\text{$\operatorname*{grad}$}_{n}g\left(  \Sigma_{s_{1}},\ldots,\Sigma_{s_{n}%
}\right)  \right\vert \mathcal{F}_{s}\right] \nonumber\\
&  =\mathbb{E}\left[  \left.  Q_{s_{n}}//_{s_{n}}^{-1}%
\text{$\operatorname*{grad}$}_{n}f\left(  \Sigma_{s_{1}},\ldots,\Sigma_{s_{n}%
},\Sigma_{s_{n+1}}\right)  \right\vert \mathcal{F}_{s}\right] \nonumber\\
&  \qquad+\mathbb{E}\left[  \left.  Q_{s_{n+1}}//_{s_{n+1}}^{-1}\left(
\text{$\operatorname*{grad}$}_{n+1}f\right)  (\Sigma_{s_{1}},\ldots
,\Sigma_{s_{n}},\Sigma_{s_{n+1}})\right\vert \mathcal{F}_{s}\right]
\label{e.7.62}%
\end{align}
Assembling Eqs. (\ref{e.7.59}), (\ref{e.7.60}), (\ref{e.7.61}) and
(\ref{e.7.62}) implies
\begin{align*}
g  &  (\Sigma_{s_{1}},\ldots,\Sigma_{s_{n}})\\
&  =\mathbb{E}\left[  F\left(  \Sigma\right)  \right]  +\int_{0}^{1}\sum
_{i=1}^{n}\left\langle Q_{s}^{-1}\mathbb{E}\left[  \left.  1_{s\leq s_{i}%
}Q_{s_{i}}//_{s_{i}}^{-1}\text{$\operatorname*{grad}$}_{i}f\left(
\Sigma_{s_{1}},\ldots,\Sigma_{s_{n}},\Sigma_{s_{n+1}}\right)  \right\vert
\mathcal{F}_{s}\right]  ,db_{s}\right\rangle \\
&  \qquad+\int_{0}^{1}\left\langle Q_{s}^{-1}\mathbb{E}\left[  \left.
1_{s\leq s_{n}}Q_{s_{n+1}}//_{s_{n+1}}^{-1}\left(  \text{$\operatorname*{grad}%
$}_{n+1}f\right)  (\Sigma_{s_{1}},\ldots,\Sigma_{s_{n}},\Sigma_{s_{n+1}%
})\right\vert \mathcal{F}_{s}\right]  ,db_{s}\right\rangle
\end{align*}
which combined with Eq. (\ref{e.7.58}) shows
\begin{multline*}
F\left(  \Sigma\right)  =\mathbb{E}\left[  F\left(  \Sigma\right)  \right] \\
+\int_{0}^{1}\left\langle Q_{s}^{-1}\mathbb{E}\left[  \left.  \sum_{i=1}%
^{n+1}1_{s\leq s_{i}}Q_{s_{i}}//_{s_{i}}^{-1}\text{$\operatorname*{grad}$}%
_{i}f\left(  \Sigma_{s_{1}},\ldots,\Sigma_{s_{n}},\Sigma_{s_{n+1}}\right)
\right\vert \mathcal{F}_{s}\right]  ,db_{s}\right\rangle .
\end{multline*}
This completes the induction argument and hence the proof.
\end{proof}

\begin{proposition}
\label{p.7.48}Equation (\ref{e.7.51}) may also be written as%
\begin{equation}
F\left(  \Sigma\right)  =\mathbb{E}\left[  F\left(  \Sigma\right)  \right]
+\int_{0}^{1}\left\langle \mathbb{E}\left[  \left.  \xi_{s}-\frac{1}{2}%
\int_{s}^{1}Q_{s}^{-1}Q_{r}\operatorname{Ric}_{//_{r}}\xi_{r}dr\right\vert
\mathcal{F}_{s}\right]  ,db_{s}\right\rangle . \label{e.7.63}%
\end{equation}
where
\[
\xi_{s}:=//_{s}^{-1}\frac{d}{ds}\left(  DF\right)  _{s}.
\]

\end{proposition}

\begin{proof}
Let $v_{i}:=//_{s_{i}}^{-1}\operatorname*{grad}_{i}f\left(  \Sigma_{s_{1}%
},\ldots,\Sigma_{s_{n}}\right)  ,$ so that
\[
\xi_{s}:=//_{s}^{-1}\frac{d}{ds}\left(  DF\right)  _{s}=\sum_{i=1}%
^{n}1_{s<s_{i}}v_{i},
\]
and let
\[
\alpha_{s}:=\sum_{i=1}^{n}1_{s\leq s_{i}}Q_{s}^{-1}Q_{s_{i}}//_{s_{i}}%
^{-1}\text{$\operatorname*{grad}$}_{i}f\left(  \Sigma_{s_{1}},\ldots
,\Sigma_{s_{n}}\right)  =\sum_{i=1}^{n}1_{s\leq s_{i}}Q_{s}^{-1}Q_{s_{i}}%
v_{i}.
\]
Then the Lebesgue-Stieljtes measure associate to $\xi_{s}$ is%
\[
d\xi_{s}=-\sum_{i=1}^{n}\delta_{s_{i}}\left(  ds\right)  v_{i}%
\]
and therefore%
\[
\alpha_{s}=-Q_{s}^{-1}\int_{s}^{1}Q_{r}d\xi_{r}=-\int_{s}^{1}Q_{s}^{-1}%
Q_{r}d\xi_{r}.
\]
So by integration by parts we have, for $s\notin\left\{  0,s_{1},\dots
,s_{n},1\right\}  ,$%
\begin{align*}
\alpha_{s}  &  =-\int_{s}^{1}Q_{s}^{-1}Q_{r}d\xi_{r}=-\left[  Q_{s}^{-1}%
Q_{r}\xi_{r}\right]  |_{r=s}^{r=1}+\int_{s}^{1}Q_{s}^{-1}\left[  \frac{d}%
{dr}Q_{r}\right]  \xi_{r}\\
&  =\xi_{s}-\frac{1}{2}\int_{s}^{1}Q_{s}^{-1}Q_{r}\operatorname*{Ric}%
\nolimits_{//_{r}}\xi_{r}%
\end{align*}
where we have used $\xi_{1}=0.$ This completes the proof since from Eqs.
(\ref{e.7.51}) and (\ref{e.7.52}),
\[
F\left(  \Sigma\right)  =\mathbb{E}\left[  F\left(  \Sigma\right)  \right]
+\int_{0}^{1}\left\langle E\left[  \left.  \alpha_{s}\right\vert
\mathcal{F}_{s}\right]  ,db_{s}\right\rangle .
\]

\end{proof}

\begin{corollary}
\label{c.7.49} Let $F$ be a smooth cylinder function, then there is a
predictable, piecewise continuously differentiable Cameron-Martin vector field
$X$ such that $F=\mathbb{E}\left[  F\right]  +D^{\ast}X.$
\end{corollary}

\begin{proof}
Just follow the proof of Lemma \ref{l.7.42} using Theorem \ref{t.7.47} in
place of Corollary \ref{c.7.20}.
\end{proof}

\subsubsection{The equivalence of integration by parts and the representation
formula\label{s.7.7.1}}

\begin{corollary}
\label{c.7.50}The representation formula in Theorem \ref{t.7.47} may be used
to prove the integration by parts Theorem \ref{t.7.32} in the case $F$ is a
cylinder function.
\end{corollary}

\begin{proof}
Let $F$ be a cylinder function, $a_{s}$ be as in Eq. (\ref{e.7.52}), $h$ be an
adapted Cameron-Martin process and $k_{s}:=\left(  Q_{s}^{\operatorname*{tr}%
}\right)  ^{-1}h_{s}.$ Then, by the product rule and Eq. (\ref{e.7.39}),%
\[
h_{s}^{\prime}+\frac{1}{2}\operatorname{Ric}_{//_{s}}h_{s}=\left(  \frac
{d}{ds}+\frac{1}{2}\operatorname{Ric}_{//_{s}}\right)  Q_{s}%
^{\operatorname*{tr}}k_{s}=Q_{s}^{\operatorname*{tr}}k_{s}^{\prime}.
\]
Hence,%
\begin{align*}
\mathbb{E}  &  \left[  F\int_{0}^{1}\langle h_{s}^{\prime}+\frac{1}%
{2}\operatorname{Ric}_{//_{s}}h_{s},db_{s}\rangle\right] \\
&  =\mathbb{E}\left[  \left(  \mathbb{E}F+\int_{0}^{1}\left\langle
a_{s},db_{s}\right\rangle \right)  \int_{0}^{1}\langle Q_{s}%
^{\operatorname*{tr}}k_{s}^{\prime},db_{s}\rangle\right] \\
&  =\mathbb{E}\left[  \int_{0}^{1}\langle Q_{s}^{\operatorname*{tr}}%
k_{s}^{\prime},a_{s}\rangle ds\right] \\
&  =\mathbb{E}\left[  \int_{0}^{1}\langle Q_{s}^{\operatorname*{tr}}%
k_{s}^{\prime},\sum_{i=1}^{n}1_{s\leq s_{i}}Q_{s}^{-1}Q_{s_{i}}//_{s_{i}}%
^{-1}\text{$\operatorname*{grad}$}_{i}f\left(  \Sigma_{s_{1}},\ldots
,\Sigma_{s_{n}}\right)  \rangle ds\right] \\
&  =\mathbb{E}\left[  \int_{0}^{1}\langle k_{s}^{\prime},\sum_{i=1}%
^{n}1_{s\leq s_{i}}Q_{s_{i}}//_{s_{i}}^{-1}\text{$\operatorname*{grad}$}%
_{i}f\left(  \Sigma_{s_{1}},\ldots,\Sigma_{s_{n}}\right)  \rangle ds\right] \\
&  =\mathbb{E}\left[  \sum_{i=1}^{n}\langle k_{s_{i}},Q_{s_{i}}//_{s_{i}}%
^{-1}\text{$\operatorname*{grad}$}_{i}f\left(  \Sigma_{s_{1}},\ldots
,\Sigma_{s_{n}}\right)  \rangle\right] \\
&  =\mathbb{E}\left[  \sum_{i=1}^{n}\langle//_{s_{i}}h_{s_{i}}%
,\text{$\operatorname*{grad}$}_{i}f\left(  \Sigma_{s_{1}},\ldots,\Sigma
_{s_{n}}\right)  \rangle\right]  =\mathbb{E}\left[  X^{h}F\right]  .
\end{align*}

\end{proof}

Conversely we may give a proof of Theorem \ref{t.7.47} which is based on the
integration by parts Theorem \ref{t.7.32}.

\begin{theorem}
[Representation Formula]\label{t.7.51}Suppose $F$ is a cylinder function on
$W\left(  M\right)  $ as in Eq. (\ref{e.7.2}) and $\xi_{s}:=//_{s}^{-1}%
\frac{d}{ds}\left(  DF\right)  _{s},$ then%
\begin{equation}
F=\mathbb{E}F+\int_{0}^{1}\left\langle \mathbb{E}\left[  \left.  \xi_{s}%
-\frac{1}{2}\int_{s}^{1}Q_{s}^{-1}Q_{r}\operatorname{Ric}_{//_{r}}\xi
_{r}dr\right\vert \mathcal{F}_{s}\right]  ,db_{s}\right\rangle .
\label{e.7.64}%
\end{equation}
where $Q_{s}$ is the solution to Eq. (\ref{e.6.1}).
\end{theorem}

\begin{proof}
Let $h\in\mathcal{X}_{a}$ be a predictable adapted Cameron-Martin valued
process such that $\mathbb{E}\int_{0}^{1}\left\vert h_{s}^{\prime}\right\vert
^{2}ds<\infty.$ By the martingale representation property in Corollary
\ref{c.7.20},
\begin{equation}
F=\mathbb{E}F+\int_{0}^{1}\langle a,db\rangle\label{e.7.65}%
\end{equation}
for some predictable process $a$ such that $\mathbb{E}\int_{0}^{1}\left\vert
a_{s}\right\vert ^{2}ds<\infty.$ Then from Corollary \ref{c.7.35} and the
It\^{o} isometry property,%
\begin{align}
\mathbb{E}\left[  X^{Q^{\mathrm{tr}}h}F\right]   &  =\mathbb{E}\left[
F\cdot\left(  X^{Q^{\mathrm{tr}}h}\right)  ^{\ast}1\right]  =\mathbb{E}\left[
F\cdot\int_{0}^{1}\langle Q^{\mathrm{tr}}h^{\prime},db\rangle\right]
\nonumber\\
&  =\mathbb{E}\left[  \int_{0}^{1}\langle Q_{s}^{\mathrm{tr}}h_{s}^{\prime
},a_{s}\rangle ds\right]  =\mathbb{E}\left[  \int_{0}^{1}\langle h_{s}%
^{\prime},Q_{s}a_{s}\rangle ds\right]  . \label{e.7.66}%
\end{align}
On the other hand we may compute $\mathbb{E}\left[  X^{Q^{\mathrm{tr}}%
h}F\right]  $ as:
\begin{align}
\mathbb{E}\left[  X^{Q^{\mathrm{tr}}h}F\right]   &  =\mathbb{E}\left[  \langle
DF,//Q^{\mathrm{tr}}h\rangle_{H}\right]  =\mathbb{E}\int_{0}^{1}\langle\xi
_{s},\frac{d}{ds}\left(  Q^{\mathrm{tr}}h\right)  _{s}\rangle ds\nonumber\\
&  =\mathbb{E}\int_{0}^{1}\left\langle \xi_{s},Q_{s}^{\mathrm{tr}}%
h_{s}^{\prime}-\frac{1}{2}\operatorname{Ric}_{//_{s}}Q_{s}^{\mathrm{tr}}%
h_{s}\right\rangle ds \label{e.7.67}%
\end{align}
where we have used Eq. (\ref{e.7.39}) in the last equality. We will now
rewrite the right side of Eq. (\ref{e.7.67}) so that it has the same form as
Eq. (\ref{e.7.66}) To do this let $\rho_{s}:=\frac{1}{2}\operatorname*{Ric}%
_{//_{s}}$ and notice that%
\begin{align*}
\int_{0}^{1}\langle\xi_{s},\rho_{s}Q_{s}^{\mathrm{tr}}h_{s}\rangle ds  &
=\int_{0}^{1}\left\langle Q_{s}\rho_{s}^{\ast}\xi_{s},\left(  \int_{0}%
^{s}h_{r}^{\prime}dr\right)  \right\rangle ds\\
&  =\int drds1_{0\leq r\leq s\leq1}\langle Q_{s}\rho_{s}^{\ast}\xi_{s}%
,h_{r}^{\prime}\rangle=\int_{0}^{1}\left\langle \int_{s}^{1}Q_{r}\rho
_{r}^{\ast}\xi_{r}dr,h_{s}^{\prime}\right\rangle ds
\end{align*}
wherein the last equality we have interchanged the role of $r$ and $s.$ Using
this result back in Eq. (\ref{e.7.67}) implies%
\begin{equation}
\mathbb{E}\left[  X^{Q^{\mathrm{tr}}h}F\right]  =\mathbb{E}\int_{0}%
^{1}\left\langle Q_{s}\xi_{s}-\int_{s}^{1}Q_{r}\rho_{r}^{\ast}\xi_{r}%
dr,h_{s}^{\prime}\right\rangle ds. \label{e.7.68}%
\end{equation}
and comparing this with Eq. (\ref{e.7.66}) shows%
\begin{equation}
\mathbb{E}\int_{0}^{1}\left\langle Q_{s}a_{s}-Q_{s}\xi_{s}+\int_{s}^{1}%
Q_{r}\rho_{r}^{\ast}\xi_{r}dr,h_{s}^{\prime}\right\rangle ds=0 \label{e.7.69}%
\end{equation}
for all $h\in\mathcal{X}_{a}.$

Up to now we have only used $F\in\mathcal{D}\left(  D\right)  $ and not the
fact that $F$ is a cylinder function. We will use this hypothesis now. From
the easy part of Theorem \ref{t.7.47} we know that $a_{s}$ satisfies the
additional properties of being 1) bounded, 2) zero if $s\geq s_{n}$ and most
importantly 3) $s\rightarrow a_{s}$ is continuous off the partition set,
$\{s_{1},\ldots,s_{n}\}.$

Fix $\tau\in(0,1)\setminus\{s_{1},\ldots,s_{n}\},$ $v\in T_{o}M$ and let $G$
be a bounded $\mathcal{F}_{\tau}$ -- measurable function. For $n\in\mathbb{N}$
let
\[
l_{n}\left(  s\right)  :=\int_{0}^{s}n1_{\tau\leq r\leq\tau+\frac{1}{n}}dr.
\]
Replacing $h$ in Eq. (\ref{e.7.69}) by $h_{n}\left(  s\right)  :=G\cdot
l_{n}\left(  s\right)  v$ and then passing to the limit as $n\rightarrow
\infty,$ implies
\begin{align*}
0  &  =\lim_{n\rightarrow\infty}\mathbb{E}\int_{0}^{1}\left\langle Q_{s}%
a_{s}-Q_{s}\xi_{s}+\int_{s}^{1}Q_{r}\rho_{r}^{\ast}\xi_{r}dr,h_{n}^{\prime
}\left(  s\right)  \right\rangle ds\\
&  =\mathbb{E}\left[  G\left\langle Q_{\tau}a_{\tau}-Q_{\tau}\xi_{\tau}%
+\int_{\tau}^{1}Q_{r}\rho_{r}^{\ast}\xi_{r}dr,v\right\rangle \right]
\end{align*}
and since $G$ and $v$ were arbitrary we conclude from this equation that%
\[
\mathbb{E}\left[  \left.  Q_{\tau}\xi_{\tau}-\int_{\tau}^{1}Q_{r}\rho_{r}%
\xi_{r}dr\right\vert \mathcal{F}_{\tau}\right]  =Q_{\tau}a_{\tau}.
\]
Thus for all but finitely many $s\in\lbrack0,1],$%
\begin{align*}
a_{s}  &  =Q_{s}^{-1}\mathbb{E}\left[  \left.  Q_{s}\xi_{s}-\int_{s}^{1}%
Q_{r}\rho_{r}\xi_{r}dr\right\vert \mathcal{F}_{s}\right] \\
&  =\mathbb{E}\left[  \left.  \xi_{s}-\frac{1}{2}\int_{s}^{1}Q_{s}^{-1}%
Q_{r}\operatorname{Ric}_{//_{r}}\xi_{r}dr\right\vert \mathcal{F}_{s}\right]  .
\end{align*}
Combining this with Eq. (\ref{e.7.65}) proves Eq. (\ref{e.7.64}).
\end{proof}

\subsection{Logarithmic-Sobolev Inequality for $W\left(  M\right)
$\label{s.7.8}}

The next theorem is the \textquotedblleft curved\textquotedblright%
\ generalization of Theorem \ref{t.7.24}.

\begin{theorem}
[Hsu's Logarithmic Sobolev Inequality]\label{t.7.52}Let $M$ be a compact
Riemannian manifold, then for all $F\in\mathcal{D}\left(  \bar{D}\right)  $
\begin{align}
\mathbb{E}\left[  F^{2}\log F^{2}\right]  \leq &  \mathbb{E}F^{2}\cdot
\log\mathbb{E}F^{2}\nonumber\\
&  +2\mathbb{E}\int_{0}^{1}\left\vert //_{s}^{-1}\left(  DF\right)
_{s}^{\prime}-\frac{1}{2}\int_{s}^{1}Q_{s}^{-1}Q_{r}\operatorname{Ric}%
_{//_{r}}//_{r}^{-1}\left(  DF\right)  _{r}^{\prime}dr\right\vert ^{2}ds,
\label{e.7.70}%
\end{align}
where $\left(  DF\right)  _{s}^{\prime}:=\frac{d}{ds}\left(  DF\right)  _{s}.$
Moreover, there is a constant $C=C\left(  \operatorname*{Ric}\right)  $ such
that%
\begin{equation}
\mathbb{E}\left[  F^{2}\log F^{2}\right]  \leq C\mathbb{E}\left[  \langle
DF,DF\rangle_{H\left(  T_{o}M\right)  }\right]  +\mathbb{E}F^{2}\cdot
\log\mathbb{E}F^{2}. \label{e.7.71}%
\end{equation}

\end{theorem}

\begin{proof}
The proof we give here follows the paper of Capitaine, Hsu and Ledoux
\cite{CHL97}. We begin in the same way as the proof of Theorem \ref{t.7.24}.
Let $F\in\mathcal{F}C^{1}\left(  W\left(  M\right)  \right)  ,$ $\varepsilon
>0,$ $H_{\varepsilon}:=F^{2}+\varepsilon\in\mathcal{D}\left(  \bar{D}\right)
$ and%
\[
a_{s}:=\mathbb{E}\left[  \left.  \xi_{s}-\frac{1}{2}\int_{s}^{1}Q_{s}%
^{-1}Q_{r}\operatorname{Ric}_{//_{r}}\xi_{r}dr\right\vert \mathcal{F}%
_{s}\right]
\]
where
\[
\xi_{s}=//_{s}^{-1}\frac{d}{ds}\left(  DH_{\varepsilon}\right)  _{s}%
=2F\cdot//_{s}^{-1}\frac{d}{ds}\left(  DF\right)  _{s}.
\]
Then by Theorem \ref{t.7.47},%
\[
H_{\varepsilon}=\mathbb{E}H_{\varepsilon}+\int_{0}^{1}\langle a,db\rangle.
\]
The same proof used to derive Eq. (\ref{e.7.23}) shows%
\begin{align*}
\mathbb{E}\left[  \phi\left(  H_{\varepsilon}\right)  \right]   &
=\mathbb{E}\left[  \phi\left(  M_{1}\right)  \right]  =\phi\left(
\mathbb{E}M_{1}\right)  +\frac{1}{2}\mathbb{E}\left[  \int_{0}^{1}\frac
{1}{M_{s}}\left\vert a_{s}\right\vert ^{2}ds\right] \\
&  =\phi\left(  \mathbb{E}H_{\varepsilon}\right)  +\frac{1}{2}\mathbb{E}%
\left[  \int_{0}^{1}\frac{1}{\mathbb{E}\left[  H_{\varepsilon}|\mathcal{F}%
_{s}\right]  }\left\vert a_{s}\right\vert ^{2}ds\right]  .
\end{align*}
By the Cauchy-Schwarz inequality and the contractive properties of conditional
expectations,%
\begin{align*}
\left\vert a_{s}\right\vert ^{2}  &  =\left\vert \mathbb{E}\left[  2F\left.
\left\{  //_{s}^{-1}\left(  DF\right)  _{s}^{\prime}-\frac{1}{2}\int_{s}%
^{1}Q_{s}^{-1}Q_{r}\operatorname{Ric}_{//_{r}}//_{r}^{-1}\left(  DF\right)
_{r}^{\prime}dr\right\}  \right\vert \mathcal{F}_{s}\right]  \right\vert
^{2}\\
&  \leq4\mathbb{E}\left[  F^{2}|\mathcal{F}_{s}\right]  \cdot\mathbb{E}\left[
\left.  \left\vert //_{s}^{-1}\left(  DF\right)  _{s}^{\prime}-\frac{1}{2}%
\int_{s}^{1}Q_{s}^{-1}Q_{r}\operatorname{Ric}_{//_{r}}//_{r}^{-1}\left(
DF\right)  _{r}^{\prime}dr\right\vert ^{2}\right\vert \mathcal{F}_{s}\right]
\end{align*}
Combining the last two equations along with Eq. (\ref{e.7.24}) implies%
\begin{align*}
\mathbb{E}\phi\left(  H_{\varepsilon}\right)  \leq &  \phi\left(
\mathbb{E}H_{\varepsilon}\right) \\
&  +2\mathbb{E}\int_{0}^{1}\mathbb{E}\left[  \left.  \left\vert //_{s}%
^{-1}\left(  DF\right)  _{s}^{\prime}-\frac{1}{2}\int_{s}^{1}Q_{s}^{-1}%
Q_{r}\operatorname{Ric}_{//_{r}}//_{r}^{-1}\left(  DF\right)  _{r}^{\prime
}dr\right\vert ^{2}\right\vert \mathcal{F}_{s}\right]  ds\\
=  &  \phi\left(  \mathbb{E}H_{\varepsilon}\right) \\
&  +2\mathbb{E}\int_{0}^{1}\left\vert //_{s}^{-1}\left(  DF\right)
_{s}^{\prime}-\frac{1}{2}\int_{s}^{1}Q_{s}^{-1}Q_{r}\operatorname{Ric}%
_{//_{r}}//_{r}^{-1}\left(  DF\right)  _{r}^{\prime}dr\right\vert ^{2}ds.
\end{align*}
We may now let $\varepsilon\downarrow0$ in this inequality to learn Eq.
(\ref{e.7.70}) holds for all $F\in\mathcal{F}C^{1}\left(  W\right)  .$ By
compactness of $M,$ $\operatorname{Ric}_{m}$ is bounded on $M$ and so by
simple Gronwall type estimates on $Q$ and $Q^{-1},$ there is a non-random
constant $K<\infty$ such that%
\[
\left\Vert Q_{s}^{-1}Q_{r}\operatorname{Ric}_{//_{r}}\right\Vert _{op}\leq
K\text{ for all }r,s.
\]
Therefore,%
\begin{align*}
&  \left\vert //_{s}^{-1}\left(  DF\right)  _{s}^{\prime}-\frac{1}{2}\int
_{s}^{1}Q_{s}^{-1}Q_{r}\operatorname{Ric}_{//_{r}}//_{r}^{-1}\left(
DF\right)  _{r}^{\prime}dr\right\vert ^{2}\\
&  \qquad\leq\left[  \left\vert \left(  DF\right)  _{s}^{\prime}\right\vert
+\frac{1}{2}K\int_{0}^{1}\left\vert \left(  DF\right)  _{s}^{\prime
}\right\vert ds\right]  ^{2}\\
&  \qquad\leq2\left\vert \left(  DF\right)  _{s}^{\prime}\right\vert
^{2}+\frac{1}{2}K^{2}\left[  \int_{0}^{1}\left\vert \left(  DF\right)
_{s}^{\prime}\right\vert ds\right]  ^{2}\\
&  \qquad\leq2\left\vert \left(  DF\right)  _{s}^{\prime}\right\vert
^{2}+\frac{1}{2}K^{2}\int_{0}^{1}\left\vert \left(  DF\right)  _{s}^{\prime
}\right\vert ^{2}ds
\end{align*}
and hence%
\begin{align*}
2\mathbb{E}  &  \int_{0}^{1}\left\vert \left(  DF\right)  _{s}^{\prime}%
-\frac{1}{2}\int_{s}^{1}Q_{s}^{-1}Q_{r}\operatorname{Ric}_{//_{r}}\left(
DF\right)  _{r}^{\prime}dr\right\vert ^{2}ds\\
&  \leq\left(  4+K^{2}\right)  \int_{0}^{1}\left\vert \left(  DF\right)
_{s}^{\prime}\right\vert ^{2}ds.
\end{align*}
Combining this estimate with Eq. (\ref{e.7.70}) implies Eq. (\ref{e.7.71})
holds with $C=\left(  4+K^{2}\right)  .$ Again, since $\mathcal{F}C^{1}\left(
W\right)  $ is a core for $\bar{D},$ standard limiting arguments show that Eq.
(\ref{e.7.70}) and Eq. (\ref{e.7.71}) are valid for all $F\in\mathcal{D}%
\left(  \bar{D}\right)  .$
\end{proof}

Theorem \ref{t.7.52} was first proved by Hsu \cite{Hsu3} with an independent
proof given shortly thereafter by Aida and Elworthy \cite{AE}. Hsu's original
proof relied on a Markov dependence version of a standard additivity property
for logarithmic Sobolev inequalities and makes key use of Corollary
\ref{c.7.37}. On the other hand Aida and Elworthy show, using the projection
construction of Brownian motion, the logarithmic Sobolev inequality on $W(M)$
is a consequence of Gross' \cite{Gr4} original logarithmic Sobolev inequality
on the classical Wiener space $W(\mathbb{R}^{N}),$ see Theorem \ref{t.7.24}.
In Aida's and Elworthy's proof, Theorem \ref{t.5.43} plays an important role.

\subsection{More References}

Many people have now proved some version of integration by parts for path and
loop spaces in one context or another, see for example
\cite{Bismut84a,CM2,Cross96,CM1,CM2,CM3,D5,D6,D9,ES1,ES2,FM2,Fr,Leandre93b,MM4,No3,S3,S4,S5,Shig97a,HUZ02}%
. We have followed Bismut in these notes who proved integration by parts
formulas for cylinder functions depending on one time. However, as is pointed
out by Leandre and Malliavin and Fang, Bismut's technique works with out any
essential change for arbitrary cylinder functions. In \cite{D5,D6}, the flow
associated to a general class of vector fields on paths and loop spaces of a
manifold were constructed. The reader is also referred to the texts
\cite{El-LJ-Li,Hsu2003i,Stroock2001} and the related articles
\cite{Fang95a,Fang95b,Cr-Fang96,Fang98,Fang99a,Fang99b,Fang99c,Cruz-M99,Cruz-M00,CFM00,Cruz-M01,CF01,Cruz-M02,xdl00}%
.

Many of the results in this section extend to pinned Wiener measure on loop
spaces, see \cite{D6} for example. Loop spaces are more interesting than path
spaces since they have nontrivial topology, The issue of the spectral gap and
logarithmic Sobolev inequalities for general loop spaces is still an open
problem. In \cite{Gr6}, Gross has prove a logarithmic Sobolev inequality on
Loop groups with an added \textquotedblleft potential term\textquotedblright%
\ for a special geometry on loop groups. Here Gross uses pinned Wiener measure
as the reference measure. In Driver and Lohrenz \cite{DL}, it is shown that a
logarithmic Sobolev inequality \textbf{without} a potential term does hold on
the Loop group provided one replace pinned Wiener measure by a
\textquotedblleft heat kernel\textquotedblright\ measure. The
quasi-invarariance properties of the heat kernel measure on loop groups was
first established in \cite{Driver97b,Driver98b}. For more results on heat
kernel measures on the loop groups see for example,
\cite{Driver01a,Driver00a,Carson97,Carson99,Fang99a,Fang99b,Inahama01}.

The question as to when or if the potential is needed in Gross's setting for
logarithmic Sobolev inequalities is still an open question, but see Gong,
R\"{o}ckner and Wu \cite{GRW01} for a positive result in this direction.
Eberle \cite{Eberle01,Eberle02,Eberle03a,Eberle03b} has provided examples of
Riemannian manifolds where the spectral gap inequality fails in the loop space
setting. The reader is referred to \cite{Driver03,Driver04} and the references
therein for some more perspective on the stochastic analysis on loop spaces.

\newpage

\section{Malliavin's Methods for Hypoelliptic Operators\label{s.8}}

In this section we will be concerned with determining smoothness properties of
the $\mathrm{Law}\left(  \Sigma_{t}\right)  $ where $\Sigma_{t}$ denotes the
solution to Eq. (\ref{e.5.1}) with $\Sigma_{0}=o$ and $\beta=B$ is an
$\mathbb{R}^{n}$ -- valued Brownian motion. Unlike the previous sections in
these notes, the map $\mathbf{X}\left(  m\right)  :\mathbb{R}^{n}\rightarrow
T_{m}M$ is \textbf{not }assumed to be surjective. Equivalently put, the
diffusion generator $L:=\frac{1}{2}\sum_{i=1}^{n}X_{i}^{2}+X_{0}$ is no longer
assumed to be elliptic. However we will always be assuming that the vector
fields $\left\{  X_{i}\right\}  _{i=0}^{n}$ satisfy H\"{o}rmander's restricted
bracket condition at $o\in M$ as in Definition \ref{d.8.1}. Let $\mathcal{K}%
_{1}:=\left\{  X_{1},\dots,X_{n}\right\}  $ and $\mathcal{K}_{l}$ be defined
inductively by
\[
\mathcal{K}_{l+1}=\left\{  [X_{i},K]:K\in\mathcal{K}_{l}\right\}
\cup\mathcal{K}_{l}.
\]
For example%
\begin{align*}
\mathcal{K}_{2}= &  \left\{  X_{1},\dots,X_{n}\right\}  \cup\left\{  \lbrack
X_{j},X_{i}]:i,j=1,\dots,n\right\}  \text{ and}\\
\mathcal{K}_{3}= &  \left\{  X_{1},\dots,X_{n}\right\}  \cup\left\{  \lbrack
X_{j},X_{i}]:i,j=1,\dots,n\right\}  \\
&  \qquad\cup\left\{  \left[  X_{k},[X_{j},X_{i}]\right]  :i,j,k=1,\dots
,n\right\}  \text{ etc.}%
\end{align*}

\begin{definition}
\label{d.8.1}The collection of vector fields, $\left\{  X_{i}\right\}
_{i=0}^{n}\subset\Gamma\left(  TM\right)  ,$ satisfies \textbf{H\"{o}rmander's
restricted bracket condition at }$m\in M$ if there exist $l\in\mathbb{N}$ such
that%
\[
\mathrm{span}(\left\{  K(m):K\in\mathcal{K}_{l}\right\}  )=T_{m}M.
\]

\end{definition}

Under this condition it follows from a classical theorem of H\"{o}rmander that
solutions to the heat equation $\partial_{t}u=Lu$ are necessarily smooth.
Since the fundamental solution to this equation at $o\in M$ is the law of the
process $\Sigma_{t},$ it follows that the $\mathrm{Law}\left(  \Sigma
_{t}\right)  $ is absolutely continuous relative to the volume measure
$\lambda$ on $M$ and its Radon-Nikodym derivative is a smooth function on $M.$
Malliavin, in his 1976 pioneering paper \cite{Malliavin78}, gave a
probabilistic proof of this fact. Malliavin's original paper was followed by
an avalanche of papers carrying out and extending Malliavin's program
including the fundamental works of Stroock \cite{St1a,St1b,St2}, Kusuoka and
Stroock \cite{KS1,KS2,KS3}, and Bismut \cite{Bismut84a}. See also
\cite{Bell5,Bell5a,BH,IW2,Malliavin97,Rao95,Norris86b,Nu1,S1,S2,Wa1} and the
references therein. The purpose of this section is to briefly explain
(omitting some details) Malliavin methods.

\subsection{Malliavin's Ideas in Finite Dimensions\label{s.8.1}}

To understand Malliavin's methods it is best to begin with a finite
dimensional analogue.

\begin{theorem}
[Malliavin's Ideas in Finite Dimensions]\label{t.8.2}Let $W\mathbb{=R}^{N},$
$\mu$ be the Gaussian measure on $W$ defined by%
\[
d\mu\left(  x\right)  :=\left(  2\pi\right)  ^{-N/2}e^{-\frac{1}{2}\left\vert
x\right\vert ^{2}}dm\left(  x\right)  .
\]
Further suppose $F:W\rightarrow\mathbb{R}^{d}$ (think $F=\Sigma_{t})$ is a
function satisfying:
\end{theorem}

\begin{enumerate}
\item $F$ is smooth and all of its partial derivatives are in
\[
L^{\infty-}\left(  \mu\right)  :=\cap_{1\leq p<\infty}L^{p}(W,\mu).
\]

\item $F$ is a submersion or equivalently assume the \textquotedblleft
Malliavin\textquotedblright\ matrix
\[
C(\omega):=DF(\omega)DF(\omega)^{\ast}%
\]
is invertible for all $\omega\in W.$

\item Let
\[
\Delta(\omega):=\det C(\omega)=\det(DF(\omega)DF(\omega)^{\ast})
\]
and assume $\Delta^{-1}\in L^{\infty-}\left(  \mu\right)  .$
\end{enumerate}

Then the law $(\mu_{F}=F_{\ast}\mu=\mu\circ F^{-1})$ of $F$ is absolutely
continuous relative to Lebesgue measure, $\lambda,$ on $\mathbb{R}^{d}$ and
the Radon-Nikodym derivative, $\rho:=d\mu_{F}/d\lambda,$ is smooth.

\begin{proof}
For each vector field $Y\in\Gamma\left(  T\mathbb{R}^{d}\right)  ,$ define%
\begin{equation}
\mathbb{Y}(\omega)=DF(\omega)^{\ast}C(\omega)^{-1}Y\left(  F\left(
\omega\right)  \right)  \label{e.8.1}%
\end{equation}
--- a smooth vector field on $W$ such that $DF(\omega)\mathbb{Y}%
(\omega)=Y\left(  F\left(  \omega\right)  \right)  $ or in more geometric
notation,%
\begin{equation}
F_{\ast}\mathbb{Y}(\omega)=Y\left(  F\left(  \omega\right)  \right)  .
\label{e.8.2}%
\end{equation}
For the purposes of this proof, it is sufficient to restrict our attention to
the case where $Y$ is a constant vector field.

Explicit computations using the chain rule and Cramer's rule for computing
$C(\omega)^{-1}$ shows that $D^{k}\mathbb{Y}$ may be expressed as a polynomial
in $\Delta^{-1}$ and $D^{\ell}F$ for $\ell=0,1,2\ldots,k.$ In particular
$D^{k}\mathbb{Y}$ is in $L^{\infty-}\left(  \mu\right)  .$ Suppose
$f,g:W\rightarrow\mathbb{R}$ are $C^{1}$ functions such that $f,g,$ and their
first order derivatives are in $L^{\infty-}\left(  \mu\right)  .$ Then by a
standard truncation argument and integration by parts, one shows that
\[
\int_{W}(\mathbb{Y}f)g\,d\mu=\int_{W}f(\mathbb{Y}^{\ast}g)\,d\mu,
\]
where
\[
\mathbb{Y}^{\ast}=-\mathbb{Y}+\delta(\mathbb{Y})\text{ and }\delta
(\mathbb{Y})(\omega):=-\operatorname*{div}(\mathbb{Y})(\omega)+\mathbb{Y}%
(\omega)\cdot\omega.
\]

Suppose that $\phi\in C_{c}^{\infty}(\mathbb{R}^{d})$ and $Y_{i}\in
\mathbb{R}^{d}\subset\Gamma\left(  \mathbb{R}^{d}\right)  ,$ then from Eq.
(\ref{e.8.2}) and induction,
\[
(Y_{1}Y_{2}\cdots Y_{k}\phi)(F(\omega))\,=(\mathbb{Y}_{1}\mathbb{Y}_{2}%
\cdots\mathbb{Y}_{k}(\phi\circ F))(\omega)
\]
and therefore,
\begin{align}
\int_{\mathbb{R}^{d}}(Y_{1}Y_{2}\cdots Y_{k}\phi)d\mu_{F}  &  =\int_{W}%
(Y_{1}Y_{2}\cdots Y_{k}\phi)(F(\omega))\,d\mu(\omega)\nonumber\\
&  =\int_{W}(\mathbb{Y}_{1}\mathbb{Y}_{2}\cdots\mathbb{Y}_{k}(\phi\circ
F))(\omega)\,d\mu(\omega)\nonumber\\
&  =\int_{W}\phi(F(\omega))\cdot(\mathbb{Y}_{k}^{\ast}\mathbb{Y}_{k-1}^{\ast
}\cdots\mathbb{Y}_{1}^{\ast}1)(\omega)\,d\mu(\omega). \label{e.8.3}%
\end{align}
By the remarks in the previous paragraph, $(\mathbb{Y}_{k}^{\ast}%
\mathbb{Y}_{k-1}^{\ast}\cdots\mathbb{Y}_{1}^{\ast}1)\in L^{\infty-}\left(
\mu\right)  $ which along with Eq. (\ref{e.8.3}) shows%
\[
\left\vert \int_{\mathbb{R}^{d}}(Y_{1}Y_{2}\cdots Y_{k}\phi)d\mu
_{F}\right\vert \leq C\left\Vert \phi\right\Vert _{L^{\infty}\left(
\mathbb{R}^{d}\right)  },
\]
where $C=\left\Vert \mathbb{Y}_{k}^{\ast}\mathbb{Y}_{k-1}^{\ast}%
\cdots\mathbb{Y}_{1}^{\ast}1\right\Vert _{L^{1}\left(  \mu\right)  }<\infty.$
It now follows from Sobolev imbedding theorems or simple Fourier analysis that
$\mu_{F}\ll\lambda$ and that $\rho:=d\mu_{F}/d\lambda$ is a smooth function.
\end{proof}

The remainder of Section \ref{s.8} will be devoted to an infinite dimensional
analogue of Theorem \ref{t.8.2} (see Theorem \ref{t.8.9}) where $\mathbb{R}%
^{d}$ is replaced by a manifold $M^{d},$%
\[
W:=\left\{  \omega\in C\left(  [0,\infty),\mathbb{R}^{n}\right)
:\omega\left(  0\right)  =0\right\}  ,
\]
$\mu$ is taken to be Wiener measure on $W,$ $B_{t}:W\rightarrow\mathbb{R}^{n}$
be defined by $B_{t}\left(  \omega\right)  =\omega_{t}$ and $F:=\Sigma
_{t}:W\left(  \mathbb{R}^{n}\right)  \rightarrow M$ is a solution to Eq.
(\ref{e.5.1}) with $\Sigma_{0}=o\in M$ and $\beta=B.$ Recall that $\mu$ is the
unique measure on $\mathcal{F}:=\sigma\left(  B_{t}:t\in\lbrack0,\infty
)\right)  $ such that $\left\{  B_{t}\right\}  _{t\geq0}$ is a Brownian
motion. I am now using $t$ as the dominant parameter rather than $s$ to be in
better agreement with the literature on this subject.

\subsection{Smoothness of Densities for H\"{o}rmander Type Diffusions
\label{s.8.2}}

For simplicity of the exposition, it will be assumed that $M$ is a compact
Riemannian manifold. However this can and should be relaxed. For example most
everything we are going to say would work if $M$ is an imbedded submanifold in
$\mathbb{R}^{N}$ and the vector fields $\left\{  X_{i}\right\}  _{i=0}^{n}$
are the restrictions of smooth vector fields on $\mathbb{R}^{N}$ whose partial
derivatives to any order greater than $0$ are all bounded.

\begin{remark}
\label{r.8.3}The choice of Riemannian metric here is somewhat arbitrary and is
an artifact of the method to be described below. It is the author's belief
that this issue has still not been adequately addressed in the literature.
\end{remark}

To abbreviate the notation, let
\[
H=\left\{  h\in W:\langle h,h\rangle_{H}:=\int_{0}^{\infty}\left\vert \dot
{h}\left(  t\right)  \right\vert ^{2}dt<\infty\right\}
\]
and $D\Sigma_{t}:H\rightarrow T_{\Sigma_{t}}M$ be defined by $\left(
D\Sigma_{t}\right)  h:=\partial_{h}T_{t}^{B}\left(  o\right)  $ as defined
Theorem \ref{t.7.26}. Recall from Theorem \ref{t.7.26} that%
\begin{equation}
\left(  D\Sigma_{t}\right)  h:=Z_{t}\int_{0}^{t}Z_{\tau}^{-1}\mathbf{X}\left(
\Sigma_{\tau}\right)  \dot{h}_{\tau}d\tau=//_{t}z_{t}\int_{0}^{t}z_{\tau}%
^{-1}//_{\tau}^{-1}\mathbf{X}\left(  \Sigma_{\tau}\right)  \dot{h}_{\tau}%
d\tau, \label{e.8.4}%
\end{equation}
where $\dot{h}_{\tau}:=\frac{d}{d\tau}h_{\tau},$ $Z_{t}:=\left(  T_{t}%
^{B}\right)  _{\ast o}:T_{o}M\rightarrow T_{\Sigma_{t}}M,$ $//_{t}$ is
stochastic parallel translation along $\Sigma$ and $z_{t}:=//_{t}^{-1}Z_{t}.$
In the sequel, adjoints will be denote by either \textquotedblleft\ $^{\ast}$
\textquotedblright\ or \textquotedblleft\ $^{\mathrm{tr}}$ \textquotedblright%
\ with the former being used if an infinite dimensional space is involved and
the latter if all spaces involved are finite dimensional.

\begin{definition}
[Reduced Malliavin Covariance]\label{d.8.4}The $\operatorname*{End}\left(
T_{o}M\right)  $ -- valued random variable,%
\begin{align}
\bar{C}_{t}  &  :=\int_{0}^{t}Z_{\tau}^{-1}\mathbf{X}\left(  \Sigma_{\tau
}\right)  \mathbf{X}\left(  \Sigma_{\tau}\right)  ^{\mathrm{tr}}\left(
Z_{\tau}^{-1}\right)  ^{\mathrm{tr}}d\tau\label{e.8.5}\\
&  =\int_{0}^{t}z_{\tau}^{-1}//_{\tau}^{-1}\mathbf{X}\left(  \Sigma_{\tau
}\right)  \mathbf{X}\left(  \Sigma_{\tau}\right)  ^{\mathrm{tr}}//_{\tau
}\left(  z_{\tau}^{-1}\right)  ^{\mathrm{tr}}d\tau, \label{e.8.6}%
\end{align}
will be called the \textbf{reduced Malliavin covariance} \textbf{matrix.}
\end{definition}

\begin{theorem}
\label{t.8.5}The adjoint, $\left(  D\Sigma_{t}\right)  ^{\ast}:T_{\Sigma_{t}%
}M\rightarrow H,$ of the map $D\Sigma_{t}$ is determined by%
\begin{equation}
\frac{d}{d\tau}\left[  \left(  D\Sigma_{t}\right)  ^{\ast}//_{t}v\right]
_{\tau}=1_{\tau\leq t}\mathbf{X}\left(  \Sigma_{\tau}\right)  ^{\mathrm{tr}%
}//_{\tau}\left(  z_{t}z_{\tau}^{-1}\right)  ^{\mathrm{tr}}v \label{e.8.7}%
\end{equation}
for all $v\in T_{o}M.$ The Malliavin covariance matrix $C_{t}:=D\Sigma
_{t}\left(  D\Sigma_{t}\right)  ^{\ast}:T_{\Sigma_{t}}M\rightarrow
T_{\Sigma_{t}}M$ is given by $C_{t}=Z_{t}\bar{C}_{t}Z_{t}^{\mathrm{tr}}$ or
equivalently%
\begin{equation}
C_{t}=D\Sigma_{t}\left(  D\Sigma_{t}\right)  ^{\ast}=//_{t}z_{t}\bar{C}%
_{t}z_{t}^{\mathrm{tr}}//_{t}^{-1}. \label{e.8.8}%
\end{equation}

\end{theorem}

\begin{proof}
Using Eq. (\ref{e.8.4}),%
\begin{align}
\left\langle D\Sigma_{t}h,//_{t}v\right\rangle _{T_{\Sigma_{t}}M}  &
=\left\langle Z_{t}\int_{0}^{t}Z_{\tau}^{-1}\mathbf{X}\left(  \Sigma_{\tau
}\right)  \dot{h}_{\tau}d\tau,//_{t}v\right\rangle _{T_{\Sigma_{t}}%
M}\nonumber\\
&  =\left\langle //_{t}z_{t}\int_{0}^{t}z_{\tau}^{-1}//_{\tau}^{-1}%
\mathbf{X}\left(  \Sigma_{\tau}\right)  \dot{h}_{\tau}d\tau,//_{t}%
v\right\rangle _{T_{\Sigma_{t}}M}\nonumber\\
&  =\int_{0}^{t}\left\langle z_{t}z_{\tau}^{-1}//_{\tau}^{-1}\mathbf{X}\left(
\Sigma_{\tau}\right)  \dot{h}_{\tau},v\right\rangle _{T_{o}M}d\tau\nonumber\\
&  =\int_{0}^{t}\left\langle \dot{h}_{\tau},\mathbf{X}\left(  \Sigma_{\tau
}\right)  ^{\mathrm{tr}}//_{\tau}\left(  z_{t}z_{\tau}^{-1}\right)
^{\mathrm{tr}}v\right\rangle _{\mathbb{R}^{n}}d\tau\label{e.8.9}%
\end{align}
which implies Eq. (\ref{e.8.7}). Combining Eqs. (\ref{e.8.4}) and
(\ref{e.8.7}), using
\[
Z_{\tau}^{\mathrm{tr}}=\left(  //_{\tau}z_{\tau}\right)  ^{\mathrm{tr}%
}=z_{\tau}^{\mathrm{tr}}//_{\tau}^{\mathrm{tr}}=z_{\tau}^{\mathrm{tr}}%
//_{\tau}^{-1},
\]
shows%
\begin{align*}
D\Sigma_{t}\left(  D\Sigma_{t}\right)  ^{\ast}//_{t}v  &  =Z_{t}\int_{0}%
^{t}Z_{\tau}^{-1}\mathbf{X}\left(  \Sigma_{\tau}\right)  \mathbf{X}\left(
\Sigma_{\tau}\right)  ^{\mathrm{tr}}//_{\tau}\left(  z_{t}z_{\tau}%
^{-1}\right)  ^{\mathrm{tr}}vd\tau\\
&  =Z_{t}\int_{0}^{t}Z_{\tau}^{-1}\mathbf{X}\left(  \Sigma_{\tau}\right)
\mathbf{X}\left(  \Sigma_{\tau}\right)  ^{\mathrm{tr}}\left(  Z_{\tau}%
^{-1}\right)  ^{\mathrm{tr}}Z_{t}^{\mathrm{tr}}//_{t}vd\tau.
\end{align*}
Therefore,
\[
C_{t}=Z_{t}\bar{C}_{t}Z_{t}^{\mathrm{tr}}=//_{t}z_{t}\bar{C}_{t}%
z_{t}^{\mathrm{tr}}//_{t}^{-1}%
\]
from which Eq. (\ref{e.8.8}) follows.
\end{proof}

The next crucial theorem is at the heart of Malliavin's method and constitutes
the deepest part of the theory. The proof of this theorem will be postponed
until Section \ref{s.8.4} below.

\begin{theorem}
[Non-degeneracy of $\bar{C}_{t}$]\label{t.8.6}Let $\bar{\Delta}_{t}%
:=\det\left(  \bar{C}_{t}\right)  .$ If H\"{o}rmander's restricted bracket
condition at $o\in M$ holds then $\bar{\Delta}_{t}>0$ a.e. (i.e. $\bar{C}_{t}$
is invertible a.e.) and moreover $\bar{\Delta}_{t}^{-1}\in L^{\infty-}\left(
\mu\right)  .$
\end{theorem}

Following the general strategy outlined in Theorem \ref{t.8.2}, given a vector
field $Y\in\Gamma\left(  TM\right)  $ we wish to lift it via the map
$\Sigma_{t}:W\rightarrow M$ to a vector field $\mathbb{Y}^{t}$ on $W:=W\left(
\mathbb{R}^{n}\right)  .$ According to the prescription used in Eq.
(\ref{e.8.1}) in Theorem \ref{t.8.2},
%Perhaps we should be lifting $//_{t}v$ rather than an element $Y\inG(TM)?$%
\begin{equation}
\mathbb{Y}^{t}:=\left(  D\Sigma_{t}\right)  ^{\ast}\left(  D\Sigma_{t}\left(
D\Sigma_{t}\right)  ^{\ast}\right)  ^{-1}Y\left(  \Sigma_{t}\right)  =\left(
D\Sigma_{t}\right)  ^{\ast}C_{t}^{-1}Y\left(  \Sigma_{t}\right)  \in H.
\label{e.8.10}%
\end{equation}
From Eq. (\ref{e.8.8})
\[
C_{t}^{-1}=//_{t}\left(  z_{t}^{\mathrm{tr}}\right)  ^{-1}\bar{C}_{t}%
^{-1}z_{t}^{-1}//_{t}^{-1}%
\]
and combining this with Eq. (\ref{e.8.10}), using Eq. (\ref{e.8.7}), implies%
\begin{align*}
\frac{d}{d\tau}\mathbb{Y}_{\tau}^{t}  &  =1_{\tau\leq t}\frac{d}{d\tau}\left[
\left(  D\Sigma_{t}\right)  //_{t}\left(  z_{t}^{\mathrm{tr}}\right)
^{-1}\bar{C}_{t}^{-1}z_{t}^{-1}//_{t}^{-1}Y\left(  \Sigma_{t}\right)  \right]
_{\tau}\\
&  =1_{\tau\leq t}\mathbf{X}\left(  \Sigma_{\tau}\right)  ^{\mathrm{tr}%
}//_{\tau}\left(  z_{t}z_{\tau}^{-1}\right)  ^{\mathrm{tr}}\left(
z_{t}^{\mathrm{tr}}\right)  ^{-1}\bar{C}_{t}^{-1}z_{t}^{-1}//_{t}^{-1}Y\left(
\Sigma_{t}\right) \\
&  =1_{\tau\leq t}\mathbf{X}\left(  \Sigma_{\tau}\right)  ^{\mathrm{tr}%
}//_{\tau}\left(  z_{\tau}^{-1}\right)  ^{\mathrm{tr}}\bar{C}_{t}^{-1}%
Z_{t}^{-1}Y\left(  \Sigma_{t}\right) \\
&  =1_{\tau\leq t}\mathbf{X}\left(  \Sigma_{\tau}\right)  ^{\mathrm{tr}%
}\left(  Z_{\tau}^{-1}\right)  ^{\mathrm{tr}}\bar{C}_{t}^{-1}Z_{t}%
^{-1}Y\left(  \Sigma_{t}\right)  .
\end{align*}
Hence, the formula for $\mathbb{Y}^{t}$ in Eq. (\ref{e.8.10}) may be
explicitly written as%
\begin{equation}
\mathbb{Y}_{s}^{t}=\left[  \int_{0}^{s\wedge t}\left(  Z_{\tau}^{-1}%
\mathbf{X}\left(  \Sigma_{\tau}\right)  \right)  ^{\mathrm{tr}}d\tau\right]
\bar{C}_{t}^{-1}Z_{t}^{-1}Y\left(  \Sigma_{t}\right)  . \label{e.8.11}%
\end{equation}
The reader should observe that the process $s\rightarrow\mathbb{Y}_{s}^{t}$ is
non-adapted since~$\bar{C}_{t}^{-1}Z_{t}^{-1}Y\left(  \Sigma_{t}\right)  $
depends on the entire path of $\Sigma$ up to time $t.$

\begin{theorem}
\label{t.8.7}Let $Y\in\Gamma\left(  TM\right)  $ and $\mathbb{Y}^{t}$ be the
non-adpated Cameron-Martin process defined in Eq. (\ref{e.8.11}). Then
$\mathbb{Y}^{t}$ is \textquotedblleft Malliavin smooth,\textquotedblright%
\ i.e. $\mathbb{Y}^{t}$ is $H$ -- differentiable (in the sense of Theorem
\ref{t.7.14}) to all orders with all differentials being in $L^{\infty
-}\left(  \mu\right)  ,$ see Nualart \cite{Nu1} for more precise definitions.
Moreover if $f\in C^{\infty}\left(  M\right)  ,$ then $f\left(  \Sigma
_{t}\right)  $ is Malliavin smooth and
\begin{equation}
\langle\bar{D}\left[  f\left(  \Sigma_{t}\right)  \right]  ,\mathbb{Y}%
^{t}\rangle_{H}=Yf\left(  \Sigma_{t}\right)  \label{e.8.12}%
\end{equation}
where $\bar{D}$ is the closure of the gradient operator defined in Corollary
\ref{c.7.16}.
\end{theorem}

\begin{proof}
We only sketch the proof here and refer the reader to
\cite{Norris86b,Bell5a,Nu1} with regard to some of the technical details which
are omitted below. Let $\left\{  e_{i}\right\}  _{i=1}^{d}$ be an orthonormal
basis for $T_{o}M,$ then
\begin{equation}
\mathbb{Y}_{s}^{t}=\sum_{i=1}^{d}\left\langle e_{i},\bar{C}_{t}^{-1}Z_{t}%
^{-1}Y\left(  \Sigma_{t}\right)  \right\rangle \int_{0}^{s}\left(  Z_{\tau
}^{-1}\mathbf{X}\left(  \Sigma_{\tau}\right)  \right)  ^{\mathrm{tr}}%
e_{i}d\tau=\sum_{i=1}^{d}a_{i}h_{s}^{i} \label{e.8.13}%
\end{equation}
where
\[
a_{i}:=\left\langle e_{i},\bar{C}_{t}^{-1}Z_{t}^{-1}Y\left(  \Sigma
_{t}\right)  \right\rangle \text{ and }h_{s}^{i}:=\int_{0}^{s\wedge t}\left(
Z_{\tau}^{-1}\mathbf{X}\left(  \Sigma_{\tau}\right)  \right)  ^{\mathrm{tr}%
}e_{i}d\tau.
\]
It is well known that solutions to stochastic differential equations with
smooth coefficients are Malliavin smooth from which it follows that $h^{i},$
$Z_{t}^{-1}Y\left(  \Sigma_{t}\right)  ,$ and $\bar{C}_{t}$ are Malliavin
smooth. It also follows from the general theory, under the conclusion of
Theorem \ref{t.8.6}, that $\bar{C}_{t}^{-1}$ is Malliavin smooth and hence so
are each the functions $a_{i}$ for $i=1,\dots d.$ Therefore, $\mathbb{Y}%
^{t}=\sum_{i=1}^{d}a_{i}h^{i}$ is Malliavin smooth as well and in particular
$\mathbb{Y}^{t}\in\mathcal{D}\left(  D^{\ast}\right)  .$ It now only remains
to verify Eq. (\ref{e.8.12}).

Let $h$ be a non-random element of $H.$ Then from Theorems \ref{t.7.14},
\ref{t.7.15}, \ref{t.7.26} and the chain rule for Wiener calculus,%
\begin{align*}
\mathbb{E}\left[  f\left(  \Sigma_{t}\right)  \cdot D^{\ast}h\right]   &
=\mathbb{E}\left[  \partial_{h}\left[  f\left(  \Sigma_{t}\right)  \right]
\right]  =\mathbb{E}\left[  df\left(  D\Sigma_{t}h\right)  \right] \\
&  =\mathbb{E}\left[  df\left(  Z_{t}\int_{0}^{t}Z_{\tau}^{-1}\mathbf{X}%
\left(  \Sigma_{\tau}\right)  \dot{h}_{\tau}d\tau\right)  \right] \\
&  =\mathbb{E}\left[  \left\langle \vec{\nabla}f\left(  \Sigma_{t}\right)
,Z_{t}\int_{0}^{t}Z_{\tau}^{-1}\mathbf{X}\left(  \Sigma_{\tau}\right)  \dot
{h}_{\tau}d\tau\right\rangle _{T_{\Sigma_{t}}M}\right] \\
&  =\mathbb{E}\left[  \int_{0}^{t}\left\langle \mathbf{X}\left(  \Sigma_{\tau
}\right)  ^{\mathrm{tr}}\left(  Z_{\tau}^{-1}\right)  ^{\mathrm{tr}}%
Z_{t}^{\mathrm{tr}}\vec{\nabla}f\left(  \Sigma_{t}\right)  \vec{\nabla
}f\left(  \Sigma_{t}\right)  ,\dot{h}_{\tau}\right\rangle _{\mathbb{R}^{n}%
}d\tau\right]
\end{align*}
from which we conclude that $f\left(  \Sigma_{t}\right)  \in\mathcal{D}\left(
D^{\ast\ast}\right)  =\mathcal{D}\left(  \bar{D}\right)  $ and
\[
\left(  \bar{D}\left[  f\left(  \Sigma_{t}\right)  \right]  \right)  _{s}%
=\int_{0}^{s\wedge t}\mathbf{X}\left(  \Sigma_{\tau}\right)  ^{\mathrm{tr}%
}\left(  Z_{\tau}^{-1}\right)  ^{\mathrm{tr}}Z_{t}^{\mathrm{tr}}\vec{\nabla
}f\left(  \Sigma_{t}\right)  d\tau.
\]
From this formula and the definition of $\mathbb{Y}^{t}$ it follows that
\begin{align*}
\langle &  \bar{D}\left[  f\left(  \Sigma_{t}\right)  \right]  ,\mathbb{Y}%
^{t}\rangle_{H}\\
&  =\int_{0}^{t}\left\langle \mathbf{X}\left(  \Sigma_{\tau}\right)
^{\mathrm{tr}}\left(  Z_{\tau}^{-1}\right)  ^{\mathrm{tr}}Z_{t}^{\mathrm{tr}%
}\vec{\nabla}f\left(  \Sigma_{t}\right)  ,\mathbf{X}\left(  \Sigma_{\tau
}\right)  ^{\mathrm{tr}}\left(  Z_{\tau}^{-1}\right)  ^{\mathrm{tr}}\bar
{C}_{t}^{-1}Z_{t}^{-1}Y\left(  \Sigma_{t}\right)  \right\rangle d\tau\\
&  =\left\langle \vec{\nabla}f\left(  \Sigma_{t}\right)  ,Z_{t}\left(
\int_{0}^{t}Z_{\tau}^{-1}\mathbf{X}\left(  \Sigma_{\tau}\right)  \left(
Z_{\tau}^{-1}\mathbf{X}\left(  \Sigma_{\tau}\right)  \right)  ^{\mathrm{tr}%
}d\tau\right)  \bar{C}_{t}^{-1}Z_{t}^{-1}Y\left(  \Sigma_{t}\right)
\right\rangle \\
&  =\left\langle \vec{\nabla}f\left(  \Sigma_{t}\right)  ,Z_{t}\bar{C}_{t}%
\bar{C}_{t}^{-1}Z_{t}^{-1}Y\left(  \Sigma_{t}\right)  \right\rangle
=\left\langle \vec{\nabla}f\left(  \Sigma_{t}\right)  ,Y\left(  \Sigma
_{t}\right)  \right\rangle \\
&  =\left(  Yf\right)  \left(  \Sigma_{t}\right)  .
\end{align*}

\end{proof}

\begin{notation}
\label{n.8.8}Let $\mathbb{Y}^{t}$ act on Malliavin smooth functions by the
formula, $\mathbb{Y}^{t}F:=\left\langle \bar{D}F,\mathbb{Y}^{t}\right\rangle
_{H}$ and let $\left(  \mathbb{Y}^{t}\right)  ^{\ast}$ denote the
$L^{2}\left(  \mu\right)  $ -- adjoint of $\mathbb{Y}^{t}.$
\end{notation}

With this notation, Theorem \ref{t.8.7} asserts that
\begin{equation}
\mathbb{Y}^{t}\left[  f\left(  \Sigma_{t}\right)  \right]  =\left(  Yf\right)
\left(  \Sigma_{t}\right)  . \label{e.8.14}%
\end{equation}
Now suppose $F,G:W\rightarrow\mathbb{R}$ are Malliavin smooth functions, then%
\begin{align*}
\mathbb{E}\left[  \mathbb{Y}^{t}F\cdot G+F\cdot\mathbb{Y}^{t}G\right]   &
=\mathbb{E}\left[  \mathbb{Y}^{t}\left[  FG\right]  \right]  =\mathbb{E}%
\left[  \left\langle \bar{D}\left[  FG\right]  ,\mathbb{Y}^{t}\right\rangle
_{H}\right] \\
&  =\mathbb{E}\left[  F\cdot GD^{\ast}\mathbb{Y}^{t}\right]
\end{align*}
from which it follows that $G\in\mathcal{D}\left(  \left(  \mathbb{Y}%
^{t}\right)  ^{\ast}\right)  $ and%
\begin{equation}
\left(  \mathbb{Y}^{t}\right)  ^{\ast}G=-\mathbb{Y}^{t}G+GD^{\ast}%
\mathbb{Y}^{t}. \label{e.8.15}%
\end{equation}
From the general theory (see \cite{Nu1} for example), $D^{\ast}U$ is Malliavin
smooth if $U$ is Malliavin smooth. In particular $\left(  \mathbb{Y}%
^{t}\right)  ^{\ast}G$ is Malliavin smooth if $G$ is Malliavin smooth.

\begin{theorem}
[Smoothness of Densities]\label{t.8.9}Assume the restricted H\"{o}rmander
condition holds at $o\in M$ (see Definition \ref{d.8.1}) and suppose $f\in
C^{\infty}\left(  M\right)  $ and $\left\{  Y_{i}\right\}  _{i=1}^{k}%
\subset\Gamma\left(  TM\right)  .$ Then%
\begin{align}
\mathbb{E}\left[  \left(  Y_{1}\dots Y_{k}f\right)  \left(  \Sigma_{t}\right)
\right]   &  =\mathbb{E}\left[  \mathbb{Y}_{1}^{t}\mathbb{\dots Y}_{k}%
^{t}\left[  f\left(  \Sigma_{t}\right)  \right]  \right] \nonumber\\
&  =\mathbb{E}\left[  \left[  f\left(  \Sigma_{t}\right)  \right]  \left(
\mathbb{Y}_{k}^{t}\right)  ^{\ast}\mathbb{\dots}\left(  \mathbb{Y}_{1}%
^{t}\right)  ^{\ast}1\right]  . \label{e.8.16}%
\end{align}
Moreover, the law of $\Sigma_{t}$ is smooth.
\end{theorem}

\begin{proof}
By an induction argument using Eq. (\ref{e.8.14}),
\[
\mathbb{Y}_{1}^{t}\mathbb{\dots Y}_{k}^{t}\left[  f\left(  \Sigma_{t}\right)
\right]  =\left(  Y_{1}\dots Y_{k}f\right)  \left(  \Sigma_{t}\right)
\]
from which Eq. (\ref{e.8.16}) is a simple consequence. As has already been
observed, $\left(  \mathbb{Y}_{k}^{t}\right)  ^{\ast}\mathbb{\dots}\left(
\mathbb{Y}_{1}^{t}\right)  ^{\ast}1$ is Malliavin smooth and in particular
$\left(  \mathbb{Y}_{k}^{t}\right)  ^{\ast}\mathbb{\dots}\left(
\mathbb{Y}_{1}^{t}\right)  ^{\ast}1\in L^{1}\left(  \mu\right)  .$ Therefore
it follows from Eq. (\ref{e.8.16}) that
\begin{equation}
\left\vert \mathbb{E}\left[  \left(  Y_{1}\dots Y_{k}f\right)  \left(
\Sigma_{t}\right)  \right]  \right\vert \leq\left\Vert \left(  \mathbb{Y}%
_{k}^{t}\right)  ^{\ast}\mathbb{\dots}\left(  \mathbb{Y}_{1}^{t}\right)
^{\ast}1\right\Vert _{L^{1}\left(  \mu\right)  }\left\Vert f\right\Vert
_{\infty}. \label{e.8.17}%
\end{equation}
Since the argument used in the proof of Theorem \ref{t.8.2} after Eq.
(\ref{e.8.16}) is local in nature, it follows from Eq. (\ref{e.8.17}) that the
$\mathrm{Law}(\Sigma_{t})$ has a smooth density relative to any smooth measure
on $M$ and in particular the Riemannian volume measure.
\end{proof}

\subsection{The Invertibility of $\bar{C}_{t}$ in the Elliptic
Case\label{s.8.3}}

As a warm-up to the proof of the full version of Theorem \ref{t.8.6} let us
first consider the special case where $\mathbf{X}\left(  m\right)
:\mathbb{R}^{n}\rightarrow T_{m}M$ is surjective for all $m\in M.$ Since $M$
is compact this will imply there exists and $\varepsilon>0$ such that
\[
\mathbf{X}(m)\mathbf{X}^{\mathrm{tr}}(m)\geq\varepsilon I_{T_{m}M}\text{ for
all }m\in M.
\]

\begin{notation}
\label{n.8.10}We will write $f\left(  \varepsilon\right)  =O\left(
\varepsilon^{\infty-}\right)  $ if, for all $p<\infty,$%
\[
\lim_{\varepsilon\downarrow0}\frac{\left\vert f\left(  \varepsilon\right)
\right\vert }{\varepsilon^{p}}=0.
\]

\end{notation}

\begin{proposition}
[Elliptic Case]\label{p.8.11}Suppose there is an $\varepsilon>0$ such that
\[
\mathbf{X}(m)\mathbf{X}^{\mathrm{tr}}(m)\geq\varepsilon I_{T_{m}M}%
\]
for all $m\in M,$ then $\left[  \det\left(  \bar{C}_{t}\right)  \right]
^{-1}\in L^{\infty-}\left(  \mu\right)  .$
\end{proposition}

\begin{proof}
Let $\delta\in(0,1)$ and
\begin{equation}
T_{\delta}:=\inf\left\{  t>0:\left\vert z_{t}-I_{T_{o}M}\right\vert
>\delta\right\}  \label{e.8.18}%
\end{equation}
where, as usual,%
\[
z_{t}:=//_{t}^{-1}Z_{t}=//_{t}^{-1}\left(  T_{t}^{B}\right)  _{\ast o}.
\]
Since for all $a\in T_{o}M,$%
\begin{align*}
\langle Z_{\tau}^{-1}\mathbf{X}(\Sigma_{\tau})  &  \mathbf{X}^{\mathrm{tr}%
}(\Sigma_{\tau})\left(  Z_{\tau}^{\mathrm{tr}}\right)  ^{-1}a,a\rangle\\
&  =\left\langle \mathbf{X}(\Sigma_{\tau})\mathbf{X}^{\mathrm{tr}}%
(\Sigma_{\tau})\left(  Z_{\tau}^{\mathrm{tr}}\right)  ^{-1}a,\left(  Z_{\tau
}^{\mathrm{tr}}\right)  ^{-1}a\right\rangle \\
&  \geq\varepsilon\left\langle \left(  Z_{\tau}^{\mathrm{tr}}\right)
^{-1}a,\left(  Z_{\tau}^{\mathrm{tr}}\right)  ^{-1}a\right\rangle
=\varepsilon\left\langle a,Z_{\tau}^{\mathrm{tr}}\left(  Z_{\tau}%
^{\mathrm{tr}}\right)  ^{-1}a\right\rangle ,
\end{align*}
we have%
\begin{align*}
Z_{\tau}^{-1}\mathbf{X}  &  (\Sigma_{\tau})\mathbf{X}^{\mathrm{tr}}%
(\Sigma_{\tau})\left(  Z_{\tau}^{\mathrm{tr}}\right)  ^{-1}\\
&  \geq\varepsilon Z_{\tau}^{\mathrm{tr}}\left(  Z_{\tau}^{\mathrm{tr}%
}\right)  ^{-1}=\varepsilon z_{t}^{\mathrm{tr}}//_{t}^{\mathrm{tr}}\left(
//_{t}^{\mathrm{tr}}\right)  ^{-1}\left(  z_{t}^{\mathrm{tr}}\right)
^{-1}=\varepsilon z_{t}^{\mathrm{tr}}\left(  z_{t}^{\mathrm{tr}}\right)
^{-1}.
\end{align*}
Hence%
\begin{align*}
\bar{C}_{t}  &  =\int_{0}^{t}Z_{\tau}^{-1}\mathbf{X}(\Sigma_{\tau}%
)\mathbf{X}^{\mathrm{tr}}(\Sigma_{\tau})\left(  Z_{\tau}^{\mathrm{tr}}\right)
^{-1}d\tau\\
&  \geq\varepsilon\int_{0}^{t}Z_{\tau}^{-1}\left(  Z_{\tau}^{\mathrm{tr}%
}\right)  ^{-1}d\tau\geq\varepsilon\int_{0}^{t\wedge T_{\delta}}z_{\tau
}^{\mathrm{tr}}\left(  z_{\tau}^{\mathrm{tr}}\right)  ^{-1}d\tau
\end{align*}
and therefore,%
\[
\bar{\Delta}_{t}=\det\left(  \bar{C}_{t}\right)  \geq\varepsilon^{n}%
\det\left(  \int_{0}^{t\wedge T_{\delta}}z_{\tau}^{\mathrm{tr}}\left(
z_{\tau}^{\mathrm{tr}}\right)  ^{-1}d\tau\right)  .
\]
By choosing $\delta>0$ sufficiently small we may arrange that
\[
\left\Vert z_{\tau}^{\mathrm{tr}}\left(  z_{\tau}^{\mathrm{tr}}\right)
^{-1}-I\right\Vert \leq1/2
\]
for all $\tau\leq t\wedge T_{\delta}$ in which case%
\[
\int_{0}^{t\wedge T_{\delta}}z_{\tau}^{\mathrm{tr}}\left(  z_{\tau
}^{\mathrm{tr}}\right)  ^{-1}d\tau\geq\frac{1}{2}t\wedge T_{\delta}\cdot Id
\]
and hence $\bar{\Delta}_{t}=\det\left(  \bar{C}_{t}\right)  \geq
\varepsilon^{n}\frac{1}{2}t\wedge T_{\delta}.$ From this it follows that
\[
\mathbb{E}\left[  \bar{\Delta}_{t}^{-p}\right]  \leq2^{p}\varepsilon
^{-np}\mathbb{E}\left(  \left(  \frac{1}{t\wedge T_{\delta}}\right)
^{p}\right)  .
\]
Now%
\begin{align*}
\mathbb{E}\left(  \left(  \frac{1}{t\wedge T_{\delta}}\right)  ^{p}\right)
&  =\mathbb{E}\left(  -\int_{t\wedge T_{\delta}}^{\infty}\frac{d}{d\tau}%
\tau^{-p}d\tau\right)  =\mathbb{E}\left(  p\int_{0}^{\infty}1_{t\wedge
T_{\delta}\leq\tau}\cdot\tau^{-p-1}d\tau\right) \\
&  =p\int_{0}^{\infty}\tau^{-p-1}\mu\left(  t\wedge T_{\delta}\leq\tau\right)
d\tau
\end{align*}
which will be finite for all $p>1$ iff $\mu\left(  t\wedge T_{\delta}\leq
\tau\right)  =\mu\left(  T_{\delta}\leq\tau\right)  =O(\tau^{k})$ as
$\tau\downarrow0$ for all $k>0.$

By Chebyschev's inequalities and Eq. (\ref{e.9.10}) of Proposition \ref{p.9.5}
below,%
\begin{equation}
\mu\left(  T_{\delta}\leq\tau\right)  =\mu\left(  \sup_{s\leq\tau}\left\vert
z_{s}-I\right\vert >\delta\right)  \leq\delta^{-p}\mathbb{E}\left[
\sup_{s\leq\tau}\left\vert z_{s}-I\right\vert ^{p}\right]  =O(\tau^{p/2}).
\label{e.8.19}%
\end{equation}
Since $p\geq2$ was arbitrary it follows that $\mu(T_{\delta}\leq\tau)=O\left(
\tau^{\infty-}\right)  $ which completes the proof.
\end{proof}

\subsection{Proof of Theorem \ref{t.8.6}\label{s.8.4}}

\begin{notation}
\label{n.8.12}Let $S:=\left\{  v\in T_{o}M:\langle v,v\rangle=1\right\}  ,$
i.e. $S$ is the unit sphere in $T_{o}M.$
\end{notation}

\begin{proof}
(\emph{Proof of Theorem \ref{t.8.6}.}) To show $\bar{C}_{t}^{-1}\in
L^{\infty-}\left(  \mu\right)  $ it suffices to shows
\[
\mu(\inf_{v\in S}\langle\bar{C}_{t}v,v\rangle<\varepsilon)=O(\varepsilon
^{\infty-}).
\]
To verify this claim, notice that $\lambda_{0}:=\inf_{v\in S}\langle\bar
{C}_{t}v,v\rangle$ is the smallest eigenvalue of $\bar{C}_{t}.$ Since
$\det\bar{C}_{t}$ is the product of the eigenvalues of $\bar{C}_{t}$ it
follows that $\bar{\Delta}_{t}:=\det\bar{C}_{t}\geq\lambda_{0}^{n}$ and so
$\left\{  \det\bar{C}_{t}<\varepsilon^{n}\right\}  \subset\left\{  \lambda
_{0}<\varepsilon\right\}  $ and hence%
\[
\mu\left(  \det\bar{C}_{t}<\varepsilon^{n}\right)  \leq\mu\left(  \lambda
_{0}<\varepsilon\right)  =O(\varepsilon^{\infty-}).
\]
By replacing $\varepsilon$ by $\varepsilon^{1/n}$ above this implies
$\mu\left(  \bar{\Delta}_{t}<\varepsilon\right)  =O(\varepsilon^{\infty-}).$
From this estimate it then follows that
\begin{align*}
\mathbb{E}\left[  \bar{\Delta}_{t}^{-q}\right]   &  =\mathbb{E}\int
_{\bar{\Delta}_{t}}^{\infty}q\tau^{-q-1}d\tau=q\mathbb{E}\int_{0}^{\infty
}1_{\bar{\Delta}_{t}\leq\tau}~\tau^{-q-1}d\tau\\
&  =q\int_{0}^{\infty}\mu(\bar{\Delta}_{t}\leq\tau)~\tau^{-q-1}d\tau=q\int
_{0}^{\infty}O(\tau^{p})~\tau^{-q-1}d\tau
\end{align*}
which is seen to be finite by taking $p\geq q+1.$

More generally if $T$ is any stopping time with $T\leq t,$ since $\langle
\bar{C}_{T}v,v\rangle\leq\langle\bar{C}_{t}v,v\rangle$ for all $v\in S$ it
suffices to prove%
\begin{equation}
\mu\left(  \inf_{v\in S}\langle\bar{C}_{T}v,v\rangle<\varepsilon\right)
=O(\varepsilon^{\infty-}). \label{e.8.20}%
\end{equation}
According to Lemma \ref{l.8.13} and Proposition \ref{p.8.15} below, Eq.
(\ref{e.8.20}) holds with%
\begin{equation}
T=T_{\delta}:=\inf\left\{  t>0:\max\left\{  \left\vert z_{t}-I_{T_{o}%
M}\right\vert ,\mathrm{dist}(\Sigma_{t},\Sigma_{0})\right\}  >\delta\right\}
\label{e.8.21}%
\end{equation}
provided $\delta>0$ is chosen sufficiently small.
\end{proof}

The rest of this section is now devoted to the proof of Lemma \ref{l.8.13} and
Proposition \ref{p.8.15} below. In what follows we will make repeated use of
the identity,%
\begin{equation}
\langle\bar{C}_{T}v,v\rangle=\sum_{i=1}^{n}\int_{0}^{T}\left\langle Z_{\tau
}^{-1}X_{i}(\Sigma_{\tau}),v\right\rangle ^{2}d\tau. \label{e.8.22}%
\end{equation}
To prove this, let $\left\{  e_{i}\right\}  _{i=1}^{n}$ be the standard basis
for $\mathbb{R}^{n}.$ Then
\begin{align*}
Z_{\tau}^{-1}\mathbf{X}(\Sigma_{\tau})\mathbf{X}^{\mathrm{tr}}(\Sigma_{\tau
})\left(  Z_{\tau}^{\mathrm{tr}}\right)  ^{-1}v  &  =\sum_{i=1}^{n}Z_{\tau
}^{-1}\mathbf{X}(\Sigma_{\tau})e_{i}~\langle e_{i},\mathbf{X}^{\mathrm{tr}%
}(\Sigma_{\tau})\left(  Z_{\tau}^{\mathrm{tr}}\right)  ^{-1}v\rangle\\
&  =\sum_{i=1}^{n}\langle Z_{\tau}^{-1}X_{i}(\Sigma_{\tau}),v\rangle~Z_{\tau
}^{-1}X_{i}(\Sigma_{\tau})
\end{align*}
so that
\[
\left\langle Z_{\tau}^{-1}\mathbf{X}(\Sigma_{\tau})\mathbf{X}^{\mathrm{tr}%
}(\Sigma_{\tau})\left(  Z_{\tau}^{\mathrm{tr}}\right)  ^{-1}v,v\right\rangle
=\sum_{i=1}^{n}\left\langle Z_{\tau}^{-1}X_{i}(\Sigma_{\tau}),v\right\rangle
^{2}%
\]
which upon integrating on $\tau$ gives Eq. (\ref{e.8.22}).

In the proofs below, there will always be an implied sum on repeated indices.

\begin{lemma}
[Compactness Argument]\label{l.8.13}Let $T_{\delta}$ be as in Eq.
(\ref{e.8.21}) and suppose for all $v\in S$ there exists $i\in\{1,\dots,n\}$
and an open neighborhood $N\subset_{o}S$ of $v$ such that
\begin{equation}
\sup_{u\in N}\mu\left(  \int_{0}^{T_{\delta}}\left\langle Z_{\tau}^{-1}%
X_{i}(\Sigma_{\tau}),u\right\rangle ^{2}d\tau<\varepsilon\right)  =O\left(
\varepsilon^{\infty-}\right)  , \label{e.8.23}%
\end{equation}
then Eq. (\ref{e.8.20}) holds provided $\delta>0$ is sufficiently small.
\end{lemma}

\begin{proof}
By compactness of $S,$ it follows from Eq. (\ref{e.8.23}) that%
\begin{equation}
\sup_{u\in S}\mu\left(  \int_{0}^{T_{\delta}}\left\langle Z_{\tau}^{-1}%
X_{i}(\Sigma_{\tau}),u\right\rangle ^{2}d\tau<\varepsilon\right)  =O\left(
\varepsilon^{\infty-}\right)  . \label{e.8.24}%
\end{equation}
For $w\in T_{o}M,$ let $\partial_{w}$ denote the directional derivative acting
on functions $f\left(  v\right)  $ with $v\in T_{o}M.$ Because for all
$v,w\in\mathbb{R}^{n}$ with $\left\vert v\right\vert \leq1$ and $\left\vert
w\right\vert \leq1$ (using Eq. (\ref{e.8.22})),%
\begin{align*}
\left\vert \partial_{w}\left\langle \bar{C}_{T_{\delta}}v,v\right\rangle
\right\vert \leq &  2\sum_{i=1}^{n}\int_{0}^{T_{\delta}}\left\vert
\left\langle Z_{\tau}^{-1}X_{i}(\Sigma_{\tau}),v\right\rangle \left\langle
Z_{\tau}^{-1}X_{i}(\Sigma_{\tau}),w\right\rangle \right\vert d\tau\\
\leq &  2\sum_{i=1}^{n}\int_{0}^{T_{\delta}}\left\vert Z_{\tau}^{-1}%
X_{i}(\Sigma_{\tau})\right\vert _{\operatorname*{Hom}\left(  \mathbb{R}%
^{n},T_{o}M\right)  }^{2}d\tau\\
&  \qquad=2\sum_{i=1}^{n}\int_{0}^{T_{\delta}}\left\vert z_{\tau}^{-1}%
//_{\tau}^{-1}X_{i}(\Sigma_{\tau})\right\vert _{\operatorname*{Hom}\left(
\mathbb{R}^{n},T_{o}M\right)  }^{2}d\tau,
\end{align*}
by choosing $\delta>0$ in Eq. (\ref{e.8.21}) sufficiently small we may assume
there is a non-random constant $\theta<\infty$ such that
\[
\sup_{\left\vert v\right\vert ,\left\vert w\right\vert \leq1}\left\vert
\partial_{w}\left\langle \bar{C}_{T_{\delta}}v,v\right\rangle \right\vert
\leq\theta<\infty.
\]
With this choice of $\delta,$ if $v,w\in S$ satisfy $\left\vert v-w\right\vert
<\theta/\varepsilon$ then%
\begin{equation}
\left\vert \left\langle \bar{C}_{T_{\delta}}v,v\right\rangle -\left\langle
\bar{C}_{T_{\delta}}w,w\right\rangle \right\vert <\varepsilon. \label{e.8.25}%
\end{equation}
There exists $D<\infty$ satisfying: for any $\varepsilon>0,$ there is an open
cover of $S$ with at most $D\cdot\left(  \theta/\varepsilon\right)  ^{n}$
balls of the form $B(v_{j},\varepsilon/\theta).$ From Eq. (\ref{e.8.25}), for
any $v\in S$ there exists $j$ such that $v\in B(v_{j},\varepsilon/\theta)\cap
S$ and
\[
\left\vert \left\langle \bar{C}_{T_{\delta}}v,v\right\rangle -\left\langle
\bar{C}_{T_{\delta}}v_{j},v_{j}\right\rangle \right\vert <\varepsilon.
\]
So if $\inf_{v\in S}\left\langle \bar{C}_{T_{\delta}}v,v\right\rangle
<\varepsilon$ then $\min_{j}\left\langle \bar{C}_{T_{\delta}}v_{j}%
,v_{j}\right\rangle <2\varepsilon,$ i.e.%
\[
\left\{  \inf_{v\in S}\left\langle \bar{C}_{T_{\delta}}v,v\right\rangle
<\varepsilon\right\}  \subset\left\{  \min_{j}\left\langle \bar{C}_{T_{\delta
}}v_{j},v_{j}\right\rangle <2\varepsilon\right\}  \subset\bigcup_{j}\left\{
\left\langle \bar{C}_{T_{\delta}}v_{j},v_{j}\right\rangle <2\varepsilon
\right\}  .
\]
Therefore,%
\begin{align*}
\mu\left(  \inf_{v\in S}\left\langle \bar{C}_{T_{\delta}}v,v\right\rangle
<\varepsilon\right)   &  \leq\sum_{j}\mu\left(  \left\langle \bar
{C}_{T_{\delta}}v_{j},v_{j}\right\rangle <2\varepsilon\right) \\
&  \leq D\cdot\left(  \theta/\varepsilon\right)  ^{n}\cdot\sup_{v\in S}%
\mu\left(  \left\langle \bar{C}_{T_{\delta}}v,v\right\rangle <2\varepsilon
\right) \\
&  \leq D\cdot\left(  \theta/\varepsilon\right)  ^{n}O(\varepsilon^{\infty
-})=O(\varepsilon^{\infty-}).
\end{align*}

\end{proof}

The following important proposition is the Stochastic version of Theorem
\ref{t.4.9}. It gives the first hint that H\"{o}rmander's condition in
Definition \ref{d.8.1} is relevant to showing $\bar{\Delta}_{t}^{-1}\in
L^{\infty-}\left(  \mu\right)  $ or equivalently that $\bar{C}_{t\ }^{-1}\in
L^{\infty-}\left(  \mu\right)  .$

\begin{proposition}
[The appearance of commutators]\label{p.8.14}Let $W\in\Gamma\left(  TM\right)
,$ then%
\begin{equation}
\delta\left[  Z_{s}^{-1}W(\Sigma_{s})\right]  =Z_{s}^{-1}[X_{0},W](\Sigma
_{s})ds+Z_{s}^{-1}\sum_{i=1}^{n}[X_{i},W](\Sigma_{s})\delta B_{s}^{i}.
\label{e.8.26}%
\end{equation}
This may also be written in It\^{o} form as%
\begin{align}
d\left[  Z_{s}^{-1}W(\Sigma_{s})\right]   &  =Z_{s}^{-1}[X_{i},W](\Sigma
_{s})dB_{s}^{i}\nonumber\\
&  +\left\{  Z_{s}^{-1}[X_{0},W](\Sigma_{s})+\frac{1}{2}\sum_{i=1}^{n}%
Z_{s}^{-1}\left(  L_{X_{i}}^{2}W\right)  (\Sigma_{s})\right\}  ds,
\label{e.8.27}%
\end{align}
where $L_{X}W:=\left[  X,W\right]  $ as in Theorem \ref{t.4.9}.
\end{proposition}

\begin{proof}
Write $W\left(  \Sigma_{s}\right)  =Z_{s}w_{s},$ i.e. let $w_{s}:=Z_{s}%
^{-1}W(\Sigma_{s}).$ By Proposition \ref{p.5.36} and Theorem \ref{t.5.41},%
\begin{align*}
\nabla_{\delta\Sigma_{s}}W  &  =\delta^{\nabla}\left[  W\left(  \Sigma
_{s}\right)  \right]  =\delta^{\nabla}\left[  Z_{s}w_{s}\right]  =\left(
\delta^{\nabla}Z_{s}\right)  w_{s}+Z_{s}\delta w_{s}\\
&  =\left(  \nabla_{Z_{s}w_{s}}\mathbf{X}\right)  \delta B_{s}+\left(
\nabla_{Z_{s}w_{s}}X_{0}\right)  ds+Z_{s}\delta w_{s}.
\end{align*}
Therefore, using the fact that $\nabla$ has zero torsion (see Proposition
\ref{p.3.36}),%
\begin{align*}
\delta w_{s}  &  =Z_{s}^{-1}\left[  \nabla_{\delta\Sigma_{s}}W-\left(
\nabla_{Z_{s}w_{s}}\mathbf{X}\right)  \delta B_{s}+\left(  \nabla_{Z_{s}w_{s}%
}X_{0}\right)  ds\right] \\
&  =Z_{s}^{-1}\left[  \nabla_{\mathbf{X}\left(  \Sigma_{s}\right)  \delta
B_{s}+X_{0}\left(  \Sigma_{s}\right)  ds}W-\left(  \nabla_{W\left(  \Sigma
_{s}\right)  }\mathbf{X}\right)  \delta B_{s}+\left(  \nabla_{W\left(
\Sigma_{s}\right)  }X_{0}\right)  ds\right] \\
&  =Z_{s}^{-1}\left[  \left(  \nabla_{X_{i}\left(  \Sigma_{s}\right)
}W-\nabla_{W\left(  \Sigma_{s}\right)  }X_{i}\right)  \delta B_{s}^{i}+\left(
\nabla_{X_{0}\left(  \Sigma_{s}\right)  }W-\nabla_{W\left(  \Sigma_{s}\right)
}X_{0}\right)  ds\right] \\
&  =Z_{s}^{-1}\left(  \left[  X_{i},W\right]  \left(  \Sigma_{s}\right)
\delta B_{s}^{i}+[X_{0},W](\Sigma_{s})ds\right)
\end{align*}
which proves Eq. (\ref{e.8.26}).

Applying Eq. (\ref{e.8.26}) with $W$ replaced by $[X_{i},W]$ implies%
\[
d\left[  Z_{s}^{-1}[X_{i},W](\Sigma_{s})\right]  =Z_{s}^{-1}[X_{j}%
,[X_{i},W]](\Sigma_{s})dB_{s}^{j}+d\left[  BV\right]  ,
\]
where $BV$ denotes process of bounded variation. Hence%
\begin{align*}
Z_{s}^{-1}[X_{i},W](\Sigma_{s})\delta B_{s}^{i}  &  =Z_{s}^{-1}[X_{i}%
,W](\Sigma_{s})dB_{s}^{i}+\frac{1}{2}d\left\{  Z_{s}^{-1}[X_{i},W](\Sigma
_{s})\right\}  dB_{s}^{i}\\
&  =Z_{s}^{-1}[X_{i},W](\Sigma_{s})dB_{s}^{i}+\frac{1}{2}Z_{s}^{-1}%
[X_{j},[X_{i},W]](\Sigma_{s})dB_{s}^{j}dB_{s}^{i}\\
&  =Z_{s}^{-1}[X_{i},W](\Sigma_{s})dB_{s}^{i}+\frac{1}{2}Z_{s}^{-1}%
[X_{i},[X_{i},W]](\Sigma_{s})ds
\end{align*}
which combined with Eq. (\ref{e.8.26}) proves Eq. (\ref{e.8.27}).
\end{proof}

\begin{proposition}
\label{p.8.15}Let $T_{\delta}$ be as in Eq. (\ref{e.8.21}). If H\"{o}rmander's
restricted bracket condition holds at $o\in M$ and $v\in S$ is given, there
exists $i\in\left\{  1,2,\dots,n\right\}  $ and an open neighborhood
$U\subset_{o}S$ of $v$ such that%
\[
\sup_{u\in U}\mu\left(  \int_{0}^{T_{\delta}}\left\langle Z_{\tau}^{-1}%
X_{i}(\Sigma_{\tau}),u\right\rangle ^{2}d\tau\leq\varepsilon\right)  =O\left(
\varepsilon^{\infty-}\right)  .
\]

\end{proposition}

\begin{proof}
The proof given here will follow Norris \cite{Norris86b}. H\"{o}rmander's
condition implies there exist $l\in\mathbb{N}$ and $\beta>0$ such that%
\[
\frac{1}{\left\vert \mathcal{K}_{l}\right\vert }\sum_{K\in\mathcal{K}_{l}%
}K(o)K(o)^{\mathrm{tr}}\geq3\beta I
\]
or equivalently put for all $v\in S,$%
\[
3\beta\leq\frac{1}{\left\vert \mathcal{K}_{l}\right\vert }\sum_{K\in
\mathcal{K}_{l}}\left\langle K(o),v\right\rangle ^{2}\leq\max_{K\in
\mathcal{K}_{l}}\left\langle K(o),v\right\rangle ^{2}.
\]
By choosing $\delta>0$ in Eq. (\ref{e.8.21}) sufficiently small we may assume
that
\[
\max_{K\in\mathcal{K}_{l}}\inf_{\tau\leq T_{\delta}}\left\langle Z_{\tau}%
^{-1}K(\Sigma_{\tau}),v\right\rangle ^{2}\geq2\beta\text{ for all }v\in S.
\]
Fix a $v\in S$ and $K\in\mathcal{K}_{l}$ such that
\[
\inf_{\tau\leq T_{\delta}}\left\langle Z_{\tau}^{-1}K(\Sigma_{\tau
}),v\right\rangle ^{2}\geq2\beta
\]
and choose an open neighborhood $U\subset S$ of $v$ such that
\[
\inf_{\tau\leq T_{\delta}}\left\langle Z_{\tau}^{-1}K(\Sigma_{\tau
}),u\right\rangle ^{2}\geq\beta\text{ for all }u\in U.
\]
Then, using Eq. (\ref{e.8.19}),
\begin{align}
\sup_{u\in U}  &  \mu\left(  \int_{0}^{T_{\delta}}\left\langle Z_{\tau}%
^{-1}K(\Sigma_{\tau}),u\right\rangle ^{2}d\tau\leq\varepsilon\right)
\nonumber\\
&  \leq\mu\left(  \int_{0}^{T_{\delta}}\beta dt\leq\varepsilon\right)
=\mu\left(  T_{\delta}\leq\varepsilon/\beta\right)  =O\left(  \varepsilon
^{\infty-}\right)  . \label{e.8.28}%
\end{align}

Write $K=L_{X_{i_{r}}}\dots L_{X_{i_{2}}X_{i_{1}}}$ with $r\leq l.$ If it
happens that $r=1$ then Eq. (\ref{e.8.28}) becomes%
\[
\sup_{u\in U}\mu\left(  \left\langle \bar{C}_{T_{\delta}}u,u\right\rangle
\leq\varepsilon\right)  \leq\sup_{u\in U}\mu\left(  \int_{0}^{T_{\delta}%
}\left\langle Z_{\tau}^{-1}X_{i_{1}}(\Sigma_{\tau}),u\right\rangle ^{2}%
dt\leq\varepsilon\right)  =O\left(  \varepsilon^{\infty-}\right)
\]
and we are done. So now suppose $r>1$ and set
\[
K_{j}=L_{X_{i_{j}}}\dots L_{X_{i_{2}}}X_{i_{1}}\text{ for }j=1,2,\dots,r
\]
so that $K_{r}=K.$ We will now show by (decreasing) induction on $j$ that%
\begin{equation}
\sup_{u\in U}\mu\left(  \int_{0}^{T_{\delta}}\left\langle Z_{\tau}^{-1}%
K_{j}(\Sigma_{\tau}),u\right\rangle ^{2}dt\leq\varepsilon\right)  =O\left(
\varepsilon^{\infty-}\right)  . \label{e.8.29}%
\end{equation}
From Proposition \ref{p.8.14} we have%
\begin{align*}
d\left[  Z_{t}^{-1}K_{j-1}(\Sigma_{t})\right]  =Z_{t}^{-1}  &  [X_{i}%
,K_{j-1}](\Sigma_{t})dB^{i}(t)\\
&  +\left\{  Z_{t}^{-1}[X_{0},K_{j-1}](\Sigma_{t})+\frac{1}{2}Z_{t}%
^{-1}\left(  L_{X_{i}}^{2}K_{j-1}\right)  (\Sigma_{t})\right\}  dt
\end{align*}
which upon integrating on $t$ gives%
\begin{align*}
\left\langle Z_{t}^{-1}K_{j-1}(\Sigma_{t}),u\right\rangle  &  =\left\langle
K_{j-1}(\Sigma_{0}),u\right\rangle +\int_{0}^{t}\left\langle Z_{\tau}%
^{-1}[X_{i},K_{j-1}](\Sigma_{\tau}),u\right\rangle dB_{\tau}^{i}\\
&  +\int_{0}^{t}\left\langle Z_{\tau}^{-1}[X_{0},K_{j-1}](\Sigma_{\tau}%
)+\frac{1}{2}Z_{t}^{-1}\left(  L_{X_{i}}^{2}K_{j-1}\right)  (\Sigma_{\tau
}),u\right\rangle d\tau.
\end{align*}
Applying Proposition \ref{p.9.13} of the appendix with $T=T_{\delta,}$%
\begin{align*}
Y_{t}  &  :=\left\langle Z_{t}^{-1}K_{j-1}(\Sigma_{t}),u\right\rangle ,\text{
}y=\left\langle K_{j-1}(\Sigma_{0}),u\right\rangle ,\\
M_{t}  &  =\int_{0}^{t}\left\langle Z_{\tau}^{-1}[X_{i},K_{j-1}](\Sigma_{\tau
}),u\right\rangle dB_{\tau}^{i}\text{ and }\\
A_{t}  &  :=\int_{0}^{t}\left\langle Z_{\tau}^{-1}[X_{0},K_{j-1}](\Sigma
_{\tau})+\frac{1}{2}Z_{\tau}^{-1}\left(  L_{X_{i}}^{2}K_{j-1}\right)
(\Sigma_{\tau}),u\right\rangle dt
\end{align*}
implies%
\begin{equation}
\sup_{u\in U}\mu\left(  \Omega_{1}\left(  u\right)  \cap\Omega_{2}\left(
u\right)  \right)  =O\left(  \varepsilon^{\infty-}\right)  , \label{e.8.30}%
\end{equation}
where%
\begin{align*}
\Omega_{1}\left(  u\right)   &  :=\left\{  \int_{0}^{T_{\delta}}\left\langle
Z_{t}^{-1}K_{j-1}(\Sigma_{t}),u\right\rangle ^{2}dt<\varepsilon^{q}\right\}
,\\
\Omega_{2}\left(  u\right)   &  :=\left\{  \int_{0}^{T_{\delta}}\sum_{i=1}%
^{n}\left\langle Z_{\tau}^{-1}[X_{i},K_{j-1}](\Sigma_{\tau}),u\right\rangle
^{2}d\tau\geq\varepsilon\right\}
\end{align*}
and $q>4.$ Since
\begin{align*}
\sup_{u\in U}\mu\left(  \left[  \Omega_{2}\left(  u\right)  \right]
^{c}\right)   &  =\sup_{u\in U}\mu\left(  \int_{0}^{T_{\delta}}\sum_{i=1}%
^{n}\left\langle Z_{\tau}^{-1}[X_{i},K_{j-1}](\Sigma_{\tau}),u\right\rangle
^{2}d\tau<\varepsilon\right) \\
&  \leq\sup_{u\in U}\mu\left(  \int_{0}^{T_{\delta}}\left\langle Z_{\tau}%
^{-1}K_{j}(\Sigma_{\tau}),u\right\rangle ^{2}d\tau<\varepsilon\right)
\end{align*}
we may applying the induction hypothesis to learn,%
\begin{equation}
\sup_{u\in U}\mu\left(  \left[  \Omega_{2}\left(  u\right)  \right]
^{c}\right)  =O\left(  \varepsilon^{\infty-}\right)  . \label{e.8.31}%
\end{equation}
It now follows from Eqs. (\ref{e.8.30}) and (\ref{e.8.31}) that%
\begin{align*}
\sup_{u\in U}\mu(\Omega_{1}\left(  u\right)  )  &  \leq\sup_{u\in U}\mu
(\Omega_{1}\left(  u\right)  \cap\Omega_{2}\left(  u\right)  )+\sup_{u\in
U}\mu(\Omega_{1}\left(  u\right)  \cap\left[  \Omega_{2}\left(  u\right)
\right]  ^{c})\\
&  \leq\sup_{u\in U}\mu(\Omega_{1}\left(  u\right)  \cap\Omega_{2}\left(
u\right)  )+\sup_{u\in U}\mu(\left[  \Omega_{2}\left(  u\right)  \right]
^{c})\\
&  =O\left(  \varepsilon^{\infty-}\right)  +O\left(  \varepsilon^{\infty
-}\right)  =O\left(  \varepsilon^{\infty-}\right)  ,
\end{align*}
i.e.
\[
\sup_{u\in U}\mu\left(  \int_{0}^{T_{\delta}}\left\langle Z_{t}^{-1}%
K_{j-1}(\Sigma_{t}),u\right\rangle ^{2}dt<\varepsilon^{q}\right)  =O\left(
\varepsilon^{\infty-}\right)  .
\]
Replacing $\varepsilon$ by $\varepsilon^{1/q}$ in the previous equation, using
$O\left(  \left(  \varepsilon^{1/q}\right)  ^{\infty-}\right)  =O\left(
\varepsilon^{\infty-}\right)  ,$ completes the induction argument and hence
the proof.
\end{proof}

\subsection{More References}

The literature on the \textquotedblleft Malliavin calculus\textquotedblright%
\ is very extensive and I will not make any attempt at summarizing it here.
Let me just add to references already mentioned the articles in
\cite{Taniguchi83,Imkeller97,Schiltz98a} which carry out Malliavin's method in
the geometric context of these notes. Also see \cite{Picard02} for another
method which works if H\"{o}rmander's bracket condition holds at level $2,$
namely when
\[
\mathrm{span}(\left\{  K(m):K\in\mathcal{K}_{2}\right\}  )=T_{m}M\text{ for
all }m\in M,
\]
see Definition \ref{d.8.1}. The reader should also be aware of the deep
results of Ben Arous and Leandre in
\cite{BenArous89a,BenArous89b,B-AL91a,B-AL91b,Leand92a}.

\section{Appendix: Martingale and SDE Estimates\label{s.9}}

In this appendix $\left\{  B_{t}:t\geq0\right\}  $ will denote and
$\mathbb{R}^{n}$ -- valued Brownian motion, $\left\{  \beta_{t}:t\geq
0\right\}  $ will be a one dimensional Brownian motion and, unlike in the
text, we will use the more standard letter $P$ rather than $\mu$ to denote the
underlying probability measure.

\begin{notation}
\label{n.9.1}When $M_{t}$ is a martingale and $A_{t}$ is a process of bounded
variation let $\langle M\rangle_{t}$ be the quadratic variation of $M$ and
$\left\vert A\right\vert _{t}$ be the total variation of $A$ up to time $t.$
\end{notation}

\subsection{Estimates of Wiener Functionals Associated to SDE's\label{s.9.1}}

\begin{proposition}
\label{p.9.2}Suppose $p\in\lbrack2,\infty),$ $\alpha_{\tau}$ and $A_{\tau}$
are predictable $\mathbb{R}^{d}$ and $\operatorname*{Hom}\left(
\mathbb{R}^{n},\mathbb{R}^{d}\right)  $ -- valued processes respectively and%
\begin{equation}
Y_{t}:=\int_{0}^{t}A_{\tau}dB_{\tau}+\int_{0}^{t}\alpha_{\tau}d\tau.
\label{e.9.1}%
\end{equation}
Then, letting $Y_{t}^{\ast}:=\sup_{\tau\leq t}\left\vert Y_{\tau}\right\vert
,$%
\begin{equation}
\mathbb{E}\left(  Y_{t}^{\ast}\right)  ^{p}\leq C_{p}\left\{  \mathbb{E}%
\left(  \int_{0}^{t}\left\vert A_{\tau}\right\vert ^{2}d\tau\right)
^{p/2}+\mathbb{E}\left(  \int_{0}^{t}\left\vert \alpha_{\tau}\right\vert
d\tau\right)  ^{p}\right\}  \label{e.9.2}%
\end{equation}
where
\[
\left\vert A\right\vert ^{2}=\operatorname*{tr}\left(  AA^{\ast}\right)
=\sum_{i=1}^{n}\left(  AA^{\ast}\right)  _{ii}=\sum_{i,j}A_{ij}A_{ij}%
=\operatorname*{tr}\left(  A^{\ast}A\right)  .
\]

\end{proposition}

\begin{proof}
We may assume the right side of Eq. (\ref{e.9.2}) is finite for otherwise
there is nothing to prove. For the moment further assume $\alpha\equiv0.$ By a
standard limiting argument involving stopping times we may further assume
there is a non-random constant $C<\infty$ such that
\[
Y_{T}^{\ast}+\int_{0}^{T}\left\vert A_{\tau}\right\vert ^{2}d\tau\leq C.
\]

Let $f(y)=\left\vert y\right\vert ^{p}$ and $\hat{y}:=y/\left\vert
y\right\vert $ for $y\in\mathbb{R}^{d}.$ Then, for $a,b\in\mathbb{R}^{d},$%
\[
\partial_{a}f(y)=p\left\vert y\right\vert ^{p-1}\hat{y}\cdot a=p\left\vert
y\right\vert ^{p-2}y\cdot a
\]
and
\begin{align*}
\partial_{b}\partial_{a}f(y)  &  =p\left(  p-2\right)  \left\vert y\right\vert
^{p-4}\left(  y\cdot a\right)  \left(  y\cdot b\right)  +p\left\vert
y\right\vert ^{p-2}b\cdot a\\
&  =p\left\vert y\right\vert ^{p-2}\left[  \left(  p-2\right)  \left(  \hat
{y}\cdot a\right)  \left(  \hat{y}\cdot b\right)  +b\cdot a\right]  .
\end{align*}
So by It\^{o}'s formula%
\begin{align*}
d\left\vert Y_{t}\right\vert ^{p}  &  =d\left[  f(Y_{t})\right] \\
&  =p\left\vert Y_{t}\right\vert ^{p-1}\hat{Y}_{t}\cdot dY_{t}+\frac{p}%
{2}\left\vert Y_{t}\right\vert ^{p-2}\left[  \left(  p-2\right)  \left(
\hat{Y}_{t}\cdot dY_{t}\right)  \left(  \hat{Y}_{t}\cdot dY_{t}\right)
+dY_{t}\cdot dY_{t}\right]  .
\end{align*}
Taking expectations of this formula (using $Y$ is a martingale) then gives%
\begin{equation}
\mathbb{E}\left\vert Y_{t}\right\vert ^{p}=\frac{p}{2}\int_{0}^{t}%
\mathbb{E}\left(  \left\vert Y\right\vert ^{p-2}\left[  \left(  p-2\right)
\left(  \hat{Y}\cdot dY\right)  \left(  \hat{Y}\cdot dY\right)  +dY\cdot
dY\right]  \right)  . \label{e.9.3}%
\end{equation}
Using $dY=AdB,$%
\[
dY\cdot dY=Ae_{i}\cdot Ae_{j}dB^{i}dB^{j}=e_{i}\cdot A^{\ast}Ae_{i}%
dt=\operatorname*{tr}(A^{\ast}A)dt=\left\vert A\right\vert ^{2}dt
\]
and%
\begin{align*}
\left(  \hat{Y}\cdot dY\right)  ^{2}  &  =\left(  \hat{Y}\cdot Ae_{i}\right)
\left(  \hat{Y}\cdot Ae_{j}\right)  dB^{i}dB^{j}=\left(  A^{\ast}\hat{Y}\cdot
e_{i}\right)  \left(  A^{\ast}\hat{Y}\cdot e_{i}\right)  dt\\
&  =\left(  A^{\ast}\hat{Y}\cdot A^{\ast}\hat{Y}\right)  dt=\left(  AA^{\ast
}\hat{Y}\cdot\hat{Y}\right)  dt\leq\left\vert A\right\vert ^{2}dt.
\end{align*}
Putting these results back into Eq. (\ref{e.9.3}) implies%
\[
\mathbb{E}\left\vert Y_{t}\right\vert ^{p}\leq\frac{p}{2}(p-1)\int_{0}%
^{t}\mathbb{E}\left(  \left\vert Y_{\tau}\right\vert ^{p-2}\left\vert A_{\tau
}\right\vert ^{2}\right)  d\tau.
\]
By Doob's inequality there is a constant $C_{p}$ (for example $C_{p}=\left[
\frac{p}{p-1}\right]  ^{p}$ will work) such that
\[
\mathbb{E}\left\vert Y_{t}^{\ast}\right\vert ^{p}\leq C_{p}\mathbb{E}%
\left\vert Y_{t}\right\vert ^{p}.
\]
Combining the last two displayed equations implies%
\begin{equation}
\mathbb{E}\left\vert Y_{t}^{\ast}\right\vert ^{p}\leq C\int_{0}^{t}%
\mathbb{E}\left(  \left\vert Y_{\tau}\right\vert ^{p-2}\left\vert A_{\tau
}\right\vert ^{2}\right)  d\tau\leq C\mathbb{E}\left(  \left\vert Y_{t}^{\ast
}\right\vert ^{p-2}\int_{0}^{t}\left\vert A_{\tau}\right\vert ^{2}%
d\tau\right)  . \label{e.9.4}%
\end{equation}
Now applying H\"{o}lder's inequality to the result, with exponents $q=p\left(
p-2\right)  ^{-1}$ and conjugate exponent $q^{\prime}=p/2$ gives%
\[
\mathbb{E}\left\vert Y_{t}^{\ast}\right\vert ^{p}\leq C\left[  \mathbb{E}%
\left\vert Y_{t}^{\ast}\right\vert ^{p}\right]  ^{\frac{p-2}{p}}\left[
\mathbb{E}\left(  \int_{0}^{t}\left\vert A_{\tau}\right\vert ^{2}d\tau\right)
^{p/2}\right]  ^{2/p}%
\]
or equivalently, using $1-\left(  p-2\right)  /p=2/p,$%
\[
\left(  \mathbb{E}\left\vert Y_{t}^{\ast}\right\vert ^{p}\right)  ^{2/p}\leq
C\left[  \mathbb{E}\left(  \int_{0}^{t}\left\vert A_{\tau}\right\vert
^{2}d\tau\right)  ^{p/2}\right]  ^{2/p}.
\]
Taking the $2/p$ roots of this equation then shows
\begin{equation}
\mathbb{E}\left\vert Y_{t}^{\ast}\right\vert ^{p}\leq C\mathbb{E}\left(
\int_{0}^{t}\left\vert A_{\tau}\right\vert ^{2}d\tau\right)  ^{p/2}.
\label{e.9.5}%
\end{equation}

The general case now follows, since when $Y$ is given as in Eq. (\ref{e.9.1})
we have%
\[
Y_{t}^{\ast}\leq\left(  \int_{0}^{\cdot}A_{\tau}dB_{\tau}\right)  _{t}^{\ast
}+\int_{0}^{t}\left\vert \alpha_{\tau}\right\vert d\tau
\]
so that
\begin{align*}
\left\Vert Y_{t}^{\ast}\right\Vert _{p}  &  \leq\left\Vert \left(  \int
_{0}^{\cdot}A_{\tau}dB_{\tau}\right)  _{t}^{\ast}\right\Vert _{p}+\left\Vert
\int_{0}^{t}\left\vert \alpha_{\tau}\right\vert d\tau\right\Vert _{p}\\
&  \leq C\left[  \mathbb{E}\left(  \int_{0}^{t}\left\vert A_{\tau}\right\vert
^{2}d\tau\right)  ^{p/2}\right]  ^{1/p}+\left[  \mathbb{E}\left(  \int_{0}%
^{t}\left\vert \alpha_{\tau}\right\vert d\tau\right)  ^{p}\right]  ^{1/p}%
\end{align*}
and taking the $p^{\text{th}}$ -- power of this equation proves Eq.
(\ref{e.9.2}).
\end{proof}

\begin{remark}
\label{r.9.3}A slightly different application of H\"{o}lder's inequality to
the right side of Eq. (\ref{e.9.4}) gives%
\begin{align*}
\mathbb{E}\left\vert Y_{t}^{\ast}\right\vert ^{p}  &  \leq C\left(  \int
_{0}^{t}\mathbb{E}\left[  \left\vert Y_{t}^{\ast}\right\vert ^{p-2}\left\vert
A_{\tau}\right\vert ^{2}\right]  d\tau\right)  \leq C\left(  \int_{0}%
^{t}\left[  \mathbb{E}\left\vert Y_{t}^{\ast}\right\vert ^{p}\right]
^{\frac{p-2}{p}}\left[  \mathbb{E}\left\vert A_{\tau}\right\vert ^{p}\right]
^{2/p}d\tau\right) \\
&  =\left[  \mathbb{E}\left\vert Y_{t}^{\ast}\right\vert ^{p}\right]
^{\frac{p-2}{p}}C\int_{0}^{t}\left[  \mathbb{E}\left\vert A_{\tau}\right\vert
^{p}\right]  ^{2/p}d\tau
\end{align*}
which leads to the estimate%
\[
\mathbb{E}\left\vert Y_{t}^{\ast}\right\vert ^{p}\leq C\left(  \int_{0}%
^{t}\left[  \mathbb{E}\left\vert A_{\tau}\right\vert ^{p}\right]  ^{2/p}%
d\tau\right)  ^{p/2}.
\]

\end{remark}

Here are some applications of Proposition \ref{p.9.2}.

\begin{proposition}
\label{p.9.4}Let $\{X_{i}\}_{i=0}^{n}$ be a collection of smooth vector fields
on $\mathbb{R}^{N}$ for which $D^{k}X_{i}$ is bounded for all $k\geq1$ and
suppose $\Sigma_{t}$ denotes the solution to Eq. (\ref{e.5.1}) with
$\Sigma_{0}=x\in M:=\mathbb{R}^{N}$ and $\beta=B.$ Then for all $T<\infty$ and
$p<\infty,$%
\begin{equation}
\mathbb{E}\left(  \Sigma_{T}^{\ast}\right)  ^{p}:=\mathbb{E}\left[
\sup_{t\leq T}\left\vert \Sigma_{t}\right\vert ^{p}\right]  <\infty.
\label{e.9.6}%
\end{equation}

\end{proposition}

\begin{proof}
Since
\begin{align*}
X_{i}(\Sigma_{t})\delta B^{i}(t)  &  =X_{i}(\Sigma_{t})dB^{i}(t)+\frac{1}%
{2}d\left[  X_{i}(\Sigma_{t})\right]  \cdot dB^{i}(t)\\
&  =X_{i}(\Sigma_{t})dB^{i}(t)+\frac{1}{2}\left(  \partial_{X_{i}(\Sigma_{t}%
)}X_{i}\right)  (\Sigma_{t})dt,
\end{align*}
the It\^{o} form of Eq. (\ref{e.5.1}) is%
\[
\delta\Sigma_{t}=\left[  X_{0}(\Sigma_{t})+\frac{1}{2}\left(  \partial
_{X_{i}(\Sigma_{t})}X_{i}\right)  (\Sigma_{t})\right]  dt+X_{i}(\Sigma
_{t})dB^{i}(t)\text{ with }\Sigma_{0}=x,
\]
or equivalently,%
\[
\Sigma_{t}=x+\int_{0}^{t}X_{i}(\Sigma_{\tau})dB_{\tau}^{i}+\int_{0}^{t}\left[
X_{0}(\Sigma_{\tau})+\frac{1}{2}\left(  \partial_{X_{i}(\Sigma_{\tau})}%
X_{i}\right)  (\Sigma_{\tau})\right]  d\tau.
\]
By Proposition \ref{p.9.2},%
\begin{align}
\mathbb{E}\left\vert \Sigma_{t}\right\vert ^{p}  &  \leq\mathbb{E}\left(
\Sigma_{t}^{\ast}\right)  ^{p}\leq C_{p}\left\vert x\right\vert ^{p}%
+C_{p}\mathbb{E}\left(  \int_{0}^{t}\left\vert \mathbf{X}(\Sigma_{\tau
})\right\vert ^{2}d\tau\right)  ^{p/2}\nonumber\\
&  +C_{p}\mathbb{E}\left(  \int_{0}^{t}\left\vert X_{0}(\Sigma_{\tau}%
)+\frac{1}{2}\left(  \partial_{X_{i}(\Sigma_{\tau})}X_{i}\right)
(\Sigma_{\tau})\right\vert d\tau\right)  ^{p}. \label{e.9.7}%
\end{align}
Using the bounds on the derivatives of $X$ we learn%
\begin{align*}
&  \left\vert \mathbf{X}(\Sigma_{\tau})\right\vert ^{2}\leq C\left(
1+\left\vert \Sigma_{\tau}\right\vert ^{2}\right)  \text{ and }\\
&  \left\vert X_{0}(\Sigma_{\tau})+\frac{1}{2}\left(  \partial_{X_{i}%
(\Sigma_{\tau})}X_{i}\right)  (\Sigma_{\tau})\right\vert \leq C\left(
1+\left\vert \Sigma_{\tau}\right\vert \right)
\end{align*}
which combined with Eq. (\ref{e.9.7}) gives the estimate%
\begin{align*}
\mathbb{E}\left\vert \Sigma_{t}\right\vert ^{p}  &  \leq\mathbb{E}\left(
\Sigma_{t}^{\ast}\right)  ^{p}\\
&  \leq C_{p}\left\vert x\right\vert ^{p}+C_{p}\mathbb{E}\left(  \int_{0}%
^{t}C\left(  1+\left\vert \Sigma_{\tau}\right\vert ^{2}\right)  d\tau\right)
^{p/2}+C_{p}\mathbb{E}\left(  \int_{0}^{t}C\left(  1+\left\vert \Sigma_{\tau
}\right\vert \right)  d\tau\right)  ^{p}.
\end{align*}
Now assuming $t\leq T<\infty,$ we have by Jensen's (or H\"{o}lder's)
inequality that
\begin{align*}
\mathbb{E}\left\vert \Sigma_{t}\right\vert ^{p}\leq &  \mathbb{E}\left(
\Sigma_{t}^{\ast}\right)  ^{p}\\
\leq &  C\left\vert x\right\vert ^{p}+Ct^{p/2}\mathbb{E}\int_{0}^{t}\left(
1+\left\vert \Sigma_{\tau}\right\vert ^{2}\right)  ^{p/2}\frac{d\tau}{t}\\
&  \qquad+Ct^{p}\mathbb{E}\int_{0}^{t}\left(  1+\left\vert \Sigma_{\tau
}\right\vert \right)  ^{p}\frac{d\tau}{t}\\
\leq &  C\left\vert x\right\vert ^{p}+CT^{\left(  p/2-1\right)  }%
\mathbb{E}\int_{0}^{t}\left(  1+\left\vert \Sigma_{\tau}\right\vert
^{2}\right)  ^{p/2}d\tau\\
&  \qquad+CT^{\left(  p-1\right)  }\mathbb{E}\int_{0}^{t}\left(  1+\left\vert
\Sigma_{\tau}\right\vert \right)  ^{p}d\tau
\end{align*}
from which it follows%
\begin{equation}
\mathbb{E}\left\vert \Sigma_{t}\right\vert ^{p}\leq\mathbb{E}\left(
\Sigma_{t}^{\ast}\right)  ^{p}\leq C\left\vert x\right\vert ^{p}+C(T)\int
_{0}^{t}\left(  1+\mathbb{E}\left\vert \Sigma_{\tau}\right\vert ^{p}\right)
d\tau. \label{e.9.8}%
\end{equation}
An application of Gronwall's inequality now shows $\sup_{t\leq T}%
\mathbb{E}\left\vert \Sigma_{t}\right\vert ^{p}<\infty$ for all $p<\infty$ and
feeding this back into Eq. (\ref{e.9.8}) with $t=T$ proves Eq. (\ref{e.9.6}).
\end{proof}

\begin{proposition}
\label{p.9.5}Suppose $\{X_{i}\}_{i=0}^{n}$ is a collection of smooth vector
fields on $M,$ $\Sigma_{t}$ solves Eq. (\ref{e.5.1}) with $\Sigma_{0}=o\in M$
and $\beta=B,$ $z_{t}$ is the solution to Eq. (\ref{e.5.59}) (i.e.
$z_{t}:=//_{t}^{-1}T_{t\ast o}^{B})\mathbb{\ }$and further
assume\footnote{This will always be true when $M$ is compact.} there is a
constant $K<\infty$ such that $\left\Vert A\left(  m\right)  \right\Vert
_{op}\leq K<\infty$ for all $m\in M,$ where \thinspace$A\left(  m\right)
\in\operatorname*{End}\left(  T_{m}M\right)  $ is defined by%
\[
A\left(  m\right)  v:=\frac{1}{2}\left[  \nabla_{v}\left(  \sum_{i=1}%
^{n}\nabla_{X_{i}}X_{i}+X_{0}\right)  -\sum_{i=1}^{n}R^{\nabla}\left(
v,X_{i}\left(  m\right)  \right)  X_{i}\left(  m\right)  \right]
\]
and%
\[
\sum_{i=1}^{n}\left\vert \nabla_{v}X_{i}\right\vert \leq K\left\vert
v\right\vert \text{ for all }v\in TM.
\]
Then for all $p<\infty$ and $T<\infty,$
\begin{equation}
\mathbb{E}\left[  \sup_{t\leq T}\left\vert z_{t}\right\vert ^{p}\right]
<\infty\label{e.9.9}%
\end{equation}
and%
\begin{equation}
\mathbb{E}\left[  \left(  z_{\cdot}-I\right)  _{t}^{\ast p}\right]  =O\left(
t^{p/2}\right)  \text{ as }t\downarrow0. \label{e.9.10}%
\end{equation}

\end{proposition}

\begin{proof}
In what follows $C$ will denote a constant depending on $K,$ $T$ and $p.$ From
Theorem \ref{t.5.43}, we know that the integrated It\^{o} form of Eq.
(\ref{e.5.59}) is
\begin{equation}
z_{t}=I_{T_{o}M}+\int_{0}^{t}//_{\tau}^{-1}\left(  \nabla_{//_{\tau}z_{\tau
}\left(  \cdot\right)  }\mathbf{X}\right)  dB_{\tau}+\frac{1}{2}A_{//_{\tau}%
}z_{\tau}vd\tau\label{e.9.11}%
\end{equation}
where $A_{//_{t}}:=//_{t}^{-1}A\left(  \Sigma_{t}\right)  //_{t}.$ By
Proposition \ref{p.9.2} and the assumed bounds on $A$ and $\nabla_{\cdot
}\mathbf{X,}$%
\begin{align*}
\mathbb{E}\left(  z_{t}^{\ast}\right)  ^{p}\leq &  C\left\vert I\right\vert
^{p}+C\mathbb{E}\left(  \int_{0}^{t}\sum_{i=1}^{n}\left\vert //_{\tau}%
^{-1}\left(  \nabla_{//_{\tau}z_{\tau}\left(  \cdot\right)  }X_{i}\right)
\right\vert ^{2}d\tau\right)  ^{p/2}\\
&  \qquad+C\mathbb{E}\left(  \int_{0}^{t}\left\vert A_{//_{\tau}}z_{\tau
}\right\vert d\tau\right)  ^{p}\\
\leq &  C+C\mathbb{E}\left(  \int_{0}^{t}\left\vert z_{\tau}\right\vert
^{2}d\tau\right)  ^{p/2}+C\mathbb{E}\left(  \int_{0}^{t}\left\vert z_{\tau
}\right\vert d\tau\right)  ^{p}\\
\leq &  C+C\int_{0}^{t}\mathbb{E}\left\vert z_{\tau}\right\vert ^{p}d\tau
\end{align*}
and
\begin{align}
\mathbb{E}\left[  \left(  z_{\cdot}-I\right)  _{t}^{\ast p}\right]   &  \leq
C\mathbb{E}\left(  \int_{0}^{t}\left\vert z_{\tau}\right\vert ^{2}%
d\tau\right)  ^{p/2}+C\mathbb{E}\left(  \int_{0}^{t}\left\vert z_{\tau
}\right\vert d\tau\right)  ^{p}\nonumber\\
&  \leq C\cdot\mathbb{E}\left\vert z_{t}^{\ast}\right\vert ^{p}\cdot\left(
t^{p/2}+t^{p}\right)  \label{e.9.12}%
\end{align}
where we have made use of H\"{o}lder's (or Jensen's) inequality. Since
\begin{equation}
\mathbb{E}\left\vert z_{t}\right\vert ^{p}\leq\mathbb{E}\left(  z_{t}^{\ast
}\right)  ^{p}\leq C+C\int_{0}^{t}\mathbb{E}\left\vert z_{\tau}\right\vert
^{p}d\tau, \label{e.9.13}%
\end{equation}
Gronwall's inequality implies%
\[
\sup_{t\leq T}\mathbb{E}\left[  \left\vert z_{t}\right\vert ^{p}\right]  \leq
Ce^{CT}<\infty.
\]
Feeding the last inequality back into Eq. (\ref{e.9.13}) shows Eq.
(\ref{e.9.9}). Eq. (\ref{e.9.10}) now follows from Eq. (\ref{e.9.9}). and Eq.
(\ref{e.9.12}).
\end{proof}

\begin{exercise}
\label{exr.9.6}Show under the same hypothesis of Proposition \ref{p.9.5} that
\[
\mathbb{E}\left[  \sup_{t\leq T}\left\vert z_{t}^{-1}\right\vert ^{p}\right]
<\infty
\]
for all $p,T<\infty.$ \textbf{Hint: }Show $z_{t}^{-1}$ satisfies an equation
similar to Eq. (\ref{e.9.11}) with coefficients satisfying the same type of bounds.
\end{exercise}

\subsection{Martingale Estimates\label{s.9.2}}

This section follows the presentation in Norris \cite{Norris86b}.

\begin{lemma}
[Reflection Principle]\label{l.9.7}Let $\beta_{t}$ be a $1$ - dimensional
Brownian motion starting at $0$, $a>0$ and $T_{a}=\inf\left\{  t>0:\beta
_{t}=a\right\}  $ -- be first time $\beta_{t}$ hits height $a,$ see Figure
\ref{f.15}. Then%
\[
P(T_{a}<t)=2P(\beta_{t}>a)=\frac{2}{\sqrt{2\pi t}}\int_{a}^{\infty}%
e^{-x^{2}/2t}dx
\]
%

%TCIMACRO{\FRAME{ftbphFU}{2.3755in}{1.5995in}{0pt}{\Qcb{The first hitting time
%$T_{a}$ of level $a$ by $\beta_{t}.$}}{\Qlb{f.15}}{bmpath.eps}%
%width 2.3755in;  height 1.5995in;  depth 0pt;  original-width 2.6101in;
%original-height 2.3492in;  cropleft "0";  croptop "1";  cropright "1";
%filename 'bmpath.eps';file-properties "XNPEU";}}}%
%BeginExpansion
\begin{figure}
[ptbh]
\begin{center}
\includegraphics[
height=1.5995in,
width=2.3755in
]%
{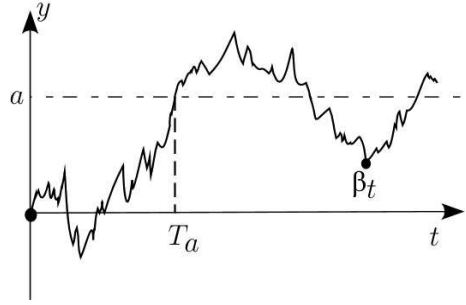}%
\caption{The first hitting time $T_{a}$ of level $a$ by $\beta_{t}.$}%
\label{f.15}%
\end{center}
\end{figure}
%EndExpansion

\end{lemma}

\begin{proof}
Since $P(\beta_{t}=a)=0,$%
\begin{align*}
P(T_{a}<t)  &  =P(T_{a}<t~\&~\beta_{t}>a)+P(T_{a}<t~\&~\beta_{t}<a)\\
&  =P(\beta_{t}>a)+P(T_{a}<t~\&~\beta_{t}<a),
\end{align*}
it suffices to prove%
\[
P(T_{a}<t~\&~\beta_{t}<a)=P(\beta_{t}>a).
\]
To do this define a new process $\tilde{\beta}_{t}$ by%
\[
\tilde{\beta}_{t}=\left\{
\begin{array}
[c]{ccc}%
\beta_{t} & \text{for} & t<T_{a}\\
2a-\beta_{t} & \text{for} & t\geq T_{a}%
\end{array}
\right.
\]
(see Figure \ref{f.16}) and notice that $\tilde{\beta}_{t}$ may also be
expressed as%
\begin{equation}
\tilde{\beta}_{t}=\beta_{t\wedge T_{a}}-1_{t\geq T_{a}}(\beta_{t}%
-\beta_{t\wedge T_{a}})=\int_{0}^{t}\left(  1_{\tau<T_{a}}-1_{\tau\geq T_{a}%
}\right)  d\beta_{\tau}. \label{e.9.14}%
\end{equation}%
%TCIMACRO{\FRAME{ftbphFU}{2.3755in}{1.6859in}{0pt}{\Qcb{The Brownian motion
%$\beta_{t}$ and its reflection $\tilde{\beta}_{t}$ about the line $y=a.$ Note
%that after time $T_{a},$ the labellings of the $\beta_{t}$ and the
%$\tilde{\beta}_{t}$ could be interchanged and the picture would still be
%possible. This should help alleviate the readers fears that Brownian motion
%has some funny asymmetry after the first hitting of level $a.$}}{\Qlb{f.16}%
%}{reflect.eps}{\special{ language "Scientific Word";  type "GRAPHIC";
%width 2.3755in;  height 1.6859in;  depth 0pt;  original-width 2.6101in;
%original-height 1.8434in;  cropleft "0";  croptop "1";  cropright "1";
%filename 'reflect.eps';file-properties "XNPEU";}}}%
%BeginExpansion
\begin{figure}
[ptbh]
\begin{center}
\includegraphics[
height=1.6859in,
width=2.3755in
]%
{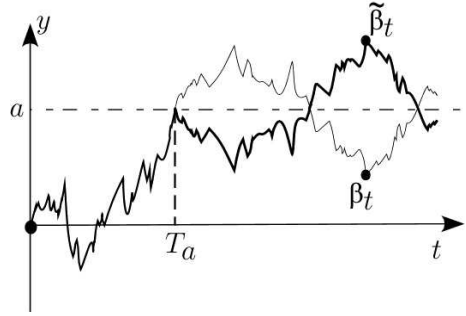}%
\caption{The Brownian motion $\beta_{t}$ and its reflection $\tilde{\beta}%
_{t}$ about the line $y=a.$ Note that after time $T_{a},$ the labellings of
the $\beta_{t}$ and the $\tilde{\beta}_{t}$ could be interchanged and the
picture would still be possible. This should help alleviate the readers fears
that Brownian motion has some funny asymmetry after the first hitting of level
$a.$}%
\label{f.16}%
\end{center}
\end{figure}
%EndExpansion
So $\tilde{\beta}_{t}=\beta_{t}$ for $t\leq T_{a}$ and $\tilde{\beta}_{t}$ is
$\beta_{t}$ reflected across the line $y=a$ for $t\geq T_{a}.$

From Eq. (\ref{e.9.14}) it follows that $\tilde{\beta}_{t}$ is a martingale
and
\[
\left(  d\tilde{\beta}_{t}\right)  ^{2}=\left(  1_{\tau<T_{a}}-1_{\tau\geq
T_{a}}\right)  ^{2}dt=dt
\]
and hence that $\tilde{\beta}_{t}$ is another Brownian motion. Since
$\tilde{\beta}_{t}$ hits level $a$ for the first time exactly when $\beta_{t}$
hits level $a,$%
\[
T_{a}=\tilde{T}_{a}:=\inf\left\{  t>0:\tilde{\beta}_{t}=a\right\}
\]
and $\left\{  \tilde{T}_{a}<t\right\}  =\left\{  T_{a}<t\right\}  .$
Furthermore (see Figure \ref{f.16}),%
\[
\left\{  T_{a}<t~\&~\beta_{t}<a\right\}  =\left\{  \tilde{T}_{a}%
<t~\&~\tilde{\beta}_{t}>a\right\}  =\left\{  \tilde{\beta}_{t}>a\right\}  .
\]
Therefore,
\[
P(T_{a}<t~\&~\beta_{t}<a)=P(\tilde{\beta}_{t}>a)=P(\beta_{t}>a)
\]
which completes the proof.
\end{proof}

\begin{remark}
\label{r.9.8}An alternate way to get a handle on the stopping time $T_{a}$ is
to compute its Laplace transform. This can be done by considering the
martingale
\[
M_{t}:=e^{\lambda\beta_{t}-\frac{1}{2}\lambda^{2}t}.
\]
Since $M_{t}$ is bounded by $e^{\lambda a}$ for $t\in\lbrack0,T_{a}]$ the
optional sampling theorem may be applied to show%
\[
e^{\lambda a}E\left[  e^{-\frac{1}{2}\lambda^{2}T_{a}}\right]  =E\left[
e^{\lambda a-\frac{1}{2}\lambda^{2}T_{a}}\right]  =EM_{T_{a}}=EM_{0}=1,
\]
i.e. this implies that $E\left[  e^{-\frac{1}{2}\lambda^{2}T_{a}}\right]
=e^{-\lambda a}.$ This is equivalent to%
\[
E\left[  e^{-\lambda T_{a}}\right]  =e^{-a\sqrt{2\lambda}}.
\]
From this point of view one would now have to invert the Laplace transform to
get the density of the law of $T_{a}.$
\end{remark}

\begin{corollary}
\label{c.9.9}Suppose now that $T=\inf\left\{  t>0:\left\vert \beta
_{t}\right\vert =a\right\}  ,$ i.e. the first time $\beta_{t}$ leaves the
strip $(-a,a).$ Then
\begin{align}
P(T  &  <t)\leq4P(\beta_{t}>a)=\frac{4}{\sqrt{2\pi t}}\int_{a}^{\infty
}e^{-x^{2}/2t}dx\nonumber\\
&  \leq\min\left(  \sqrt{\frac{8t}{\pi a^{2}}}e^{-a^{2}/2t},1\right)  .
\label{e.9.15}%
\end{align}
Notice that $P(T<t)=P(\beta_{t}^{\ast}\geq a)$ where $\beta_{t}^{\ast}%
=\max\left\{  \left\vert \beta_{\tau}\right\vert :\tau\leq t\right\}  .$ So
Eq. (\ref{e.9.15}) may be rewritten as%
\begin{equation}
P(\beta_{t}^{\ast}\geq a)\leq4P(\beta_{t}>a)\leq\min\left(  \sqrt{\frac
{8t}{\pi a^{2}}}e^{-a^{2}/2t},1\right)  \leq2e^{-a^{2}/2t}. \label{e.9.16}%
\end{equation}

\end{corollary}

\begin{proof}
By definition $T=T_{a}\wedge T_{-a}$ so that $\left\{  T<t\right\}  =\left\{
T_{a}<t\right\}  \cup\left\{  T_{-a}<t\right\}  $ and therefore%
\begin{align*}
P(T<t)  &  \leq P\left(  T_{a}<t\right)  +P\left(  T_{-a}<t\right) \\
&  =2P(T_{a}<t)=4P(\beta_{t}>a)=\frac{4}{\sqrt{2\pi t}}\int_{a}^{\infty
}e^{-x^{2}/2t}dx\\
&  \leq\frac{4}{\sqrt{2\pi t}}\int_{a}^{\infty}\frac{x}{a}e^{-x^{2}%
/2t}dx=\frac{4}{\sqrt{2\pi t}}\left.  \left(  -\frac{t}{a}e^{-x^{2}%
/2t}\right)  \right\vert _{a}^{\infty}=\sqrt{\frac{8t}{\pi a^{2}}}%
e^{-a^{2}/2t}.
\end{align*}
This proves everything but the very last inequality in Eq. (\ref{e.9.16}). To
prove this inequality first observe the elementary calculus inequality:
\begin{equation}
\min\left(  \frac{4}{\sqrt{2\pi}y}e^{-y^{2}/2},1\right)  \leq2e^{-y^{2}/2}.
\label{e.9.17}%
\end{equation}
Indeed Eq. (\ref{e.9.17}) holds $\frac{4}{\sqrt{2\pi}y}\leq2,$ i.e. if $y\geq
y_{0}:=2/\sqrt{2\pi}.$ The fact that Eq. (\ref{e.9.17}) holds for $y\leq
y_{0}$ follows from the following trivial inequality%
\[
1\leq1.4552\cong2e^{-\frac{1}{\pi}}=e^{-y_{0}^{2}/2}.
\]
Finally letting $y=a/\sqrt{t}$ in Eq. (\ref{e.9.17}) gives the last inequality
in Eq. (\ref{e.9.16}).
\end{proof}

\begin{theorem}
\label{t.9.10}Let $N$ be a continuous martingale such that $N_{0}=0$ and $T$
be a stopping time. Then for all $\varepsilon,\delta>0,$%
\[
P\left(  \langle N\rangle_{T}<\varepsilon~\&~N_{T}^{\ast}\geq\delta\right)
\leq P(\beta_{\varepsilon}^{\ast}\geq\delta)\leq2e^{-\delta^{2}/2\varepsilon
}.
\]

\end{theorem}

\begin{proof}
By the Dambis, Dubins \& Schwarz's theorem (see p.174 of \cite{KS91}) we may
write $N_{t}=\beta_{\langle N\rangle_{t}}$ where $\beta$ is a Brownian motion
(on a possibly \textquotedblleft augmented\textquotedblright\ probability
space). Therefore%
\[
\left\{  \langle N\rangle_{T}<\varepsilon~\&~N_{T}^{\ast}\geq\delta\right\}
\subset\left\{  \beta_{\varepsilon}^{\ast}\geq\delta\right\}
\]
and hence from Eq. (\ref{e.9.16}),%
\[
P\left(  \langle N\rangle_{T}<\varepsilon~\&~N_{T}^{\ast}\geq\delta\right)
\leq P(\beta_{\varepsilon}^{\ast}\geq\delta)\leq2e^{-\delta^{2}/2\varepsilon
}.
\]

\end{proof}

\begin{theorem}
\label{t.9.11}Suppose that $Y_{t}=M_{t}+A_{t}$ where $M_{t}$ is a martingale
and $A_{t}$ is a process of bounded variation which satisfy: $M_{0}=A_{0}=0,$
$\left\vert A\right\vert _{t}\leq ct$ and $\langle M\rangle_{t}\leq ct$ for
some constant $c<\infty.$ If $T_{a}:=\inf\left\{  t>0:\left\vert
Y_{t}\right\vert =a\right\}  $ and $t<a/2c,$ then
\[
P(Y_{t}^{\ast}\geq a)=P(T_{a}\leq t)\leq\frac{4}{\sqrt{\pi a}}\exp\left(
-\frac{a^{2}}{8ct}\right)  \text{ }%
\]

\end{theorem}

\begin{proof}
Since
\[
Y_{t}^{\ast}\leq M_{t}^{\ast}+A_{t}^{\ast}\leq M_{t}^{\ast}+\left\vert
A\right\vert _{t}\leq M_{t}^{\ast}+ct
\]
it follows that
\[
\left\{  Y_{t}^{\ast}\geq a\right\}  \subset\left\{  M_{t}^{\ast}\geq
a/2\right\}  \cup\left\{  ct\geq a/2\right\}  =\left\{  M_{t}^{\ast}\geq
a/2\right\}
\]
when $t<a/2c.$ Again by the Dambis, Dubins and Schwarz's theorem (see p.174 of
\cite{KS91}), we may write $M_{t}=\beta_{\langle M\rangle_{t}}$ where $\beta$
is a Brownian motion on a possibly augmented probability space. Since
\[
M_{t}^{\ast}=\max_{\tau\leq\langle M\rangle_{t}}\left\vert \beta_{\tau
}\right\vert \leq\max_{\tau\leq ct}\left\vert \beta_{\tau}\right\vert
=\beta_{ct}^{\ast}%
\]
we learn%
\begin{align*}
P(Y_{t}^{\ast}  &  \geq a)\leq P\left(  M_{t}^{\ast}\geq a/2\right)  \leq
P\left(  \beta_{ct}^{\ast}\geq a/2\right) \\
&  \leq\sqrt{\frac{8ct}{\pi\left(  a/2\right)  ^{2}}}e^{-\left(  a/2\right)
^{2}/2ct}=\sqrt{\frac{8ct}{\pi\left(  a/2\right)  ^{2}}}e^{-\left(
a/2\right)  ^{2}/2ct}\\
&  \leq\sqrt{\frac{8c\left(  a/2c\right)  }{\pi\left(  a/2\right)  ^{2}}%
}e^{-\left(  a/2\right)  ^{2}/2ct}=\frac{4}{\sqrt{\pi a}}\exp\left(
-\frac{a^{2}}{8ct}\right)
\end{align*}
wherein the last inequality we have used the restriction $t<a/2c.$
\end{proof}

\begin{lemma}
\label{l.9.12}If $f:$ $[0,\infty)\rightarrow\mathbb{R}$ is a locally
absolutely continuous function such that $f\left(  0\right)  =0,$ then%
\[
\left\vert f(t)\right\vert \leq\sqrt{2\left\Vert \dot{f}\right\Vert
_{L^{\infty}\left(  [0,t]\right)  }\left\Vert f\right\Vert _{L^{1}([0,t])}%
}~\forall~t\geq0.
\]

\end{lemma}

\begin{proof}
By the fundamental theorem of calculus,%
\[
f^{2}(t)=2\int_{0}^{t}f(\tau)\dot{f}(\tau)d\tau\leq2\left\Vert \dot
{f}\right\Vert _{L^{\infty}\left(  [0,t]\right)  }\left\Vert f\right\Vert
_{L^{1}([0,t])}.
\]

\end{proof}

We are now ready for a key result needed in the probabilistic proof of
H\"{o}rmander's theorem. Loosely speaking it states that if $Y$ is a Brownian
semi-martinagale, then it can happen \textbf{only} with small probability that
the $L^{2}$ -- norm of $Y$ is small while the quadratic variation of $Y$ is
relatively large.

\begin{proposition}
[A key martingale inequality]\label{p.9.13}Let $T$ be a stopping time bounded
by $t_{0}<\infty,$ $Y=y+M+A$ where $M$ is a continuous martingale and $A$ is a
process of bounded variation such that $M_{0}=A_{0}=0.$ Further assume, on the
set $\left\{  t\leq T\right\}  ,$ that $\langle M\rangle_{t}$ and $\left\vert
A\right\vert _{t}$ are absolutely continuous functions and there exists finite
positive constants, $c_{1}$ and $c_{2},$ such that%
\[
\frac{d\langle M\rangle_{t}}{dt}\leq c_{1}\text{ and }\frac{d\left\vert
A\right\vert _{t}}{dt}\leq c_{2}.
\]
Then for all $\nu>0$ and $q>\nu+4$ there exists constants $c=c(t_{0}%
,q,\nu,c_{1},c_{2})>0$ and $\varepsilon_{0}=\varepsilon_{0}(t_{0},q,\nu
,c_{1},c_{2})>0$ such that%
\begin{equation}
P\left(  \int_{0}^{T}Y_{t}^{2}dt<\varepsilon^{q},~\langle Y\rangle_{T}=\langle
M\rangle_{T}\geq\varepsilon\right)  \leq2\exp\left(  -\frac{1}{2c_{1}%
\varepsilon^{\nu}}\right)  =O\left(  \varepsilon^{-\infty}\right)
\label{e.9.18}%
\end{equation}
for all $\varepsilon\in(0,\varepsilon_{0}].$
\end{proposition}

\begin{proof}
Let $q_{0}=\frac{q-\nu}{2}$ (so that $q_{0}\in\left(  2,q/2\right)  ),$
$N:=\int_{0}^{\cdot}YdM$ and%
\begin{equation}
C_{\varepsilon}:=\left\{  \langle N\rangle_{T}\leq c_{1}\varepsilon
^{q},~\text{ }N_{T}^{\ast}\geq\varepsilon^{q_{0}}\right\}  . \label{e.9.19}%
\end{equation}

We will show shortly that for $\varepsilon$ sufficiently small,%
\begin{equation}
B_{\varepsilon}:=\left\{  \int_{0}^{T}Y_{t}^{2}dt<\varepsilon^{q},~\langle
Y\rangle_{T}\geq\varepsilon\right\}  \subset C_{\varepsilon}. \label{e.9.20}%
\end{equation}
By an application of Theorem \ref{t.9.10},%
\[
P(C_{\varepsilon})\leq2\exp\left(  -\frac{\varepsilon^{2q_{0}}}{2c_{1}%
\varepsilon^{q}}\right)  =2\exp\left(  -\frac{1}{2c_{1}\varepsilon^{v}%
}\right)
\]
and so assuming the validity of Eq. (\ref{e.9.20}),%
\begin{equation}
P\left(  \int_{0}^{T}Y_{t}^{2}dt<\varepsilon^{q},~\langle Y\rangle_{T}%
\geq\varepsilon\right)  \leq P(C_{\varepsilon})\leq2\exp\left(  -\frac
{1}{2c_{1}\varepsilon^{v}}\right)  \label{e.9.21}%
\end{equation}
which proves Eq. (\ref{e.9.18}). So to finish the proof it only remains to
verify Eq. (\ref{e.9.20}) which will be done by showing $B_{\varepsilon}\cap
C_{\varepsilon}^{c}=\emptyset.$

For the rest of the proof, it will be assumed that we are on the set
$B_{\varepsilon}\cap C_{\varepsilon}^{c}.$ Since $\langle N\rangle_{T}%
=\int_{0}^{T}\left\vert Y_{t}\right\vert ^{2}d\langle M\rangle_{t},$ we have%
\begin{equation}
B_{\varepsilon}\cap C_{\varepsilon}^{c}=\left\{  \int_{0}^{T}Y_{t}%
^{2}dt<\varepsilon^{q},~\langle Y\rangle_{T}\geq\varepsilon,~\int_{0}%
^{T}\left\vert Y_{t}\right\vert ^{2}d\langle M\rangle_{t}>c_{1}\varepsilon
^{q},~\text{ }N_{T}^{\ast}<\varepsilon^{q_{0}}\right\}  . \label{e.9.22}%
\end{equation}
From Lemma \ref{l.9.12} with $f(t)=\langle Y\rangle_{t}$ and the assumption
that $d\langle Y\rangle_{t}/dt\leq c_{1},$%
\begin{equation}
\langle Y\rangle_{T}\leq\sqrt{2\left\Vert \dot{f}\right\Vert _{L^{\infty
}\left(  [0,T]\right)  }\left\Vert f\right\Vert _{L^{1}([0,T])}}\leq
\sqrt{2c_{1}\int_{0}^{T}\langle Y\rangle_{t}dt}. \label{e.9.23}%
\end{equation}
By It\^{o}'s formula, the quadratic variation, $\langle Y\rangle_{t},$ of $Y$
satisfies%
\begin{equation}
\langle Y\rangle_{t}=Y_{t}^{2}-y^{2}-2\int_{0}^{t}YdY\leq Y_{t}^{2}%
+2\left\vert \int_{0}^{t}YdY\right\vert \label{e.9.24}%
\end{equation}
and on the set $\left\{  t\leq T\right\}  \cap B_{\varepsilon}\cap
C_{\varepsilon}^{c},$%
\begin{align}
\left\vert \int_{0}^{t}YdY\right\vert  &  =\left\vert \int_{0}^{t}YdM+\int
_{0}^{t}YdA\right\vert \leq\left\vert N_{t}\right\vert +\int_{0}^{t}\left\vert
Y\right\vert dA\nonumber\\
&  \leq N_{T}^{\ast}+c_{2}\int_{0}^{T}\left\vert Y_{\tau}\right\vert d\tau
\leq\varepsilon^{q_{0}}+c_{2}T^{1/2}\sqrt{\int_{0}^{T}Y_{\tau}^{2}d\tau
}\nonumber\\
&  \leq\varepsilon^{q_{0}}+c_{2}t_{0}^{1/2}\varepsilon^{q}. \label{e.9.25}%
\end{align}
Combining Eqs. (\ref{e.9.24}) and (\ref{e.9.25}) shows, on the set $\left\{
t\leq T\right\}  \cap B_{\varepsilon}\cap C_{\varepsilon}^{c}$ that%
\[
\langle Y\rangle_{t}\leq Y_{t}^{2}+2\left[  \varepsilon^{q_{0}}+c_{2}%
t_{0}^{1/2}\varepsilon^{q}\right]
\]
and using this in Eq. (\ref{e.9.23}) implies%
\begin{align}
\langle Y\rangle_{T}  &  \leq\sqrt{2c_{1}\int_{0}^{T}\left(  Y_{t}%
^{2}+2\left[  \varepsilon^{q_{0}}+c_{2}t_{0}^{1/2}\varepsilon^{q}\right]
\right)  dt}\nonumber\\
&  \leq\sqrt{2c_{1}\left[  \varepsilon^{q}+2\left[  \varepsilon^{q_{0}}%
+c_{2}t_{0}^{1/2}\varepsilon^{q}\right]  t_{0}\right]  }=O\left(
\varepsilon^{q_{0}}\right)  =o\left(  \varepsilon\right)  . \label{e.9.26}%
\end{align}
Hence we may choose $\varepsilon_{0}=\varepsilon_{0}\left(  c_{1},c_{2}%
,t_{o},q,\nu\right)  >0$ such that for $\varepsilon\leq\varepsilon_{0}$ we
have
\[
\sqrt{2c_{1}\left(  \varepsilon^{q}+2\varepsilon^{q_{0}}t_{0}+2c_{2}%
t_{0}^{3/2}\varepsilon^{q/2}\right)  }<\varepsilon
\]
and hence on $B_{\varepsilon}\cap C_{\varepsilon}^{c}$ we learn $\varepsilon
\leq\langle Y\rangle_{T}<\varepsilon$ which is absurd. So we must conclude
that $B_{\varepsilon}\cap C_{\varepsilon}^{c}=\emptyset.$
\end{proof}

\newpage

%\bibliographystyle{amsplain}
%\bibliography{eth,newbib,density2}

\providecommand{\bysame}{\leavevmode\hbox
to3em{\hrulefill}\thinspace}
\providecommand{\MR}{\relax\ifhmode\unskip\space\fi MR }
%\MRhref is called by the amsart/book/proc definition of \MR.
\providecommand{\MRhref}[2]{  \href{http://www.ams.org/mathscinet-getitem?mr=#1}{#2}
} \providecommand{\href}[2]{#2}

\end{document}